\documentclass[mnsc,nonblindrev]{informs3} 

\pdfoutput=1

\usepackage{epsfig}
\usepackage{centernot}
\usepackage[colorlinks,citecolor=blue,urlcolor=blue,filecolor=blue,hypertexnames=false]{hyperref} 
\usepackage{url}
\usepackage{multirow}
\usepackage{graphicx}
\usepackage{enumitem}
\usepackage{lscape}
\usepackage{dsfont}
\usepackage[ruled,vlined]{algorithm2e}
\usepackage{algorithmic}
\usepackage{placeins}
\usepackage{mathtools,amsmath,accents}
\usepackage{bm}

\usepackage{etoolbox}
\newcommand{\smalldisplayskips}{%
  \setlength{\abovedisplayskip}{3pt}%
  \setlength{\belowdisplayskip}{3pt}%
  \setlength{\abovedisplayshortskip}{3pt}%
  \setlength{\belowdisplayshortskip}{3pt}}
\appto{\normalsize}{\smalldisplayskips}
\appto{\small}{\smalldisplayskips}
\appto{\footnotesize}{\smalldisplayskips}

\makeatletter
\newcommand{\substackalign}[1]{%
  \vcenter{%
    \Let@ \restore@math@cr \default@tag
    \baselineskip\fontdimen10 \scriptfont\tw@
    \advance\baselineskip\fontdimen12 \scriptfont\tw@
    \lineskip\thr@@\fontdimen8 \scriptfont\thr@@
    \lineskiplimit\lineskip
    \ialign{\hfil$\m@th\scriptstyle##$&$\m@th\scriptstyle{}##$\hfil\crcr
      #1\crcr
    }%
  }%
}
\makeatother

\makeatletter
\DeclareRobustCommand{\cev}[1]{%
  \mathpalette\do@cev{#1}%
}
\newcommand{\do@cev}[2]{%
  \fix@cev{#1}{+}%
  \reflectbox{$\m@th#1\vec{\reflectbox{$\fix@cev{#1}{-}\m@th#1#2\fix@cev{#1}{+}$}}$}%
  \fix@cev{#1}{-}%
}
\newcommand{\fix@cev}[2]{%
  \ifx#1\displaystyle
    \mkern#23mu
  \else
    \ifx#1\textstyle
      \mkern#23mu
    \else
      \ifx#1\scriptstyle
        \mkern#22mu
      \else
        \mkern#22mu
      \fi
    \fi
  \fi
}
  
\makeatother

\usepackage{etoolbox} 

\makeatletter
\patchcmd{\@algocf@start}{-1.5em}{0em}{}{} 
\makeatother

\newcommand{\abovebelowset}[3]{\:
\overset{
    \text{\raisebox{1.2ex}{\smash{\scalebox{1.0}{$#1$}}}}%
}
{
\underset{
    \text{\raisebox{1.3ex}{\smash{\scalebox{1.0}{$#2$}}}}%
}
{
    \text{\raisebox{0.2ex}{\smash{$#3$}}}
}
}
\:}

\newcommand{\sspar}{{\normalfont{\text{par}}}}
\newcommand{\ssapx}{{\normalfont{\text{apx}}}}
\newcommand{\ssmin}{{\normalfont{\text{min}}}}
\newcommand{\ssmax}{{\normalfont{\text{max}}}}
\newcommand{\MT}{{\normalfont{\text{MT}}}}

\newcommand{\LP}{{\normalfont{\text{LP}}}}
\newcommand{\sssub}{{\normalfont{\text{sub}}}}
\newcommand{\sstheo}{{\normalfont{\text{theo}}}}

\newcommand{\ssstep}{{\normalfont{\text{step}}}}
\newcommand{\ssUB}{{\normalfont{\text{ub}}}}
\newcommand{\ssLB}{{\normalfont{\text{lb}}}}

\newcommand{\LPTAG}[1]{{\normalfont{(\hyperref[eqn:mt-tf-lp-any-index]{$\LP$-{#1}})}}}

\newcommand{\CPDLPDUAL}[1]{{\normalfont{(\hyperref[eqn:mt-tf-algo-lp]{\MT_{\sspar}^{*({#1})}})}}}
\newcommand{\CPDLPPRIMAL}[1]{{\normalfont{(\hyperref[eqn:mt-tf-algo-lp-dual]{\MT_{\sspar}^{({#1})}})}}}

\usepackage{amsfonts,mathrsfs,bbm,tikz}
\usetikzlibrary{calc}

\newcommand{\Td}{\mathrm{d}}

\newcommand{\BA}{\mathbf{A}}
\newcommand{\BB}{\mathbf{B}}

\newcommand{\BD}{\mathbf{D}}

\newcommand{\BH}{\mathbf{H}}
\newcommand{\BI}{\mathbf{I}}

\newcommand{\BO}{\mathbf{O}}
\newcommand{\BP}{\mathbf{P}}
\newcommand{\BQ}{\mathbf{Q}}
\newcommand{\BR}{\mathbf{R}}

\newcommand{\BT}{\mathbf{T}}
\newcommand{\BU}{\mathbf{U}}
\newcommand{\BV}{\mathbf{V}}
\newcommand{\BW}{\mathbf{W}}

\newcommand{\BIa}{{\boldsymbol{a}}}
\newcommand{\BIb}{{\boldsymbol{b}}}
\newcommand{\BIc}{{\boldsymbol{c}}}

\newcommand{\BIe}{{\boldsymbol{e}}}

\newcommand{\BIg}{{\boldsymbol{g}}}

\newcommand{\BIk}{{\boldsymbol{k}}}

\newcommand{\BIr}{{\boldsymbol{r}}}
\newcommand{\BIs}{{\boldsymbol{s}}}

\newcommand{\BIu}{{\boldsymbol{u}}}
\newcommand{\BIv}{{\boldsymbol{v}}}
\newcommand{\BIw}{{\boldsymbol{w}}}
\newcommand{\BIx}{{\boldsymbol{x}}}
\newcommand{\BIy}{{\boldsymbol{y}}}
\newcommand{\BIz}{{\boldsymbol{z}}}

\newcommand{\BIX}{{\boldsymbol{X}}}

\newcommand{\BIZ}{{\boldsymbol{Z}}}

\newcommand{\CB}{\mathcal{B}}
\newcommand{\CC}{\mathcal{C}}

\newcommand{\CF}{\mathcal{F}}
\newcommand{\CG}{\mathcal{G}}

\newcommand{\CI}{\mathcal{I}}

\newcommand{\CK}{\mathcal{K}}

\newcommand{\CP}{\mathcal{P}}

\newcommand{\CV}{\mathcal{V}}

\newcommand{\CX}{\mathcal{X}}

\newcommand{\FC}{\mathfrak{C}}
\newcommand{\FD}{\mathfrak{D}}

\newcommand{\FF}{\mathfrak{F}}

\newcommand{\FK}{\mathfrak{K}}

\newcommand{\FM}{\mathfrak{M}}

\newcommand{\FS}{\mathfrak{S}}

\newcommand{\BTheta}{{\boldsymbol{\Theta}}}

\newcommand{\Btheta}{{\boldsymbol{\theta}}}

\newcommand{\Blambda}{{\boldsymbol{\lambda}}}
\newcommand{\Bmu}{{\boldsymbol{\mu}}}
\newcommand{\Bnu}{{\boldsymbol{\nu}}}
\newcommand{\Bxi}{{\boldsymbol{\xi}}}

\newcommand{\Brho}{{\boldsymbol{\rho}}}

\newcommand{\Bchi}{{\boldsymbol{\chi}}}
\newcommand{\Bpsi}{{\boldsymbol{\psi}}}
\newcommand{\Bomega}{{\boldsymbol{\omega}}}

\newcommand{\R}{\mathbb{R}} 
\newcommand{\N}{\mathbb{N}} 
\newcommand{\vecone}{\mathbf{1}} 
\newcommand{\veczero}{\mathbf{0}} 
\newcommand{\INDI}{\mathbbm{1}} 
\newcommand{\TRANSP}{\top} 
\newcommand{\DIFF}{\Td} 
\newcommand{\DIFFX}[1]{\,\Td{#1}} 
\newcommand{\DIFFM}[2]{\,{#1}({#2})}
\newcommand{\diag}{\mathrm{diag}} 
\newcommand{\PROB}{\mathbb{P}} 
\newcommand{\EXP}{\mathbb{E}} 
\newcommand{\clos}{\mathrm{cl}} 
\newcommand{\inter}{\mathrm{int}} 
\newcommand{\relint}{\mathrm{relint}} 
\newcommand{\support}{\mathrm{supp}} 
\newcommand{\conv}{\mathrm{conv}} 
\newcommand{\aff}{\mathrm{aff}} 
\newcommand{\cone}{\mathrm{cone}} 
\newcommand{\volume}{\mathrm{vol}} 
\newcommand{\epigraph}{\mathrm{epi}} 
\newcommand{\hypograph}{\mathrm{hyp}} 
\newcommand{\tr}{\mathrm{tr}} 
\newcommand{\Lebesgue}{\mathscr{L}}

\DeclareFontFamily{U}{mathx}{\hyphenchar\font45}
\DeclareFontShape{U}{mathx}{m}{n}{
      <5> <6> <7> <8> <9> <10>
      <10.95> <12> <14.4> <17.28> <20.74> <24.88>
      mathx10
      }{}
\DeclareSymbolFont{mathx}{U}{mathx}{m}{n}
\DeclareMathSymbol{\bigtimes}{1}{mathx}{"91}

\newcommand{\overbar}[1]{\mkern 6.5mu\overline{\mkern-5.5mu#1\mkern-2.5mu}\mkern 2.5mu}

\DeclareMathOperator*{\minimize}{\mathrm{minimize}} 
\DeclareMathOperator*{\maximize}{\mathrm{maximize}} 

\usepackage{natbib}[longnamefirst]
 \bibpunct[, ]{(}{)}{,}{a}{}{,}%
 \def\bibfont{\small}%
 \def\bibsep{\smallskipamount}%
 \def\BIBand{and}%

\TheoremsNumberedBySection  
\ECRepeatTheorems

\EquationsNumberedBySection 

\MANUSCRIPTNO{MS-OPT-2024-05071.R2}

\begin{document}


\RUNAUTHOR{Neufeld and Xiang}

\RUNTITLE{Feasible approximation for large-scale matching for teams problems}

\TITLE{Feasible approximation of matching equilibria for large-scale matching for teams problems}

\ARTICLEAUTHORS{%
    \AUTHOR{Ariel Neufeld and Qikun Xiang}
    \AFF{Division of Mathematical Sciences, Nanyang Technological University, 21 Nanyang Link, 637371 Singapore 
    \EMAIL{ariel.neufeld@ntu.edu.sg}, 
    \EMAIL{qikun.xiang@ntu.edu.sg}}
} 

\ABSTRACT{%
We propose a numerical algorithm for computing feasible and approximately optimal solutions of the matching for teams problem.
Specifically, 
we introduce the notion of approximate matching equilibrium as a feasible approximation of a matching equilibrium with relaxed rationality,
and we show that a true equilibrium is recovered in the limit of 
a sequence of approximate matching equilibria with sub-optimality approaching~0.
In our approximation scheme,
we parametrize the so-called transfer functions, 
and we show that tackling the resulting parametric primal and dual optimization problems
yields two approximate matching equilibria
as well as provable and computable lower and upper bounds 
for the optimal social welfare. 
Under a flexible Euclidean setting, 
we show that 
the approximation error of our scheme can be controlled to be arbitrarily close to~0,
we derive an explicit computational complexity bound,
and we develop an algorithm for computing approximate matching equilibria that is efficient for large-scale problems involving a large number of agent populations.
We study three problems in our numerical experiments: a retail business problem, the Wasserstein barycenter problem, and a large-scale problem involving up to 1000 agent populations. 
We show that the proposed algorithm can produce nearly optimal approximate matching equilibria to 
provide quantitative managerial insights for policymakers, 
and that the computed sub-optimality estimates are much less conservative than theoretical estimates.
}%


\renewcommand{\KEYWORDSname}{{\it Keywords\/}{\kern0.7pt}:\enskip}

\KEYWORDS{matching for teams, equilibrium computation, Wasserstein barycenter, optimal transport.}

\maketitle

%

\section{Introduction}

In this paper, we develop a numerical algorithm for computing approximately optimal solutions of a general matching problem with transferable utilities introduced by \citet{carlier2010matching},
and we apply our algorithm to study three concrete problems,
including a retail business problem in which we draw concrete managerial insights for improving the social welfare.
The matching game involves possibly infinitely many agents 
divided into $N\ge 2$ equally-sized and heterogeneous populations 
$\CX_1,\ldots,\CX_N$,
and allocates the agents into teams 
where each team consists of exactly one agent from each population.
This is carried out by letting
each agent $x_i\in\nobreak\CX_i$ choose 
$z$ out of a set $\CX_0$ 
that maximizes his/her utility independently of all other agents,
and placing agents that have chosen the same~$z$ into the same team,
where members of each team then transfer utilities with each other in a transaction.
The utility maximized by agent~$x_i$ is given by 
$\varphi_i(z)-c_i(x_i,z)$,
where 
$\varphi_i(z)$ is the net utility transferred from other members of team~$z$ and 
$c_i(x_i,z)$ represents the cost.
An equilibrium of the matching game, 
termed \textit{matching equilibrium},
refers to a state in which no agent has incentive to deviate from the existing team allocations.

When $N=2$, this is a widely adopted model for two-sided interactions in economics, typically between buyers and sellers \citep*{chiappori2010hedonic},
such as
managerial talent markets \citep*{gabaix2008why, tervio2008difference}, 
real estate markets \citep*{shapley1972assignment},
and marriage \citep*{becker1973theory, chiappori2017partner, galichon2022cupid}.
To illustrate the case with $N\ge 3$ agent populations, 
let us consider 
a retail business in a city hiring $N-1$ categories of employees.
Let $\CX_i\subset\R^2$ represent the residential locations of the employees in category~$i$ for $i=1,\ldots,N-1$, 
let $\CX_N\subset\R^2$ represent the locations of the suppliers of the business, 
and let $\CX_0\subset\R^2$ represent the locations of retail outlets. 
To open a retail outlet at $z\in\CX_0$,
a team comprising one employee from each category and a supplier needs to be formed.
In the matching game,
an employee residing at $x_i\in\CX_i$ chooses to work at outlet $z\in\CX_0$ by 
maximizing the difference between the salary $\varphi_i(z)$ 
and the commuting cost $c_i(x_i,z)$,
whereas 
each business owner with supplier at $x_N\in\CX_N$ 
chooses outlet $z\in\CX_0$ 
to minimize the sum of the total salary payout $-\varphi_N(z)$ and the restocking cost $c_N(x_N,z)$  
(one could consider $c_1,\ldots,c_N$ to be distance functions in this example).
One can similarly use this model to describe 
other decentralized team building or consensus making processes involving multiple populations of interacting agents,
such as the 
construction of universities/hospitals hiring professors/doctors with different expertise.

Let us briefly recall the mathematical formulation of 
matching for teams by \citet{carlier2010matching},
and we defer a more comprehensive discussion of the model assumptions and their economic interpretations to Section~\ref{sec:model}.
The sets $(\CX_i)_{i=1:N}$ describing agents 
are termed
\textit{type spaces}.
Rather than describing finitely many agents,
the model allows the number of agents to be infinite and 
models the $i$-th agent population by a distribution on 
the compact metric space $\CX_i$.
Because the agent populations are assumed to be equally-sized,
each distribution is normalized to~1 and described by 
a probability measure
$\mu_i\in\CP(\CX_i)$
termed the \textit{type measure}.
Note that a finite population can be described by an empirical measure and is thus subsumed under this general model.
The set $\CX_0$ describing the teams is a compact metric space termed the
\textit{quality space},
where each $z\in\CX_0$ can be interpreted as describing the quality of the good transacted by the members of team~$z$.

Subsequently, each potential state of the matching game is described by
probability measures on the products of the type spaces and the quality space,
i.e.,
$\gamma_1\in\CP(\CX_1\times\CX_0),\ldots,\gamma_N\in\CP(\CX_N\times\CX_0)$.
Due to the requirement that each team consists of exactly one agent from each population, 
the marginal of $\gamma_i$ on $\CX_i$ is required to be equal to $\mu_i$
for $i=1,\ldots,N$,
whereas the marginals of $\gamma_1,\ldots,\gamma_N$ on $\CX_0$ must all be equal to some \textit{quality measure} $\nu\in\CP(\CX_0)$.
We call any $(\gamma_i)_{i=1:N},\nu$ that satisfy these marginal constraints \textit{team allocations}.
The transfers of utilities within teams are described by 
continuous functions $\varphi_1,\ldots,\varphi_N:\CX_0\to\R$
termed the \textit{transfer functions}. 
They are required to satisfy the \textit{balance condition}
$\sum_{i=1}^{N}\varphi_i=0$,
which reflects the constraint that each team is self-financed.
They can be interpreted as a price system
that drives the agents into forming distinct teams.
In the retail example above, 
the balance condition asserts that at every outlet $z\in\CX_0$,
the total salary payout,
i.e., $-\varphi_N(z)$,
is equal to the sum of the salaries received by the employees,
i.e., $\sum_{i=1}^{N-1}\varphi_i(z)$.

Any transfer functions and team allocations $(\varphi_i,\gamma_i)_{i=1:N},\nu$
are said to constitute a \textit{matching equilibrium} if, 
for $i=1,\ldots,N$, 
the \textit{rationality condition} 
$\varphi_i^{c_i}(x_i)+\varphi_i(z)=c_i(x_i,z)$ holds for 
$\gamma_i$-almost all $(x_i,z)\in\CX_i\times\CX_0$.
Here, $\varphi_i^{c_i}(x):=\inf_{z\in\CX_0}\big\{c_i(x,z)-\varphi_i(z)\big\}$ $\forall x\in\CX_i$ is called the \textit{$c_i$-transform} of~$\varphi_i$,
and the minimization $\inf_{z\in\CX_0}\big\{c_i(x,z)-\varphi_i(z)\big\}$ characterizes the rational decision
of each agent in the $i$-th population given 
the transfer function $\varphi_i$.
Thus, at a matching equilibrium,
for $\gamma_i$-almost all $(x_i,z)\in\CX_i\times\CX_0$, 
the agent~$x_i$ is allocated to the team~$z$ that maximizes his/her utility,
and there is no incentive for any agent to deviate from the team allocations
$(\gamma_i)_{i=1:N},\nu$.

Notably, the structure of the matching for teams problem also describes a class of geometry-aware distribution aggregation problems.
In particular, in the 2-Wasserstein barycenter problem \citep{agueh2011barycenters},
$(\CX_i)_{i=1:N}$ are subsets of the same Euclidean space $\R^d$,
$\CX_0\subset\R^d$ is compact, convex, and satisfies
$\CX_0\supseteq \bigcup_{i=1}^{N}\CX_i$,
and
$c_i(\BIx,\BIz):=\lambda_i\|\BIx-\nobreak\BIz\|_2^2$
is the scaled squared Euclidean distance 
with weights $(\lambda_i)_{i=1:N}\subset(0,1]$, $\sum_{i=1}^{N}\lambda_i=\nobreak 1$.
Then, 
$\nu\in\CP(\CX_0)$ in a matching equilibrium $(\varphi_i,\gamma_i)_{i=1:N},\nu$
is called a 2-Wasserstein barycenter of $(\mu_i)_{i=1:N}$ with weights $(\lambda_i)_{i=1:N}$.
In recent years, the 2-Wasserstein barycenter has found widespread applications in fields such as
statistical inference \citep*{srivastava2015wasp,srivastava2018scalable, li2020continuous, bigot2019penalization}, 
unsupervised clustering \citep*{ye2014scaling,ye2017fast,puccetti2020on}, 
pattern recognition \citep*{tabak2022distributional}, 
texture mixing \citep*{rabin2011wasserstein},
color transfer \citep*{kuang2019sample}, 
shape interpolation \citep*{solomon2015convolutional, lindheim2023simple}, etc.
See also \citep*{carlier2024wasserstein} for the properties of \mbox{1-Wasserstein} barycenters.

Many classical studies in theoretical economics 
(see, e.g., \citep{gabaix2008why, tervio2008difference})
considered one-dimensional type spaces $(\CX_i)_{i=1:N}$
and cost functions with certain monotonicity structures that 
lead to a unique matching equilibrium involving 
\textit{positive assortative matchings}, 
that is,
agents are sorted from ``low'' to ``high'' based on their types and 
are matched monotonically to teams,
where agents with ``lower'' types are matched to the same team and agents with ``higher'' types are matched to the same team.
These assumptions are often made due to their mathematical convenience rather than practical relevance.
In this paper, we do not impose these assumptions and instead
consider the general model presented above 
with possibly multi-dimensional type spaces
and any Lipschitz continuous cost functions.
As such, 
there is little hope of analytically characterizing a matching equilibrium;
in fact,
we do not even require the matching equilibrium to be unique or \textit{pure} (i.e., each agent has a unique optimal choice of team, see \citep[Section~7]{carlier2010matching} as well as \citep{ekeland2010,pass2014multi}).

In this paper, we develop an efficient numerical algorithm 
for computing states of the matching game that are \textit{approximately at equilibrium} up to a small \textit{error tolerance} $\epsilon>0$,
especially for large-scale problems involving a large number of agent populations, i.e., when $N$ is large.
Concretely,
for any transfer functions $(\varphi_i)_{i=1:N}$ and
team allocations $(\gamma_i)_{i=1:N},\nu$,
and for any error $\epsilon>0$,
we say that 
$(\varphi_i,\gamma_i)_{i=1:N},\nu$ constitute an
\textit{$\epsilon$-matching equilibrium}
if 
$\sum_{i=1}^{N}\int_{\CX_i\times\CX_0} \big(c_i(x,z)-\varphi_i(z) - \varphi_i^{c_i}(x)\big) \DIFFM{\gamma_i}{\DIFF x,\DIFF z} \le \epsilon$ holds.
Observe that 
this is a \textit{feasible approximation} 
of matching equilibrium, 
since the team allocations $(\gamma_i)_{i=1:N},\nu$ 
and the transfer functions $(\varphi_i)_{i=1:N}$ in an $\epsilon$-matching equilibrium are required to obey 
the marginal constraints and the balance condition, respectively.
The $\epsilon$-relaxation only lies in the rationality condition,
that is,
rather than requiring agents to be perfectly rational and join the optimal teams, 
we allow a small degree of irrationality in their team choices
and instead require the total sub-optimality aggregated over all agents 
to be bounded above by the error tolerance~$\epsilon$.
We will show that an $\epsilon$-matching equilibrium 
converges to a true matching equilibrium when $\epsilon$ approaches~0 
(Theorem~\ref{thm:mt-tf-equilibrium}).
Moreover, since every approximate matching equilibrium (AME)
describes a state of the matching game which approximately maximizes social welfare (we will discuss this in detail in Section~\ref{sec:model}),
computing and examining AMEs
can provide managerial insights to 
policymakers and regulators in the following three ways.
First, 
every \mbox{$\epsilon$-matching} equilibrium naturally provides lower and upper bounds for the optimal social welfare that are at most $\epsilon$ apart.
Second, 
with an $\epsilon$-matching equilibrium,
policymakers and regulators are able to compute lower and upper bounds for the improvement of the optimal social welfare when perturbing the cost functions $(c_i)_{i=1:N}$, 
e.g., via taxation or via infrastructure investments;
see Proposition~\ref{prop:perturbation} and the discussion thereafter.
Third,
by computing 
AMEs under scenarios with varying type measures $(\mu_i)_{i=1:N}$,
policymakers and regulators 
can examine how the equilibrium would evolve when the agent types gradually change.
Concretely, in the aforementioned retail example, 
city planners are able to infer how traffic patterns react to changing residential distributions,
which can then guide their decisions to 
improve transportation efficiency and social welfare.

\noindent\textbf{Our contributions.}
First,
we introduce the notion of $\epsilon$-matching equilibrium along with its economic interpretations, 
and we examine its theoretical properties (see Theorem~\ref{thm:mt-tf-equilibrium} and Proposition~\ref{prop:perturbation}).
Second, we introduce a pair of parametric primal and dual formulations 
\eqref{eqn:mt-tf-dual} and \eqref{eqn:mt-tf-lsip}
that are linear semi-infinite programming (LSIP) problems,
and we \textit{explicitly construct
two approximate matching equilibria}
from a pair of approximate optimizers of 
\eqref{eqn:mt-tf-dual} and \eqref{eqn:mt-tf-lsip} (see Theorem~\ref{thm:mt-tf-approx}).
Third, under a Euclidean setting (Setting~\ref{sett:Euclidean})
that allows us to control the approximation error to be arbitrarily small (see Theorem~\ref{thm:error-control}),
we derive an explicit computational complexity bound for computing an 
\mbox{$\epsilon$-matching} equilibrium for any $\epsilon>0$ (see Theorem~\ref{thm:complexity-explicit}),
and we develop a practically efficient numerical algorithm (Algorithm~\ref{alg:mt-tf}) for computing approximate matching equilibria as well as 
lower and upper bounds for the optimal social welfare (see Theorem~\ref{thm:mt-tf}).
Compared to existing algorithms,
our algorithm 
is applicable to general cost functions $(c_i)_{i=1:N}$ and 
general non-discrete, non-parametric type measures $(\mu_i)_{i=1:N}$,
does not discretize $(\CX_i)_{i=0:N}$,
does not pre-specify a finite support of the quality measure $\nu$
(that is, our algorithm belongs to the so-called \textit{free-support} algorithms),
and, most importantly,
computes sub-optimality bounds that are typically much less conservative than sub-optimality bounds derived from purely theoretical analyses. 
Lastly,
we showcase our algorithm in three numerical experiments.
Specifically, 
we draw concrete managerial insights for improving social welfare 
in the aforementioned retail business problem (see Section~\ref{ssec:experiment-biz}),
we compare our algorithm with state-of-the-art \mbox{2-Wasserstein} barycenter algorithms to show its superior performance (see Section~\ref{ssec:experiment-WB}),
and we demonstrate that our algorithm is capable of solving large-scale matching for teams problems with $N=1000$ agent populations (see Section~\ref{ssec:experiment-1d}).
The e-companion (i.e., online appendices) contains supplementary remarks, discussions, as well as the proofs.

\noindent\textbf{Related work.}
It is well-known that the problem \eqref{eqn:mt-primalopt} admits an equivalent multi-marginal optimal transport (MMOT) reformulation \citep[Section~6]{carlier2010matching}. 
There exist numerous studies about the computation of MMOT and related problems. 
Many of these studies either only consider discrete measures, e.g., \citep*{benamou2015iterative, tupitsa2020multimarginal, ba2022accelerating, lin2022complexity, friesecke2022genetic, altschuler2023polynomial}, or approximate non-discrete measures via discretization, e.g., \citep*{guo2019computational, eckstein2019robust}.
Some studies develop regularization-based methods for approximating non-discrete MMOT and related problems, which typically involve solving infinite-dimensional optimization problems parametrized by deep neural networks; see, e.g., \citep*{eckstein2018robust, eckstein2019computation, aquino2019bounds, aquino2020minmax, henry2019martingale}. 
See also \citep*{cuturi2013sinkhorn, nutz2021entropic, eckstein2021quantitative} for the theoretical properties of entropic regularization and Sinkhorn's algorithm.
One downside of neural network based methods is the challenge posed by the non-convexity in the training objective, and there is hence no theoretical guarantee on the quality of these neural network based approximate solutions.
Recently, \citet*{alfonsi2021approximation} and \citet{neufeld2022v7numerical} developed approximation schemes for MMOT via relaxation of the marginal constraints into finitely many linear constraints. 
Our numerical approach, however, directly tackles the matching for teams problem 
\textit{without relying on the MMOT formulation}. 
Moreover, since the cost function of the MMOT problem induced by matching for teams has a minimum-of-sum structure, 
our approach and our results are applicable to MMOT problems that possess this cost structure (e.g., the problems studied by \citet{gangbo1998optimal} and \citet{heinich2002probleme});
see our detailed discussions in Section~\ref{apx:mmot} of the e-companion.

The \mbox{2-Wasserstein} barycenter problem has recently become a highly active research area due to its various applications discussed earlier. 
Most studies about the computation of Wasserstein barycenter focus on the case where $(\mu_i)_{i=1:N}$ are discrete measures with finite support, 
e.g., \citep*{borgwardt2020lp, anderes2016discrete, puccetti2020on, borgwardt2022integer, ge2019advances, heinemann2022randomized, lindheim2023simple, xie2020fast, yang2021fast}. 
Notably, \citet{altschuler2021wasserstein, altschuler2022wasserstein} have shown that there exists a polynomial-time algorithm for the exact computation of discrete \mbox{2-Wasserstein} barycenter in any fixed dimensions, 
and that the computation is NP-hard in the dimension of the underlying space.
\citet*{chizat2023doubly, luise2019sinkhorn}, and \citet{xie2020fast} developed regularization-based methods for approximating discrete \mbox{2-Wasserstein} barycenter. 
Moreover, there are also numerical methods for continuous $(\mu_i)_{i=1:N}$. 
Some of these studies are only applicable to elliptical distributions, 
e.g., \citep*{alvarez2016fixed, chewi2020gradient}.
Specifically, 
\citet{alvarez2016fixed} proposed a 
fixed-point algorithm for approximating the \mbox{2-Wasserstein} barycenter to high accuracy when $(\mu_i)_{i=1:N}$ are elliptical distributions,
and thus 
many existing studies focus on the elliptical case in their numerical experiments and use the algorithm of \citet{alvarez2016fixed}
to produce an accurate ground truth barycenter.
However, the ellipticity assumption is highly restrictive, 
as elliptical distributions are symmetric and unimodal, 
and are unable to capture the skewness and multi-modality
of many real-world distributions; e.g., distributions of residential locations.
Other studies have considered the case where the continuous $(\mu_i)_{i=1:N}$ can only be accessed through sampling, and developed stochastic optimization algorithms for approximating \mbox{2-Wasserstein} barycenter with fixed support; see, e.g., \citep*{staib2017parallel, krawtschenko2020distributed, zhang2023asynchronous}. 
Recently, algorithms for continuous \mbox{2-Wasserstein} barycenter based on neural network parametrization have been developed; see, e.g., \citep*{li2020continuous, fan2020scalable, korotin2021continuous, korotin2022wasserstein}.
They suffer from the aforementioned downside of neural network based methods due to the non-convexity in the training objective, posing challenges to subsequent theoretical analyses. 
In Section~\ref{ssec:experiment-WB}, we compare our numerical algorithm with 
five state-of-the-art \mbox{2-Wasserstein} barycenter algorithms in our numerical experiment to highlight its superior performance and its ability to produce accurate sub-optimality estimates.

\citet*{carlier2015numerical} proposed a numerical method for matching for teams with general cost functions $(c_i)_{i=1:N}$. 
After discretizing the underlying spaces $(\CX_0)_{i=0:N}$, they developed a linear programming approximation where the number of decision variables scales linearly with respect to the number of agent populations. 
In the \mbox{2-Wasserstein} barycenter case,
they developed a discretization-based method that utilizes non-smooth concave maximization. 
However, discretization of the type spaces reduces the types of agents from infinite to finite, 
which deviates from the key model assumption that there are possibly infinitely many agent types,
and discretization of the quality space restricts the traded goods to possess only the finitely many pre-specified qualities.
Hence, the solutions of the resulting discretized problems are \textit{infeasible approximations} of the matching for teams problem in general.
On the other hand, our method does not restrict the agents to finitely many types nor restrict the traded goods to finitely many qualities.

We would like to point out some existing studies about equilibrium/optimal spatial structures described by measures that are similar to the aforementioned retail business problem. 
\citet{lucas2002on} and \citet{carlier2004structure} studied the equilibrium structure of a city by analyzing the equilibrium distribution of business and residential districts while considering the positive externality of labor. 
\citet{buttazzo2005model} and \citet{carlier2005variational} considered the optimal structure of a city rather than the equilibrium structure, when taking the congestion effect into account. 
\citet*{besbes2021surge} modeled the equilibrium in the interaction between drivers and customers in a ride-hailing platform. 
Apart from these studies, the matching for teams problem is connected to 
the Cournot--Nash equilibrium problem \citep*{blanchet2016optimal, blanchet2016existence}, 
generalized barycenter problems \citep*{tanguy2024computing},
federated learning \citep*{farnia2022collaborative}, 
and adversarial machine learning \citep*{trillos2023multimarginal}.

\noindent\textbf{Notions and notations.}
We denote vectors and vector-valued functions by boldface symbols.
For $n\in\N$, we denote the all-zero vector in $\R^n$ by $\veczero_n$, and $\veczero$ is used when the dimension is unambiguous. 
We denote $\langle\BIx,\BIy\rangle:=\BIx^\TRANSP\BIy$, denote by $\|\cdot\|_p$ the $p$-norm of a vector for $p\in[1,\infty]$,
and denote $a\vee b:=\max\{a,b\}$, $a\wedge b:=\min\{a,b\}$, $(a)^+:= a\vee 0$ $\forall a,b\in\R$. 
A bounded subset of a Euclidean space is called a polytope if it is the intersection of finitely many closed half-spaces. 
For a subset $A$ of a Euclidean space, $\aff(A)$, $\conv(A)$, $\cone(A)$, $\clos(A)$, $\inter(A)$, $\relint(A)$
denote the affine hull, convex hull, conic hull, closure, interior, relative interior of $A$, respectively.
For any Polish metric space~$\CX$,
$d_{\CX}$~denotes the metric on $\CX$, 
$\CB(\CX)$ denotes the Borel subsets of $\CX$,
$\CP(\CX)$ denotes the Borel probability measures on $\CX$, 
$\delta_x$ denotes the Dirac measure at any $x\in\CX$,
$\support(\mu)$ denotes the support of any $\mu\in\CP(\CX)$,
and $\CC(\CX)$ denotes the continuous functions on $\CX$. 
For $\mu\in\CP(\CX)$, $\nu\in\CP(\CX')$ 
on two Polish metric spaces $\CX,\CX'$, 
$\Gamma(\mu,\nu):=\big\{\gamma\in\CP(\CX\times\CX'):$ the marginals of $\gamma$ on $\CX$ and $\CX'$ are $\mu$ and $\nu\big\}$
is the set of couplings of $\mu$ and $\nu$.
Lastly, for $\mu,\nu\in\CP(\CX)$, let $W_1(\mu,\nu)$ denote their Wasserstein distance of order~1,
where $W_1(\mu,\nu):=\inf_{\gamma\in\Gamma(\mu,\nu)}\big\{\int_{\CX\times\CX}d_{\CX}(x,x')\DIFFM{\gamma}{\DIFF x,\DIFF x'}\big\}$.

\section{Model assumptions and approximate matching equilibrium (AME)}\label{sec:model}
To motivate the matching for teams model and our numerical method for computing AMEs,
let us discuss in this section the practical relevance of the model assumptions,
as well as the properties and economic interpretations of 
AMEs.

\subsection{Discussions about the model assumptions}\label{ssec:model-assumptions}

\textbf{Indivisibility of traded goods.} 
A distinct feature of the matching for teams model is that 
each team consists of exactly one agent from each population, 
and thus the agents'
transfers of utilities after forming a team describes the trading of an indivisible good, e.g., house, artwork, contractual obligation, etc.
Specifically, in the retail business problem the traded goods are employment contracts.
This means that the matching for teams model is not suitable for describing the trading of divisible goods such as crude oil and electricity.
See also the discussions of \citet{henry1970indivisibilites, shapley1972assignment} about the economic implications and the plausibility of the indivisibility assumption.
The mathematical implication of this assumption is that 
the agent populations 
(or at least the sub-populations that choose to participate in trading) 
are required to have the same size, 
which is why each population is normalized to~1 and modeled by a probability measure. 
In fact, with a slight extension to allow the cost functions $(c_i)_{i=1:N}$ to take $\infty$ as value,
this model can also describe matching games with unequally-sized agent populations where agents have the option to not participate.
This can be achieved through augmenting each space with a ``dummy point'' representing non-participation;
see \citep[Section~2.1]{chiappori2010hedonic}
for details.

\textbf{The type spaces.} 
We have adopted the assumptions of \citep{carlier2010matching} that 
$(\CX_i)_{i=1:N}$ are compact metric spaces.
In our concrete Euclidean setting in Setting~\ref{sett:Euclidean}, 
$(\CX_i)_{i=1:N}$ are assumed to be compact subsets of possibly multi-dimensional Euclidean spaces formed by unions of simplices.
Since each $\CX_i$ in Setting~\ref{sett:Euclidean} possesses the flexibility to approximate sophisticated shapes,
our setting is significantly more general than the one-dimensional settings considered by many theoretical studies, as mentioned in the introduction.
In the retail business problem, for example,
we set the type space $\CX_i\subset\R^2$ for $i=1,\ldots,N-1$ 
to be the union of triangles in a triangular mesh which can represent  disconnected residential districts with complex shapes.

\textbf{The quality space.}
We adopt the same assumptions for $\CX_0$ and $(\CX_i)_{i=1:N}$.
As an example, one can let $\CX_0$ 
be a finite set in a problem instance where the traded goods have finitely many 
distinct types of qualities, e.g., the brand and model of a car.
Alternatively, 
$\CX_0$ can be multi-dimensional and have sophisticated shapes
characterized by the union of simplices,
e.g., when the qualities of traded goods correspond to the geographic coordinates
or some continuous variables such as
the durability of the material and the time of delivery.

\textbf{The cost functions.}
As a standing assumption,
we require each cost function $c_i:\CX_i\times\CX_0\to\R$ to be \mbox{$L_{c}$-Lipschitz} continuous,
i.e., 
$\big|c_i(x,z)-c_i(x',z')\big|\le L_{c}\big(d_{\CX_i}(x,x')+d_{\CX_0}(z,z')\big)$ 
$\forall (x,z),(x',z')\in\CX_i\times\CX_0$.
The purpose of this assumption is to control the approximation error of our numerical method.
Note that 
since $(\CX_i)_{i=0:N}$ are allowed to be disconnected under
our concrete Euclidean setting (Setting~\ref{sett:Euclidean}), 
it is possible to define $c_i$ to be Lipschitz continuous on each of the disjoint ``pieces'' of the domain.
Thus,
the Lipschitz continuity of $c_i$ is not a restrictive assumption,
and such $c_i$ can (at least approximately) represent real costs in matching for teams problems such as,
for example, spatial distance (as in the retail business problem),
time duration,
or dissimilarity between the preferred quality and the actual quality. 
Moreover, $c_i$ can be a weighted sum of costs arising from various sources.

\subsection{The notion of approximate matching equilibrium}\label{ssec:model-AME}

The crucial property of matching equilibria that 
motivates our notion of $\epsilon$-matching equilibrium as well as our numerical approach is the following characterization of matching equilibria by \citep{carlier2010matching} 
via the solutions of
a pair of convex optimization problems dual to each other:
\begin{align}
    \begin{split}
        \minimize_{(\gamma_i),\,\nu}\quad & \sum_{i=1}^N\int_{\CX_i\times\CX_0} c_i \DIFFX{\gamma_i} \\
        \mathrm{subject~to} \quad & \nu\in\CP(\CX_0),\qquad \gamma_i\in \Gamma(\mu_i,\nu)\quad \forall i\in\{1,\ldots,N\},
    \end{split}
    \tag{$\MT$}
    \label{eqn:mt-primalopt}\\
    \begin{split}
        \maximize_{(\varphi_i)} \quad & \sum_{i=1}^N \int_{\CX_i}\varphi_i^{c_i} \DIFFX{\mu_i} \\
        \mathrm{subject~to}\quad & \sum_{i=1}^N \varphi_i=0,\qquad \varphi_i\in\CC(\CX_0) \quad \forall i\in\{1,\ldots,N\}.
    \end{split}
    \tag{$\MT^*$}
    \label{eqn:mt-dualopt}
\end{align}
\citet[Theorem~4]{carlier2010matching} showed that 
$(\varphi_i,\gamma_i)_{i=1:N},\nu$ is a matching equilibrium
\textit{if and only if} 
$(\gamma_i)_{i=1:N},\nu$ 
is an optimizer of \eqref{eqn:mt-primalopt}
and 
$(\varphi_i)_{i=1:N}$ 
is an optimizer of \eqref{eqn:mt-dualopt}.
In particular, strong duality holds,
i.e., \eqref{eqn:mt-primalopt} and \eqref{eqn:mt-dualopt} have identical optimal values,
which we denote by $\alpha_{\MT}^\star$ in the rest of the paper.
Since \citep[Proposition~2]{carlier2010matching} guarantees the existence of matching equilibrium,
we have $\alpha_{\MT}^\star\in\R$.

Observe that 
for any $\epsilon>0$,
our definition in the introduction says that
$(\varphi_i,\gamma_i)_{i=1:N},\nu$ constitute an $\epsilon$-matching equilibrium
if and only if
$(\gamma_i)_{i=1:N},\nu$ is \textit{feasible} for \eqref{eqn:mt-primalopt},
$(\varphi_i)_{i=1:N}$ is \textit{feasible} for \eqref{eqn:mt-dualopt},
and their \textit{optimality gap}
$\big(\sum_{i=1}^N\int_{\CX_i\times\CX_0} c_i \DIFFX{\gamma_i}\big) 
- \big(\sum_{i=1}^N\int_{\CX_i}\varphi_i^{c_i} \DIFFX{\mu_i}\big)$
is at most~$\epsilon$.
Moreover,
observe that every $\epsilon$-matching equilibrium
$(\varphi_i,\gamma_i)_{i=1:N},\nu$
also constitutes an $\epsilon$-maximizer of 
$-\big(\sum_{i=1}^N\int_{\CX_i\times\CX_0} c_i \DIFFX{\gamma_i}\big) = \sum_{i=1}^N\int_{\CX_i\times\CX_0} \big(\varphi_i(z) - c_i(x,z)\big) \DIFFM{\gamma_i}{\DIFF x,\DIFF z}$,
which can be interpreted as 
the social welfare
because it corresponds to the aggregated utilities of all agents.
Consequently,
$-\alpha_{\MT}^\star$ can be interpreted as the optimal social welfare, and
$-\big(\sum_{i=1}^N\int_{\CX_i\times\CX_0} c_i \DIFFX{\gamma_i}\big)$
and 
$-\big(\sum_{i=1}^N\int_{\CX_i}\varphi_i^{c_i} \DIFFX{\mu_i}\big)$
are lower and upper bounds for
$-\alpha_{\MT}^\star$ 
that are at most~$\epsilon$ apart.
Here, the feasibility of $(\gamma_i)_{i=1:N},\nu$
and $(\varphi_i)_{i=1:N}$
for \eqref{eqn:mt-primalopt} and \eqref{eqn:mt-dualopt}
is critical for guaranteeing the validity of these bounds,
which is why 
our notion of $\epsilon$-matching equilibrium constitutes feasible approximations of 
\eqref{eqn:mt-primalopt} and \eqref{eqn:mt-dualopt}.

Most importantly, 
the following theorem shows that
a sequence of AMEs with errors approaching~0
admits subsequential limits which are true matching equilibria,
up to a minor modification.

\begin{theorem}[Convergence of AMEs]\label{thm:mt-tf-equilibrium}
\!Let $(\epsilon^{(l)})_{l\in\N}\hspace{-1pt}\subset\hspace{-1pt} (0,\infty)$ satisfy $\lim_{l\to\infty}\epsilon^{(l)}=0$.
For each $l\in\N$,
let $({\varphi}^{(l)}_i,{\gamma}_i^{(l)})_{i=1:N},\allowbreak{\nu}^{(l)}$ be an $\epsilon^{(l)}$-matching equilibrium, 
and let $\tilde{\varphi}_i^{(l)}(z):=\inf_{x\in\CX_i}\big\{c_i(x,z)-{\varphi_i^{(l)}}^{c_i}(x)\big\}$
$\forall z\in\CX_0$, 
$\forall i\in\{1,\ldots,N-1\}$,
$\tilde{\varphi}_N^{(l)}:= -\sum_{i=1}^{N-1}\tilde{\varphi}_i^{(l)}$.
Then,
the following statements hold.
\begin{enumerate}[label=(\roman*),leftmargin=20pt]
\item\label{thms:mt-tf-equilibrium-invariance}
There exist $(\kappa_{i}^{(l)})_{i=1:N,\,l\in\N}\subset\nobreak\R$ 
with ${\sum_{i=1}^N\kappa_{i}^{(l)}=0}$ $\forall l\in\N$ 
such that $\big(\tilde{\varphi}_i^{(l)}+\kappa_{i}^{(l)}\big)_{i=1:N,\,l\in\N}$ 
are uniformly bounded and 
$\big(\tilde{\varphi}^{(l)}_i+\kappa^{(l)}_i,{\gamma}_i^{(l)})_{i=1:N},\allowbreak{\nu}^{(l)}$ is an $\epsilon^{(l)}$-matching equilibrium for all $l\in\nobreak\N$.
In particular, one may choose
$\kappa^{(l)}_i:=-\min_{z\in\CX_0}\big\{\tilde{\varphi}_i^{(l)}(z)\big\}$
$\forall i\in\{1,\ldots, N-1\}$,
$\kappa^{(l)}_{N}:=-\sum_{i=1}^{N-1}\kappa^{(l)}_i$.

\item\label{thms:mt-tf-equilibrium-compactness}
Every subsequence of 
$\big(\big(\tilde{\varphi}^{(l)}_i+\kappa^{(l)}_i,{\gamma}_i^{(l)})_{i=1:N},\allowbreak{\nu}^{(l)}\big)_{l\in\N}$
in \ref{thms:mt-tf-equilibrium-invariance}
admits a further subsequence
$\big(\big(\tilde{\varphi}^{(l_t)}_i+\nobreak\kappa^{(l_t)}_i,{\gamma}_i^{(l_t)})_{i=1:N},\allowbreak{\nu}^{(l_t)}\big)_{t\in\N}$
where
$({\nu}^{(l_t)})_{t\in\N}$ converges in $(\CP(\CX_0),W_1)$ to some~${\nu}^{(\infty)}$, 
and for $i=1,\ldots,N$, 
$\big(\tilde{\varphi}_i^{(l_t)}+\kappa_i^{(l_t)}\big)_{t\in\N}$ converges uniformly to some $\tilde{\varphi}_i^{(\infty)}\in\CC(\CX_0)$, 
$({\gamma}_i^{(l_t)})_{t\in\N}$ converges in $(\CP({\CX_i\times\CX_0}),W_1)$ to some ${\gamma}^{(\infty)}_i$,
and 
$(\tilde{\varphi}^{(\infty)}_i,{\gamma}_i^{(\infty)})_{i=1:N},{\nu}^{(\infty)}$ is a matching equilibrium.
\end{enumerate}
\end{theorem}

The modification of $(\varphi_i^{(l)})_{i=1:N-1}$ in 
Theorem~\ref{thm:mt-tf-equilibrium} is known as $c_i$-convexification \citep[Proposition~5.8]{villani2008optimal} and its purpose is to guarantee the $L_c$-Lipschitz continuity of
$(\tilde{\varphi}_i^{(l)})_{i=1:N-1}$; 
see also the proof of \citep[Proposition~1]{carlier2010matching}.
Moreover, the addition of the constants $(\kappa_i^{(l)})_{i=1:N}$
stems from the invariance of \eqref{eqn:mt-dualopt} when shifting $(\varphi_i)_{i=1:N}$
by constants that sum to~0.
Note that Theorem~\ref{thm:mt-tf-equilibrium}
results from a compactness argument, and it is possible that there are multiple subsequential limits since we do not assume 
the uniqueness of matching equilibrium.
We remark that under specific settings
such as the \mbox{2-Wasserstein} barycenter problem with appropriate curvature assumptions (see, e.g., \citep[Section~3.3.1]{chewi2020gradient}),
it is possible to characterize the convergence rate to the true matching equilibrium.

Besides the property above,
we show that an $\epsilon$-matching equilibrium
provides lower and upper bounds for the improvement of social welfare when the cost functions $(c_i)_{i=1:N}$ are perturbed.

\begin{proposition}\label{prop:perturbation}
    Let $(\varphi_i,\gamma_i)_{i=1:N},\nu$ be
    an \mbox{$\epsilon$-matching} equilibrium for $\epsilon>0$,
    let $\tilde{c}_i:\CX_i\times\CX_0\to\nobreak\R$ be continuous 
    for $i=1,\ldots,N$,
    and let $\tilde{\alpha}_{\MT}^\star$ denote the optimal value 
    of \eqref{eqn:mt-primalopt} with the cost functions 
    $(c_i)_{i=1:N}$ replaced by $(\tilde{c}_i)_{i=1:N}$.
    Then, it holds that 
    $\big(\sum_{i=1}^{N}\int_{\CX_i\times\CX_0}(c_i-\nobreak\tilde{c}_i)\DIFFX{\gamma_i}\big)-\epsilon
    \le \alpha_{\MT}^\star-\tilde{\alpha}_{\MT}^\star
    \le \big(\sum_{i=1}^{N}\int_{\CX_i}(\varphi_i^{c_i}-\varphi_i^{\tilde{c}_i})\DIFFX{\mu_i}\big)+\epsilon$.
\end{proposition}

Consequently, 
policymakers can use an \mbox{$\epsilon$-matching} equilibrium $(\varphi_i,\gamma_i)_{i=1:N},\nu$ to compute these bounds in order to guide their decisions,
for example, in infrastructure planning to find an economical way 
of improving social welfare through targeted reduction of the costs 
$(c_i)_{i=1:N}$.

\section{Parametric approximation of matching equilibrium}\label{sec:parametric}
In this section, we present how we can obtain approximate matching equilibria 
via parametric approximations of 
\eqref{eqn:mt-primalopt} and \eqref{eqn:mt-dualopt}.
Our parametric formulation is motivated by the observation that 
\eqref{eqn:mt-primalopt} and \eqref{eqn:mt-dualopt} are infinite dimensional 
linear optimization problems.
In particular, if we introduce 
continuous functions $(\psi_i)_{i=1:N}$ as
auxiliary decision variables in place of $(\varphi_i^{c_i})_{i=1:N}$ in \eqref{eqn:mt-dualopt}
and impose the constraints 
$\psi_i(x)+\varphi_i(z)\le c_i(x,z)$ $\forall(x,z)\in\CX_i\times\CX_0$, $\forall i\in\{1,\ldots,N\}$,
we get a maximization problem over 
$\big((\psi_i,\varphi_i)\in\CC(\CX_i)\times\CC(\CX_0)\big)_{i=1:N}$
with a linear objective and infinitely many linear constraints.
This motivates us to parametrize each $\CC(\CX_i)$ by the linear combinations of a set $\CG_i=(g_{i,j})_{j=1:m_i}\subset \CC(\CX_i)$ of ${m_i\in\N}$ continuous functions for $i=0,1,\ldots,N$. 
We call $(\CG_i)_{i=0:N}$ test functions whose precise choices will be specified later under a concrete setting in Section~\ref{ssec:setting-Euclidean}. 
For notational simplicity, let us denote $\BIg_i(x):=\nobreak(g_{i,1}(x),\ldots,g_{i,m_i}(x))^\TRANSP\in\R^{m_i}$ $\forall x\in\CX_i$, 
for $i=0,1,\ldots,N$, 
and denote 
$\bar{\BIg}_i:=\nobreak\big(\int_{\CX_i}g_{i,1}\DIFFX{\mu_i},\ldots,\allowbreak\int_{\CX_i}g_{i,m_i}\DIFFX{\mu_i}\big)^\TRANSP\in\R^{m_i}$ 
for $i=1,\ldots,N$.
Moreover, we introduce the following \textit{moment-based relation}.
\begin{definition}[Moment-based relation $\abovebelowset{\CG}{\varsigma}{\sim}$]\label{def:momentset}%
Let $\CX$ be a compact metric space. 
For ${\varsigma\ge0}$, a finite set $\CG\subset\CC(\CX)$, and for all $\mu,\nu\in\CP(\CX)$, 
we denote 
$\mu\abovebelowset{\CG}{\varsigma}{\sim}\nu$ 
$\Leftrightarrow$
$\sum_{g\in\CG}\big|\int_{\CX}g\DIFFX{\mu}-\int_{\CX}g\DIFFX{\nu}\big|\le \varsigma$.
When $\varsigma=\nobreak0$,
we omit $\varsigma$ in the notation and write $\overset{\CG}{\sim}$ 
for notational simplicity.
\end{definition}

With these notions, let us present the parametric formulation 
\eqref{eqn:mt-tf-dual} and \eqref{eqn:mt-tf-lsip} which are 
linear semi-infinite programming (LSIP) problem dual to each other:
\begin{align}
\begin{split}
    \minimize_{(\theta_i)} \quad & \sum_{i=1}^{N} \int_{\CX_i\times\CX_0}c_i\DIFFX{\theta_i}\\
    \mathrm{subject~to}\quad & \theta_i\in\Gamma(\bar{\mu}_i,\bar{\nu}_i),\quad \bar{\mu}_i\overset{\CG_i}{\sim}\mu_i,\quad \bar{\nu}_i\overset{\CG_0}{\sim}\bar{\nu}_1 \quad\forall i\in\{1,\ldots,N\},
\end{split}
\tag{$\MT_{\sspar}$}\label{eqn:mt-tf-dual}\\
\begin{split}
    \maximize_{(y_{i,0},\BIy_i,\BIw_i)}\quad & \sum_{i=1}^N y_{i,0}+\langle\bar{\BIg}_i,\BIy_i\rangle\\
    \mathrm{subject~to}\quad & y_{i,0}+\langle\BIg_i(x),\BIy_i\rangle+\langle\BIg_0(z),\BIw_i\rangle  \le c_i(x,z) \quad \forall (x,z)\in\CX_i\times\CX_0,\; \forall i\in\{1,\ldots,N\},\\
    & \sum_{i=1}^N\BIw_i=\veczero_{m_0}, \qquad y_{i,0}\in\R,\; \BIy_i\in\R^{m_i},\; \BIw_i\in\R^{m_0}  \quad \forall i\in\{1,\ldots,N\}.
\end{split}
\tag{$\MT^*_{\sspar}$}\label{eqn:mt-tf-lsip}
\end{align}

The first important property of \eqref{eqn:mt-tf-dual} and \eqref{eqn:mt-tf-lsip} is that strong duality holds.
\begin{theorem}[Strong duality]\label{thm:mt-tf-duality}
The optimal values of 
\eqref{eqn:mt-tf-dual} and \eqref{eqn:mt-tf-lsip} 
are identical.
\end{theorem}
Thus, let us denote the common optimal value of \eqref{eqn:mt-tf-dual} and \eqref{eqn:mt-tf-lsip}
by $\alpha_{\sspar}^{\star}$ in the rest of the paper.
On the one hand, 
since $\mu_i\overset{\CG_i}{\sim}\mu_i$ for $i=1,\ldots,N$
and $\nu\overset{\CG_0}{\sim}\nu$ $\forall\nu\in\CP(\CX_0)$,
\eqref{eqn:mt-tf-dual} relaxes the marginal constraints in \eqref{eqn:mt-primalopt}
via the moment-based relation, and we thus get the lower bound $\alpha_{\sspar}^{\star}\le \alpha^{\star}_{\MT}$.
On the other hand, observe that
every feasible solution $(y_{i,0},\BIy_i,\BIw_i)_{i=1:N}$ of \eqref{eqn:mt-tf-lsip}
directly yields transfer functions $\varphi_i(\,\cdot\,):=\langle\BIg_0(\,\cdot\,),\BIw_i\rangle$ for $i=1,\ldots,N$ that are feasible for \eqref{eqn:mt-dualopt}.

For any $\epsilon_{\sspar}>0$, $\varsigma\ge 0$, 
we call 
$(\theta_i)_{i=1:N}$,
$(y_{i,0},\BIy_i,\BIw_i)_{i=1:N}$
a pair of 
$(\epsilon_{\sspar},\varsigma)$-optimizers of 
\eqref{eqn:mt-tf-dual} and \eqref{eqn:mt-tf-lsip}
if 
$(y_{i,0},\BIy_i,\BIw_i)_{i=1:N}$ is feasible for \eqref{eqn:mt-tf-lsip},
$\theta_i\in\Gamma(\bar{\mu}_i,\bar{\nu}_i)$ with 
$\bar{\mu}_i\abovebelowset{\CG_i}{\varsigma}{\sim}\mu_i$, 
$\bar{\nu}_i\abovebelowset{\CG_0}{\varsigma}{\sim}\bar{\nu}_1$
for $i=1,\ldots,N$, 
and $\!\big(\!\sum_{i=1}^N\int_{\CX_i\times\CX_0}c_i\DIFFX{\theta_i}\hspace{-1pt}\big)\hspace{-1pt} - \hspace{-1pt}\big(\!\sum_{i=1}^N y_{i,0}+\langle\bar{\BIg}_i,\BIy_i\rangle\hspace{-1pt}\big)\! \le \epsilon_{\sspar}$.
Here, we allow the moment-based constraints 
$\bar{\mu}_i\overset{\CG_i}{\sim}\mu_i$, $\bar{\nu}_i\overset{\CG_0}{\sim}\bar{\nu}_1$ $\forall i\in\{1,\ldots,N\}$ in \eqref{eqn:mt-tf-dual} to be violated by
a small error~$\varsigma$
to permit numerical imprecision during computation.
We will present in Section~\ref{sec:numerics} methods to compute a pair of 
\mbox{$(\epsilon_{\sspar},\varsigma)$-optimizers} of 
\eqref{eqn:mt-tf-dual} and \eqref{eqn:mt-tf-lsip}
under a concrete Euclidean setting.
With a pair of $(\epsilon_{\sspar},\varsigma)$-optimizers,
the following theorem demonstrates the construction of two approximate matching equilibria,
where the approximation error depends on 
$\epsilon_{\sspar}$, $\varsigma$, and the choices of $(\CG_i)_{i=0:N}$.

\begin{theorem}[Construction of AMEs]\label{thm:mt-tf-approx}%
Let 
$(\hat{\theta}_i)_{i=1:N}$,
$(\hat{y}_{i,0},\hat{\BIy}_i,\hat{\BIw}_i)_{i=1:N}$
be a pair of 
$(\epsilon_{\sspar},\varsigma)$-optimizers of 
\eqref{eqn:mt-tf-dual} and \eqref{eqn:mt-tf-lsip}
for $\epsilon_{\sspar}>0$, $\varsigma\ge 0$.
For $i=1,\ldots,N$, 
let $\hat{\mu}_i$ and $\hat{\nu}_i$ 
denote the marginals of $\hat{\theta}_i$ on $\CX_i$ and $\CX_0$, respectively.
Moreover, 
let $\overline{\epsilon}_i(\varsigma)\ge \sup\big\{W_1(\nu,\nu'):\nu,\nu'\in\nobreak\CP(\CX_i),\;\nu\abovebelowset{\CG_i}{\varsigma}{\sim}\nu'\big\}$
for $i=0,1,\ldots,N$,
and define 
$\epsilon_{\ssapx}:=
\epsilon_{\sspar}
+L_c\big(N\overline{\epsilon}_0(\varsigma)+\sum_{i=1}^N\overline{\epsilon}_i(\varsigma)\big)$.
We construct $(\hat{\varphi}_i,\hat{\gamma}_i,\tilde{\gamma}_i)_{i=1:N},\hat{\nu},\tilde{\nu}$ via a probability space $(\Omega,\CF,\PROB)$ as follows.\footnote{See Section~\ref{apx:couplings} of the e-companion for the detailed construction of the random variables 
$Z$, $(Z_i,X_i,\bar{X}_i)_{i=1:N}$, and $\bar{Z}$ on a sufficiently rich probability space
via repeated applications of \citep[Lemma~2.11]{dudley1983invariance}.}
\begin{itemize}[leftmargin=12pt]

    \item Define $\hat{\nu}:=\hat{\nu}_{1}$, and define $Z:\Omega\to\CX_0$ to be a random variable with law~$\hat{\nu}$ (i.e., $\PROB\circ Z^{-1}=\hat{\nu}$).

    \item For $i=1,\ldots,N$, 
    define
    $\hat{\varphi}_i(z):=\langle\BIg_0(z),\hat{\BIw}_i\rangle$ $\forall z\in\CX_0$,
    define random variables $Z_i:\Omega\to\CX_0$,
    $X_i:\Omega\to\CX_i$,
    $\bar{X}_i:\Omega\to\CX_i$ such that
    $\EXP\big[d_{\CX_0}(Z,Z_i)\big]\le \overline{\epsilon}_0(\varsigma)$,
    $\EXP\big[d_{\CX_i}(X_i,\bar{X}_i)\big]\le \overline{\epsilon}_i(\varsigma)$,
    the law of $\bar{X}_i$ is $\mu_i$,
    and 
    the law of $(X_i,Z_i):\Omega\to\CX_i\times\CX_0$ is $\hat{\theta}_i$.
    Let $\hat{\gamma}_i\in\CP(\CX_i\times\CX_0)$ be the law of $(\bar{X}_i,Z):\Omega\to\CX_i\times\CX_0$.
    
    \item Let the random variable $\bar{Z}:\Omega\to\CX_0$ 
    satisfy $\bar{Z}\in\argmin_{z\in\CX_0}\big\{\sum_{i=1}^{N}c_i(\bar{X}_i,z)\big\}$
    $\PROB$-almost surely.
    Let $\tilde{\nu}$ be the law of $\bar{Z}$,
    and let $\tilde{\gamma}_i\in\CP(\CX_i\times\CX_0)$ be the law of $(\bar{X}_i,\bar{Z}):\Omega\to\CX_i\times\CX_0$ for $i=1,\ldots,N$.
\end{itemize}
Then, $(\hat{\varphi}_i,\hat{\gamma}_i)_{i=1:N},\hat{\nu}$ and $(\hat{\varphi}_i,\tilde{\gamma}_i)_{i=1:N},\tilde{\nu}$ are both $\epsilon_{\ssapx}$-matching equilibria.
\end{theorem}

Note that the AMEs constructed by Theorem~\ref{thm:mt-tf-approx} provide provable and computable lower bounds  
$-\big(\sum_{i=1}^N\int_{\CX_i\times\CX_0} c_i \DIFFX{\hat{\gamma}_i}\big)$,
$-\big(\sum_{i=1}^N\int_{\CX_i\times\CX_0} c_i \DIFFX{\tilde{\gamma}_i}\big)$
and upper bound
$-\big(\sum_{i=1}^N\int_{\CX_i}\hat{\varphi}_i^{c_i} \DIFFX{\mu_i}\big)$
for the optimal social welfare
$-\alpha_{\MT}^\star$, which is a crucial advantage of our parametric formulation over methods based on discretization.
Moreover, the linearity of our parametrization possesses advantages over non-linear parametrizations such as neural networks 
since it allows us to establish the explicit computational complexity bounds in Theorem~\ref{thm:complexity-explicit},
as well as to develop a provably convergent and efficient algorithm utilizing state-of-the-art linear programming (LP) solvers in Section~\ref{ssec:efficient-algo}.

\section{Numerical method}\label{sec:numerics}
\vspace{6pt}

\subsection{A concrete Euclidean setting with simplicial complexes}\label{ssec:setting-Euclidean}

In this section, we work under a concrete and flexible Euclidean setting 
where $(\CX_i)_{i=0:N}$ are subsets of Euclidean spaces
and the test functions in $(\CG_i)_{i=0:N}$ are continuous piece-wise affine (CPWA).
To motivate our setting, recall that $\CG_i$ 
in \eqref{eqn:mt-tf-lsip} 
serves as the basis functions in a finite-dimensional parametrization of $\CC(\CX_i)$.
Therefore, 
in order for \eqref{eqn:mt-tf-lsip} to be an accurate approximation of \eqref{eqn:mt-dualopt},
it is desirable that the linear combinations of $\CG_i$ are flexible to approximate any $\psi_i\in\CC(\CX_i)$.

To illustrate our particular choice of $\CG_i$, let us first consider the two-dimensional case, i.e., when $\CX_i\subset\R^2$.
We triangulate $\CX_i$ into a fine triangular mesh,
that is, we cover $\CX_i$ by a collection of small triangles that are disjoint up to their shared edges and vertices,
where the mesh size, i.e., its granularity, corresponds to the length of the longest edge among these triangles.
Subsequently, for any $\psi_i\in\CC(\CX_i)$,
there exists a unique continuous function $\hat{\psi}_i$ 
which is affine on each triangle
and satisfies $\hat{\psi}_i(\BIv)=\psi_i(\BIv)$ for every vertex $\BIv$;
that is, $\hat{\psi}_i$ is a CPWA interpolation of $\psi_i$ 
with respect to this triangular mesh.
To choose $\CG_i$ such that its linear combinations can represent any CPWA function with respect to the triangular mesh, 
it suffices to set each $g\in\CG_i$ to be CPWA with $g(\BIv)=1$ for exactly one vertex $\BIv$
and $g(\BIv')=\nobreak0$ for any other vertex $\BIv'$.
It intuitively follows that as the triangular mesh gets finer, i.e., as its mesh size gets closer to~0,
the above CPWA interpolation of any $\psi_i\in\CC(\CX_i)$ becomes a better approximation.
Indeed, as we will show later in Theorem~\ref{thm:error-control},
this choice of $\CG_i$ provides an explicit expression for  $\overline{\epsilon}_i(\,\cdot\,)$ in Theorem~\ref{thm:mt-tf-approx}
that allows us to control the approximation error $\epsilon_{\ssapx}$ to be arbitrarily close to~0.
Moreover, we will present an explicit computational complexity bound in Theorem~\ref{thm:complexity-explicit}
and a practically efficient algorithm in Section~\ref{ssec:efficient-algo}
under this setting.

Setting~\ref{sett:Euclidean} below provides the rigorous formulation of the informal description of $\CG_i$ above, and extends it to any number of dimensions.
It requires some heavy notions related to the studies of simplicial complexes, see, e.g.,
\citep{munkres1984elements},
which we will introduce here.
For any $d\in\N$ and any $n\in\nobreak\{0,1,\ldots,d\}$,
an \textit{$n$-simplex} in $\R^d$ is the convex hull of $n+\nobreak1$ affinely independent points in $\R^d$.
For any simplex $K$,
we let $\dim(K)$ denote its \textit{dimension},
i.e., $\dim(K)=n$ if $K$ is an $n$-simplex,
let $V(K)\subset\R^d$ denote its finite set of \textit{extreme points},
and let $\FF(K)$ denote its finite set of 
\textit{faces};
note that $\emptyset\in\FF(K)$
and we let $\FF(\emptyset)=\{\emptyset\}$.
A finite set $\FK$  is called a \textit{simplicial complex} in $\R^d$
if every non-empty $K\in\nobreak\FK$ is a simplex in $\R^d$,
and $\FK$ satisfies
$\FF(K)\subseteq\FK$ $\forall K\in\FK$
as well as
$K\cap K'\in\FF(K)\cap\nobreak\FF(K')$ $\forall K,K'\in\nobreak\FK$.
Moreover, 
we let $V(\FK):=\bigcup_{K\in\FK}V(K)$ denote the finite set of all extreme points in $\FK$, 
let $\FM(\FK):=\big\{K\in\nobreak\FK:K\notin\nobreak \FF(K')$ $\forall K'\in\FK\setminus\nobreak\{K\}\big\}$ 
denote the so-called \textit{maximal simplices} in~$\FK$ which are not faces of another simplex in $\FK$,
and
let $\eta(\FK):=\max_{K\in\FK}\big\{\max_{\BIv,\BIv'\in V(K)}\big\{\|\BIv-\BIv'\|_2\big\}\big\}$ denote
the \textit{mesh size} of $\FK$.
One checks that every $\BIx\in\bigcup_{K\in\FK}K$ can be uniquely expressed as 
$\BIx=\nobreak\sum_{\BIv\in V(K)}\lambda_{\BIv}\BIv$
where $K\in\FK$, 
$(\lambda_\BIv)_{\BIv\in V(K)}\subset(0,1]$ with $\sum_{\BIv\in V(K)}\lambda_{\BIv}=\nobreak1$.
Subsequently, the \textit{barycentric coordinate} of~$\BIx$ with respect to each $\BIv\in\nobreak V(\FK)$ is denoted by $\lambda_{\FK,\BIv}(\BIx)$
and defined by
$\lambda_{\FK,\BIv}(\BIx):=\lambda_{\BIv}$ 
for each $\BIv\in V(K)$,
$\lambda_{\FK,\BIv}(\BIx):=\nobreak 0$ 
for each $\BIv\in V(\FK)\setminus V(K)$.\footnote{See, e.g., \citep[p.10]{munkres1984elements} for the definition of barycentric coordinates. 
Moreover, it follows from \citep[Proposition~3.7(iii)]{neufeld2022v7numerical} that $(\lambda_{\FK,\BIv_i})_{i=0:n}$ are continuous.}
A~simplicial complex $\FS$ is called a \textit{subdivision} of another simplicial complex $\FK$ if each $S\in\FS$ is a subset of some $K\in\FK$
and that each $K\in\FK$ can be expressed as a finite union of some simplices in~$\FS$.

With these notions and notations, let us now introduce our concrete Euclidean setting.

\begin{setting}\label{sett:Euclidean}
    For $i=0,1,\ldots,N$, $\CX_i:=\bigcup_{K\in\FK_i}K\subset\R^{d_i}$
    is equipped with the norm $\|\cdot\|_2$, where
    $d_i\in\N$ and $\FK_i$ is a simplicial complex in $\R^{d_i}$.
    Let $\FS_i$ be a subdivision of $\FK_i$, 
    let $m_i:=\nobreak|V(\FS_i)|-\nobreak1$, 
    and let $(\BIv_{i,j})_{j=0:m_i}$ be an enumeration of $V(\FS_i)$.
    Define $g_{i,j}:=\lambda_{\FS_i,\BIv_{i,j}}$ for $j=\nobreak0,1,\ldots,m_i$,
    define $\CG_i:=\nobreak(g_{i,j})_{j=1:m_i}$,
    $\BIg_i(x):=(g_{i,1}(x),\ldots,g_{i,m_i}(x))^\TRANSP\in\R^{m_i}$ $\forall x\in\CX_i$,
    and
    define $\overline{\epsilon}_i(\varsigma):=2\eta(\FS_i)+\varsigma\max_{\BIx,\BIx'\in\CX_i}\big\{\|\BIx-\BIx'\|_2\big\}$
    $\forall \varsigma\ge 0$.
    Moreover, for $i=1,\ldots,N$,
    $c_i:\CX_i\times\CX_0\to\R$ is $L_{c}$-Lipschitz continuous for $L_{c}>0$, and
    $\mu_i\in\CP(\CX_i)$ satisfies the property:
    for every $K\in\FM(\FK_i)$, 
    there exists an affine function $p_{i,K}:K\to\R_+$
    not equal to the zero function, such that 
    $\mu_i(E\cap K)\ge \int_{E}p_{i,K}\DIFFX{\Lebesgue_K}$
    $\forall E\in\CB(\CX_i)$,
    where 
    $\Lebesgue_K$ denotes the $\dim(K)$-dimensional Lebesgue measure on $K$
    (when $\dim(K)=\nobreak0$ and $K=\{\BIv\}$,
    we abuse the notation and denote $\Lebesgue_K:=\delta_{\BIv}$).
    As a particular case, consider the following class of $(\mu_i)_{i=1:N}$:\widowpenalties-8000
    \begin{enumerate}[label=(P),beginpenalty=10000]
        \item\label{settc:Euclidean-particular}%
        for $i=1,\ldots,N$ and any $K\in\FK_i$, 
        there is an affine function
        $p_{i,K}:K\to\R_+$ satisfying
        $\mu_i\big(E \,\cap\nobreak\, \relint(K)\big)=\int_{E}p_{i,K}\DIFFX{\Lebesgue_K}$
        $\forall E\in\CB(\CX_i)$, where $p_{i,K}\ne 0$ $\forall K\in\FM(\FK_i)$.
    \end{enumerate}
    Lastly, let
    $\bar{\BIg}_i:=\big(\int_{\CX_i}g_{i,1}\DIFFX{\mu_i},\ldots,\allowbreak\int_{\CX_i}g_{i,m_i}\DIFFX{\mu_i}\big)^\TRANSP\in\R^{m_i}$
    for $i=1,\ldots,N$.
\end{setting}

Note that Setting~\ref{sett:Euclidean}
does not restrict the cost functions $(c_i)_{i=1:N}$ apart from their Lipschitz continuity,
and only requires each type measure $\mu_i$ to dominate a measure on $\CX_i$ with piece-wise affine density.
In particular, an admissible choice of $\mu_i$ under
Setting~\ref{sett:Euclidean}\ref{settc:Euclidean-particular} 
is an absolutely continuous measure with a continuous piece-wise affine density function,
which can be skewed, multi-modal, and can approximate any probability measure on $\CX_i$ in the topology of weak convergence.
Thus, Setting~\ref{sett:Euclidean} provides a flexible model for 
matching for teams problems in practice;
see also our discussions in Section~\ref{ssec:model-assumptions}.

We show in the following theorem that the approximation error in Theorem~\ref{thm:mt-tf-approx} can be controlled to be arbitrarily small under Setting~\ref{sett:Euclidean}.

\begin{theorem}[Controlling the error]\label{thm:error-control}
    Under Setting~\ref{sett:Euclidean},
    the following statements hold.
    \begin{enumerate}[label=(\roman*),beginpenalty=10000]
        \item\label{thms:error-control-W1}
        For $i=0,1,\ldots,N$ and for any subdivision $\FS_i$ of $\FK_i$,
        it holds for
        $\CG_i$ and $\overline{\epsilon}_i(\,\cdot\,)$ defined in Setting~\ref{sett:Euclidean} that
        $\overline{\epsilon}_i(\varsigma)\ge \sup\big\{W_1(\nu,\nu'):\nu,\nu'\in\CP(\CX_i),\;\nu\abovebelowset{\CG_i}{\varsigma}{\sim}\nu'\big\}$ 
        $\forall \varsigma\ge 0$.
        
        \item\label{thms:error-control-arbitrary}
        For any $\epsilon\in(0,1)$,
        the subdivisions $(\FS_i)_{i=0:N}$ of $(\FK_i)_{i=0:N}$ in Setting~\ref{sett:Euclidean}
        can be explicitly constructed 
        to satisfy $\eta(\FS_i)\le \frac{\epsilon}{16NL_c}$ 
        and 
        $m_i:=|V(\FS_i)|-1
        \le \big(\frac{64ND_id_i(L_c\vee 1)}{\epsilon}\big)^{d_i} |\FM(\FK_i)|(d_i+\nobreak1)$
        for $i=\nobreak0,1,\ldots,N$,
        where 
        $D_i:=\max_{\BIx,\BIx'\in\CX_i}\big\{\|\BIx-\BIx'\|_2\big\}$ 
        for $i=\nobreak0,1,\ldots,N$;
        see the explicit procedure $\mathtt{InitTestFuncs}$
        and Lemma~\ref{lem:proc-init-test-funcs} 
        in the e-companion.
        Subsequently, 
        applying Theorem~\ref{thm:mt-tf-approx}
        with respect to $\epsilon_{\sspar}\leftarrow \nobreak\frac{\epsilon}{2}$,
        $\varsigma\leftarrow\nobreak \frac{\epsilon}{8NL_c}\Big(\max_{0\le i\le N}\big\{\max_{\BIx,\BIx'\in\CX_i}\big\{\|\BIx-\BIx'\|_2\big\}\big\}\Big)^{-1}$,
        along with $(\overline{\epsilon}_i(\varsigma),\CG_i)_{i=0:N}$ defined with respect to $(\FS_i)_{i=0:N}$
        yield two $\epsilon$-matching equilibria
        $(\hat{\varphi}_i,\hat{\gamma}_i)_{i=1:N},\hat{\nu}$ and $(\hat{\varphi}_i,\tilde{\gamma}_i)_{i=1:N},\tilde{\nu}$.
    \end{enumerate}
\end{theorem}

In particular, 
Theorem~\ref{thm:error-control}\ref{thms:error-control-W1} 
quantifies the dependence of the error $\epsilon_{\ssapx}$ in Theorem~\ref{thm:mt-tf-approx}
on the choice of the test functions $(\CG_i)_{i=0:N}$
via the mesh sizes $(\eta(\FS_i))_{i=0:N}$.
We remark that 
it is in general not possible to obtain 
a true matching equilibrium via Theorem~\ref{thm:mt-tf-approx}
unless $(\CX_i)_{i=0:N}$ are all finite, 
due to the violation of the marginal constraints in \eqref{eqn:mt-primalopt}.

\subsection{Explicit computational complexity bound}\label{ssec:complexity}

In this subsection, we work under Setting~\ref{sett:Euclidean}, and
we present an algorithm that uses Theorem~\ref{thm:mt-tf-approx} and Theorem~\ref{thm:error-control} as a template 
to compute an \mbox{$\epsilon$-matching} equilibrium
for any $\epsilon>0$.
Crucially, this algorithm yields an explicit complexity bound 
under the arithmetic model (also known as the real random-access machine),
where all input, output, and intermediate values are assumed to be stored as exact real numbers,
and its complexity is characterized by the arithmetic operations incurred, i.e., addition, subtraction, multiplication, division, and comparison of real numbers.

Due to its lengthiness, we only present here
a sketch of the overall algorithm in Algorithm~\ref{alg:complexity-overall}
as well as a brief overview of the procedures referenced in Algorithm~\ref{alg:complexity-overall},
and we refer the reader to Section~\ref{apx:complexity} of the e-companion 
for the detailed procedures.
In Line~\ref{alglin:complexity-overall-init-test-funcs-forloop},
$\mathtt{InitTestFuncs}$
constructs a subdivision $\FS_i$ of $\FK_i$ with 
$\eta(\FS_i)\le \frac{\epsilon}{16NL_c}$,
and returns the constructed test functions $\CG_i$
as well as an upper bound 
$\overline{D}_i\ge \max_{\BIx,\BIx'\in\CX_i}\big\{\|\BIx-\nobreak\BIx'\|_2\big\}$;
see Lemma~\ref{lem:proc-init-test-funcs}.
In Line~\ref{alglin:complexity-overall-init-type-measure-forloop},
$\mathtt{InitTypeMeasure}$ computes 
the objective vector $\bar{\BIg}_i$ in \eqref{eqn:mt-tf-lsip}
by evaluating the integral $\int_{\CX_i}g_{i,j}\DIFFX{\mu_i}$ for $j=\nobreak0,1,\ldots,m_i$;
see Lemma~\ref{lem:proc-init-type-measure}.
Subsequently,
$\mathtt{SolveLSIPDual}$ and $\mathtt{SolveLSIPPrimal}$ in 
Lines~\ref{alglin:complexity-overall-solve-LSIP-dual}--\ref{alglin:complexity-overall-solve-LSIP-primal}
utilize the volumetric center algorithm of \citet{vaidya1996new}
and the LP algorithm of \citet{vandenbrand2020deterministic}
to compute a pair of 
$(\epsilon_{\sspar},\varsigma)$-optimizers of 
\eqref{eqn:mt-tf-dual} and \eqref{eqn:mt-tf-lsip}
with respect to 
$\epsilon_{\sspar}\leftarrow \frac{\epsilon}{2}$ and
$\varsigma\leftarrow\nobreak \frac{\epsilon}{8N(L_c\vee 1)\max_{0\le i\le N}\{\overline{D}_i\}}$;
see Theorem~\ref{thm:proc-solve-LSIP}.
Lastly, 
$\mathtt{ConstructTransFuncs}$ and
$\mathtt{ConstructAllocs}$ in
Lines~\ref{alglin:complexity-overall-construct-trans-funcs}--\ref{alglin:complexity-overall-construct-allocs}
carry out the construction of an 
$\epsilon$-matching equilibrium $(\hat{\varphi}_i,\hat{\gamma}_i)_{i=1:N},\hat{\nu}$
via Theorem~\ref{thm:mt-tf-approx};
see Lemma~\ref{lem:proc-construct-trans-funcs} and Lemma~\ref{lem:proc-construct-allocs}.
The detailed complexity results are presented in the following theorem.

\IncMargin{3pt}
\begin{algorithm}[t]
    \DontPrintSemicolon
    \caption{\textbf{Computation of AME with explicit complexity bound}}\label{alg:complexity-overall}
    \begin{footnotesize}
    \KwIn{$N\ge2$, 
    $(\CX_i)_{i=0:N}$, 
    $(\mu_i,c_i)_{i=1:N}$
    (see Setting~\ref{sett:Euclidean}), 
    error tolerance $\epsilon\in(0,1)$,
    Lipschitz constant $L_c>0$.}



    \nl\label{alglin:complexity-overall-init-test-funcs-forloop}%
    \textbf{for} $i=0,1,\ldots,N$ \textbf{do} 
    $(\CG_i,\overline{D}_i)\leftarrow \mathtt{InitTestFuncs}\big(\CX_i,\frac{\epsilon}{16N(L_c\vee 1)}\big)$.


    \nl\label{alglin:complexity-overall-init-type-measure-forloop}%
    \textbf{for} $i=1,\ldots,N$ \textbf{do}
    $(\bar{\BIg}_i,\mathtt{Aux}_i)\leftarrow \mathtt{InitTypeMeasure}(\CX_i,\mu_i,\CG_i)$.



    \nl\label{alglin:complexity-overall-solve-LSIP-dual}%
    $\big((\hat{y}_{i,0},\hat{\BIy}_i,\hat{\BIw}_i)_{i=1:N},\widehat{\CI}\big)\leftarrow \mathtt{SolveLSIPDual}\big(N,(\CG_i)_{i=0:N},(\bar{\BIg}_i,c_i)_{i=1:N},\frac{\epsilon}{4},L_c\big)$.

    \nl\label{alglin:complexity-overall-solve-LSIP-primal}%
    $(\hat{\theta}_{i})_{i=1:N}\leftarrow\mathtt{SolveLSIPPrimal}\Big(N,(\CG_i)_{i=0:N},(\bar{\BIg}_i,c_i)_{i=1:N},\widehat{\CI},\frac{\epsilon}{8N(L_c\vee 1)\max_{0\le i\le N}\{\overline{D}_i\}}\Big)$.

    \nl\label{alglin:complexity-overall-construct-trans-funcs}%
    $(\hat{\varphi}_i)_{i=1:N}\leftarrow \mathtt{ConstructTransFuncs}\big(\CG_0,(\hat{\BIw}_i)_{i=1:N}\big)$.

    \nl\label{alglin:complexity-overall-construct-allocs}%
    $\big((\hat{\gamma}_i)_{i=1:N},\hat{\nu}\big)\leftarrow \mathtt{ConstructAllocs}\big(N,(\CX_i,\CG_i)_{i=0:N},(\mu_i,\mathtt{Aux}_i,\hat{\theta}_i)_{i=1:N}\big)$.
    
    \KwOut{$\big((\hat{\varphi}_i,\hat{\gamma}_i)_{i=1:N},\hat{\nu}\big)$.}
    \end{footnotesize}
\end{algorithm}

\begin{theorem}[Computational complexity of Algorithm~\ref{alg:complexity-overall}]\label{thm:complexity-explicit}
    Under Setting~\ref{sett:Euclidean},
    let us denote $d:=\max_{0\le i\le N}\{d_i\}$,
    $J:=\max_{0\le i\le N}\big\{|\FM(\FK_i)|\big\}$,
    $D:=\max_{0\le i\le N}\big\{\!\max_{\BIx,\BIx'\in\CX_i}\big\{\|\BIx-\BIx'\|_2\big\}\!\big\}\vee\nobreak 1$,
    $\overline{L}_{c}:=L_{c}\vee 1$,
    $I_{\ssmin}:=\min_{1\le i\le N}\big\{\!\min_{K\in\FM(\FK_i)}\big\{\int_{K}p_{i,K}\DIFFX{\Lebesgue_K}\big\}\!\big\}>0$,
    and denote the computational complexity of multiplying two 
    $n$-by-$n$ matrices by $O(n^{\omega})$.\footnote{For example, with the standard matrix multiplication procedure, the computational complexity of this operation is $O(n^3)$. However, it is known that $\omega<2.38$; see, e.g., \citep{coppersmith1990matrix}.}
    Let us make the following assumptions.
    \begin{itemize}
        \item For $i=1,\ldots,N$, evaluating 
        $c_i(\BIx,\BIz)$ at any $(\BIx,\BIz)\in\CX_i\times\CX_0$ incurs 
        $T_{\mathtt{Eval}}$ arithmetic operations.

        \item For $i=1,\ldots,N$ and any $n$-simplex $C\subseteq K$ for some $K\in \FK_i$ (where $n\in\{0,1,\ldots,d_i\}$),
        \begin{itemize}
            \item evaluating 
            $\big(\int_{\CX_i}\lambda_{\FF(C),\BIv}\INDI_{\relint(C)}\DIFFX{\mu_i}\big)_{\BIv\in V(C)}\subset\R_{+}$ incurs $T_{\mathtt{Int}}$ arithmetic operations, and

            \item if $\mu_i\big(\relint(C)\big)>0$, 
            then generating a random sample from 
            $\tilde{\mu}_{i,C}\in\CP(\CX_i)$ defined by $\frac{\DIFF\tilde{\mu}_{i,C}}{\DIFF\mu_i}:=\frac{\INDI_{\relint(C)}}{\mu_i(\relint(C))}$ incurs $T_{\mathtt{Samp}}$ arithmetic operations.
        \end{itemize}
    \end{itemize}
    \begin{enumerate}[label=(\roman*),beginpenalty=10000]
        \item\label{thms:complexity-explicit-general}
        Then, for any $\epsilon\in(0,1)$,
        Algorithm~\ref{alg:complexity-overall} incurs
        $O\big(\big(\frac{64NDd\overline{L}_{c}}{\epsilon}\big)^{(\omega+1)d}N^{\omega+1}J^{\omega+1}d^{\omega+2}\log\big(\frac{NJDd\overline{L}_{c}}{I_{\ssmin}\epsilon}\big)
        +\big(\frac{64NDd\overline{L}_{c}}{\epsilon}\big)^{3d}N^{3}J^{3}d^{4}\log\big(\frac{NJDd\overline{L}_{c}}{I_{\ssmin}\epsilon}\big)T_{\mathtt{Eval}}
        +\big(\frac{128NDd\overline{L}_{c}}{\epsilon}\big)^{d}NJT_{\mathtt{Int}}\big)$
        arithmetic operations
        to compute 
        an $\epsilon$-matching equilibrium $(\hat{\varphi}_i,\hat{\gamma}_i)_{i=1:N},\hat{\nu}$.
        Subsequently, 
        generating a random sample from $\hat{\nu}$ incurs 
        $O\big(\big(\frac{64NDd\overline{L}_{c}}{\epsilon}\big)^{d}Jd\big)$
        arithmetic operations,
        and for $i=1,\ldots,N$,
        evaluating $\hat{\varphi}_i(\BIz)$ at any $\BIz\in\CX_0$ incurs 
        $O\big(\big(\frac{64NDd\overline{L}_{c}}{\epsilon}\big)^{d}Jd^{\omega}\big)$
        arithmetic operations,
        and 
        generating a random sample from $\hat{\gamma}_i$ incurs
        $O\big(\big(\frac{128NDd\overline{L}_{c}}{\epsilon}\big)^{d}J+T_{\mathtt{Samp}}\big)$
        arithmetic operations;
        see the procedures
        $\mathtt{EvalTransFunc}$,
        $\mathtt{SampQual}$, and 
        $\mathtt{SampAlloc}$ in Section~\ref{sapx:complexity-AME} of the e-companion.
        \item\label{thms:complexity-explicit-particular-case}
        Under the particular class of
        $(\mu_i)_{i=1:N}$ in Setting~\ref{sett:Euclidean}\ref{settc:Euclidean-particular},
        statement~\ref{thms:complexity-explicit-general}
        holds with respect to 
        $T_{\mathtt{Int}}=O\big(d^{\omega+1}+\log(J)\big)$ and 
        $T_{\mathtt{Samp}}=O\big(d^{\omega+1}+\log(J)\big)$.
    \end{enumerate}
\end{theorem}

We remark that the constant terms omitted by the big-$O$ notations in Theorem~\ref{thm:complexity-explicit} 
are absolute constants and do not depend on any input or parameter of Algorithm~\ref{alg:complexity-overall}.
If we hold $D,J,d,\allowbreak \overline{L}_{c},\allowbreak I_{\ssmin},T_{\mathtt{Eval}}$ fixed and treat them as constants,
then the number of arithmetic operations incurred by Algorithm~\ref{alg:complexity-overall} in
Theorem~\ref{thm:complexity-explicit}\ref{thms:complexity-explicit-particular-case}
simplifies to 
$O\big(N^{(\omega+1)(d+1)}\epsilon^{-(\omega+1)d}\log\big(\frac{N}{\epsilon}\big)\big)$,
which is polynomial in $N$ and $\frac{1}{\epsilon}$.
With the upper bound $\omega<2.38$,
this complexity is 
$O\big(N^{6.76}\epsilon^{-3.38}\log\big(\frac{N}{\epsilon}\big)\big)$
when $d=1$,
$O\big(N^{10.14}\epsilon^{-6.76}\log\big(\frac{N}{\epsilon}\big)\big)$
when $d=2$, and
$O\big(N^{13.52}\epsilon^{-10.14}\log\big(\frac{N}{\epsilon}\big)\big)$
when $d=3$.

\subsection{A practically efficient algorithm}\label{ssec:efficient-algo}
Although Algorithm~\ref{alg:complexity-overall}
admits an explicit complexity bound, 
it utilizes the volumetric center algorithm of \citet{vaidya1996new} for solving the LSIP problem \eqref{eqn:mt-tf-lsip},
which is outperformed in practice by 
cutting-plane type algorithms utilizing state-of-the-art linear programming (LP) solvers.
This dichotomy between theoretical and practical performance of optimization algorithms can also be seen in, e.g., \citep{kahale2017superreplication, altschuler2021wasserstein, altschuler2023polynomial}.
Modern LP solvers such as the Gurobi optimizer \citep{gurobi}
use highly-efficient implementations of the simplex algorithm or the interior point algorithm 
which can effectively exploit sparsity patterns.
For this reason,
we will develop in this subsection a practically efficient algorithm
which utilizes modern LP solvers
for computing $\epsilon$-matching equilibrium for any $\epsilon>0$.

\begin{algorithm}[t]
\caption{\textbf{Cutting-plane discretization algorithm}}\label{alg:cp-tf}
\begin{footnotesize}

\KwIn{$N\ge2$, 
    $(\CX_i,\CG_i)_{i=0:N}$, 
    $(\mu_i,c_i,\bar{\BIg}_i)_{i=1:N}$
    (see Setting~\ref{sett:Euclidean}), 
    $\big(\CK_i^{(0)}\big)_{i=1:N}$,
    error tolerance $\epsilon_{\sspar}>0$.}


\nl\label{alglin:cp-tf-global-tol}%
Choose a sequence $(\tau^{(r)})_{r\in\N_0}\subset (0,\infty)$
with $\limsup_{r\to\infty}\tau^{(r)}<\frac{\epsilon_{\sspar}}{N}$.

\nl%
\For{$r=0,1,2,\ldots$}{

    \nl\label{alglin:cp-tf-lp}%
    Solve the LP s 
    \eqref{eqn:mt-tf-algo-lp} 
    and  
    \eqref{eqn:mt-tf-algo-lp-dual}, 
    denote their common optimal value by $\alpha^{(r)}$, 
    and denote the computed optimizers 
    of \eqref{eqn:mt-tf-algo-lp} 
    and  
    \eqref{eqn:mt-tf-algo-lp-dual}
    by 
    $\big(y_{i,0}^{(r)},\BIy_i^{(r)},\BIw_i^{(r)}\big)_{i=1:N}$ 
    and 
    $\big(\theta^{(r)}_{i,\BIx,\BIz}\big)_{(\BIx,\BIz)\in\CK_i^{(r)},i=1:N}$, $\Bxi^{(r)}$, 
    respectively.

    \nl\label{alglin:cp-tf-inner-forloop}%
    \For{$i=1,\ldots,N$}{

        \nl\label{alglin:cp-tf-global}%
        Approximately solve
        $\displaystyle\min_{\BIx\in\CX_i,\,\BIz\in\CX_0}\big\{c_i(\BIx,\BIz)-\langle\BIg_i(\BIx),\BIy_i^{(r)}\rangle-\langle\BIg_0(\BIz),\BIw_i^{(r)}\rangle\big\}$
        to compute $(\tilde{\BIx}_i^{(r)},\tilde{\BIz}_i^{(r)})\in\allowbreak\CX_i\times\CX_0$
        and $\underline{\beta}_i^{(r)}\in\R$ satisfying 
        $\underline{\beta}_i^{(r)}
        \le \displaystyle\min_{\BIx\in\CX_i,\,\BIz\in\CX_0}\big\{c_i(\BIx,\BIz)-\langle\BIg_i(\BIx),\BIy_i^{(r)}\rangle-\langle\BIg_0(\BIz),\BIw_i^{(r)}\rangle\big\}
        \le \allowbreak c_i(\tilde{\BIx}_i^{(r)},\tilde{\BIz}_i^{(r)})\allowbreak-\langle\BIg_i(\tilde{\BIx}_i^{(r)}),\BIy_i^{(r)}\rangle-\langle\BIg_0(\tilde{\BIz}_i^{(r)}),\BIw_i^{(r)}\rangle
        \le \underline{\beta}_i^{(r)} + \tau^{(r)}$.

        \nl\label{alglin:cp-tf-aggregate}%
        Let 
        $\widetilde{\CK}_i^{(r)}\subseteq\CX_i\times\CX_0$ 
        be a finite set such that 
        $\big(\tilde{\BIx}_i^{(r)},\tilde{\BIz}_i^{(r)}\big)
        \in\widetilde{\CK}_i^{(r)}$, 
        and update 
        $\CK_i^{(r+1)}\leftarrow \CK^{(r)}_i\cup\widetilde{\CK}_i^{(r)}$.

    }

    \nl\label{alglin:cp-tf-termination}%
    \textbf{if}
    $\sum_{i=1}^N(y_{i,0}^{(r)}-\underline{\beta}_i^{(r)})\le\epsilon_{\sspar}$ \textbf{then} skip to Line~\ref{alglin:cp-tf-bounds}, \textbf{else} continue to the next iteration.
}

\nl\label{alglin:cp-tf-bounds}%
$\alpha_{\sspar}^{\ssLB}\leftarrow\alpha^{(r)}-\big(\sum_{i=1}^N y_{i,0}^{(r)}-\underline{\beta}_i^{(r)}\big)$,
$\alpha_{\sspar}^{\ssUB}\leftarrow\alpha^{(r)}$.

\nl\label{alglin:cp-tf-primal-dual}%
\textbf{for} $i=1,\ldots,N$ \textbf{do}
$\hat{\theta}_i\leftarrow\sum_{(\BIx,\BIz)\in\CK_i^{(r)}}\theta^{(r)}_{i,\BIx,\BIz}\delta_{(\BIx,\BIz)}$,
$\hat{y}_{i,0}\leftarrow \underline{\beta}_i^{(r)}$, 
$\hat{\BIy}_i\leftarrow\BIy_i^{(r)}$, 
$\hat{\BIw}_i\leftarrow\BIw_i^{(r)}$.
    


\KwOut{$(\hat{\theta}_i)_{i=1:N}$,
$(\hat{y}_{i,0},\hat{\BIy}_i,\hat{\BIw}_i)_{i=1:N}$, $\alpha_{\sspar}^{\ssLB}$, 
$\alpha_{\sspar}^{\ssUB}$.} 
\end{footnotesize}
\end{algorithm}

Let us first develop a cutting-plane discretization algorithm inspired by the Conceptual Algorithm~11.4.1 in \citep{goberna1998linear} for solving the LSIP problems \eqref{eqn:mt-tf-dual} and \eqref{eqn:mt-tf-lsip} under Setting~\ref{sett:Euclidean}.
In iteration $r\in\N_0$ of the algorithm, we replace the semi-infinite constraint 
$y_{i,0}
+\langle\BIg_i(\BIx),\BIy_i\rangle
+\langle\BIg_0(\BIz),\BIw_i\rangle  
\le c_i(\BIx,\BIz)$ 
$\forall (\BIx,\BIz)\in\CX_i\times\CX_0$, $\forall i\in\{1,\ldots,N\}$
in \eqref{eqn:mt-tf-lsip} with finitely many constraints characterized by finite sets $\CK_1^{(r)}\subseteq\CX_1\times\CX_0,\ldots,\CK_N^{(r)}\subseteq\CX_N\times\CX_0$ to relax \eqref{eqn:mt-tf-lsip} by the following LP problem \eqref{eqn:mt-tf-algo-lp} and its dual \eqref{eqn:mt-tf-algo-lp-dual}:
\begin{align}
\maximize_{(y_{i,0},\BIy_i,\BIw_i)}\quad & \sum_{i=1}^N y_{i,0}+\langle\bar{\BIg}_i,\BIy_i\rangle 
\tag{$\MT^{*(r)}_{\sspar}$}
\label{eqn:mt-tf-algo-lp}\\*
\mathrm{subject~to}\quad & 
y_{i,0}
+\langle\BIg_i(\BIx),\BIy_i\rangle
+\langle\BIg_0(\BIz),\BIw_i\rangle  
\le c_i(\BIx,\BIz) 
\quad\forall (\BIx,\BIz)\in\CK_i^{(r)},\; \forall i\in\{1,\ldots,N\},\nonumber\\*
& \sum_{i=1}^N\BIw_i=\veczero_{m_0}, 
\qquad y_{i,0}\in\R,\; \BIy_i\in\R^{m_i},\; \BIw_i\in\R^{m_0} 
\quad \forall i\in\{1,\ldots,N\},\nonumber\\
\minimize_{(\theta_{i,\BIx,\BIz}),\,\Bxi}\quad 
& \sum_{i=1}^N\sum_{(\BIx,\BIz)\in\CK_i^{(r)}}\theta_{i,\BIx,\BIz}c_i(\BIx,\BIz)\label{eqn:mt-tf-algo-lp-dual}
\tag{$\MT^{(r)}_{\sspar}$}\\*
\mathrm{subject~to}\quad 
& \sum_{(\BIx,\BIz)\in\CK_i^{(r)}}\theta_{i,\BIx,\BIz}=1,\; 
\sum_{(\BIx,\BIz)\in\CK_i^{(r)}}\theta_{i,\BIx,\BIz}\BIg_i(\BIx)=\bar{\BIg}_i, \; 
\sum_{(\BIx,\BIz)\in\CK_i^{(r)}}\theta_{i,\BIx,\BIz}\BIg_0(\BIz)=\Bxi 
\quad \forall i\in\{1,\ldots,N\}, \nonumber \\*
& \Bxi\in\R^{m_0}, \qquad  \theta_{i,\BIx,\BIz}\in\R_+  \quad \forall (\BIx,\BIz)\in\CK_i^{(r)},\; \forall i\in\{1,\ldots,N\}.  \nonumber
\end{align}
We begin with an initial LP relaxation $\CPDLPDUAL{0}$,
and we iteratively add more constraints, i.e., by adding elements to $\big(\CK_i^{(r)}\big)_{i=1:N}$ to ``cut'' the feasible set of \eqref{eqn:mt-tf-algo-lp}, 
until the constraint violation falls below the pre-specified threshold~$\epsilon_{\sspar}$.
Our cutting-plane discretization algorithm is presented in Algorithm~\ref{alg:cp-tf} and its properties are presented in Proposition~\ref{prop:cp-tf-properties}. 
We present in Sections~\ref{sapx:efficient-algo-LP}--\ref{sapx:efficient-algo-global} additional remarks about Algorithm~\ref{alg:cp-tf},
including how the LP problem in Line~\ref{alglin:cp-tf-lp}
and the global minimization problem in Line~\ref{alglin:cp-tf-global} are solved.

\begin{proposition}[Properties of Algorithm~\ref{alg:cp-tf}]\label{prop:cp-tf-properties}\vspace{-4pt}
The following statement holds.
\begin{enumerate}[label=(\roman*),beginpenalty=10000]
    \item\label{props:cp-tf-boundedness}%
    Under Setting~\ref{sett:Euclidean}, there exist finite sets
    $\CK_i^{(0)}\subseteq \CX_i\times\CX_0$ $\forall i\in\{1,\ldots,N\}$
    such that 
    the set of optimizers of $\CPDLPDUAL{0}$ is bounded.
    In particular,
    one may choose 
    $\CK_i^{(0)}:=\big\{(\BIv_{i,j},\BIv_{0,0}):j\in\nobreak\{0,1,\ldots,m_i\}\big\}
    \cup \big\{(\BIv_{i,0},\BIv_{0,l}):l\in\{1,\ldots,m_0\}\big\}$
    for $i=1,\ldots,N$.

\end{enumerate}
In the following, let us work under Setting~\ref{sett:Euclidean} and let
$\big(\CK_i^{(0)}\big)_{i=1:N}$
be such that the set of optimizers of $\CPDLPDUAL{0}$ is bounded.
Then, for any $\epsilon_{\sspar}>0$, the following statements hold.
\begin{enumerate}[label=(\roman*),topsep=0pt,parsep=0pt,partopsep=0pt,itemsep=-1pt,beginpenalty=10000]
    \setcounter{enumi}{1}
    \item\label{props:cp-tf-termination} %
    Algorithm~\ref{alg:cp-tf} terminates after finitely many iterations. 
    
    \item\label{props:cp-tf-bounds} 
    $\alpha_{\sspar}^{\ssLB}\le\alpha_{\sspar}^{\star}\le\alpha_{\sspar}^{\ssUB}$ where $\alpha_{\sspar}^{\ssUB}-\alpha_{\sspar}^{\ssLB}\le\epsilon_{\sspar}$.

    \item\label{props:cp-tf-dual}%
    $(\hat{y}_{i,0},\hat{\BIy}_i,\hat{\BIw}_i)_{i=1:N}$ 
    is an 
    $\epsilon_{\sspar}$-optimizer of \eqref{eqn:mt-tf-lsip} 
    and 
    $\sum_{i=1}^N\hat{y}_{i,0}+\langle\bar{\BIg}_i,\hat{\BIy}_i\rangle=\alpha_{\sspar}^{\ssLB}$. 
    
    \item\label{props:cp-tf-primal}%
    $(\hat{\theta}_i)_{i=1:N}$ is an 
    $\epsilon_{\sspar}$-optimizer of \eqref{eqn:mt-tf-dual} where
    $\sum_{i=1}^N\int_{\CX_i\times\CX_0}c_i\DIFFX{\hat{\theta}_i}=\alpha_{\sspar}^{\ssUB}$
    and 
    each $\hat{\theta}_i$ has finite support.
\end{enumerate}
Consequently, $(\hat{\theta}_i)_{i=1:N}$, $(\hat{y}_{i,0},\hat{\BIy}_i,\hat{\BIw}_i)_{i=1:N}$
is a pair of $(\epsilon_{\sspar},0)$-optimizers of \eqref{eqn:mt-tf-dual} and \eqref{eqn:mt-tf-lsip}.
\end{proposition}
Note that unlike Theorem~\ref{thm:complexity-explicit}, Proposition~\ref{prop:cp-tf-properties} only shows the finite termination of Algorithm~\ref{alg:cp-tf} without an explicit computational complexity.
Despite this, Algorithm~\ref{alg:cp-tf}
is highly efficient in practice thanks to the advancement of modern LP solvers, as shown by our experiments in Section~\ref{sec:experiments}.

Next,
we work under the following additional assumptions 
for $(\CX_i,\mu_i)_{i=1:N}$
which allow tractable construction of $W_1$-optimal couplings 
via classical results that we recall in Proposition~\ref{prop:OT-construction}.

\begin{setting}\label{sett:typemeasure}
For $i=1,\ldots,N$, $\CX_i$ and $\mu_i$ satisfy one of the following conditions:
\begin{enumerate}[label=\normalfont{(A\arabic*)},leftmargin=28pt,topsep=0pt,parsep=0pt,partopsep=0pt,itemsep=-2pt,beginpenalty=10000]
    \item\label{setts:typemeasure-discrete}
    $\mu_i$ is supported on finitely many points, i.e., every $K\in\FM(\FK_i)$ consists of a single point;
    
    \item\label{setts:typemeasure-continuous}
    $\mu_i$ is absolutely continuous with respect to the Lebesgue measure on~$\CX_i$;
    
    \item\label{setts:typemeasure-1d}
    $\CX_i\subset\R$, i.e., $d_i=1$.
\end{enumerate}
\end{setting}

\begin{algorithm}[t]
\caption{\textbf{Construction of AME}}\label{alg:mt-tf}
\begin{footnotesize}
\KwIn{$N\ge2$, 
    $(\CX_i,\CG_i)_{i=0:N}$, 
    $(\mu_i,c_i,\bar{\BIg}_i)_{i=1:N}$
    (see Setting~\ref{sett:Euclidean}), 
    $\big(\CK_i^{(0)}\big)_{i=1:N}$,
    error tolerance $\epsilon_{\sspar}>0$.}


\nl\label{alglin:mt-tf-cpalgo}%
$\big((\hat{\theta}_i)_{i=1:N}$,
$(\hat{y}_{i,0},\hat{\BIy}_i,\hat{\BIw}_i)_{i=1:N}$,
$\alpha_{\sspar}^{\ssLB}$, 
$\alpha_{\sspar}^{\ssUB}\big)\leftarrow$
the output of Algorithm~\ref{alg:cp-tf}.

\nl\For{$i=1,\ldots,N$}{
    \nl\label{alglin:mt-tf-dualsol}%
    $\hat{\varphi}_i(\,\cdot\,)\leftarrow\langle\BIg_0(\,\cdot\,),\hat{\BIw}_i\rangle$.

    \nl\label{alglin:mt-tf-marginals}%
    $\hat{\mu}_i\leftarrow $ the marginal of $\hat{\theta}_i$ on $\CX_{i}$,
    $\hat{\nu}_i\leftarrow $ the marginal of $\hat{\theta}_i$ on $\CX_{0}$.
    
}

\nl\label{alglin:mt-tf-primal1}%
$\hat{\nu}\leftarrow \hat{\nu}_{1}$.
Let $(\Omega,\CF,\PROB)$ be a probability space 
and let 
$\BIZ:\Omega\to\CX_{0}$ be a random variable with law $\hat{\nu}$.

\nl\label{alglin:mt-tf-bindingloop}%
\For{$i=1,\ldots,N$}{

    \nl\label{alglin:mt-tf-binding1}%
    Construct $\BIZ_i:\Omega\to\CX_0$ with law $\hat{\nu}_i$ such that 
    $\EXP\big[\|\BIZ-\BIZ_i\|_2\big]=W_1(\hat{\nu},\hat{\nu}_i)$
    (e.g., via Proposition~\ref{prop:OT-construction}).

    \nl\label{alglin:mt-tf-binding2}%
    Define $\BIX_i:\Omega\to\CX_i$ 
    such that 
    $\PROB[\BIX_i\in E|\BIZ_i=\BIz]
    =\frac{\hat{\theta}_i(E\times\{\BIz\})}{\hat{\nu}_i(\{\BIz\})}$ 
    $\forall \BIz\in\support(\hat{\nu}_i)$, 
    $\forall E\in\CB(\CX_i)$.

    \nl\label{alglin:mt-tf-reassembly}%
    Construct 
    $\overbar{\BIX}_i:\Omega\to\CX_i$ with law $\mu_i$ such that 
    $\EXP\big[\|\BIX_i-\overbar{\BIX}_i\|_2\big]=W_1(\hat{\mu}_i,\mu_i)$
    (e.g., via Proposition~\ref{prop:OT-construction}).

}

\nl\label{alglin:mt-tf-primal2}%
Let $\overbar{\BIZ}:\Omega\to\CX_0$ satisfy 
$\overbar{\BIZ}\in\argmin_{\BIz\in\CX_0}\big\{\sum_{i=1}^{N}c_i(\overbar{\BIX}_i,\BIz)\big\}$
$\PROB$-almost surely.
$\tilde{\nu}\leftarrow$ the law of $\overbar{\BIZ}$.

\nl\label{alglin:mt-tf-primalcoup}%
\textbf{for} $i=1,\ldots,N$ \textbf{do} 
$\hat{\gamma}_i\leftarrow$ 
the law of 
$(\overbar{\BIX}_i,\BIZ)$, 
$\hat{\alpha}_i\leftarrow\EXP\big[c_i(\overbar{\BIX}_i,\BIZ)\big]$. 
$\tilde{\gamma}_i\leftarrow$ 
the law of 
$(\overbar{\BIX}_i,\overbar{\BIZ})$, 
$\tilde{\alpha}_i\leftarrow\EXP\big[c_i(\overbar{\BIX}_i,\overbar{\BIZ})\big]$.



\nl\label{alglin:mt-tf-bounds}%
$\hat{\alpha}_{\MT}^{\ssLB}\leftarrow 
\alpha_{\sspar}^{\ssLB}$, 
$\hat{\alpha}_{\MT}^{\ssUB}\leftarrow 
\sum_{i=1}^N\hat{\alpha}_i$, 
$\tilde{\alpha}_{\MT}^{\ssUB}\leftarrow 
\sum_{i=1}^N\tilde{\alpha}_i$. 
$\hat{\epsilon}_{\sssub}\leftarrow 
\hat{\alpha}_{\MT}^{\ssUB}-\hat{\alpha}_{\MT}^{\ssLB}$, 
$\tilde{\epsilon}_{\sssub}\leftarrow 
\tilde{\alpha}_{\MT}^{\ssUB}-\hat{\alpha}_{\MT}^{\ssLB}$. 

\KwOut{
$(\hat{\varphi}_i,\hat{\gamma}_i)_{i=1:N},\hat{\nu}$,
$(\hat{\varphi}_i,\tilde{\gamma}_i)_{i=1:N},\tilde{\nu}$,
$\hat{\alpha}_{\MT}^{\ssLB}$,
$\hat{\alpha}_{\MT}^{\ssUB}$,
$\tilde{\alpha}_{\MT}^{\ssUB}$,
$\hat{\epsilon}_{\sssub}$,
$\tilde{\epsilon}_{\sssub}$.}

\end{footnotesize}
\end{algorithm}

Under Settings~\ref{sett:Euclidean} and \ref{sett:typemeasure},
Algorithm~\ref{alg:mt-tf} constructs two 
AMEs via Theorem~\ref{thm:mt-tf-approx}
based on the output of Algorithm~\ref{alg:cp-tf}.
The properties of Algorithm~\ref{alg:mt-tf} are presented in Theorem~\ref{thm:mt-tf}; see also Section~\ref{sapx:efficient-algo-coupling} for additional remarks related to Algorithm~\ref{alg:mt-tf}.

\vspace{-4pt}
\begin{theorem}[Properties of Algorithm~\ref{alg:mt-tf}]\label{thm:mt-tf}%
Under Settings~\ref{sett:Euclidean} and \ref{sett:typemeasure}, 
let $\epsilon_{\sspar}>0$ be arbitrary, 
and let 
$\epsilon_{\sstheo}:=\epsilon_{\sspar}+2L_c\big(N\eta(\FS_{0})+\sum_{i=1}^N\eta(\FS_i)\big)$. 
Then, the following statements hold.
\begin{enumerate}[label=(\roman*),itemsep=-2pt,beginpenalty=10000]
\item\label{thms:mt-tf-algo-equilibrium1}%
$(\hat{\varphi}_i,\hat{\gamma}_i)_{i=1:N},\hat{\nu}$ constitute an $\hat{\epsilon}_{\sssub}$-matching equilibrium.

\item\label{thms:mt-tf-algo-equilibrium2}%
$(\hat{\varphi}_i,\tilde{\gamma}_i)_{i=1:N},\tilde{\nu}$ constitute an $\tilde{\epsilon}_{\sssub}$-matching equilibrium.

\item\label{thms:mt-tf-algo-bound}%
$\hat{\alpha}_{\MT}^{\ssLB}\le\alpha_{\MT}^{\star}\le\tilde{\alpha}_{\MT}^{\ssUB}\le\hat{\alpha}_{\MT}^{\ssUB}$ and $\tilde{\epsilon}_{\sssub}\le\hat{\epsilon}_{\sssub}\le\epsilon_{\sstheo}$.

\item\label{thms:mt-tf-algo-control}%
For any $\epsilon>0$
and any $\epsilon_{\sspar}\in(0,\epsilon)$, 
the subdivisions $(\FS_i)_{i=0:N}$ of $(\FK_i)_{i=0:N}$ in Setting~\ref{sett:Euclidean}
can be explicitly constructed to satisfy 
$\eta(\FS_i)\le \frac{\epsilon-\epsilon_{\sspar}}{4NL_c}$ $\forall 0\le i\le N$;
see the explicit procedure 
$\mathtt{InitTestFuncs}$ and Lemma~\ref{lem:proc-init-test-funcs}
in the e-companion.
Subsequently, Algorithm~\ref{alg:mt-tf} with respect to 
the constructed
$(\CG_i)_{i=0:N}$ yields $\tilde{\epsilon}_{\sssub}\le\hat{\epsilon}_{\sssub}\le\epsilon$.
\end{enumerate}
\end{theorem}

From a theoretical perspective, for any given $\epsilon>0$
and $\epsilon_{\sspar}\in(0,\epsilon)$, 
Theorem~\ref{thm:mt-tf}\ref{thms:mt-tf-algo-control} provides 
\textit{explicit} choices of 
$(\CG_i)_{i=0:N}$
such that 
$(\hat{\varphi}_i,\hat{\gamma}_i)_{i=1:N},\hat{\nu}$ and 
$(\hat{\varphi}_i,\tilde{\gamma}_i)_{i=1:N},\tilde{\nu}$ 
computed by Algorithm~\ref{alg:mt-tf} are \mbox{$\epsilon$-matching} equilibria.
However, in practice, 
one often specifies $(\CG_i)_{i=0:N}$ and $\epsilon_{\sspar}>0$ directly, 
and then uses $\hat{\epsilon}_{\sssub}$ and $\tilde{\epsilon}_{\sssub}$ in the output of Algorithm~\ref{alg:mt-tf} to estimate the sub-optimality of the computed AMEs. 
Note that
$\epsilon_{\sstheo}$ in Theorem~\ref{thm:mt-tf} 
is a theoretical upper bound for $\hat{\epsilon}_{\sssub}$ and $\tilde{\epsilon}_{\sssub}$
derived from Theorem~\ref{thm:mt-tf-approx} using the upper estimates 
$\big(2\eta(\FS_i)\big)_{i=0:N}$ of 
$\big(\sup\big\{W_1(\mu,\mu'):\mu,\mu'\in\CP(\CX_i),\;\mu\overset{\CG_i}{\sim}\mu'\big\}\big)_{i=0:N}$.
Therefore, we call $\epsilon_{\sstheo}$ an a priori upper bound for 
$\hat{\epsilon}_{\sssub}$ and $\tilde{\epsilon}_{\sssub}$ 
since it can be computed independent of Algorithm~\ref{alg:mt-tf}. 
In practice, $\hat{\epsilon}_{\sssub}$ and $\tilde{\epsilon}_{\sssub}$ are typically much less conservative than $\epsilon_{\sstheo}$, as demonstrated by the numerical experiments in Section~\ref{sec:experiments}.

We would like to remark that even when $\CX_0$ is an infinite set,
$\hat{\nu}$ computed by Algorithm~\ref{alg:mt-tf} is a discrete probability measure (see Proposition~\ref{prop:cp-tf-properties}\ref{props:cp-tf-primal}),
whereas $\tilde{\nu}$ is typically non-discrete. 
However, we emphasize that we do not pre-specify the finite support of $\hat{\nu}$ in Algorithm~\ref{alg:mt-tf}.

\section{Numerical experiments}\label{sec:experiments}
In this section, we 
demonstrate our numerical algorithm (i.e., Algorithm~\ref{alg:mt-tf})
in three experiments
to show how it computes approximate matching equilibria that provide managerial insights, 
to show its superior performance compared to five state-of-the-art \mbox{2-Wasserstein} barycenter algorithms,
and to show its empirical scalability to a large number of agent populations (up to $N=1000$).%
\footnote{The code used in the experiments is available at: \url{https://github.com/qikunxiang/MatchingForTeams}; part of our code utilizes the Computational Geometry Algorithms Library by \citet{CGAL} and the Gurobi optimizer by \citet{gurobi}.}

\subsection{Experiment~1: retail business problem}\label{ssec:experiment-biz}
\begin{figure}[t]
    \begin{center}
        \includegraphics[width=1.0\linewidth]{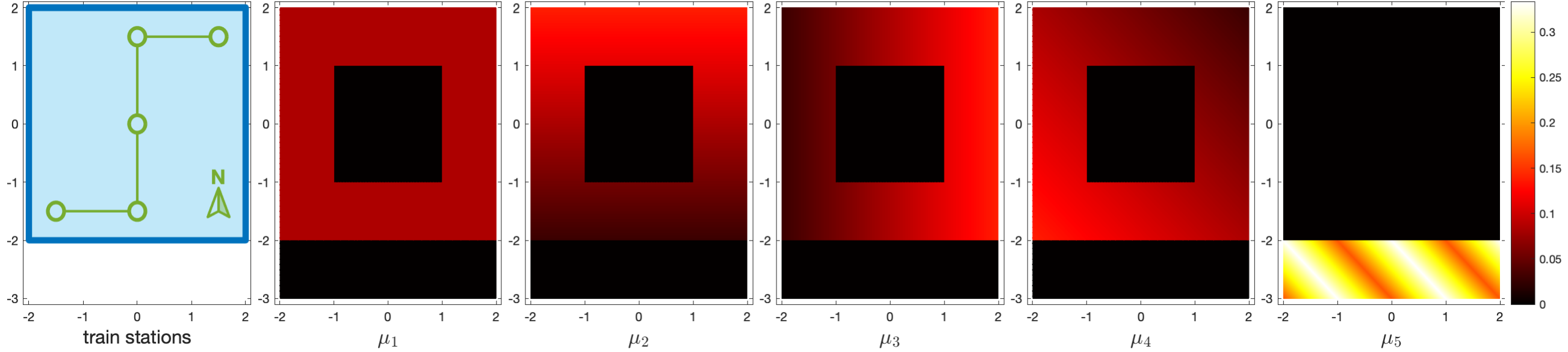}
    \end{center}
    \caption{Experiment~1 -- The railway line, the locations of the train stations in the city, and the probability density functions of $\mu_1,\ldots,\mu_N$.}\label{fig:exp-biz-density}
\end{figure}

In our first experiment, we study the retail business problem presented in the introduction.
Let us consider a square-shaped city and a business hiring 4 categories of employees, that is, 
$N=5$ agent populations are involved.
The cost for each employee to travel on foot between $\BIx\in\R^2$ and $\BIz\in\R^2$ is the scaled city block distance $c_{\mathrm{walk}}\|\BIx-\BIz\|_1$ for $c_{\mathrm{walk}}>0$.
There is a railway line through the city with 5 train stations at $\BIu_1,\ldots,\BIu_5\in\R^2$ (see the leftmost panel of Figure~\ref{fig:exp-biz-density}), where the travel cost between consecutive stations is $c_{\mathrm{train}}>0$. 
Therefore, for $i=1,\ldots,4$, we define the commuting cost of each category~$i$ employee by
\begin{align*}
    c_i(\BIx,\BIz):=\min_{1\le j\le 5,\,1\le j'\le 5}\big\{c_{\mathrm{walk}}\|\BIx-\BIu_j\|_1+c_{\mathrm{walk}}\|\BIz-\BIu_{j'}\|_1+c_{\mathrm{train}}|j-{j'}|\big\} \wedge c_{\mathrm{walk}}\|\BIx-\BIz\|_1.
\end{align*}
Moreover, we let the restocking cost of retail outlets be 
$c_5(\BIx,\BIz):=c_{\mathrm{restock}}{\|\BIx-\BIz\|_1}$ for $c_{\mathrm{restock}}>\nobreak0$.
The detailed setting of Experiment~1, including the MIP formulation 
of the global minimization problem in Line~\ref{alglin:cp-tf-global} of Algorithm~\ref{alg:cp-tf}, 
is presented in Section~\ref{sapx:experiments-biz},
where we also show that 
Setting~\ref{sett:Euclidean} is satisfied
with respect to $L_{c}\leftarrow \nobreak\sqrt{2}(c_{\mathrm{walk}}\vee c_{\mathrm{restock}})$.

Figure~\ref{fig:exp-biz-density} shows the shape of the type spaces $(\CX_i)_{i=1:N}$, 
and shows the probability density functions of the type measures $(\mu_i)_{i=1:N}$ as color plots.
When computing AMEs by Algorithm~\ref{alg:mt-tf},
we fix ${\epsilon_{\sspar}=2\times 10^{-4}}$, and
use 5 combinations of test functions $(\CG_i)_{i=0:N}$ 
constructed via 
Setting~\ref{sett:Euclidean} with shrinking mesh sizes $(\eta(\FS_i))_{i=0:N}$.
The resulting number of decision variables in \eqref{eqn:mt-tf-lsip}, 
i.e., $q:=\sum_{i=1}^{N}(1+m_i+m_0)$, is between 234 and 35409. 
The upper bounds $\hat{\alpha}_{\MT}^{\ssUB}$, 
$\tilde{\alpha}_{\MT}^{\ssUB}$
are computed via Monte Carlo integration with $10^4$ independent samples. 
Moreover, each Monte Carlo integration is repeated 100 times to examine the Monte Carlo error. 

\begin{figure}[t]
    \begin{center}
        \includegraphics[width=0.40\linewidth]{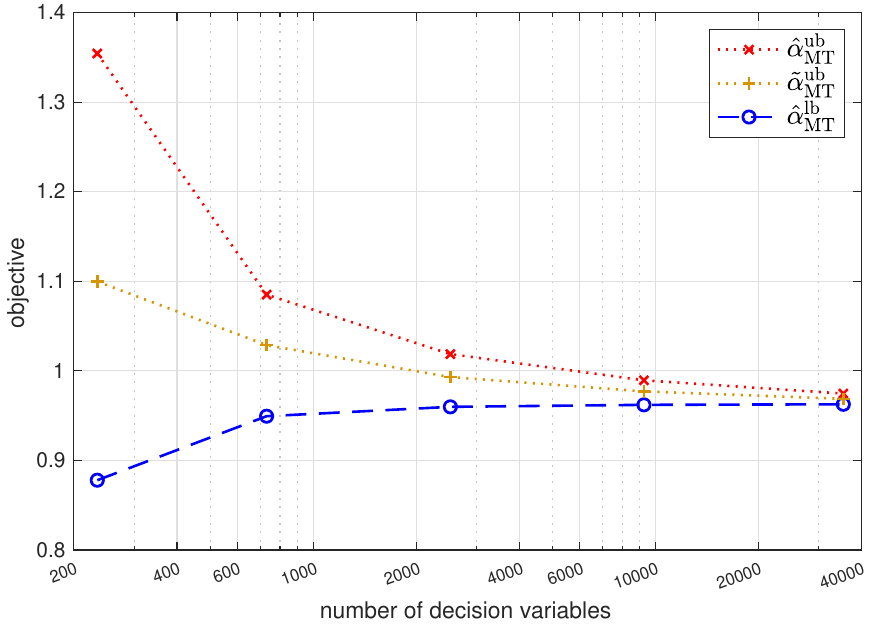}
        ~
        \includegraphics[width=0.40\linewidth]{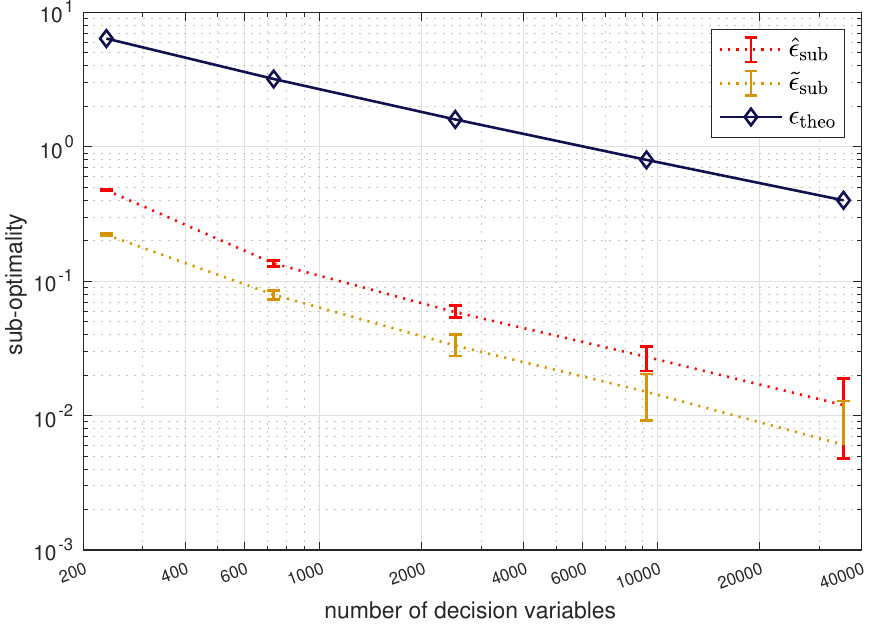} 
    \end{center}
    \caption{Experiment~1 -- The objective bounds and sub-optimality estimates computed by Algorithm~\ref{alg:mt-tf}.}\label{fig:exp-biz-bounds}
\end{figure}

The left panel of Figure~\ref{fig:exp-biz-bounds} shows the bounds 
$\hat{\alpha}_{\MT}^{\ssUB}$, 
$\tilde{\alpha}_{\MT}^{\ssUB}$, 
$\hat{\alpha}_{\MT}^{\ssLB}$ for the optimal value $\alpha_{\MT}^\star$ of \eqref{eqn:mt-primalopt} computed by Algorithm~\ref{alg:mt-tf}. 
The horizontal axis shows the number of decision variables $q$ in \eqref{eqn:mt-tf-lsip}.
It can be seen that $\tilde{\alpha}_{\MT}^{\ssUB}$ is considerably smaller than $\hat{\alpha}_{\MT}^{\ssUB}$. 
The differences between the upper bounds $\hat{\alpha}_{\MT}^{\ssUB}$, $\tilde{\alpha}_{\MT}^{\ssUB}$ and the lower bound $\hat{\alpha}_{\MT}^{\ssLB}$ are large when $q=234$, 
and decrease as $q$ increases. 
The right panel of Figure~\ref{fig:exp-biz-bounds} shows the sub-optimality estimates $\hat{\epsilon}_{\sssub}$ and $\tilde{\epsilon}_{\sssub}$ computed by Algorithm~\ref{alg:mt-tf} and their a priori upper bound $\epsilon_{\sstheo}$ in Theorem~\ref{thm:mt-tf}.
We show the sub-optimality estimates on the log-scale and use error bars to indicate the Monte Carlo errors in the computation of the upper bounds 
$\hat{\alpha}_{\MT}^{\ssUB}$, 
$\tilde{\alpha}_{\MT}^{\ssUB}$. 
Note that the lower branch of an error bar is not shown if it reaches below~0.
The results show that $\epsilon_{\sstheo}$ is around 13 to 33 times larger than $\hat{\epsilon}_{\sssub}$ and around 29 to 66 times larger than $\tilde{\epsilon}_{\sssub}$.
This demonstrates that Algorithm~\ref{alg:mt-tf} produces not only two AMEs
$(\hat{\varphi}_i,\hat{\gamma}_i)_{i=1:N},\hat{\nu}$ and 
$(\hat{\varphi}_i,\tilde{\gamma}_i)_{i=1:N},\tilde{\nu}$
but also sub-optimality estimates $\hat{\epsilon}_{\sssub}$, $\tilde{\epsilon}_{\sssub}$ of the AMEs that are much less conservative than suggested by purely theoretical analyses, as discussed in Section~\ref{ssec:efficient-algo}. 
Specifically, when $q=\nobreak35409$, $\epsilon_{\sstheo}=0.4002$
provides a highly conservative a priori sub-optimality estimate,
whereas the values $\hat{\epsilon}_{\sssub}=0.0120$ and $\tilde{\epsilon}_{\sssub}=0.0061$ 
show that $(\hat{\varphi}_i,\hat{\gamma}_i)_{i=1:N},\hat{\nu}$ and 
$(\hat{\varphi}_i,\tilde{\gamma}_i)_{i=1:N},\tilde{\nu}$ are close to true matching equilibria. 

\begin{figure}[t]
    \begin{center}
        \includegraphics[width=1.0\linewidth]{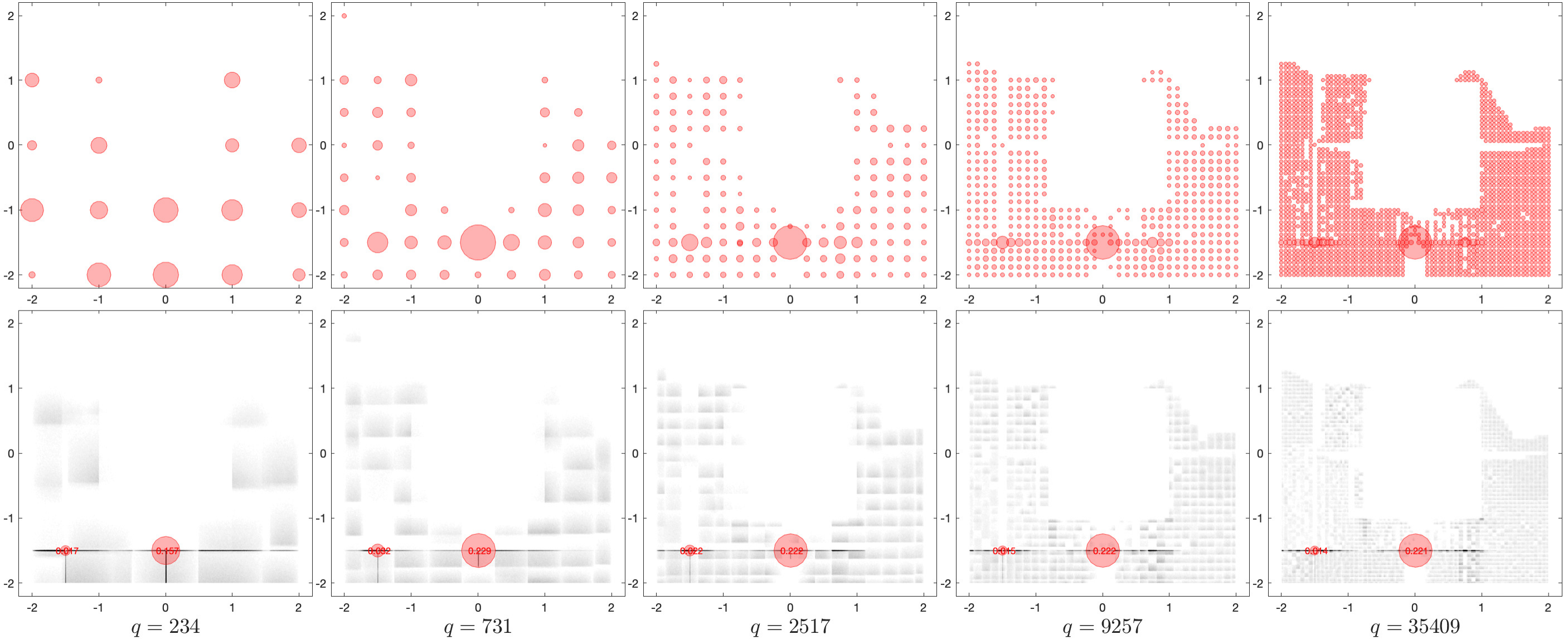}%
    \end{center}%
    \caption{Experiment~1 -- The probability measures $\hat{\nu}$ and $\tilde{\nu}$ computed by Algorithm~\ref{alg:mt-tf}.}\label{fig:exp-biz-hist}
\end{figure}

Finally, Figures~\ref{fig:exp-biz-hist}--\ref{fig:exp-biz-tf} illustrate $\hat{\nu}$, $\tilde{\nu}$, $(\tilde{\gamma}_i,\hat{\varphi}_i)_{i=1:N}$ computed by Algorithm~\ref{alg:mt-tf}.
The top row of Figure~\ref{fig:exp-biz-hist} shows $\hat{\nu}$ as bubble plots, where the locations 
and sizes
of the red circles represent 
the locations and probabilities of the atoms in $\hat{\nu}$.
The bottom row of Figure~\ref{fig:exp-biz-hist} shows $\tilde{\nu}$ as grayscale color plots superimposed with bubble plots. The reason for this choice is that for every combination of test functions, $\tilde{\nu}$ is a mixed probability measure containing a discrete component with two atoms and a non-discrete component. 
The probabilities of the two atoms in $\tilde{\nu}$ are shown as text in the bubble plot, whereas the histograms of the non-discrete part of $\tilde{\nu}$ are shown as grayscale color plots;
observe that some probability mass in $\tilde{\nu}$ is concentrated on a horizontal line. 
Figure~\ref{fig:exp-biz-coup} shows 1000 coupled samples from the computed couplings $(\tilde{\gamma}_i)_{i=1:N}$
when $q=35409$, 
where the black and red dots represent the locations of $\BIx$ and $\BIz$ in the samples, respectively, and the blue lines connecting the dots indicate the coupling between the locations. 
Samples from $(\hat{\gamma}_i)_{i=1:N}$ look similar to those from $(\tilde{\gamma}_i)_{i=1:N}$ and are thus omitted here. 
These coupled samples illustrate how the employees and the suppliers choose the outlets.
Figure~\ref{fig:exp-biz-tf} shows the computed transfer functions $(\hat{\varphi}_i)_{i=1:N}$ in color plots.
Recall that $\sum_{i=1}^{5}\hat{\varphi}_i=\nobreak0$ by the balance condition. 
The following insights can be drawn from Figures~\ref{fig:exp-biz-hist}--\ref{fig:exp-biz-tf}. 
\begin{itemize}[leftmargin=10pt]
    \item As the granularity of our approximation scheme 
    (e.g., as expressed by the number of decision variables $q$ in \eqref{eqn:mt-tf-lsip}) increases, 
    the structures of both the discrete quality measure $\hat{\nu}$ and the non-discrete quality measure $\tilde{\nu}$ increase in sophistication. 
    
    \item $(\hat{\varphi}_i,\hat{\gamma}_i)_{i=1:N},\hat{\nu}$ can be interpreted as an AME in which the retail outlets are only located at finitely many locations. 
    When $q=35409$, the two atoms in $\hat{\nu}$ with the largest probabilities are $(0, -1.5)^\TRANSP$ and $(-1.5, -1.5)^\TRANSP$ that correspond to two train stations in the south. 
    
    \item $(\hat{\varphi}_i,\tilde{\gamma}_i)_{i=1:N},\tilde{\nu}$ represents an AME in which retail outlets are distributed over uncountably many locations. 
    It shows some highly non-trivial features including the presence of two atoms at $(0, -1.5)^\TRANSP$ and $(-1.5, -1.5)^\TRANSP$, 
    positive probability concentrated on a horizontal line, 
    and an absolutely continuous component. 
    Moreover, Figure~\ref{fig:exp-biz-bounds} shows that $\tilde{\nu}$ has considerably smaller sub-optimality than~$\hat{\nu}$. 
    This demonstrates that, when the business is close to an equilibrium state, 
    around 22.1\% of outlets will be located at the train station at $(0, -1.5)^\TRANSP$, 
    around 1.4\% of outlets will be located at the train station at $(-1.5, -1.5)^\TRANSP$, 
    and the remaining outlets will not be concentrated at specific locations but dispersed into a continuum of locations, where a considerable portion will be dispersed along a horizontal line passing through the two aforementioned train stations. 
    
\begin{figure}[t]
    \begin{center}
        \includegraphics[width=1.0\linewidth]{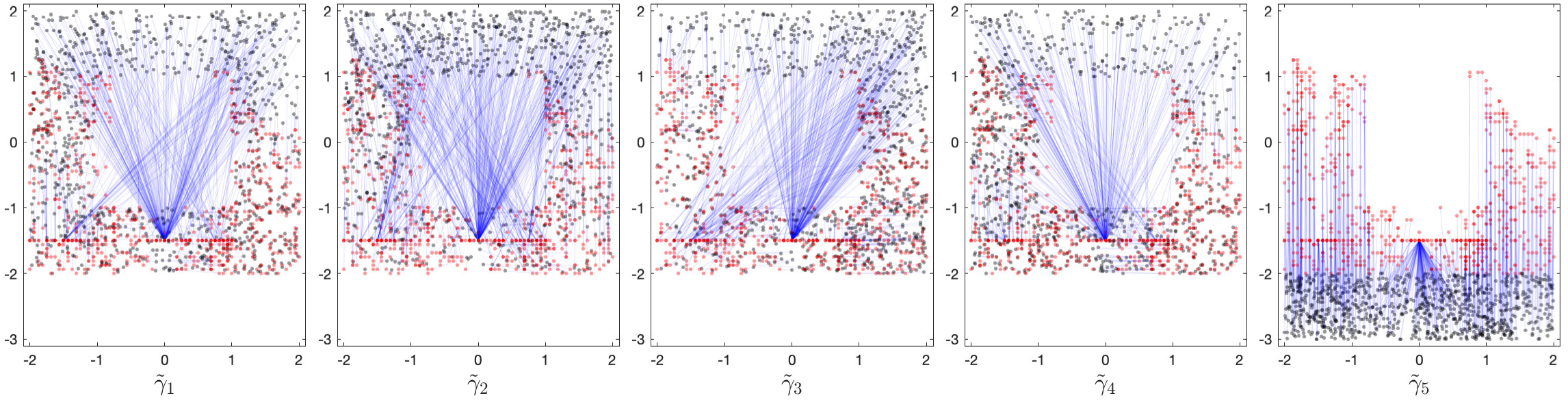}%
    \end{center}%
    \caption{Experiment~1 -- 1000 coupled samples from $(\tilde{\gamma}_i)_{i=1:N}$ computed by Algorithm~\ref{alg:mt-tf}.}\label{fig:exp-biz-coup}
\end{figure}

\begin{figure}[t]
    \begin{center}
        \includegraphics[width=1.0\linewidth]{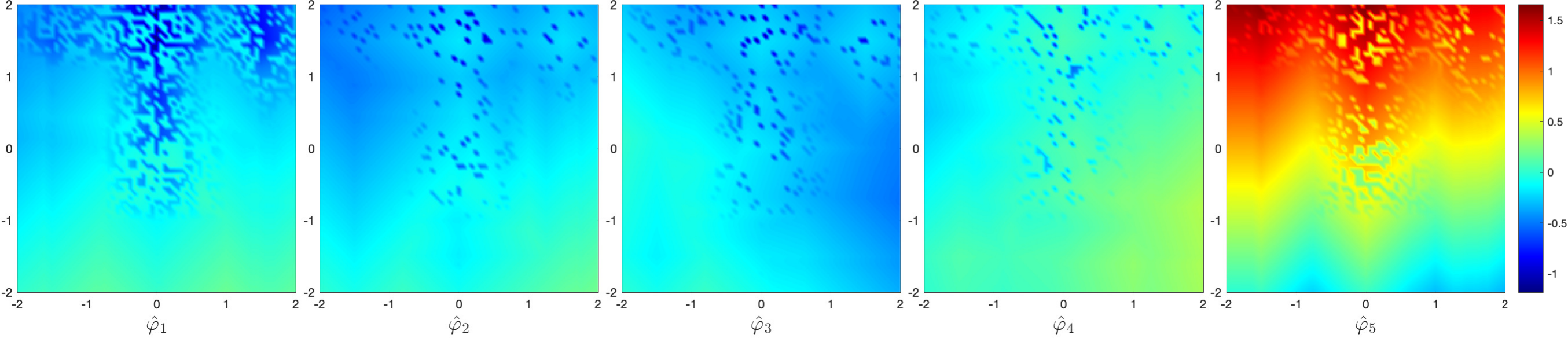}%
    \end{center}%
    \caption{Experiment~1 -- The approximately optimal transfer functions $(\hat{\varphi}_i)_{i=1:N}$ computed by Algorithm~\ref{alg:mt-tf}.}\label{fig:exp-biz-tf}
\end{figure}%
    
    \item The approximately optimal couplings $\tilde{\gamma}_1,\tilde{\gamma}_2,\tilde{\gamma}_3,\tilde{\gamma}_4$ demonstrate how the presence of the train in the city impacts the workplace choices of the employees. 
    Since the train connects the northern and southern parts of the city,
    many northern residents will commute to the southern part by train. 
    On the other hand, some residents in the western and eastern parts of the city that live far from train stations will work at nearby retail outlets. 
    
    \item The approximately optimal transfer functions 
    $\hat{\varphi}_1,\hat{\varphi}_2,\hat{\varphi}_3,\hat{\varphi}_4$ can be interpreted as salary maps indicating the amount of salary paid out to each category of employees at each location, 
    up to adding constants to $(\hat{\varphi}_i)_{i=1:5}$ that sum up to~0. 
    On the other hand, $\hat{\varphi}_5$ corresponds to the negative of the total salary paid out to the employees at each location. 
    Due to all the suppliers being located in the south, 
    retail outlets located in the northern part of the city need to spend more on restocking and thus will pay lower salary. 
    Besides this observation, it can be seen from Figure~\ref{fig:exp-biz-tf} that the approximately optimal transfer functions exhibit complex and heterogeneous patterns, and that 
    there are irregular patterns that lie outside $\support(\hat{\nu})$.
    The irregular patterns can be removed by applying $c_i$-convexification 
    described in Theorem~\ref{thm:mt-tf-equilibrium}.
\end{itemize}

Using the computed AMEs, the policymakers (e.g., city planners) can analyze the effects of transportation infrastructures, e.g., highways, railways, on shaping the geographical structure of retail outlets at equilibrium. 
For example, in the problem instance studied here, the experimental results above offer the following recommendations.
\begin{itemize}[leftmargin=10pt, beginpenalty=10000]
    \item The railway operators can use the computed approximately optimal couplings $(\tilde{\gamma}_{i})_{i=1:N}$ to gauge the demand for the train services at each train station. Therefore, the railway operations can be planned accordingly to cater to the pattern of commuters and avoid congestion. 
    
    \item A second railway line could be constructed to connect residents in the western and eastern parts of the city to the southern part due to the proximity of the southern part to the suppliers. 
    Its effects on the optimal social welfare 
    can be analyzed by
    Proposition~\ref{prop:perturbation} with respect to
    the modified cost functions 
    and the computed AMEs 
    $(\hat{\varphi}_i,\hat{\gamma}_i)_{i=1:N},\hat{\nu}$ and 
    $(\hat{\varphi}_i,\tilde{\gamma}_i)_{i=1:N},\tilde{\nu}$.
    
    \item Since an AME only reflects an approximately optimal structure of the business close to equilibrium, 
    the actual structure may differ from a computed equilibrium. 
    Therefore, policymakers can utilize the computed AME to implement additional policies to incentivize and facilitate the shift towards the optimal structure in order to improve the overall economic efficiency.
\end{itemize}

\subsection{Experiment~2: \mbox{2-Wasserstein} barycenter}\label{ssec:experiment-WB}

Numerical experiments about the computation of \mbox{2-Wasserstein} barycenters in existing studies mostly focus on the restrictive case where 
$(\mu_i)_{i=1:N}$ belong to a family of elliptical distributions,
as we have discussed in the introduction.
To the best of our knowledge, 
the only study in which the ground truth barycenter is accessible in the non-elliptical case is by \citet[\mbox{Lemmas~2--3}]{korotin2022wasserstein},
where they fix a probability measure~$\bar{\mu}$ and generate $N$ deformations $(\mu_i)_{i=1:N}$ such that $\bar{\mu}$ is the \mbox{2-Wasserstein} barycenter of $(\mu_i)_{i=1:N}$.
However, their deformations appear to be restrictive (see \citep[Figure~4]{korotin2022wasserstein}), possibly due to the high-dimensionality of their setting.
Moreover, 
their criterion to evaluate their algorithm only considers the mean and covariance discrepancies and neglects more nuanced differences.
Hence, due to the lack of accessible ground truth barycenters, there is a lack of quantitative analyses about the empirical approximation errors of \mbox{2-Wasserstein} barycenter algorithms for general non-parametric probability measures. 

\begin{figure}[t]
    \begin{center}
        \includegraphics[width=\linewidth]{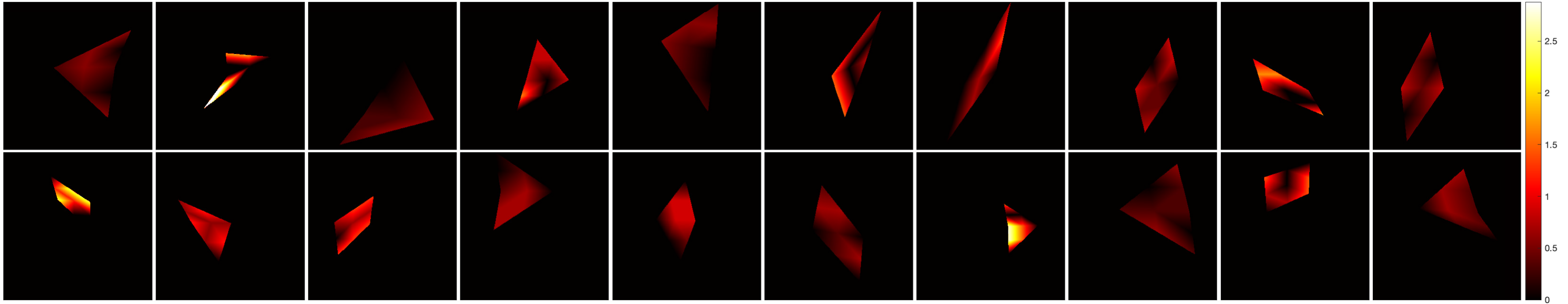}
    \end{center}
    \caption{Experiment~2 -- The probability density functions of $\mu_1,\ldots,\mu_{20}$.}\label{fig:exp-WB-density}
\end{figure}

In this experiment, we approximately compute the equally weighted \mbox{2-Wasserstein} barycenter of $N=20$ general non-parametric probability measures $(\mu_i)_{i=1:N}$ on $\R^2$, 
i.e., we solve the problem $\inf_{\nu\in\CP(\CX_0)}\big\{\frac{1}{N}\sum_{i=1}^N W_2(\mu_i,\nu)^2\big\}$ where $W_2(\mu_i,\nu):=\inf_{\gamma_i\in\Gamma(\mu_i,\nu)}\big\{\int_{\CX_i\times\CX_0}\|\BIx-\BIz\|_2^2\DIFFM{\gamma_i}{\DIFF\BIx,\DIFF\BIz}\big\}^{\frac{1}{2}}$.
Figure~\ref{fig:exp-WB-density} shows the probability density functions of $(\mu_i)_{i=1:N}$ as color plots.
Since the \mbox{2-Wasserstein} barycenter is guaranteed to be concentrated on the Minkowski sum $\sum_{i=1}^N\frac{1}{N}\conv(\CX_i)$ (see, e.g., \citep[Section~2.2]{carlier2015numerical}), we set $\CX_0:=\sum_{i=1}^N\frac{1}{N}\conv(\CX_i)$.
The detailed setting of Experiment~2,
including the detailed implementation of the global minimization problem in Line~\ref{alglin:cp-tf-global} of Algorithm~\ref{alg:cp-tf},
is presented in Section~\ref{sapx:experiments-WB},
where we also show that Setting~\ref{sett:Euclidean}
is satisfied with respect to 
$L_{c}\leftarrow \frac{4}{N}\max_{0\le i\le N}\!\big\{\!\max_{\BIx\in\CX_i}\!\big\{\|\BIx\|_2\big\}\!\big\}$.
In Algorithm~\ref{alg:mt-tf}, we fix $\epsilon_{\sspar}=2\times 10^{-4}$ and use 6 combinations of the test functions $(\CG_i)_{i=0:N}$ constructed via Setting~\ref{sett:Euclidean} with shrinking mesh sizes 
$(\eta(\FS_i))_{i=0:N}$.
The resulting number 
of decision variables in \eqref{eqn:mt-tf-lsip},
i.e., $q:=\sum_{i=1}^{N}(1+m_i+m_0)$, 
is between 2479 and 262895.
The upper bounds $\hat{\alpha}_{\MT}^{\ssUB}$ and $\tilde{\alpha}_{\MT}^{\ssUB}$ are computed via Monte Carlo integration using $10^7$ independent samples with 100 repetitions to examine the Monte Carlo error. 

\begin{figure}[t]
    \begin{center}
        \includegraphics[width=0.40\linewidth]{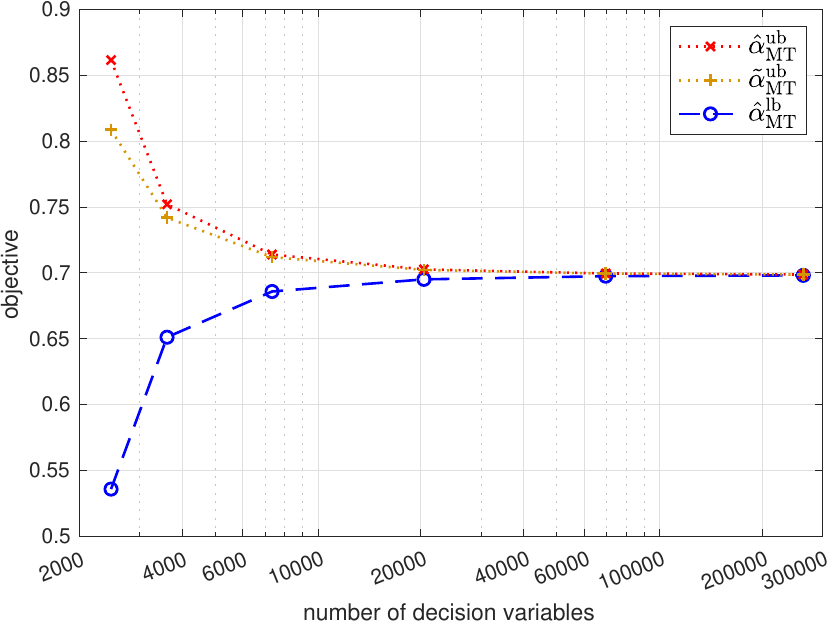}
        \includegraphics[width=0.40\linewidth]{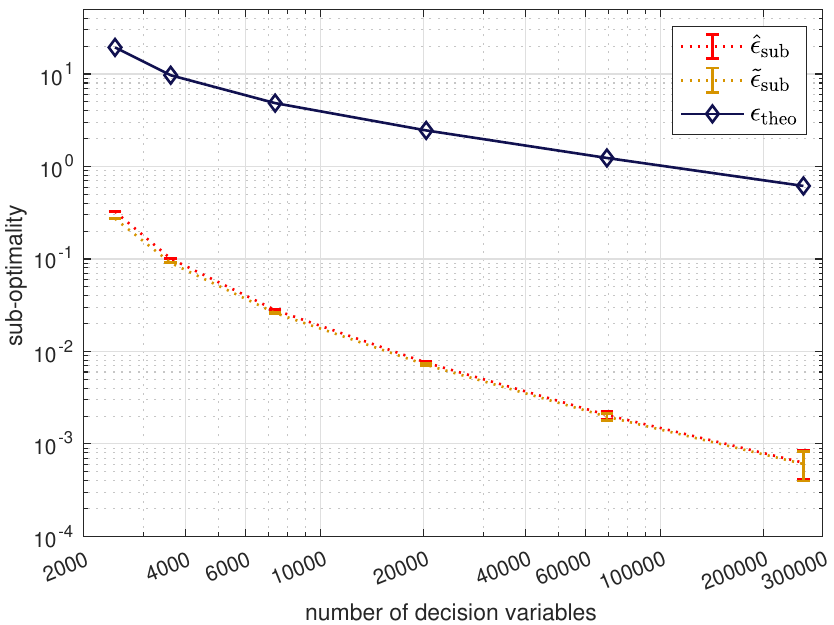}
    \end{center}
    \caption{Experiment~2 -- The objective bounds and sub-optimality estimates computed by Algorithm~\ref{alg:mt-tf}.}\label{fig:exp-WB-bounds}
\end{figure}

\begin{figure}[t]
    \begin{center}
        \includegraphics[width=1.0\linewidth]{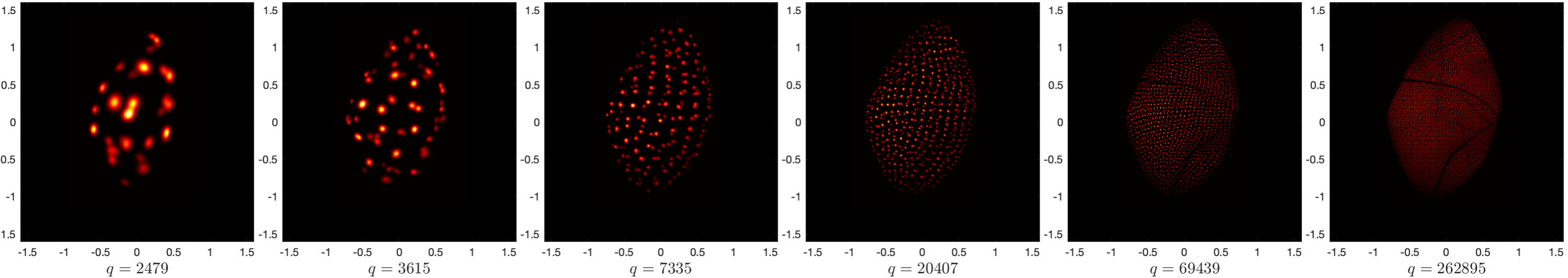}
    \end{center}
    \caption{Experiment~2  -- Histograms of the approximate \mbox{2-Wasserstein} barycenter $\tilde{\nu}$.}\label{fig:exp-WB-hist}
\end{figure}

The left panel of Figure~\ref{fig:exp-WB-bounds} shows the bounds 
$\hat{\alpha}_{\MT}^{\ssLB}$,
$\hat{\alpha}_{\MT}^{\ssUB}$, 
$\tilde{\alpha}_{\MT}^{\ssUB}$ 
for $\inf_{\nu\in\CP(\CX_0)}\big\{\frac{1}{N}\sum_{i=1}^NW_2(\mu_i,\nu)^2\big\}$,
i.e., the optimal value of the \mbox{2-Wasserstein} barycenter problem, computed by Algorithm~\ref{alg:mt-tf}.
As~$q$~increases, 
the differences between the upper bounds 
$\hat{\alpha}_{\MT}^{\ssUB}$, $\tilde{\alpha}_{\MT}^{\ssUB}$
and the lower bound $\hat{\alpha}_{\MT}^{\ssLB}$ decrease.
It can be seen that $\hat{\alpha}_{\MT}^{\ssUB}$ and $\tilde{\alpha}_{\MT}^{\ssUB}$ are comparable, with $\tilde{\alpha}_{\MT}^{\ssUB}$ being only slightly smaller than $\hat{\alpha}_{\MT}^{\ssUB}$, and the difference becomes negligible for large $q$. 
The right panel of Figure~\ref{fig:exp-WB-bounds} shows the sub-optimality estimates $\hat{\epsilon}_{\sssub}$, $\tilde{\epsilon}_{\sssub}$ computed by Algorithm~\ref{alg:mt-tf} 
as well as their a priori upper bound $\epsilon_{\sstheo}$ in 
Theorem~\ref{thm:mt-tf}.
These values are shown on the log-scale with error bars indicating the Monte Carlo errors when computing $\hat{\alpha}_{\MT}^{\ssUB}$ and $\tilde{\alpha}_{\MT}^{\ssUB}$. 
We observe that $\epsilon_{\sstheo}$ is around 60 to 1000 times larger than $\hat{\epsilon}_{\sssub}$, $\tilde{\epsilon}_{\sssub}$, 
which shows that the computed sub-optimality estimates 
are much less conservative than suggested by an a priori theoretical analysis. 
Figure~\ref{fig:exp-WB-hist} shows the histograms of the approximate \mbox{2-Wasserstein} barycenter $\tilde{\nu}$. 
One can observe that $\tilde{\nu}$ is a continuous probability measure that approximates the true \mbox{2-Wasserstein} barycenter by a finite number of ``blobs''.
When $q=262895$, we have $\tilde{\epsilon}_{\sssub}=6.0732\times 10^{-4}$, which indicates that the computed approximate \mbox{2-Wasserstein} barycenter $\tilde{\nu}$ is close to the true \mbox{2-Wasserstein} barycenter.

Moreover, we use the tightest lower bound $\hat{\alpha}_{\MT}^{\ssLB}=0.698157$ computed by Algorithm~\ref{alg:mt-tf} to compare our algorithm with five state-of-the-art \mbox{2-Wasserstein} barycenter algorithms by
\citet{staib2017parallel}, \citet{fan2020scalable}, \citet{korotin2021continuous}, \citet{korotin2022wasserstein}, and \citet{neufeld2022v7numerical}.
Specifically, \citet{staib2017parallel} approximately compute a \textit{fixed support} barycenter, where they fix a finite number of atoms in $\CX_0$ beforehand and optimize over the probabilities of the atoms via stochastic projected subgradient ascent. 
\mbox{\citet{fan2020scalable}} parametrize the problem by 
generative neural networks (GNNs) and input convex neural networks (ICNNs), and subsequently train these neural networks via a stochastic min-max-min scheme.
\citet{korotin2021continuous} adopt a parametrization with ICNNs but avoid the min-max-min structure in \citep{fan2020scalable} by enforcing the optimality conditions via regularization. 
\citet{korotin2022wasserstein} combine the conceptual fixed-point algorithm of \citet{alvarez2016fixed} with GNNs to develop a neural networks based fixed-point algorithm.
Lastly, \citet{neufeld2022v7numerical} tackle the multi-marginal optimal transport formulation via a cutting-plane algorithm.
To conduct a meaningful comparison of these algorithms and ours, we compute a discrete probability measure $\hat{\nu}\in\CP(\CX_0)$ via each algorithm, and evaluate its objective $\frac{1}{N}\sum_{i=1}^NW_2(\mu_i,\hat{\nu})^2$ and its sub-optimality $\big(\frac{1}{N}\sum_{i=1}^NW_2(\mu_i,\hat{\nu})^2\big)-\hat{\alpha}_{\MT}^{\ssLB}$.
\citet{staib2017parallel} compute a \textit{fixed support} \mbox{2-Wasserstein} barycenter which is already discrete.
Since the algorithms of \citet{fan2020scalable}, \citet{korotin2021continuous}, and \citet{korotin2022wasserstein} train generative models that can produce an arbitrary number of independent samples from the approximate barycenter, we take the empirical measure formed by 10000 independent samples produced from each of these models.
In the algorithm of \citet{neufeld2022v7numerical}, we take the discrete probability measure $\hat{\nu}$ from \citep[Proposition~5.3(ii)]{neufeld2022v7numerical}.
As for our algorithm, we take the discrete probability measure $\hat{\nu}$ computed by Algorithm~\ref{alg:mt-tf} for evaluation.
Section~\ref{sapx:experiments-WB} provides additional details about the settings of these algorithms for comparison.

\begin{table}[t]
\caption{Experiment~2 -- Comparison of \mbox{2-Wasserstein} barycenter algorithms.}\label{tab:experiment-WB-comparison}
\begin{footnotesize}
    \begin{center}
        \begin{tabular}{|r|r|r|}
        \hline
        Algorithm & ~~Objective & ~~Sub-optimality \\
        \hline
        \citet{staib2017parallel} &         $0.735321$     &    $3.7164\times 10^{-2}$ \\
        \citet{fan2020scalable}   &         $0.699104$     &    $9.4682\times 10^{-4}$ \\
        \citet{korotin2021continuous}  &    $0.705472$     &    $7.3150\times 10^{-3}$ \\
        \citet{korotin2022wasserstein}  &   $0.698804$     &    $6.4663\times 10^{-4}$ \\
        \citet{neufeld2022v7numerical} &    $0.698849$     &    $6.9213\times 10^{-4}$ \\
        \textbf{our algorithm}    &         $0.698537$     &    $3.7984\times 10^{-4}$ \\
        \hline 
        \end{tabular}
    \end{center}
\end{footnotesize}
\end{table}

Table~\ref{tab:experiment-WB-comparison} shows the objective and sub-optimality values of the discrete probability measures from the six algorithms.
Overall, our algorithm has produced the approximate \mbox{2-Wasserstein} barycenter with the lowest sub-optimality. 
The approximate \mbox{2-Wasserstein} barycenters produced by the algorithms of \citet{fan2020scalable}, \citet{korotin2022wasserstein}, and \citet{neufeld2022v7numerical} also have low sub-optimality.
On the other hand, despite that the algorithm of \citet{korotin2021continuous} utilizes ICNNs, its resulting sub-optimality is much worse.
One possible explanation is that the large number of probability measures ($N=20$) have rendered the regularization-based technique in \citep{korotin2021continuous} ineffective. 
The fixed support algorithm of \citet{staib2017parallel} performed the worst overall.
This demonstrates the need for so-called \textit{free support} \mbox{2-Wasserstein} barycenter algorithms.
The algorithm of \citet{neufeld2022v7numerical} has achieved low sub-optimality.
Besides, their algorithm also produced a lower bound $0.696569$ for the optimal value, which is smaller than the lower bound $\hat{\alpha}_{\MT}^{\ssLB}=0.698157$ computed by our algorithm.
However,
the computational efficiency of their algorithm is considerably worse than our algorithm.

\subsection{Experiment~3: one-dimensional type spaces}\label{ssec:experiment-1d}

In the third numerical experiment, we examine the scalability of Algorithm~\ref{alg:mt-tf} in terms of how its empirical running time changes with the number of agent populations (i.e., $N$) in the matching for teams problem. 
To that end, let us study the following matching for teams problem with one-dimensional type spaces, i.e., $\CX_1,\ldots,\CX_N\subset\R$ and a two-dimensional quality space, i.e., $\CX_0\subset\R^2$.

\begin{example}\label{exp:experiment-1d}%
The matching for teams problem in Experiment~3 is specified as follows.
\begin{itemize}[leftmargin=10pt]
\item For $i=1,\ldots,N$, the type space $\CX_i=[0,1]\subset\R$ 
represents a scalar-valued preference variable of category~$i$ agents that is between 0 and 1. 

\item Each good is characterized by two quality variables $z_1,z_2\in\R_{+}$ whose sum is at most~1,
i.e., 
$\CX_0=\big\{(z_1,z_2)^\TRANSP:z_1+z_2 \le 1,\; z_1,z_2\in\R_+\big\}\subset\R^2$.

\item For $i=1,\ldots,N$, 
the type measure $\mu_i\in\CP(\CX_i)$
representing the distribution of the preference variable within category~$i$
is absolutely continuous with respect to the Lebesgue measure on $\CX_i$, and its probability density function is a continuous piece-wise affine function on $\CX_i$. 

\item For $i=1,\ldots,N$, the cost function $c_i$ is given by 
$c_i(x,\BIz):=\frac{1}{N}\big(\big(|x-\langle\BIs_i,\BIz\rangle|\wedge \kappa_{i,2}\big)-\kappa_{i,1}\big)^+$ $\forall x\in\CX_i$, $\forall\BIz\in\CX_0$.
An agent evaluates a good with quality $\BIz\in\CX_0$
by comparing the assessment $\langle\BIs_i,\BIz\rangle$
with his/her preference $x\in\CX_i$,
where $\BIs_i\in\R^2$, $\|\BIs_i\|_2=\nobreak1$.
The cost is~0 if 
$\big|x-\langle\BIs_i,\BIz\rangle\big|$
is below the first threshold $\kappa_{i,1}>0$,
and the cost grows linearly when 
$\big|x-\langle\BIs_i,\BIz\rangle\big|$ is between $\kappa_{i,1}$
and the second threshold $\kappa_{i,2}>\kappa_{i,1}$,
until it reaches $\kappa_{i,2}-\kappa_{i,1}$ and then remains constant.
The factor $\frac{1}{N}$ in $c_i$ guarantees that 
Setting~\ref{sett:Euclidean} is satisfied with respect to 
$L_{c}\leftarrow \frac{1}{N}$ for any $N\ge 2$,
which yields an a priori sub-optimality bound $\epsilon_{\sstheo}$ in Theorem~\ref{thm:mt-tf} that is constant for all values of $N$. 
\end{itemize}%
\end{example}

We generate 10 problem instances of Example~\ref{exp:experiment-1d}, where $(\BIs_i,\kappa_{i,1},\kappa_{i,2})_{i=1:N}$ and the probability density function of $(\mu_i)_{i=1:N}$ are independently randomly generated. 
The test functions $(\CG_i)_{i=0:N}$ are constructed via Setting~\ref{sett:Euclidean} such that $|\CG_0|=560$, $|\CG_i|=49$ for $i=1,\ldots,N$.
The global minimization problem in Line~\ref{alglin:cp-tf-global} of Algorithm~\ref{alg:cp-tf} 
is formulated into a mixed-integer programming (MIP) problem and solved by the MIP solver of the Gurobi optimizer \citep{gurobi}; see Section~\ref{sapx:experiments-1d} 
for the concrete formulation.
We would like to remark that these global minimization problems 
can be solved in parallel. 
However, we solve them sequentially in our implementation except for the $N=1000$ case in order not to over-complicate the running time analysis. 

We fix $\epsilon_{\sspar}=5\times 10^{-5}$, 
run Algorithm~\ref{alg:mt-tf} for the 10 randomly generated problem instances,
and record the computed values of $\hat{\epsilon}_{\sssub}$ as well as the running time of Algorithm~\ref{alg:cp-tf} for $N=4,8,16,\allowbreak 32,64,128,1000$ agent populations. 
In particular, the largest instances where $N=1000$ are considerably larger than problem instances considered in existing studies, which typically have $N\le 20$, see, e.g., \citep{carlier2015numerical, anderes2016discrete, srivastava2018scalable, altschuler2021wasserstein, tanguy2024computing}.
$\hat{\epsilon}_{\sssub}$~is computed via Monte Carlo integration using $10^7$ independent samples. 
We only examine the values of $\hat{\epsilon}_{\sssub}$ here because $\tilde{\epsilon}_{\sssub}\le \hat{\epsilon}_{\sssub}$ and the computation of $\tilde{\epsilon}_{\sssub}$ requires solving a global minimization problem $\min_{\BIz\in\CX_0}\big\{\sum_{i=1}^Nc_i(x_i,\BIz)\big\}$ 
that is more computationally costly.
In addition, we only examine the running time of Line~\ref{alglin:mt-tf-cpalgo} in Algorithm~\ref{alg:mt-tf}, as the running time of the rest of Algorithm~\ref{alg:mt-tf} consists mostly of time spent computing $\hat{\alpha}_{\MT}^{\ssUB}$ via Monte Carlo integration in Line~\ref{alglin:mt-tf-primalcoup}, which is negligible.
Note that 
the running time in the $N=1000$ case has benefited from running 10 MIP solvers in parallel for the for-loop in Line~\ref{alglin:cp-tf-inner-forloop} of Algorithm~\ref{alg:cp-tf}.

\begin{table}[t]
\caption{Experiment~3 -- The computed sub-optimality estimate $\hat{\epsilon}_{\sssub}$ and running time ($\epsilon_{\sstheo}=0.1293$).\protect\footnotemark}\label{tab:experiment-1d-2d}
\begin{footnotesize}
\begin{center}
    \begin{tabular}{|r|rr|rr|rr|rr|}%
    \hline
    \multicolumn{1}{|c|}{} & 
    \multicolumn{1}{c}{Avg.\ $\hat{\epsilon}_{\sssub}$} & 
    \multicolumn{1}{c|}{Max.\ $\hat{\epsilon}_{\sssub}$} & 
    \multicolumn{1}{c}{Avg.\ LP} & 
    \multicolumn{1}{c|}{Max.\ LP} & 
    \multicolumn{1}{c}{Avg.\ MIP} & 
    \multicolumn{1}{c|}{Max.\ MIP} & 
    \multicolumn{1}{c}{Avg.\ total} & 
    \multicolumn{1}{c|}{Max.\ total}  \\
    \multicolumn{1}{|r|}{$N$} & 
    \multicolumn{1}{c}{[$\times10^{-4}$]} & 
    \multicolumn{1}{c|}{[$\times10^{-4}$]} &  
    \multicolumn{1}{c}{time \![$\mathrm{sec.}/N$]} & 
    \multicolumn{1}{c|}{\!\!\!time \![$\mathrm{sec.}/N$]} &  
    \multicolumn{1}{c}{time [$\mathrm{sec.}/N$]} & 
    \multicolumn{1}{c|}{\!\!\!time \![$\mathrm{sec.}/N$]} &  
    \multicolumn{1}{c}{time \![$\mathrm{sec.}/N$]} & 
    \multicolumn{1}{c|}{\!\!\!time \![$\mathrm{sec.}/N$]} \\
    \hline
    4 &    9.566 &   21.241 &   5.55 &   8.14 & 107.00 & 139.14 & 112.62 & 145.88 \\
    8 &    8.859 &   13.107 &   6.17 &   7.95 & 121.21 & 150.60 & 127.43 & 157.88 \\
    16 &    6.798 &    8.892 &   6.13 &   7.47 & 131.73 & 169.47 & 137.90 & 176.99 \\
    32 &    5.191 &    7.271 &   6.28 &   8.64 & 138.07 & 178.90 & 144.38 & 187.58 \\
    64 &    4.930 &    6.625 &   7.25 &   8.36 & 147.14 & 160.70 & 154.43 & 168.64 \\
    128 &    4.361 &    5.001 &  11.10 &  14.46 & 173.69 & 217.46 & 184.83 & 231.09 \\
    \hline
    1000 &    4.167 &    4.672 &  42.78 &  49.58 &  *53.93 &  *64.64 &  *96.75 & *114.26 \\
    \hline 
    \end{tabular} 
\end{center}
\end{footnotesize}
\end{table}

\footnotetext{In the $N=1000$ case, we solve 10 global minimization problems in Line~\ref{alglin:cp-tf-global} of Algorithm~\ref{alg:cp-tf} in parallel, and we thus use * to indicate running time in the table that has benefited from parallel computation.}

Columns~2 and 3 of Table~\ref{tab:experiment-1d-2d} show the average and maximum values of the computed sub-optimality estimate~$\hat{\epsilon}_{\sssub}$.
Observe that the values of $\hat{\epsilon}_{\sssub}$ 
are about two orders of magnitude smaller than the a priori upper bound $\epsilon_{\sstheo}=0.1293$. 
The rest of Table~\ref{tab:experiment-1d-2d} shows the average and maximum running time of the LP solver (Line~\ref{alglin:cp-tf-lp}), the MIP solver (Line~\ref{alglin:cp-tf-global}), and the entire Algorithm~\ref{alg:cp-tf},
which are divided by $N$ for better interpretability.
We have observed that the total running time of Algorithm~\ref{alg:cp-tf} is dominated by the MIP solver,
and that the running time seems to be increasing at a faster than linear but slower than quadratic rate with respect to~$N$.
Lastly, we highlight that in the $N=1000$ case we are able to drastically improve the computational efficiency by using 10 parallel processes when solving the MIP problems.

\section{Conclusion}\label{sec:conclusion}

We have introduced the notion of $\epsilon$-matching equilibrium 
as a feasible approximation of a matching equilibrium,
which corresponds to a plausible state of the matching for teams game 
where agents are approximately rational and the aggregated sub-optimality of their decisions is at most~$\epsilon$.
We have shown that each $\epsilon$-matching equilibrium provides lower and upper bounds for the optimal social welfare, and that computing and examining an $\epsilon$-matching equilibrium can aid policymakers and regulators 
in improving the overall economic efficiency.
Under a flexible Euclidean setting,
we have derived an explicit computational complexity bound for the computation of $\epsilon$-matching equilibria,
and we have developed an algorithm 
for computing $\epsilon$-matching equilibria
which is highly efficient and scalable to problem instances with a large number of agent populations.
In our numerical experiments,
we have demonstrated the practical scalability of our algorithm,
and that is capable of computing nearly optimal approximate matching equilibria 
that lead to managerial insights.
In particular, our algorithm outperformed state-of-the-art algorithms for the Wasserstein barycenter problem in our comparative study.


%
%

%

\ACKNOWLEDGMENT{
    AN and QX gratefully acknowledge the financial support by the MOE AcRF Tier~2 Grant \textit{MOE-T2EP20222-0013}.
}

\setlength{\bibsep}{1pt}

\bibliographystyle{informs2014} 


\ECSwitch
%
%
\ECHead{Appendices}
In this e-companion, also known as the online appendices, we present
technical theoretical results that are omitted from the main paper,
additional remarks and discussions about the numerical experiments,
as well as the proofs of all the theoretical results.
The e-companion is organized as follows.
Section~\ref{apx:couplings} presents 
the detailed procedure to construct the random variables in Theorem~\ref{thm:mt-tf-approx}
on a sufficiently rich probability space.
Section~\ref{apx:complexity} 
contains the explicit procedures referenced in Algorithm~\ref{alg:complexity-overall}
as well as their computational complexity analyses.
In Section~\ref{apx:efficient-algo},
we present additional details and remarks about 
Algorithm~\ref{alg:cp-tf} and Algorithm~\ref{alg:mt-tf}.
Section~\ref{apx:mmot} discusses the connection between the matching for teams problem and the multi-marginal optimal transport problem.
Section~\ref{apx:experiments} 
explains the specific settings of our numerical experiments.
In particular, it includes the detailed procedures that we carried out when comparing Algorithm~\ref{alg:mt-tf} with five state-of-the-art \mbox{2-Wasserstein} barycenter algorithms.
Finally,
Section~\ref{apx:proof} contains the proofs of the theoretical results in the paper and the proofs of the auxiliary theoretical results in this e-companion.

\renewcommand{\theequation}{\thesection.\arabic{equation}}
\numberwithin{equation}{section}

\section{Details about the construction of the random variables in Theorem~\ref{thm:mt-tf-approx}}\label{apx:couplings}

In this section, let us present the detailed construction of the random variables 
$Z$, $(Z_i,X_i,\bar{X}_i)_{i=1:N}$, and $\bar{Z}$ that appeared in Theorem~\ref{thm:mt-tf-approx} on a sufficiently rich probability space
$(\Omega,\CF,\PROB)$.
To that end, 
let $(\hat{\theta}_i,\hat{\mu}_i,\hat{\nu}_i)_{i=1:N}$ and $\big(\overline{\epsilon}_i(\varsigma)\big)_{i=0:N}$ 
be given by the statement of Theorem~\ref{thm:mt-tf-approx}.
Moreover, let 
$\Omega:=[0,1]$, $\CF:=\CB([0,1])$, 
and let $\PROB$ be the Lebesgue measure on $[0,1]$.
It follows from classical results such as \citep[Lemma~1]{ECskorohod1976on}
that there exist $(3N+1)$ independent random variables $U_0,U_1,\ldots,U_{3N}:\Omega\to[0,1]$ where each $U_i$ is uniformly distributed on $[0,1]$,
i.e., $\PROB\circ U_i^{-1}=\PROB$ $\forall 0\le\nobreak i\le\nobreak 3N$.

Subsequently, we 
let $Y(\omega):=0$ $\forall \omega\in\Omega$,
and
construct the random variable $Z:\Omega\to\CX_0$ with law $\hat{\nu}:=\hat{\nu}_{1}$
by \citep[Lemma~2.11]{ECdudley1983invariance} with respect to 
$S\leftarrow\{0\}$,
$T\leftarrow\CX_0$,
$Q\leftarrow \delta_{0}\otimes \hat{\nu}$ (i.e., the product measure of $\delta_{0}$ and $\hat{\nu}$),
$\mu\leftarrow \delta_{0}$,
$X\leftarrow Y$,
$U\leftarrow U_0$,
$R\leftarrow [0,1]$,
$Y\leftarrow Z$.

Next, for $i=1,\ldots,N$, 
let $\eta_i\in\Gamma(\hat{\nu},\hat{\nu}_i)$
satisfy 
$\int_{\CX_0\times\CX_0}d_{\CX_0}(z,z')\DIFFM{\eta_i}{\DIFF z,\DIFF z'}\le \overline{\epsilon}_0(\varsigma)$,
and let $\rho_i\in\nobreak\Gamma(\hat{\mu}_i,\mu_i)$
satisfy 
$\int_{\CX_i\times\CX_i}d_{\CX_i}(x,x')\DIFFM{\rho_i}{\DIFF x,\DIFF x'}\le\overline{\epsilon}_i(\varsigma)$.
The existence of $\eta_i$ and $\rho_i$ is guaranteed by the assumption that 
$\hat{\nu}_i\abovebelowset{\CG_0}{\varsigma}{\sim}\hat{\nu}_{1}$,
$\hat{\mu}_i\abovebelowset{\CG_i}{\varsigma}{\sim}\mu_i$,
as well as 
$\overline{\epsilon}_0(\varsigma)\ge \sup\big\{W_1(\nu,\nu'):\nu,\nu'\in\CP(\CX_0),\;\nu\abovebelowset{\CG_0}{\varsigma}{\sim}\nu'\big\}$,
$\overline{\epsilon}_i(\varsigma)\ge \sup\big\{W_1(\nu,\nu'):\nu,\nu'\in\CP(\CX_i),\;\nu\abovebelowset{\CG_i}{\varsigma}{\sim}\nu'\big\}$.
Moreover, let $\tilde{\theta}_i\in\CP(\CX_0\times\CX_i)$ be defined by $\tilde{\theta}_i(A\times\nobreak B):=\hat{\theta}_i(B\times\nobreak A)$ 
$\forall A\in\CB(\CX_0)$, $\forall B\in\CB(\CX_i)$.
Then, we construct the random variable $Z_i:\Omega\to\CX_0$
by \citep[Lemma~2.11]{ECdudley1983invariance} with respect to 
$S\leftarrow\CX_0$,
$T\leftarrow\CX_0$,
$Q\leftarrow \eta_i$,
$\mu\leftarrow \hat{\nu}$,
$X\leftarrow Z$,
$U\leftarrow U_{i}$,
$R\leftarrow\nobreak [0,1]$,
$Y\leftarrow Z_i$,
construct the random variable
$X_i:\Omega\to\CX_i$ 
by \citep[Lemma~2.11]{ECdudley1983invariance} with respect to 
$S\leftarrow\CX_0$,
$T\leftarrow\CX_i$,
$Q\leftarrow \tilde{\theta}_i$,
$\mu\leftarrow \hat{\nu}_i$,
$X\leftarrow Z_i$,
$U\leftarrow U_{N+i}$,
$R\leftarrow\nobreak [0,1]$,
$Y\leftarrow X_i$,
and construct the random variable
$\bar{X}_i:\Omega\to\CX_i$
by \citep[Lemma~2.11]{ECdudley1983invariance} with respect to 
$S\leftarrow\CX_i$,
$T\leftarrow\nobreak\CX_i$,
$Q\leftarrow \rho_i$,
$\mu\leftarrow \hat{\mu}_i$,
$X\leftarrow X_i$,
$U\leftarrow U_{2N+i}$,
$R\leftarrow [0,1]$,
$Y\leftarrow \bar{X}_i$.
One checks that 
$\EXP\big[d_{\CX_0}(Z,Z_i)\big]=\int_{\CX_0\times\CX_0}d_{\CX_0}(z,z')\DIFFM{\eta_i}{\DIFF z,\DIFF z'}\le \overline{\epsilon}_0(\varsigma)$,
$\EXP\big[d_{\CX_i}(X_i,\bar{X}_i)\big]=\int_{\CX_i\times\CX_i}d_{\CX_i}(x,x')\DIFFM{\rho_i}{\DIFF x,\DIFF x'}\le\overline{\epsilon}_i(\varsigma)$,
the law of $\bar{X}_i$ is $\mu_i$,
and 
the law of $(X_i,Z_i):\Omega\to\CX_i\times\CX_0$ is $\hat{\theta}_i$.
Consequently, $(Z_i,X_i,\bar{X}_i)_{i=1:N}$ satisfy the properties in the statement of Theorem~\ref{thm:mt-tf-approx}.

Lastly,
let $\bar{z}:\CX_1\times\cdots\times\CX_N\to\CX_0$ be a Borel measurable function satisfying
\begin{align*}
    \sum_{i=1}^{N}c_i\big(x_i,\bar{z}(x_1,\ldots,x_N)\big)=\min_{z\in\CX_0}\Bigg\{\sum_{i=1}^{N}c_i(x_i,z)\Bigg\}
    \qquad \forall (x_1,\ldots,x_N)\in\CX_1\times\cdots\times\CX_N,
\end{align*}
which exists due to \citep[Proposition~7.33]{ECbertsekas1978stochastic}.
We construct the random variable $\bar{Z}:\Omega\to\CX_0$ by 
$\bar{Z}:=\bar{z}(\bar{X}_1,\ldots,\bar{X}_N)$.
This completes the construction of the random variables 
$Z$, $(Z_i,X_i,\bar{X}_i)_{i=1:N}$, and $\bar{Z}$ in Theorem~\ref{thm:mt-tf-approx}.

\section{Details of the procedures in Algorithm~\ref{alg:complexity-overall}}\label{apx:complexity}

In this section,
we present the detailed procedures in Algorithm~\ref{alg:complexity-overall}
as well as their explicit computational complexities 
under the arithmetic model.
This section is organized as follows.
We first present in Section~\ref{sapx:complexity-notations}
the notations that are used throughout this section.
Next,
since the arithmetic model assumes that all inputs, outputs, and intermediate values of the procedures are stored as exact real numbers,
we present in Section~\ref{sapx:complexity-encoding} how the various non-real valued mathematical objects throughout the procedures are encoded by collections of real numbers.
Subsequently, 
we present in Sections~\ref{sapx:complexity-init}--\ref{sapx:complexity-AME}
the detailed procedures for computing an $\epsilon$-matching equilibrium,
and we derive explicit bounds for their computational complexities.
Moreover, 
Section~\ref{sapx:complexity-AME} also contains the procedures for accessing the computed $\epsilon$-matching equilibrium.
Lastly, 
under the particular class of $(\mu_i)_{i=1:N}$ in Setting~\ref{sett:Euclidean}\ref{settc:Euclidean-particular},
Section~\ref{sapx:complexity-particular-case}
presents the procedures for the two operations in 
Theorem~\ref{thm:complexity-explicit} related to $\mu_i$,
and
derives explicit upper bounds for $T_{\mathtt{Int}}$ and~$T_{\mathtt{Samp}}$.

To facilitate the derivation of the computational complexities,
we present the number of arithmetic operations within each line of pseudocode in big-$O$ notations as a comment at the right-end of the line.
Specifically, the comment accompanying each non-loop line of pseudocode 
counts the number of arithmetic operations within each execution of that line,
whereas
the comment accompanying each loop (i.e., a line of pseudocode with the keyword \textbf{for} or \textbf{while})
counts the total number of arithmetic operations in each execution of the entire loop.
To avoid clutter, we omit the comments in lines involving $O(1)$ arithmetic operations.
We remark that the constant terms omitted by the big-$O$ notations 
in this section 
are absolute constants and do not depend on any input or any parameter
of the procedures.

\IncMargin{3pt}
\SetProcNameFnt{\tt}
\SetCommentSty{}
\SetKwHangingKw{KwInAndOut}{}

\subsection{Notations}\label{sapx:complexity-notations}
Throughout this section, 
we work under Setting~\ref{sett:Euclidean}, and
we continue to use the notations introduced in Section~\ref{ssec:setting-Euclidean},
including 
$\dim(\,\cdot\,)$,
$V(\,\cdot\,)$,
$\FF(\,\cdot\,)$,
$\FM(\,\cdot\,)$,
$\eta(\,\cdot\,)$.
In addition to the notations in Section~\ref{ssec:setting-Euclidean},
let $O(n^\omega)$ denote the number of arithmetic operations incurred when multiplying two $n$-by-$n$ matrices.
For any \mbox{$n$-simplex} $C$ in $\R^d$ where $n\in\{1,\ldots,d\}$,
let $\Lebesgue_{\aff(C)}$ denote the Lebesgue measure on the \mbox{$n$-dimensional} affine hull of $C$,
and let 
$\volume_n(C)$ denote the \textit{$n$-dimensional volume} of $C$,
that is, 
$\volume_n(C):=\int_{\aff(C)}\INDI_{C}\DIFFX{\Lebesgue_{\aff(C)}}$.
Recall the Gram determinant formula for the \mbox{$n$-dimensional} volume of a simplex: 
$\volume_n(C)=\frac{1}{n!}\det(\BV_{C}^\TRANSP\BV_{C})^{\frac{1}{2}}$,
where 
\begin{align*}
    \BV_{C}:=\left(\begin{tabular}{cccc}
        $|$ & $|$ & & $|$ \\
        $(\BIv_{1}-\BIv_{0})$ & $(\BIv_{2}-\BIv_{0})$ & $\cdots$ & $(\BIv_{n}-\BIv_{0})$ \\
        $|$ & $|$ & & $|$
    \end{tabular}\right)\in\R^{d\times n};
\end{align*}
see, e.g., \citep{ECklain2004intuitive}.
Moreover, let $\Lebesgue_{C}$ denote the \mbox{$n$-dimensional} 
Lebesgue measure on $C$,
which is the restriction of $\Lebesgue_{\aff(C)}$ to $\CB(C)$,
i.e.,
it holds that 
$\Lebesgue_C(E)=\Lebesgue_{\aff(C)}(E)$ $\forall E\in\CB(C)$.
Furthermore, for any $n\in\N$,
let $\Lebesgue_{\R^n}$ denote the $n$-dimensional Lebesgue measure on $\R^n$.
In addition, in the case where $C$ is a $0$-simplex in $\R^d$, i.e., $C$ contains a single point $\BIv\in\R^d$,
we abuse the notations and denote 
$\Lebesgue_{\aff(C)}:=\delta_{\BIv}$,
$\Lebesgue_{C}:=\delta_{\BIv}$,
and $\volume_0(C):=1$.

\subsection{Encoding mathematical objects by real numbers}\label{sapx:complexity-encoding}

The following list presents how various mathematical objects in the inputs and outputs of the procedures in this section are encoded by real numbers.
\begin{enumerate}[label=\normalfont{(E\arabic*)},beginpenalty=10000,leftmargin=26pt]
    \item\label{enc:simplex}
    Every $n$-simplex $K$ in $\R^d$ ($n\in\{0,1,\ldots,d\}$) is encoded by 
    $\big(d,n,(\BIv_l)_{l=0:n}\big)$, 
    where $d,n\in\N$ and $(\BIv_l)_{l=0:n}\subset\R^d$ is an enumeration of the $n+1$ extreme points of $K$ in $V(K)$.

    \item\label{enc:complex}
    Every simplicial complex $\FK$ in $\R^d$ is encoded by 
    $\big(d,J,(K_j)_{j=1:J}\big)$,
    where $d\in\N$, $J:=|\FM(\FK)|\in\N$, and 
    $(K_j)_{j=1:J}$ is an enumeration of $\FM(\FK)$ where each $K_j$ is encoded via \ref{enc:simplex}.

    \item\label{enc:spaces}
    For $i=0,1,\ldots,N$, 
    the encoding of the quality/type space $\CX_i$ 
    corresponds to the ``triangulation encoding'' of 
    the simplicial complex $\FK_i$ characterizing $\CX_i$ in Setting~\ref{sett:Euclidean}.
    Specifically, $\CX_i$ is encoded by 
    $\big(d_i,q_i,\allowbreak(\BIu_{i,l})_{l=0:q_i},\allowbreak J_i,\allowbreak \big(\rho_{i,j},(\kappa_{i,j,\iota})_{\iota=0:\rho_{i,j}}\big)_{j=1:J_i}\big)$,
    where $d_i\in\N$ is the dimension of $\CX_i$,
    $q_i:=|V(\FK_i)|-\nobreak1\in\nobreak\N_0$,
    $(\BIu_{i,l})_{l=0:q_i}\subset\R^{d_i}$ is an enumeration of $V(\FK_i)$,
    $J_i:=|\FM(\FK_i)|\in\N$,
    $(K_{i,j})_{j=1:J_i}$ is an enumeration of $\FM(\FK_i)$,
    $\rho_{i,j}:=\dim(K_{i,j})\in\nobreak\N_0$ $\forall j\in\{1,\ldots,J_i\}$,
    and 
    $(\kappa_{i,j,\iota})_{\iota=0:\rho_{i,j}}\subset\{0,1,\ldots,q_i\}$ 
    satisfy 
    $(\BIu_{i,\kappa_{i,j,\iota}})_{\iota=0:\rho_{i,j}}=V(K_{i,j})$ $\forall j\in\{1,\ldots,J_i\}$.

    \item\label{enc:typemeasures}
    Under Setting~\ref{sett:Euclidean},
    the type measures $(\mu_i)_{i=1:N}$ 
    are not assumed to have specific parametric forms, and are only accessed through their associated oracle procedures $\mathtt{Int}$ and $\mathtt{Samp}$.
    Specifically, 
    for any $i\in\{1,\ldots,N\}$,
    any $\hat{j}\in\{1,\ldots,J_i\}$,
    and any $n$-simplex $C$ such that 
    $n\in\nobreak\{0,1,\ldots,d_i\}$, $V(C)=\{\BIv_l\}_{l=0:n}$,
    $C\subseteq K_{i,\hat{j}}$ (see~\ref{enc:spaces}),
    $\mathtt{Int}(\CX_i,\mu_i,C,\hat{j})$
    incurs $T_{\mathtt{Int}}$ arithmetic operations and
    returns 
    $(h_l)_{l=0:n}$ where
    $h_l=\nobreak\int_{\CX_i}\lambda_{\FF(C),\BIv_l}\INDI_{\relint(C)}\DIFFX{\mu_i}\in\R_+$ $\forall l\in\{0,1,\ldots,n\}$.
    Moreover, given an infinite sequence of independent realizations $(U_t)_{t\in\N}$ of $\mathrm{Uniform}[0,1]$ random variable
    and assuming that $\mu_i\big(\relint(C)\big)>0$,
    $\mathtt{Samp}\big(\CX_i,\mu_i,C,\hat{j},(U_t)_{t\in\N}\big)$
    incurs $T_{\mathtt{Samp}}$ arithmetic operations and
    returns $\BIx\in\CX_i$ and $\overline{t}\in\N_0$,
    where $\BIx$ is a random sample from $\tilde{\mu}_{i,C}\in\CP(\CX_i)$ defined by $\frac{\DIFF\tilde{\mu}_{i,C}}{\DIFF\mu_i}:=\frac{\INDI_{\relint(C)}}{\mu_i(\relint(C))}$,
    and $\overline{t}\in\N_0$ is the number of realizations in $(U_t)_{t\in\N}$ consumed in $\mathtt{Samp}$,
    that is, $\BIx$ is a Borel measurable function of $(U_t)_{t=1:\overline{t}}$.
    We do not specify the encoding of $(\mu_i)_{i=1:N}$
    in the absence of additional assumptions.
    The particular class of $(\mu_i)_{i=1:N}$ in Setting~\ref{sett:Euclidean}\ref{settc:Euclidean-particular}
    will be encoded via \ref{enc:typemeasures-particular} introduced later in Section~\ref{sapx:complexity-particular-case}.

    \item\label{enc:costfuncs}
    Under Setting~\ref{sett:Euclidean},
    the cost functions $(c_i)_{i=1:N}$ are not assumed to have specific parametric forms,
    and we thus do not specify their encoding.
    Our only assumptions about $(c_i)_{i=1:N}$ are 
    that they are $L_c$-Lipschitz continuous for $L_{c}>0$,
    and for $i=1,\ldots,N$,
    the evaluation of $c_i(\BIx,\BIz)$ at any $(\BIx,\BIz)\in \CX_i\times\CX_0$ incurs $T_{\mathtt{Eval}}$ arithmetic operations.

    \item\label{enc:testfuncs}
    For $i=0,1,\ldots,N$,
    the encoding of the test functions $\CG_i$ on $\CX_i$ 
    corresponds to the ``triangulation encoding'' of the subdivision $\FS_i$ of $\FK_i$
    in Setting~\ref{sett:Euclidean}, similar to the encoding of $\CX_i$ in \ref{enc:spaces}.
    Specifically, $\CG_i$ is encoded by 
    $\big(d_i,m_i,(\BIv_{i,l})_{l=0:m_i},\widetilde{J}_i,\big(\varrho_{i,\tilde{j}},\allowbreak (\varkappa_{i,\tilde{j},\iota})_{\iota=0:\varrho_{i,\tilde{j}}},\allowbreak\varpi_{i,\tilde{j}}\big)_{\tilde{j}=1:\widetilde{J}_i},\overline{\eta}_i\big)$,
    where $d_i\in\N$ is the dimension of $\CX_i$,
    $m_i:=|V(\FS_i)|-\nobreak1\in\nobreak\N_0$,
    $(\BIv_{i,l})_{l=0:m_i}\subset\R^{d_i}$ is an enumeration of $V(\FS_i)$,
    $\widetilde{J}_i:=|\FM(\FS_i)|\in\nobreak\N$,
    $J_i:=|\FM(\FK_i)|\in\nobreak\N$,
    $(S_{i,\tilde{j}})_{\tilde{j}=1:\widetilde{J}_i}$ is an enumeration of $\FM(\FS_i)$,
    $(K_{i,j})_{j=1:J_i}$ is the enumeration of $\FM(\FK_i)$ in \ref{enc:spaces},
    $\varrho_{i,\tilde{j}}:=\nobreak\dim(S_{i,\tilde{j}})\in\N_0$ 
    $\forall j\in\{1,\ldots, \widetilde{J}_i\}$,
    $(\varkappa_{i,\tilde{j},\iota})_{\iota=0:\varrho_{i,\tilde{j}}}\subset\{0,1,\ldots,m_i\}$ and
    $\varpi_{i,\tilde{j}}\in\{1,\ldots,J_i\}$
    satisfy
    $\{\BIv_{i,\varkappa_{i,\tilde{j},\iota}}\}_{\iota=0:\varrho_{i,\tilde{j}}}=\nobreak V(S_{i,\tilde{j}})$, 
    $S_{i,\tilde{j}}\subseteq\nobreak K_{i,\varpi_{i,\tilde{j}}}$
    $\forall j\in\{1,\ldots, \widetilde{J}_i\}$,
    and $\overline{\eta}_i\ge \eta(\FS_i)$.

    \item\label{enc:typemeasures-aux}
    For $i=1,\ldots,N$, $\mathtt{Aux}_i$ represents the auxiliary data structure associated with the type measure $\mu_i$ that is computed by the procedure $\mathtt{InitTypeMeasure}$.
    $\mathtt{Aux}_i$ is encoded by
    $\big(\widetilde{R}_i,(F_{i,\tilde{r}},\tilde{\varpi}_{i,\tilde{r}})_{\tilde{r}=1:\widetilde{R}_i},\BP_i\big)$,
    where $\widetilde{R}_i:=|\FS_i|-1\in\N$,
    $(F_{i,\tilde{r}})_{\tilde{r}=1:\widetilde{R}_i}$
    is an enumeration of $\FS_i\setminus\nobreak\{\emptyset\}$,
    with each $F_{i,\tilde{r}}$ being a simplex encoded via \ref{enc:simplex}
    and satisfying $F_{i,\tilde{r}}\subseteq K_{i,\tilde{\varpi}_{i,\tilde{r}}}$;
    here $(K_{i,j})_{j=1:J_i}$ is the enumeration of $\FM(\FK_i)$ in \ref{enc:spaces}.
    $\BP_i\in\R_{+}^{(m_i+1)\times\widetilde{R}_i}$ 
    describes a coupling between $\mu_i$ and its discrete approximation, and
    satisfies 
    $\langle\vecone_{m_i+1},\BP_i\vecone_{\widetilde{R}_i}\rangle=1$
    as well as the properties in Lemma~\ref{lem:proc-init-type-measure}\ref{lems:proc-init-type-measure-error}.

    \item\label{enc:par-primal}
    Every possibly infeasible discrete approximate optimizer 
    $(\hat{\theta}_i)_{i=1:N}$ 
    of \eqref{eqn:mt-tf-dual}
    with $\support(\hat{\theta}_i)\subseteq V(\FS_{i})\times V(\FS_{0})$ 
    $\forall i\in\{1,\ldots,N\}$ 
    is encoded by
    $\big(n_i,\allowbreak(j_{i,\iota},l_{i,\iota},\hat{\theta}_{i,\iota})_{\iota=1:n_i}\big)_{i=1:N}$,
    where for $i=1,\ldots,N$,
    $n_i\in\N$,
    $(j_{i,\iota},l_{i,\iota})_{\iota=1:n_i}\subseteq\{0,1,\ldots,m_i\}\times\{0,1,\ldots,m_0\}$,
    $(\hat{\theta}_{i,\iota})_{\iota=1:n_i}\subset[0,1]$ satisfy 
    $\sum_{\iota=1}^{n_i}\hat{\theta}_{i,\iota}=\nobreak 1$ and
    $\hat{\theta}_i=\sum_{\iota=1}^{n_i}\hat{\theta}_{i,\iota}\delta_{(\BIv_{i,j_{i,\iota}},\BIv_{0,l_{i,\iota}})}$.

    \item\label{enc:par-dual}
    Every approximate optimizer of \eqref{eqn:mt-tf-lsip} is encoded by $(\hat{y}_{i,0},\hat{\BIy}_i,\hat{\BIw}_i)_{i=1:N}$,
    where
    $\hat{y}_{i,0}\in\R$,
    $\hat{\BIy}_i\in\nobreak\R^{m_i}$,
    $\hat{\BIw}_i\in\R^{m_0}$
    for $i=1,\ldots,N$.

    \item\label{enc:transfuncs}
    For $i=1,\ldots,N$, every transfer function $\hat{\varphi}_i:=\sum_{l=1}^{m_0}w_{i,l}g_{0,l}$ is encoded by 
    $\big(\CG_0,(w_{i,l})_{l=1:m_0}\big)$ 
    where $\CG_0$ is encoded via \ref{enc:testfuncs}.
    See the procedure $\mathtt{EvalTransFunc}$ in Section~\ref{sapx:complexity-AME} which allows one to evaluate $\hat{\varphi}_i$ at any $\BIz\in\CX_0$.

    \item\label{enc:qualitymeasure}
    Every discrete quality measure $\hat{\nu}\in\CP(\CX_0)$ in team allocations 
    $(\hat{\gamma}_i)_{i=1:N},\hat{\nu}$
    is encoded by
    $\big(n_0,d_0,(\BIz_{0,l},\hat{\nu}_{0,l})_{l=0:n_0-1}\big)$,
    where $n_0:=|\support(\hat{\nu})|\in\N$, 
    $(\BIz_{0,l})_{l=0:n_0-1}=\support(\hat{\nu})\subseteq\CX_0$,
    and
    $(\hat{\nu}_{0,l})_{l=0:n_0-1}\subset[0,1]$ satisfy $\sum_{l=0}^{n_0-1}\hat{\nu}_{0,l}=\nobreak 1$,
    $\hat{\nu}=\sum_{l=0}^{n_0-1}\hat{\nu}_{0,l}\delta_{\BIz_{0,l}}$.
    See the procedure $\mathtt{SampQual}$ in Section~\ref{sapx:complexity-AME}
    which allows one to generate any number of random samples from $\hat{\nu}$.

    \item\label{enc:allocations}
    For $i=1,\ldots,N$, every coupling $\hat{\gamma}_i\in\Gamma(\mu_i,\hat{\nu})$ in team allocations $(\hat{\gamma}_i)_{i=1:N},\hat{\nu}$ is encoded by
    $\big(\hat{\nu},\CX_i,\mu_i,\mathtt{Aux}_i,\BT_i\big)$,
    where $\hat{\nu}$ is encoded by $\big(n_0,d_0,(\BIz_{0,l},\hat{\nu}_{0,l})_{l=0:n_0-1}\big)$ via~\ref{enc:qualitymeasure},
    $\CX_i$ is encoded via~\ref{enc:spaces},
    $\mu_i$ is encoded via~\ref{enc:typemeasures},
    $\mathtt{Aux}_i$ is the auxiliary data structure encoded via \ref{enc:typemeasures-aux}, 
    $\BT_i\in\nobreak\R_{+}^{n_0\times \widetilde{R}_i}$ is a matrix with each row summing to~1.
    See the procedure $\mathtt{SampAlloc}$ in Section~\ref{sapx:complexity-AME}
    which allows one to generate any number of random samples from $\hat{\gamma}_i$.
    To make the definition of $\hat{\gamma}_i$ precise, 
    given the aforementioned encoding of $\hat{\gamma}_i$,
    let $(\Omega,\CF,\PROB)$ be a probability space,
    let $(U_t)_{t\in\N}$ be an infinite sequence of independent $\mathrm{Uniform}[0,1]$ random variables on $\Omega$,
    let $(\widetilde{\BIX}_i,\widetilde{\BIZ},\overline{t})$ be the output of 
    $\mathtt{SampAlloc}\big(\hat{\gamma}_i,(U_t)_{t\in\N}\big)$,
    and let $\hat{\gamma}_i\in\CP(\CX_i\times\CX_0)$ be defined as the law of the random variable 
    $(\widetilde{\BIX}_i,\widetilde{\BIZ}):\Omega\to\CX_i\times\CX_0$.
    Lemma~\ref{lem:proc-construct-allocs}\ref{lems:proc-construct-allocs-samp-alloc} in Section~\ref{sapx:complexity-AME} then shows that 
    $\hat{\gamma}_i\in\Gamma(\mu_i,\hat{\nu})$ indeed holds.

\end{enumerate}

\subsection{Initialization}\label{sapx:complexity-init}

\begin{procedure}[t]
    \DontPrintSemicolon
    \caption{EdgewiseDivide(){$(C, k)$ $\to$ $\FS$} \quad\citep{ECedelsbrunner1999edgewise}}
    \begin{footnotesize}
    \KwInAndOut{\!\!\!\!\!\begin{tabular}{lll}
        \textbf{Input} & 
        $C$ & an $n$-simplex in $\R^d$ encoded by $\big(d,n,(\BIv_{l})_{l=0:n}\big)$ as in \ref{enc:simplex} \\
        & $k\in\N$ & the number of pieces to divide each edge of $C$ into \\
        \textbf{Output} &
        $\FS$ & the $k$-fold edge-wise subdivision $\FS$ of $C$ encoded via \ref{enc:complex}
    \end{tabular}}
    \nl%
    \If{$n=0$}{

        \nl%
        \Return{$\FS$ encoded by $\big(d,1,(C)\big)$.} \tcp*[f]{$O(d)$}
    }

    \nl%
    $j\leftarrow 1$. 

    \nl%
    \For(\tcp*[f]{$O(k^{d+1}d^2)$}){each $(t_1,\ldots,t_n)\in\{1,\ldots,k\}^n$}{

        \nl%
        $r\leftarrow 0$.\label{alglin:edgewise-divide-construction-begin}

        \nl%
        \For(\tcp*[f]{$O(d^2)$}){$l=0,1,\ldots,n$}{
            \nl $\BIu_l\leftarrow \veczero_d$.\tcp*[f]{$O(d)$}
        }

        \nl\label{alglin:edgewise-divide-row-loop}%
        \For(\tcp*[f]{$O(kd^2)$}){$i=1,\ldots,k$}{
            \nl\label{alglin:edgewise-divide-column-loop}%
            \For(\tcp*[f]{$O(d^2)$}){$l=0,1,\ldots,n$}{
                \nl\label{alglin:edgewise-divide-check-begin}%
                \If{$l\ne 0$ and $t_l=i$}{
                    \nl%
                    $r\leftarrow r+1$. 
                }

                \nl\label{alglin:edgewise-divide-check-end}%
                $\BIu_l\leftarrow \BIu_l+\frac{\BIv_{r}}{k}$.\tcp*[f]{$O(d)$}
            }
        }
        
        \nl\label{alglin:edgewise-divide-construction-end}%
        $S_j\leftarrow \big(d,n,(\BIu_l)_{l=0:n}\big)$.\tcp*[f]{$O(d^2)$}

        \nl%
        $j\leftarrow j+1$.
    }

    \nl%
    Let $\FS$ be encoded by $\big(d,k^n,(S_j)_{j=1:k^n}\big)$ as in \ref{enc:complex}.\tcp*[f]{$O(d)$}

    \nl \Return{$\FS$.} 
    \end{footnotesize} 
\end{procedure}

Let us first present the $k$-fold edge-wise subdivision procedure
in $\EdgewiseDivide$,
which was developed by \citet{ECedelsbrunner1999edgewise}.
The properties of $\EdgewiseDivide$ are presented in 
Lemma~\ref{lem:proc-edgewise-divide} below.

\begin{lemma}\label{lem:proc-edgewise-divide}
    For any $n$-simplex $C$ in $\R^d$ and any $k\in\N$,
    $\EdgewiseDivide(C,k)$ 
    and its output $\FS$ satisfy the following statements.\widowpenalties-1000
    \begin{enumerate}[label=(\roman*),beginpenalty=10000]
        \item\label{lems:proc-edgewise-divide-complexity}
        $\EdgewiseDivide(C,k)$ incurs $O(k^{d+1}d^2)$ arithmetic operations.
        
        \item\label{lems:proc-edgewise-divide-simpcomp}
        $\FS$ is a simplicial complex with 
        $\bigcup_{S\in\FS}S=C$
        and $\volume_n(S)=\frac{\volume_n(C)}{k^n}$ $\forall S\in\FS$.

        \item\label{lems:proc-edgewise-divide-counts}
        $|\FM(\FS)|\le k^d$ 
        and $|V(\FS)|\le (d+1)k^d$.
        
        \item\label{lems:proc-edgewise-divide-barycoord}
        For any $\BIv\in V(C)$ and any $S\in\FM(\FS)$,
        there exists $\BIu\in V(S)$ with 
        $\lambda_{\FF(C),\BIv}(\BIu)\ge \frac{1}{k}$; 
        note that $\FF(C)$ is a simplicial complex.
        
        \item\label{lems:proc-edgewise-divide-meshsize}
        $\eta(\FS)\le \frac{d}{k}\max_{\BIv,\BIv'\in V(C)}\big\{\|\BIv-\BIv'\|_2\big\}$.
        
        \item\label{lems:proc-edgewise-divide-faces}
        For any non-empty $F\in\FF(C)$, 
        the output of $\EdgewiseDivide(F,k)$ is equal to $\big\{S\in\FS:S\subseteq F\big\}$.

    \end{enumerate}
\end{lemma}

\proof{Proof of Lemma~\ref{lem:proc-edgewise-divide}.}
    See Section~\ref{sapx:proof-complexity-init}.
\endproof

\begin{procedure}[t]
    \DontPrintSemicolon
    \caption{InitTestFuncs(){$(\CX_i,\eta_{\ssUB})$} $\to$ $(\CG_i,\overline{D}_i)$}
    \begin{footnotesize}
    \KwInAndOut{\!\!\!\!\!\begin{tabular}{lll}
        \textbf{Input} & 
        $\CX_i$ & the quality/type space encoded by $\big(d_i,q_i,\allowbreak(\BIu_{i,l})_{l=0:q_i},\allowbreak J_i,\allowbreak \big(\rho_{i,j},(\kappa_{i,j,\iota})_{\iota=0:\rho_{i,j}}\big)_{j=1:J_i}\big)$ as \\
        && in \ref{enc:spaces} \\
        & $\eta_{\ssUB}\in(0,1)$ & the upper bound for the mesh size when subdividing the simplicial complex $\FK_i$ \\
        \textbf{Output} &
        $\CG_i$ & the test functions on $\CX_i$ encoded via \ref{enc:testfuncs} \\
        & $\overline{D}_i\ge 1$ & an upper bound for 
        $\max_{\BIx,\BIx'\in \CX_i}\big\{\|\BIx-\BIx'\|_2\big\}$
    \end{tabular}\!\!\!\!\!\!\!\!\!\!\!\!\!\!}

    \tcc{Let us denote $d:=\max_{0\le i\le N}\{d_i\}$, $J:=\max_{0\le i\le N}\{J_i\}$, and $D:=\max_{0\le i\le N}\big\{\max_{\BIx,\BIx'\in \CX_i}\big\{\|\BIx-\BIx'\|_2\big\}\big\}\vee 1$ in the following computational complexities.}
        
    \nl\label{alglin:init-test-funcs-meshsize-init}%
    $Q_i\leftarrow \max_{1\le j\le J_i,\,0\le \iota< \iota'\le \rho_{i,j}}\big\{\|\BIu_{i,\kappa_{i,j,\iota}}-\BIu_{i,\kappa_{i,j,\iota'}}\|_2^2\big\}$.\tcp*[f]{$O(Jd^2)$}

    \nl\label{alglin:init-test-funcs-search-overlineD}%
    $\overline{D}_i\leftarrow \min\{2^q:2^{2q}\ge Q_i,\; q\in\N_0\}$.
    \tcp*[f]{$O\big(\log(D)\big)$}

    \nl\label{alglin:init-test-funcs-search-k}%
    $k_i\leftarrow \min\big\{2^p: 2^p\ge \frac{\overline{D}_id_i}{\eta_{\ssUB}},\; p\in\N_0\big\}$.
    \tcp*[f]{$O\big(\log\big(\frac{Dd}{\eta_{\ssUB}}\big)\big)$}

    \nl%
    $\widetilde{J}_i\leftarrow 0$, $m_i\leftarrow -1$.

    \nl%
    \For(\tcp*[f]{$O\big(\big(\frac{4Dd}{\eta_{\ssUB}}\big)^{2d}J^2d^3\big)$}){$j=1,\ldots,J_i$}{

        \nl\label{alglin:init-test-funcs-divide}%
        $\FC_{i,j}\leftarrow \EdgewiseDivide(K_{i,j},k_i)$, where 
        $K_{i,j}$ is encoded by $\big(d_i,\rho_{i,j},(\BIu_{i,\kappa_{i,j,\iota}})_{\iota=0:\rho_{i,j}}\big)$ as in \ref{enc:simplex},
        and let
        $\big(d_i,t_{i,j},(A_{i,j,l})_{l=1:t_{i,j}}\big)$
        be the encoding of $\FC_{i,j}$
        as in \ref{enc:complex}. 
        \tcp*[f]{$O\big(\big(\frac{4Dd}{\eta_{\ssUB}}\big)^{d+1}d^2\big)$}

        \nl\label{alglin:init-test-funcs-add-simplex-begin}%
        \For(\tcp*[f]{$O\big(\big(\frac{4Dd}{\eta_{\ssUB}}\big)^{2d}Jd^3\big)$}){$l=1,\ldots,t_{i,j}$}{

            \nl%
            $\tilde{j}\leftarrow \widetilde{J}_i+l$, $\varpi_{i,\tilde{j}}\leftarrow j$.
            
            \nl%
            Let $\big(d_i,\varrho_{i,\tilde{j}},(\BIw_{i,j,l,\iota})_{\iota=0:\varrho_{i,\tilde{j}}}\big)$ be the encoding of $A_{i,j,l}$ as in \ref{enc:simplex}. 

            \nl%
            \For(\tcp*[f]{$O\big(\big(\frac{4Dd}{\eta_{\ssUB}}\big)^{d}Jd^3\big)$}){$\iota=0,1,\ldots,\varrho_{i,\tilde{j}}$}{

                \nl\label{alglin:init-test-funcs-search-existing}%
                \uIf(\tcp*[f]{$O\big(\big(\frac{4Dd}{\eta_{\ssUB}}\big)^{d}Jd^2\big)$}){there exists $l'\in\{0,1,\ldots,m_i\}$ with $\BIw_{i,j,l,\iota}= \BIv_{i,l'}$}{

                    \nl%
                    $\varkappa_{i,\tilde{j},\iota}\leftarrow l'$.
                }
                \nl%
                \Else{
                    \nl\label{alglin:init-test-funcs-add-point-begin}%
                    $\BIv_{i,m_i+1}\leftarrow \BIw_{i,j,l,\iota}$, 
                    $\varkappa_{i,\tilde{j},\iota}\leftarrow m_i+1$.\tcp*[f]{$O(d)$}

                    \nl\label{alglin:init-test-funcs-add-point-end}%
                    $m_i\leftarrow m_i+1$.
                }
            }
        }

        \nl\label{alglin:init-test-funcs-add-simplex-end}%
        $\widetilde{J}_i\leftarrow \widetilde{J}_i+t_{i,j}$.
    }

    \nl%
    Let $\CG_i$ be encoded by $\big(d_i,m_i,(\BIv_{i,l})_{l=0:m_i},\widetilde{J}_i,\big(\varrho_{i,\tilde{j}},(\varkappa_{i,\tilde{j},\iota})_{\iota=0:\varrho_{i,\tilde{j}}},\varpi_{i,\tilde{j}}\big)_{\tilde{j}=1:\widetilde{J}_i},\eta_{\ssUB}\big)$ as in \ref{enc:testfuncs}.

    \nl%
    \Return{$(\CG_i,\overline{D}_i)$.}
    \end{footnotesize}
\end{procedure}

Using the procedure $\EdgewiseDivide$, we can subdivide the simplicial complexes $(\FK_i)_{i=0:N}$ in Setting~\ref{sett:Euclidean} 
into $(\FS_i)_{i=0:N}$ with arbitrarily small mesh sizes $(\eta(\FS_i))_{i=0:N}$.
We present the detailed procedure in $\InitTestFuncs$,
and we present the properties of $\InitTestFuncs$ in Lemma~\ref{lem:proc-init-test-funcs}.

\begin{lemma}\label{lem:proc-init-test-funcs}
    Under Setting~\ref{sett:Euclidean}, 
    let us denote $d:=\max_{0\le i\le N}\{d_i\}$,
    $J_i:=|\FM(\FK_i)|$ $\forall 0\le\nobreak i\le\nobreak N$,
    $D_i:=\nobreak\max_{\BIx,\BIx'\in \CX_i}\big\{\|\BIx-\BIx'\|_2\big\}$
    $\forall 0\le\nobreak i\le\nobreak N$,
    $J:=\max_{0\le i\le N}\{J_i\}$,
    and 
    $D:=\nobreak\max_{0\le i\le N}\{D_i\}\vee\nobreak 1$.
    Moreover, let $i\in\nobreak\{0,1,\ldots,N\}$ and $\eta_{\ssUB}\in(0,1)$ be arbitrary,
    and let $(\CG_i,\overline{D}_i)$ be
    the output of 
    $\InitTestFuncs(\CX_i,\eta_{\ssUB})$.
    Subsequently, let 
    $\CG_i$ be encoded by
    $\big(d_i,m_i,(\BIv_{i,l})_{l=0:m_i},\allowbreak \widetilde{J}_i,\allowbreak\big(\varrho_{i,\tilde{j}},\allowbreak(\varkappa_{i,\tilde{j},\iota})_{\iota=0:\varrho_{i,\tilde{j}}},\allowbreak\varpi_{i,\tilde{j}}\big)_{\tilde{j}=1:\widetilde{J}_i},\allowbreak\overline{\eta}_i\big)$ as in \ref{enc:testfuncs},
    define 
    $S_{i,\tilde{j}}:=\conv\big(\big\{\BIv_{i,\varkappa_{i,\tilde{j},\iota}}:\iota\in\{0,1,\ldots,\varrho_{i,\tilde{j}}\}\big\}\big)$ $\forall 1\le \nobreak\tilde{j}\le \nobreak\widetilde{J}_i$,
    and define
    $\FS_i:=\big\{\conv\big(\{\BIv_{i,\varkappa_{i,\tilde{j},\iota}}:\iota\in I\}\big): I \subseteq \{0,1,\ldots,\varrho_{i,\tilde{j}}\},\;
    \tilde{j}\in\{1,\ldots,\widetilde{J}_i\}\big\}$.
    Then, the following statements hold.\widowpenalty-1000
    \begin{enumerate}[label=(\roman*),beginpenalty=10000]
        \item\label{lems:proc-init-test-funcs-complexity}
        $\InitTestFuncs(\CX_i,\eta_{\ssUB})$
        incurs $O\big(\big(\frac{4Dd}{\eta_{\ssUB}}\big)^{2d}J^2d^3\big)$ arithmetic operations.
        
        \item\label{lems:proc-init-test-funcs-subdivision}
        It holds that 
        $\overline{D}_i\ge D_i$,
        $\overline{D}_i=O(D)$,
        $m_i
        \le\big(\frac{4D_id_i}{\eta_{\ssUB}}\big)^{d_i}J_i(d_i+1)
        =O\big(\big(\frac{4Dd}{\eta_{\ssUB}}\big)^{d}Jd\big)$,
        $\widetilde{J}_i=O\big(\big(\frac{4Dd}{\eta_{\ssUB}}\big)^{d}J\big)$,
        and that
        $\FS_i$ is a subdivision of $\FK_i$ 
        satisfying 
        $V(\FS_i)=\{\BIv_{i,l}\}_{l=0:m_i}$
        and $\eta(\FS_i)\le \overline{\eta}_i= \eta_{\ssUB}$.
        Moreover, for $\tilde{j}=1,\ldots,\widetilde{J}_i$,
        it holds that 
        $S_{i,\tilde{j}}\subseteq K_{i,\varpi_{i,\tilde{j}}}$,
        where $(K_{i,j})_{j=1:J_i}$ is the enumeration of $\FM(\FK_i)$ in~\ref{enc:spaces}.

        \item\label{lems:proc-init-test-funcs-LB}
        Let $\CG_i=(g_{i,j})_{j=0:m_i}$ be defined by Setting~\ref{sett:Euclidean},
        and let 
        $(p_{i,K})_{K\in\FM(\FK_i)}$ 
        satisfy the properties in Setting~\ref{sett:Euclidean}.
        Then, we have
        $\min_{0\le j\le m_i}\big\{\int_{\CX_i}g_{i,j}\DIFFX{\mu_i}\big\}\ge \frac{1}{d+2}\big(\frac{\eta_{\ssUB}}{4Dd}\big)^{d+1}\min_{K\in\FM(\FK_i)}\big\{\int_{K}p_{i,K}\DIFFX{\Lebesgue_K}\big\}$.
    \end{enumerate}
\end{lemma}

\proof{Proof of Lemma~\ref{lem:proc-init-test-funcs}.}
    See Section~\ref{sapx:proof-complexity-init}.
\endproof

\begin{procedure}[t]
    \DontPrintSemicolon
    \caption{InitTypeMeasure(){$(\CX_i,\mu_i,\CG_i)$} $\to$ $(\bar{\BIg}_i, \mathtt{Aux}_i)$}
    \begin{footnotesize}
    \KwInAndOut{\!\!\!\!\!\begin{tabular}{lll}
        \textbf{Input} & 
        $\CX_i$ & the type space encoded by $\big(d_i,q_i,\allowbreak(\BIu_{i,l})_{l=0:q_i},\allowbreak J_i,\allowbreak \big(\rho_{i,j},(\kappa_{i,j,\iota})_{\iota=0:\rho_{i,j}}\big)_{j=1:J_i}\big)$ as in \ref{enc:spaces} \\
        & $\mu_i\in\CP(\CX_i)$ & the type measure on $\CX_i$ with the associated oracle procedures $\mathtt{Int}$ and $\mathtt{Samp}$; see \ref{enc:typemeasures} \\
        & $\CG_i$ & the test functions on $\CX_i$ encoded by $\big(d_i,m_i,(\BIv_{i,l})_{l=0:m_i},\widetilde{J}_{i},\big(\varrho_{i,\tilde{j}},(\varkappa_{i,\tilde{j},\iota})_{\iota=0:\varrho_{i,\tilde{j}}},$ \\
        && $\varpi_{i,\tilde{j}}\big)_{\tilde{j}=1:\widetilde{J}_{i}},\overline{\eta}_i\big)$ as in \ref{enc:testfuncs}\\
        \textbf{Output}
        & $\bar{\BIg}_i\in\R^{m_i}$ & the vector containing the integrals of the test functions in $\CG_i$ with respect to $\mu_i$ \\
        & $\mathtt{Aux}_i$ & the auxiliary data structure encoded via \ref{enc:typemeasures-aux}
    \end{tabular}\!\!\!\!\!\!\!\!\!\!\!\!\!\!}

    \tcc{In the following computational complexities, let us denote $d:=\max_{0\le i\le N}\{d_i\}$, $\widetilde{J}:=\max_{0\le i\le N}\{\widetilde{J}_i\}$,
    and let $T_{\mathtt{Int}}$ denote the number of arithmetic operations incurred by 
    $\mathtt{Int}$.}

    \nl%
    $\mathtt{FaceSet}\leftarrow\emptyset$, $\tilde{r}\leftarrow 0$.

    \nl%
    \For(\tcp*[f]{$O\big(2^d\widetilde{J}(T_{\mathtt{Int}}+\widetilde{J}d)\big)$}){$\tilde{j}=1,\ldots,\widetilde{J}_i$}{
        
        \nl\label{alglin:init-type-measure-enumerate-face}%
        \For(\tcp*[f]{$O\big(2^d(T_{\mathtt{Int}}+\widetilde{J}d)\big)$}){each non-empty $I\subseteq \{0,1,\ldots,\varrho_{i,\tilde{j}}\}$}{

            \nl\label{alglin:init-type-measure-search-face}%
            \If(\tcp*[f]{$O\big(d\log(d)+\log(\widetilde{J})\big)$}){$\{\varkappa_{i,\tilde{j},\iota}\}_{\iota\in I} \in \mathtt{FaceSet}$}{

                \nl \textbf{continue}.
            }

            \nl%
            $\tilde{r}\leftarrow \tilde{r} + 1$. 

            \nl\label{alglin:init-type-measure-initzero}%
            $(p_{0,\tilde{r}},p_{1,\tilde{r}},\ldots, p_{m_i,\tilde{r}})^\TRANSP\leftarrow \veczero_{m_i+1}$.\tcp*[f]{$O(\widetilde{J}d)$}

            \nl\label{alglin:init-type-measure-encode-face}%
            Let the simplex $F_{i,\tilde{r}}$ be encoded by 
            $\big(d_i,|I|-1,(\BIv_{i,\varkappa_{i,\tilde{j},\iota}})_{\iota\in I}\big)$ as in \ref{enc:simplex}.\tcp*[f]{$O(d^2)$}

            \nl\label{alglin:init-type-measure-simplex-index}%
            $\tilde{\varpi}_{i,\tilde{r}}\leftarrow \varpi_{i,\tilde{j}}$.

            \nl\label{alglin:init-type-measure-insert-face}$%
            \mathtt{FaceSet}\leftarrow \mathtt{FaceSet}\cup \big\{\{\varkappa_{i,\tilde{j},\iota}\}_{\iota\in I}\big\}$.\tcp*[f]{$O\big(d\log(d)+\log(\widetilde{J})\big)$}

            \nl\label{alglin:init-type-measure-compute-integral}%
            $(p_{\varkappa_{i,\tilde{j},\iota},\tilde{r}})_{\iota\in I}\leftarrow \mathtt{Int}(\CX_i,\mu_i,F_{i,\tilde{r}},\tilde{\varpi}_{i,\tilde{r}})$.\tcp*[f]{$O(T_{\mathtt{Int}})$}
        }
    }

    \nl%
    $\widetilde{R}_{i}\leftarrow \tilde{r}$.

    \nl\label{alglin:init-type-measure-probmat}%
    $\BP_i\leftarrow(p_{l,\tilde{r}})_{l=0:m_i,\,\tilde{r}=1:\widetilde{R}_i}\in\R^{(m_i+1)\times\widetilde{R}_i}$.\tcp*[f]{$O(2^d\widetilde{J}^2d)$}

    \nl\label{alglin:init-type-measure-aggregate}%
    $(\bar{g}_{i,0},\bar{\BIg}_i^\TRANSP)^\TRANSP \leftarrow \BP_i\vecone_{\widetilde{R}_i}$.\tcp*[f]{$O(2^d\widetilde{J}^2d)$}

    \nl%
    \Return{$\big(\bar{\BIg}_i, \big(\widetilde{R}_i, (F_{i,\tilde{r}},\tilde{\varpi}_{i,\tilde{r}})_{\tilde{r}=1:\widetilde{R}_i}, \BP_i\big)\big)$.}
    \end{footnotesize}
\end{procedure}

With the test functions $\CG_i$ computed by the procedure $\InitTestFuncs$,
we can subsequently compute the objective vector 
$\bar{\BIg}_i$ in \eqref{eqn:mt-tf-lsip} 
via the procedure $\InitTypeMeasure$.
$\InitTypeMeasure$
also computes an auxiliary data structure 
$\mathtt{Aux}_i$
that is used later to construct an $\epsilon$-matching equilibrium
in the procedure $\mathtt{ConstructAllocs}$.
The properties of $\InitTypeMeasure$ are presented in 
Lemma~\ref{lem:proc-init-type-measure}.

\begin{lemma}\label{lem:proc-init-type-measure}
    Under Setting~\ref{sett:Euclidean}, 
    let us denote $d:=\max_{0\le i\le N}\{d_i\}$,
    $\widetilde{J}:=\max_{0\le i\le N}\!\big\{|\FM(\FS_i)|\big\}$, 
    and let us denote by $T_{\mathtt{Int}}$ the number of arithmetic operations incurred by the oracle procedure $\mathtt{Int}$ in \ref{enc:typemeasures}.
    Let $i\in\nobreak\{1,\ldots,N\}$ be arbitrary, let    
    $(\bar{\BIg}_i,\mathtt{Aux}_i)$ be
    the output of 
    $\InitTypeMeasure(\CX_i,\mu_i,\allowbreak\CG_i)$
    where $\mathtt{Aux}_i$ is encoded by 
    $\big(\widetilde{R}_i,(F_{i,\tilde{r}},\tilde{\varpi}_{i,\tilde{r}})_{\tilde{r}=1:\widetilde{R}_i},\BP_i\big)$ as in \ref{enc:typemeasures-aux},
    and denote $\BP_i=(p_{l,\tilde{r}})_{l=0:m_i,\,\tilde{r}=1:\widetilde{R}_i}$.
    Then, the following statements hold.
    \begin{enumerate}[label=(\roman*),beginpenalty=10000]
        \item\label{lems:proc-init-type-measure-complexity}
        $\InitTypeMeasure(\CX_i,\mu_i,\allowbreak\CG_i)$
        incurs
        $O\big(2^d\widetilde{J}(T_{\mathtt{Int}}+\widetilde{J}d)\big)$ arithmetic operations.

        \item\label{lems:proc-init-type-measure-auxiliary}
        It holds that
        $\widetilde{R}_{i}\ge m_{i}$,
        $\widetilde{R}_{i}=O(2^d\widetilde{J})$,
        $(F_{i,\tilde{r}})_{\tilde{r}=1:\widetilde{R}_{i}}$ 
        is an enumeration of $\FS_i\setminus \{\emptyset\}$
        with $F_{i,\tilde{r}}\subseteq\nobreak K_{i,\tilde{\varpi}_{i,\tilde{r}}}$ 
        $\forall \tilde{r}\in\{1,\ldots,\widetilde{R}_i\}$,
        where $(K_{i,j})_{j=1:J_i}$ is the enumeration of $\FM(\FK_i)$ in~\ref{enc:spaces}.

        \item\label{lems:proc-init-type-measure-integrals}
        It holds that $\bar{\BIg}_i=\big(\int_{\CX_i}g_{i,1}\DIFFX{\mu_i},\ldots,\int_{\CX_i}g_{i,m_i}\DIFFX{\mu_i}\big)^\TRANSP$
        and that 
        $\BP_i\vecone_{\widetilde{R}_i}=\big(\int_{\CX_i}g_{i,0}\DIFFX{\mu_i},\allowbreak\int_{\CX_i}g_{i,1}\DIFFX{\mu_i},\ldots,\allowbreak\int_{\CX_i}g_{i,m_i}\DIFFX{\mu_i}\big)^\TRANSP$.
        
        \item\label{lems:proc-init-type-measure-error}
        Let $(\Omega,\CF,\PROB)$ be a probability space, 
        let the random variables 
        $Q:\Omega\to\{0,1,\ldots,m_i\}$ and
        $S:\Omega\to\{1,\ldots,\widetilde{R}_i\}$
        satisfy
        $\PROB[Q=l,\;S=\tilde{r}]=p_{l,\tilde{r}}$
        $\forall l\in\{0,1,\ldots,m_i\}$,
        $\forall \tilde{r}\in\{1,\ldots,\widetilde{R}_i\}$,
        define
        $\BIX_i:\Omega\to\CX_i$
        by 
        $\BIX_i:=\sum_{l=0}^{m_i}\BIv_{i,l}\INDI_{\{Q=l\}}$,
        and let the random variable 
        $\overbar{\BIX}_i:\Omega\to\CX_i$ satisfy
        $\PROB[\overbar{\BIX}_i\in\nobreak E|S=\nobreak\tilde{r}]=\frac{\mu_i(E\,\cap\, \relint(F_{i,\tilde{r}}))}{\mu_i(\relint(F_{i,\tilde{r}}))}$
        $\forall E\in \CB(\CX_i)$, 
        $\forall \tilde{r}\in\big\{\tilde{r}'\in\{1,\ldots,\widetilde{R}_i\}:\PROB[S=\tilde{r}']>0\big\}$.
        Then, 
        it holds that $\PROB[S=\nobreak\tilde{r}]=\mu_i\big(\relint(F_{i,\tilde{r}})\big)$ 
        $\forall \tilde{r}\in\{1,\ldots,\widetilde{R}_i\}$,
        the law of $\overbar{\BIX}_i$ is equal to~$\mu_i$, and that
        $\EXP\big[\|\BIX_i-\nobreak\overbar{\BIX}_i\|_2\big]\le\nobreak \eta(\FS_i)$.
    \end{enumerate}
\end{lemma}

\proof{Proof of Lemma~\ref{lem:proc-init-type-measure}.}
    See Section~\ref{sapx:proof-complexity-init}.
\endproof

Now that we can explicitly compute 
$(\CG_i)_{i=0:N}$ and $(\bar{\BIg}_i)_{i=1:N}$,
the next subsection will present explicit procedures 
$\mathtt{SolveLSIPDual}$ and $\mathtt{SolveLSIPPrimal}$
for approximately solving 
\eqref{eqn:mt-tf-lsip} and \eqref{eqn:mt-tf-dual}
with $(\CG_i)_{i=0:N}$ and $(\bar{\BIg}_i)_{i=1:N}$ as inputs.

\subsection{Approximately solving \eqref{eqn:mt-tf-lsip} and \eqref{eqn:mt-tf-dual}}\label{sapx:complexity-lsip}

Our strategy for computing an approximate optimizer of \eqref{eqn:mt-tf-lsip}
is to consider the linear programming (LP) relaxation of \eqref{eqn:mt-tf-lsip}
with respect to the grids $\big(V(\FS_i)\times V(\FS_0)\big)_{i=1:N}$
and then
use the volumetric center algorithm of \citet{ECvaidya1996new} to approximately solve the LP relaxation.
Since the volumetric center algorithm of \citet{ECvaidya1996new}
requires the superlevel sets of the LP relaxation to contain open balls,
we relax the equality constraint $\sum_{i=1}^{N}\BIw_i=\veczero_{m_0}$
into an inequality constraint $\sum_{i=1}^{N}\BIw_i\ge \veczero_{m_0}$.
Moreover, for technical reasons related to the simplification of $\|\cdot\|_{\infty}$ 
upper bounds of certain superlevel sets of the LP relaxation,
we will shift the right-hand side of the inequalities by $-c_{\ssmax}$,
where $c_{\ssmax}:=\max_{1\le i\le N,\,0\le j\le m_i,\,0\le l\le m_0}\big\{c_{i}(\BIv_{i,j},\BIv_{0,l})\big\}$.
Thus, the shifted LP relaxation of \eqref{eqn:mt-tf-lsip}
with respect to the grids $\big(V(\FS_i)\times V(\FS_0)\big)_{i=1:N}$ is given by:
\begin{align}
    \begin{split}
        \maximize_{(y_{i,0},\BIy_i,\BIw_i)}\quad & \sum_{i=1}^{N} y_{i,0} + \langle\bar{\BIg}_i,\BIy_i\rangle\\
        \mathrm{subject~to}\quad & y_{i,0}+\langle\BIg_i(\BIv_{i,j}),\BIy_i\rangle + \langle\BIg_0(\BIv_{0,l}),\BIw_i\rangle \le c_i(\BIv_{i,j},\BIv_{0,l})-c_{\ssmax} \\
        & \hspace{108.2pt} \forall 0\le j\le m_i,\; \forall 0\le l\le m_0,\; \forall 1\le i\le N,\\
        & \sum_{i=1}^{N}\BIw_i\ge\veczero_{m_0},\qquad y_{i,0}\in\R, \; \BIy_i\in\R^{m_i},\; \BIw_i\in\R^{m_0} \quad \forall 1\le i\le N.
    \end{split}
    \label{eqn:mt-tf-lp-full}
\end{align}

\begin{procedure}[t]
    \DontPrintSemicolon
    \caption{GetSepHyp(){$\big(N,(\CG_i)_{i=0:N},(\bar{\BIg}_i,c_i)_{i=1:N},c_{\ssmax},\Bchi\big)$} $\to$ $(\mathtt{Inside},\mathtt{Tag},\BIa)$}
    \begin{footnotesize}
    \KwInAndOut{\!\!\!\!\!\begin{tabular}{lll}
        \textbf{Input} &
        $N\in\N\cap[2,\infty)$ & the number of agent populations \\ 
        & $\CG_i$ & the test functions on $\CX_i$ encoded by $\big(d_i,m_i,(\BIv_{i,l})_{l=0:m_i},\widetilde{J}_{i},\big(\varrho_{i,\tilde{j}},(\varkappa_{i,\tilde{j},\iota})_{\iota=0:\varrho_{i,\tilde{j}}},$ \\
        && $\varpi_{i,\tilde{j}}\big)_{\tilde{j}=1:\widetilde{J}_{i}},\overline{\eta}_i\big)$ as in \ref{enc:testfuncs}\\
        & $\bar{\BIg}_i\in\R^{m_i}$ & the vector containing the integrals of the test functions in $\CG_i$ with respect to $\mu_i$ \\
        & $c_i:\CX_i\times\CX_0\to\R$ & the cost function on $\CX_i\times\CX_0$; see \ref{enc:costfuncs} \\ 
        & $c_{\ssmax}\in\R$ & the pre-computed value of $c_{\ssmax}:=\max_{1\le i\le N,\,0\le j\le m_i,\,0\le l\le m_0}\big\{c_i(\BIv_{i,j},\BIv_{0,l})\big\}$ \\
        & $\Bchi\in\R^{q}$ & a vector in $\R^q$ where $q:=\sum_{i=1}^{N}(1+m_i+m_0)$ \\
        \textbf{Output} &
        $\mathtt{Inside}$ & a boolean value indicating whether $\Bchi$ is inside the feasible set of the problem \\
        & $\mathtt{Tag}$ & a tuple with additional information about the separating hyperplane whose value is \\
        && either $(0,0,0)$ or $(i,j,l)$ where $i\in\{1,\ldots,N\}$, $j\in\{0,1,\ldots,m_i\}$, $l\in\{0,1,\ldots,m_0\}$ \\
        & $\BIa\in\R^{q}$ & a separating hyperplane; see Lemma~\ref{lem:proc-get-sep-hyp}\ref{lems:proc-get-sep-hyp-separation}
    \end{tabular}\!\!\!\!\!\!\!\!\!\!\!\!\!\!}

    \tcc{In the following computational complexities, let us denote $m:=\max_{0\le i\le N}\{m_i\}$,
    and let $T_{\mathtt{Eval}}$ denote the number of arithmetic operations incurred when evaluating $c_i(\BIx,\BIz)$ at any $(\BIx,\BIz)\in\CX_i\times\CX_i$.}

    \nl%
    $(y_{1,0},y_{1,1},\ldots,y_{1,m_1},w_{1,1},\ldots,w_{1,m_0},\ldots,y_{N,0},y_{N,1},\ldots,y_{N,m_N},w_{N,1},\ldots,w_{N,m_0})^\TRANSP\leftarrow \Bchi$.\tcp*[f]{$O(Nm)$}

    \nl\label{alglin:get-sep-hyp-balancecut-forloop}%
    \For(\tcp*[f]{$O(Nm)$}){$l=1,\ldots,m_0$}{

        \nl\label{alglin:get-sep-hyp-balancecut-check}%
        \If(\tcp*[f]{$O(N)$}){$\sum_{i=1}^{N}w_{i,l}<0$}{

            \nl\label{alglin:get-sep-hyp-balancecut}%
            $\BIa\leftarrow (0,\veczero_{m_1}^\TRANSP,\BIe_{0,l}^\TRANSP,0,\veczero_{m_2}^\TRANSP,\BIe_{0,l}^\TRANSP,\ldots,0,\veczero_{m_N}^\TRANSP,\BIe_{0,l}^\TRANSP)^\TRANSP$, where $\BIe_{0,l}$ denotes the $l$-th standard basis vector of $\R^{m_0}$.\tcp*[f]{$O(Nm)$}

            \nl\label{alglin:get-sep-hyp-balancecut-return}%
            \Return{$\big(\mathtt{False},(0,0,0),\BIa\big)$}.
        }
    }

    \nl\label{alglin:get-sep-hyp-feascut-forloop}%
    \For(\tcp*[f]{$O(Nm^2T_{\mathtt{Eval}})$}){$i=1,\ldots,N$}{

        \nl%
        \For(\tcp*[f]{$O(m^2T_{\mathtt{Eval}}+Nm)$}){$j=0,1,\ldots,m_i$}{

            \nl\label{alglin:get-sep-hyp-feascut-forloop-innermost}%
            \For(\tcp*[f]{$O(mT_{\mathtt{Eval}}+Nm)$}){$l=0,1,\ldots,m_0$}{

                \nl\label{alglin:get-sep-hyp-feascut-check}%
                \If(\tcp*[f]{$O(T_{\mathtt{Eval}})$}){$c_i(\BIv_{i,j},\BIv_{0,l})-c_{\ssmax}-y_{i,0}-\INDI_{\{j\ne 0\}}y_{i,j}-\INDI_{\{l\ne 0\}}w_{i,l} < 0$}{

                    \nl\label{alglin:get-sep-hyp-feascutstart}%
                    $\BIa_i\leftarrow (1,\BIe_{i,j}^\TRANSP,\BIe_{0,l}^\TRANSP)^\TRANSP$, 
                    where $\BIe_{i,0}:=\veczero_{m_i}$, 
                    $\BIe_{i,j}$ denotes the $j$-th standard basis vector of $\R^{m_i}$ if $j\ne 0$,
                    $\BIe_{0,0}:=\veczero_{m_0}$,
                    $\BIe_{0,l}$ denotes the $l$-th standard basis vector of $\R^{m_0}$ if $l\ne 0$.\tcp*[f]{$O(m)$}

                    \nl%
                    $\BIa_{i'}\leftarrow \veczero_{1+m_{i'}+m_0}$ $\forall i'\in\{1,\ldots,N\}\setminus\{i\}$.\tcp*[f]{$O(Nm)$}

                    \nl\label{alglin:get-sep-hyp-feascutend}%
                    $\BIa\leftarrow -(\BIa_1^\TRANSP,\ldots,\BIa_{N}^\TRANSP)^\TRANSP$.\tcp*[f]{$O(Nm)$}

                    \nl\label{alglin:get-sep-hyp-feascut-return}%
                    \Return{$\big(\mathtt{False},(i,j,l),\BIa\big)$}.
                }
            }
        }
    }

    \nl%
    $\BIa\leftarrow (1,\bar{\BIg}_1^\TRANSP,\veczero_{m_0}^\TRANSP,\ldots,1,\bar{\BIg}_N^\TRANSP,\veczero_{m_0}^\TRANSP)^\TRANSP$.\tcp*[f]{$O(Nm)$}

    \nl\label{alglin:get-sep-hyp-objcut-return}%
    \Return{$(\mathtt{True}, (0,0,0), \BIa)$.}
    \end{footnotesize}
\end{procedure}

\begin{procedure}[t]
    \DontPrintSemicolon
    \caption{GetVolumetricQuantities(){$\big(q,k,(\BIa_j,b_j)_{j=1:k},\Bchi\big)$} $\to$ $\big(\BH,(\sigma_j)_{j=1:k}, \Bchi_{\ssstep}\big)$}
    \begin{footnotesize}
    \KwInAndOut{\!\!\!\!\!\begin{tabular}{lll}
        \textbf{Input} &
        $q\in\N$ & the dimension $\sum_{i=1}^{N}(1+m_i+m_0)$ of the decision variable of \eqref{eqn:mt-tf-lsip} \\
        & $k\in\N$ & the number of closed half-spaces characterizing a polytope; $k\ge q$ \\
        & $\BIa_j\in \R^{q}$, $b_j\in\R$ & characterization of a closed half-space $\{\Bchi'\in\R^{q}:\langle\BIa_j,\Bchi'\rangle\ge b_j\}\subset \R^{q}$; $(\BIa_j,b_j)_{j=1:k}$ \\
        && characterize the polytope $P=\{\Bchi'\in\R^{q}:\langle\BIa_j,\Bchi'\rangle\ge b_j\; \forall 1\le j\le k\}\subset \R^{q}$ \\
        & $\Bchi\in\R^{q}$ & a vector in $\inter(P)$ \\
        \textbf{Output} &
        $\BH\in\R^{q\times q}$ & the computed value of $\BH(\Bchi)$ defined in \citep[p.294]{ECvaidya1996new} \\
        & $(\sigma_j)_{j=1:k}\subset (0,\infty)$ & the computed values of $\big(\sigma_j(\Bchi)\big)_{j=1:k}$ defined in \citep[p.295]{ECvaidya1996new} \\
        & $\Bchi_{\ssstep}\in\R^{q}$ & the computed value of $\BQ(\Bchi)^{-1}\nabla F(\Bchi)$ defined in \citep[p.295--296]{ECvaidya1996new} 
    \end{tabular}\!\!\!\!\!\!\!\!\!\!\!\!\!\!}

    \nl%
    $\BH\leftarrow \sum_{j=1}^{k}\frac{\BIa_j\BIa_j^\TRANSP}{(\langle\BIa_j,\Bchi\rangle-b_j)^2}$.\tcp*[f]{$O(q^\omega)$}

    \nl%
    $\sigma_j\leftarrow \frac{\langle\BIa_j,\BH^{-1}\BIa_j\rangle}{(\langle\BIa_j,\Bchi\rangle-b_j)^2}$ $\forall j\in\{1,\ldots,k\}$.\tcp*[f]{$O(q^\omega)$}

    \nl%
    $\nabla F\leftarrow -\sum_{j=1}^{k}\frac{\sigma_j\BIa_j}{\langle\BIa_j,\Bchi\rangle-b_j}$,
    $\BQ\leftarrow\sum_{j=1}^{k}\frac{\sigma_j\BIa_j\BIa_j^\TRANSP}{(\langle\BIa_j,\Bchi\rangle-b_j)^2}$.\tcp*[f]{$O(q^\omega)$}

    \nl%
    $\Bchi_{\ssstep}\leftarrow \BQ^{-1}\nabla F$.\tcp*[f]{$O(q^\omega)$}

    \nl%
    \Return{$\big(\BH,(\sigma_j)_{j=1:k},\Bchi_{\ssstep}\big)$.}
    \end{footnotesize}
\end{procedure}

\begin{procedure}
    \DontPrintSemicolon
    \caption{SolveLSIPDual(){$\big(N,(\CG_i)_{i=0:N},(\bar{\BIg}_i,c_i)_{i=1:N},\epsilon_0,L_c\big)$}  $\to$ $\big((\hat{y}_{i,0},\hat{\BIy}_i,\hat{\BIw}_i)_{i=1:N},\widehat{\CI}\big)$}
    \begin{footnotesize}
    \KwInAndOut{\!\!\!\!\!\begin{tabular}{lll}
        \textbf{Input} &
        $N\in\N\cap[2,\infty)$ & the number of agent populations \\ 
        & $\CG_i$ & the test functions on $\CX_i$ encoded by $\big(d_i,m_i,(\BIv_{i,l})_{l=0:m_i},\widetilde{J}_{i},\big(\varrho_{i,\tilde{j}},(\varkappa_{i,\tilde{j},\iota})_{\iota=0:\varrho_{i,\tilde{j}}},$ \\
        && $\varpi_{i,\tilde{j}}\big)_{\tilde{j}=1:\widetilde{J}_{i}},\overline{\eta}_i\big)$ as in \ref{enc:testfuncs}\\
        & $\bar{\BIg}_i\in\R^{m_i}$ & the vector containing the integrals of the test functions in $\CG_i$ with respect to $\mu_i$ \\
        & $c_i:\CX_i\times\CX_0\to\R$ & the cost function on $\CX_i\times\CX_0$; see \ref{enc:costfuncs} \\ 
        & $\epsilon_0\in(0,1)$ & the tolerance value for the algorithm of \citet{ECvaidya1996new}; see Theorem~\ref{thm:proc-solve-LSIP} \\
        & $L_c>0$ & the Lipschitz constant of the cost functions $(c_i)_{i=1:N}$ \\
        \textbf{Output} &
        $(\hat{y}_{i,0},\hat{\BIy}_i,\hat{\BIw}_i)_{i=1:N}$ & an approximate optimizer of \eqref{eqn:mt-tf-lsip} encoded via \ref{enc:par-dual} \\ 
        & $\widehat{\CI}$ & a set of tuples $(i,j,l)$ where $j\in\{0,1,\ldots,m_i\}$, $l\in\{0,1,\ldots,m_0\}$, $i\in\{1,\ldots,N\}$
    \end{tabular}\!\!\!\!\!\!\!\!\!\!\!\!\!\!}

    \tcc{In the following computational complexities, let us denote $m:=\max_{0\le i\le N}\{m_i\}$,
    and let $T_{\mathtt{Eval}}$ denote the number of arithmetic operations incurred when evaluating $c_i(\BIx,\BIz)$ at any $(\BIx,\BIz)\in\CX_i\times\CX_0$.}

    \nl\label{alglin:solve-LSIP-dual-init-start}%
    $c_{\ssmax}\leftarrow \max_{1\le i\le N,\,0\le j\le m_i,\,0\le l\le m_0}\big\{c_i(\BIv_{i,j},\BIv_{0,l})\big\}$.\tcp*[f]{$O(Nm^2T_{\mathtt{Eval}})$}

    \nl%
    $c_{\ssmin}\leftarrow \min_{1\le i\le N,\,0\le j\le m_i,\,0\le l\le m_0}\big\{c_i(\BIv_{i,j},\BIv_{0,l})\big\}$.\tcp*[f]{$O(Nm^2T_{\mathtt{Eval}})$}

    \nl%
    \For(\tcp*[f]{$O(Nm)$}){$i=1,\ldots,N$}{

        \nl%
        $(\bar{g}_{i,1},\ldots,\bar{g}_{i,m_i})^\TRANSP\leftarrow \bar{\BIg}_i$, $\bar{g}_{i,0}\leftarrow 1-\langle\bar{\BIg}_i,\vecone_{m_i}\rangle$.\tcp*[f]{$O(m)$}
    }

    \nl%
    $\bar{g}_{\ssmin}\leftarrow \min_{1\le i\le N,\,0\le j\le m_i}\{\bar{g}_{i,j}\}$,
    $q\leftarrow \sum_{i=1}^{N}(1+m_i+m_0)$.\tcp*[f]{$O(Nm)$}

    \nl\label{alglin:solve-LSIP-dual-init-center}%
    $M\leftarrow (1+\bar{g}_{\ssmin}^{-1})N^2(c_{\ssmax}-c_{\ssmin}+1)$, 
    $\Bchi_0\leftarrow \veczero_{q}$, 
    $\hat{\Bchi}\leftarrow \veczero_{q}$,
    $\hat{\tau}\leftarrow -\infty$,
    $\overline{k}\leftarrow 2q$.\tcp*[f]{$O(Nm)$}

    \nl\label{alglin:solve-LSIP-dual-init-box}%
    $\BIa_{j}\leftarrow \BIe_{j}$, $b_{j}\leftarrow -M$, $\iota_{j}\leftarrow (0,0,0)$, $t_{j}\leftarrow 1$, 
    $\BIa_{q+j}\leftarrow -\BIe_{j}$, $b_{q+j}\leftarrow -M$, $\iota_{q+j}\leftarrow (0,0,0)$, $t_{q+j}\leftarrow 1$ 
    $\forall j\in\{1,\ldots,q\}$, 
    where $\BIe_j$ denotes the $j$-th standard basis vector of~$\R^{q}$.\tcp*[f]{$O(Nm^2)$}

    \nl\label{alglin:solve-LSIP-dual-init-iternum}%
    $\overline{r}\leftarrow 2\cdot 10^7 q \big(18+\min\big\{u\in\N:2^{u}\ge \frac{qNM}{\epsilon_0}\big\}\big)$. \tcp*[f]{$O\big(\log\big(\frac{Nm(c_{\ssmax}-c_{\ssmin}+1)}{\epsilon_0 \bar{g}_{\ssmin}}\big)\big)$}

    \nl\label{alglin:solve-LSIP-dual-forloop}%
    \For(\tcp*[f]{$O\big(Nm\log\big(\frac{Nm(c_{\ssmax}-c_{\ssmin}+1)}{\epsilon_0\bar{g}_{\ssmin}}\big)(N^{\omega}m^{\omega}+Nm^2T_{\mathtt{Eval}})\big)$}){$r=1,\ldots,\overline{r}$}{

        \nl\label{alglin:solve-LSIP-dual-activecuts}%
        $\{j_{l}\}_{l=1:k}\leftarrow \big\{j\in\{1,\ldots,\overline{k}\}:t_j=1\big\}$.\tcp*[f]{$O(Nm)$}

        \nl%
        $\big(\BH,(\sigma_l)_{l=1:k},\Bchi_{\ssstep}\big)\leftarrow \GetVolumetricQuantities\big(q,k,(\BIa_{j_{l}},b_{j_{l}})_{l=1:k},\Bchi_{r-1}\big)$.\tcp*[f]{$O(N^{\omega}m^{\omega})$}

        \nl%
        Let $\hat{l}\in\argmin_{1\le l\le k}\{\sigma_l\}$.\tcp*[f]{$O(Nm)$}

        \nl\label{alglin:solve-LSIP-dual-init-check-addorrem}%
        \uIf{$\sigma_{\hat{l}}\ge 10^{-7}$}{

            \nl\label{alglin:solve-LSIP-dual-add-start}%
            $(\mathtt{Inside},\mathtt{Tag},\hat{\BIa})\leftarrow \GetSepHyp\big(N,(\CG_i)_{i=0:N},(\bar{\BIg}_i,c_i)_{i=1:N},c_{\ssmax},\Bchi_{r-1}\big)$.\tcp*[f]{$O(Nm^2T_{\mathtt{Eval}})$}

            \nl\label{alglin:solve-LSIP-dual-inside-check}%
            \If(\tcp*[f]{$O(Nm)$}){$\mathtt{Inside}=\mathtt{True}$ and $\langle\hat{\BIa},\Bchi_{r-1}\rangle> \hat{\tau}$}{
                
                \nl\label{alglin:solve-LSIP-dual-inside-update}%
                $\hat{\tau}\leftarrow \langle\hat{\BIa},\Bchi_{r-1}\rangle$, $\hat{\Bchi}\leftarrow\Bchi_{r-1}$. \tcp*[f]{$O(Nm)$}
            }

            \nl\label{alglin:solve-LSIP-dual-RHS}%
            Choose $\hat{b}\in(-\infty,\langle\hat{\BIa},\Bchi_{r-1}\rangle)$
            such that
            $1.5811376\cdot 10^{-6}\le \frac{\langle\hat{\BIa},\BH^{-1}\hat{\BIa}\rangle}{(\langle\hat{\BIa},\Bchi_{r-1}\rangle-\hat{b})^2} \le 1.5811388\cdot 10^{-6}$.\tcp*[f]{$O(N^{\omega}m^{\omega})$}

            \nl\label{alglin:solve-LSIP-dual-add}%
            $\BIa_{\overline{k}+1}\leftarrow \hat{\BIa}$, $b_{\overline{k}+1}\leftarrow \hat{b}$, 
            $\iota_{\overline{k}+1}\leftarrow \mathtt{Tag}$,
            $t_{\overline{k}+1}\leftarrow 1$.\tcp*[f]{$O(Nm)$}

            \nl%
            $\overline{k}\leftarrow \overline{k}+1$, $\Bchi_{r}\leftarrow \Bchi_{r-1}$.\tcp*[f]{$O(Nm)$}

            \nl\label{alglin:solve-LSIP-dual-add-Newtonloop}%
            \For(\tcp*[f]{$O(N^{\omega}m^{\omega})$}){$s=1,\ldots,2197$}{

                \nl%
                $\big(\BH,(\sigma_l)_{l=1:k+1},\Bchi_{\ssstep}\big)\leftarrow \GetVolumetricQuantities\big(q,k+1,(\BIa_{j},b_{j})_{j\in\{j':t_{j'}=1\}},\Bchi_{r}\big)$.\tcp*[f]{$O(N^{\omega}m^{\omega})$}

                \nl$\label{alglin:solve-LSIP-dual-add-end}%
                \Bchi_{r}\leftarrow \Bchi_{r}-0.18\Bchi_{\ssstep}$.\tcp*[f]{$O(Nm)$}
            }

        }
        \nl%
        \Else{

            \nl\label{alglin:solve-LSIP-dual-rem-start}%
            $t_{j_{\hat{l}}}\leftarrow 0$, $\Bchi_{r}\leftarrow\Bchi_{r-1}$.\tcp*[f]{$O(Nm)$}

            \nl\label{alglin:solve-LSIP-dual-rem-Newtonloop}%
            \For(\tcp*[f]{$O(N^{\omega}m^{\omega})$}){$s=1,\ldots,1493$}{

                \nl%
                $\big(\BH,(\sigma_l)_{l=1:k-1},\Bchi_{\ssstep}\big)\leftarrow \GetVolumetricQuantities\big(q,k-1,(\BIa_{j},b_{j})_{j\in\{j':t_{j'}=1\}},\Bchi_{r}\big)$.\tcp*[f]{$O(N^{\omega}m^{\omega})$}

                \nl\label{alglin:solve-LSIP-dual-rem-end}%
                $\Bchi_{r}\leftarrow \Bchi_{r}-0.18\Bchi_{\ssstep}$.\tcp*[f]{$O(Nm)$}
            }

        }
    }
    
    \nl\label{alglin:solve-LSIP-dual-postproc-start}%
    $(y_{1,0},\BIy_1^\TRANSP,\BIw_1^\TRANSP,\ldots,y_{N,0},\BIy_N^\TRANSP,\BIw_N^\TRANSP)^\TRANSP\leftarrow \hat{\Bchi}$, where $y_{i,0}\in\R$, $\BIy_i\in\R^{m_i}$, $\BIw_i\in\R^{m_0}$ $\forall i\in\{1,\ldots,N\}$.\tcp*[f]{$O(Nm)$}

    \nl\label{alglin:solve-LSIP-dual-postproc-solution}%
    $\hat{y}_{i,0}\leftarrow y_{i,0}-L_c(\overline{\eta}_i+\overline{\eta}_{0})+c_{\ssmax}$,
    $\hat{\BIy}_i\leftarrow\BIy_i$ 
    $\forall i\in\{1,\ldots,N\}$,
    $\hat{\BIw}_1\leftarrow -\sum_{i=2}^{N}\BIw_{i}$,
    $\hat{\BIw}_i\leftarrow \BIw_{i}$ $\forall i\in\{2,\ldots,N\}$.\tcp*[f]{$O(Nm)$}

    \nl\label{alglin:solve-LSIP-dual-postproc-indexset}%
    $\widehat{\CI}\leftarrow \big\{\iota_{j}:j\in\{1,\ldots,\overline{k}\},\; \iota_j\ne (0,0,0),\; t_{j}=1\big\}$.\tcp*[f]{$O(Nm)$}

    \nl%
    \Return{$\big((\hat{y}_{i,0},\hat{\BIy}_i,\hat{\BIw}_i)_{i=1:N},\widehat{\CI}\big)$.}
    \end{footnotesize}
\end{procedure}

Our detailed procedure for approximately solving (\ref{eqn:mt-tf-lp-full})
is presented in $\SolveLSIPDual$.
It utilizes two auxiliary sub-procedures:
$\GetSepHyp$ is used for computing hyperplanes separating candidate points encountered in 
$\SolveLSIPDual$
and certain superlevel sets of (\ref{eqn:mt-tf-lp-full}) 
(see Lemma~\ref{lem:proc-get-sep-hyp} in Section~\ref{sapx:proof-complexity-lsip} for the properties of $\GetSepHyp$),
whereas $\GetVolumetricQuantities$ 
carries out the computation of quantities related to the volumetric center 
defined in \citep[p.294--295]{ECvaidya1996new}.
Specifically,
Lines~\ref{alglin:solve-LSIP-dual-init-start}--\ref{alglin:solve-LSIP-dual-init-iternum}
of $\SolveLSIPDual$ initialize the parameters in the volumetric center algorithm.
Lines~\ref{alglin:solve-LSIP-dual-forloop}--\ref{alglin:solve-LSIP-dual-rem-end}
of $\SolveLSIPDual$ correspond to the main iteration of the volumetric center algorithm.
We additionally maintain a list $(\iota_r)_{r=1:\overline{k}}$ throughout the iterations,
where each $\iota_r$ is a tuple $(i,j,l)$ indicating that 
the constraint $y_{i,0}+\langle\BIg_i(\BIv_{i,j}),\BIy_i\rangle + \langle\BIg_0(\BIv_{0,l}),\BIw_i\rangle \le c_i(\BIv_{i,j},\BIv_{0,l})-c_{\ssmax}$ 
of (\ref{eqn:mt-tf-lp-full}) had been violated at some point during the iterations.
Lines~\ref{alglin:solve-LSIP-dual-postproc-start}--\ref{alglin:solve-LSIP-dual-postproc-indexset} of $\SolveLSIPDual$
then construct 
an approximate optimizer 
$(\hat{y}_{i,0},\hat{\BIy}_i,\hat{\BIw}_i)_{i=1:N}$ of \eqref{eqn:mt-tf-lsip}
as well as an index set $\widehat{\CI}$ formed by aggregating tuples 
in the list $(\iota_r)_{r=1:\overline{k}}$.
The number of arithmetic operations incurred by $\SolveLSIPDual$ is shown in 
Theorem~\ref{thm:proc-solve-LSIP}\ref{thms:proc-solve-LSIP-dual-complexity},
and the feasibility of the constructed $(\hat{y}_{i,0},\hat{\BIy}_i,\hat{\BIw}_i)_{i=1:N}$ for \eqref{eqn:mt-tf-lsip}
is shown in Theorem~\ref{thm:proc-solve-LSIP}\ref{thms:proc-solve-LSIP-dual-feasibility}.

\begin{procedure}[t]
    \DontPrintSemicolon
    \caption{SolveLSIPPrimal(){$\big(N,(\CG_i)_{i=0:N},(\bar{\BIg}_i,c_i)_{i=1:N},\widehat{\CI},\varsigma\big)$}  $\to$ $(\hat{\theta}_i)_{i=1:N}$}
    \begin{footnotesize}
    \KwInAndOut{\!\!\!\!\!\begin{tabular}{lll}
        \textbf{Input} &
        $N\in\N\cap[2,\infty)$ & the number of agent populations \\ 
        & $\CG_i$ & the test functions on $\CX_i$ encoded by $\big(d_i,m_i,(\BIv_{i,l})_{l=0:m_i},\widetilde{J}_{i},\big(\varrho_{i,\tilde{j}},(\varkappa_{i,\tilde{j},\iota})_{\iota=0:\varrho_{i,\tilde{j}}},$ \\
        && $\varpi_{i,\tilde{j}}\big)_{\tilde{j}=1:\widetilde{J}_{i}},\overline{\eta}_i\big)$ as in \ref{enc:testfuncs}\\
        & $\bar{\BIg}_i\in\R^{m_i}$ & the vector containing the integrals of the test functions in $\CG_i$ with respect to $\mu_i$ \\
        & $c_i:\CX_i\times\CX_0\to\R$ & the cost function on $\CX_i\times\CX_0$; see \ref{enc:costfuncs} \\ 
        & $\widehat{\CI}$ & a set of tuples $(i,j,l)$ where $j\in\{0,1,\ldots,m_i\}$, $l\in\{0,1,\ldots,m_0\}$, $i\in\{1,\ldots,N\}$ \\
        & $\varsigma\in(0,1)$ & the tolerance value for the algorithm of \citet{ECvandenbrand2020deterministic}; see Theorem~\ref{thm:proc-solve-LSIP} \\
        \textbf{Output} &
        $(\hat{\theta}_i)_{i=1:N}$ & a possibly infeasible approximate optimizer of \eqref{eqn:mt-tf-dual} encoded via \ref{enc:par-primal} 
    \end{tabular}\!\!\!\!\!\!\!\!\!\!\!\!\!\!}

    \tcc{In the following computational complexities, let us denote 
    $m:=\max_{0\le i\le N}\{m_i\}$, 
    assume $|\widehat{\CI}|=O(Nm)$, 
    and let $T_{\mathtt{Eval}}$ denote the number of arithmetic operations incurred when evaluating $c_i(\BIx,\BIz)$ at any $(\BIx,\BIz)\in\CX_i\times\CX_0$.}

    \nl%
    $c_{\ssmax}\leftarrow \max_{1\le i\le N,\,0\le j\le m_i,\,0\le l\le m_0}\big\{c_i(\BIv_{i,j},\BIv_{0,l})\big\}$.\tcp*[f]{$O(Nm^2T_{\mathtt{Eval}})$}

    \nl%
    $c_{\ssmin}\leftarrow \min_{1\le i\le N,\,0\le j\le m_i,\,0\le l\le m_0}\big\{c_i(\BIv_{i,j},\BIv_{0,l})\big\}$.\tcp*[f]{$O(Nm^2T_{\mathtt{Eval}})$}

    \nl\label{alglin:solve-LSIP-primal-aggindex}%
    $\widehat{\CI}_{+}\leftarrow \widehat{\CI}\cup \big\{(i,j,0):j\in\{0,1,\ldots,m_i\},\; i\in\{1,\ldots,N\}\big\} \cup \big\{(i,0,l): l\in\{1,\ldots,m_0\},\; i\in\{1,\ldots,N\}\big\}$.\tcp*[f]{$O(N^2m^2)$}

    \nl%
    $q\leftarrow \sum_{i=1}^{N}(1+m_i+m_0)$, $p\leftarrow 0$.\tcp*[f]{$O(N)$} 

    \nl\label{alglin:solve-LSIP-primal-buildconstr}%
    \For(\tcp*[f]{$O(N^2m^2+NmT_{\mathtt{Eval}})$}){each $(i,j,l)\in\widehat{\CI}_{+}$}{

        \nl\label{alglin:solve-LSIP-primal-ineqconstr}%
        $\BIa_{p+1,i}\leftarrow (1,\BIe_{i,j}^\TRANSP,\BIe_{0,l}^\TRANSP)^\TRANSP$, 
        where $\BIe_{i,0}:=\veczero_{m_i}$, 
        $\BIe_{i,j}$ denotes the $j$-th standard basis vector of $\R^{m_i}$ if $j\ne 0$,
        $\BIe_{0,0}:=\veczero_{m_0}$,
        $\BIe_{0,l}$ denotes the $l$-th standard basis vector of $\R^{m_0}$ if $l\ne 0$.\tcp*[f]{$O(m)$}

        \nl%
        $\BIa_{p+1,i'}\leftarrow \veczero_{1+m_{i'}+m_0}$ $\forall i'\in\{1,\ldots,N\}\setminus\{i\}$.\tcp*[f]{$O(Nm)$}

        \nl%
        $\BIa_{p+1}\leftarrow (\BIa_{p+1,1}^\TRANSP,\ldots,\BIa_{p+1,N}^\TRANSP)^\TRANSP$.\tcp*[f]{$O(Nm)$}

        \nl%
        $\tilde{c}_{p+1}\leftarrow c_i(\BIv_{i,j},\BIv_{0,l})-c_{\ssmax}$.\tcp*[f]{$O(T_{\mathtt{Eval}})$}

        \nl\label{alglin:solve-LSIP-primal-ineqconstr-end}%
        $p\leftarrow p+1$.
    }

    \nl%
    \For(\tcp*[f]{$O(Nm^2)$}){$l=1,\ldots,m_0$}{

        \nl\label{alglin:solve-LSIP-primal-balanceconstr}%
        $\BIa_{p+l}\leftarrow -(0,\veczero_{m_1}^\TRANSP,\BIe_{0,l}^\TRANSP,0,\veczero_{m_2}^\TRANSP,\BIe_{0,l}^\TRANSP,\ldots,0,\veczero_{m_N}^\TRANSP,\BIe_{0,l}^\TRANSP)^\TRANSP$, 
        where $\BIe_{0,l}\in\R^{m_0}$ is defined as in Line~\ref{alglin:solve-LSIP-primal-ineqconstr}.\tcp*[f]{$O(Nm)$}
    }

    \nl\label{alglin:solve-LSIP-primal-aggregate}%
    $\BA\leftarrow\left(\begin{smallmatrix}
        | & & | \\
        \BIa_{1} & \cdots & \BIa_{p+m_0} \\
        | & & | 
    \end{smallmatrix}\right)\in\R^{q\times (p+m_0)}$,
    $\bar{\BIg}\leftarrow (1,\bar{\BIg}_1^\TRANSP,\veczero_{m_0}^\TRANSP,\ldots,1,\bar{\BIg}_N^\TRANSP,\veczero_{m_0}^\TRANSP)^\TRANSP$,
    $\tilde{\BIc}\leftarrow (\tilde{c}_1,\ldots,\tilde{c}_{p}, \veczero_{m_0}^\TRANSP)^\TRANSP$.\tcp*[f]{$O(N^2m^2)$}

    \nl\label{alglin:solve-LSIP-primal-solver}%
    $(\Btheta^\TRANSP,\Bxi^\TRANSP)^\TRANSP\leftarrow$ the output of \citep[Algorithm~2]{ECvandenbrand2020deterministic} with input $\big(\BA,\bar{\BIg},\tilde{\BIc},\frac{\varsigma}{4(N+1)(p+1)}\big)$,
    where $\Btheta\in\R^{p}$, $\Bxi\in\R^{m_0}$.
    \tcp*[f]{$O\big(N^{\omega}m^{\omega}\log(Nm)^2\log\big(\frac{Nm}{\varsigma}\big)\big)$}

    \nl%
    $(\theta_{i,j,l})_{(i,j,l)\in\widehat{\CI}_{+}}^{\TRANSP}\leftarrow \Btheta$.\tcp*[f]{$O(Nm)$}

    \nl%
    \For(\tcp*[f]{$O(N^2m)$}){$i=1,\ldots,N$}{

        \nl\label{alglin:solve-LSIP-primal-gatherindices}%
        $\widehat{\CI}_{i}\leftarrow \big\{(j,l):(i,j,l)\in\widehat{\CI}_{+}\big\}$.\tcp*[f]{$O(Nm)$}

        \nl\label{alglin:solve-LSIP-primal-normalize}%
        $\hat{\theta}_{i,j,l}\leftarrow \big(\sum_{(j',l')\in\widehat{\CI}_{i}}\theta_{i,j',l'}\big)^{-1}\theta_{i,j,l}$
        $\forall (j,l)\in\widehat{\CI}_i$.\tcp*[f]{$O(Nm)$}

        \nl\label{alglin:solve-LSIP-primal-buildmeas}%
        Let $\hat{\theta}_i$ be encoded by $\big(|\widehat{\CI}_i|,(j,l,\hat{\theta}_{i,j,l})_{(j,l)\in\widehat{\CI}_i}\big)$ as in \ref{enc:par-primal}.
    }

    \nl%
    \Return{$(\hat{\theta}_i)_{i=1:N}$.}
    \end{footnotesize}
\end{procedure}

Subsequently, using the index set $\widehat{\CI}$ constructed and returned by $\SolveLSIPDual$,
the procedure $\SolveLSIPPrimal$ utilizes 
\citep[Algorithm~2]{ECvandenbrand2020deterministic}
to approximately solve the dual LP problem of (\ref{eqn:mt-tf-lp-full}),
which results in a possibly infeasible approximate optimizer of \eqref{eqn:mt-tf-dual}.
The number of arithmetic operations incurred by $\SolveLSIPPrimal$
is shown in Theorem~\ref{thm:proc-solve-LSIP}\ref{thms:proc-solve-LSIP-primal-complexity},
whereas Theorem~\ref{thm:proc-solve-LSIP}\ref{thms:proc-solve-LSIP-primal-feasibility}
and Theorem~\ref{thm:proc-solve-LSIP}\ref{thms:proc-solve-LSIP-gap}
show that 
$(\hat{\theta}_i)_{i=1:N}$,
$(\hat{y}_{i,0},\hat{\BIy}_i,\hat{\BIw}_i)_{i=1:N}$
computed by $\SolveLSIPDual$ and $\SolveLSIPPrimal$
constitute a pair of 
$(\epsilon_{\sspar},\varsigma)$-optimizers 
of \eqref{eqn:mt-tf-dual} and \eqref{eqn:mt-tf-lsip},
which is required by Theorem~\ref{thm:mt-tf-approx}
for the construction of an AME.

We would like to remark that we do not use algorithms for LP problems such as 
\citep[Algorithm~2]{ECvandenbrand2020deterministic} to directly solve 
(\ref{eqn:mt-tf-lp-full}) and its dual LP problem 
due to the observation that 
(\ref{eqn:mt-tf-lp-full}) has 
$\sum_{i=1}^{N}(1+m_i+m_0)=O(Nm)$ decision variables 
and 
$m_0+\sum_{i=1}^{N}(1+m_i)(1+m_0)=O(Nm^2)$ constraints,
where $m$ denotes $\max_{0\le i\le N}\{m_i\}$.
Thus, 
when $m$ is much greater than~$N$,
the number of constraints is much greater than the number of decision variables, which motivates us to use the volumetric center algorithm 
of \citet{ECvaidya1996new}.
As explained in the proof of Theorem~\ref{thm:proc-solve-LSIP} 
in Section~\ref{sapx:proof-complexity-lsip},
the volumetric center algorithm 
of \citet{ECvaidya1996new}
produces a relaxation of (\ref{eqn:mt-tf-lp-full}) with 
$O(Nm)$ constraints whose optimal value is close to that of (\ref{eqn:mt-tf-lp-full}).
Therefore, the use of the volumetric center algorithm 
of \citet{ECvaidya1996new}
allows us to reduce the computational complexity of solving 
\eqref{eqn:mt-tf-dual} and \eqref{eqn:mt-tf-lsip},
which is the most costly part of Algorithm~\ref{alg:complexity-overall}.

Now, let us present below the theorem that establishes the properties 
of $\SolveLSIPDual$ and $\SolveLSIPPrimal$.
\begin{theorem}\label{thm:proc-solve-LSIP}
    Under Setting~\ref{sett:Euclidean},
    let $\epsilon_{0}\in(0,1)$, $\varsigma\in(0,1)$ be arbitrary.
    For $i=0,1,\ldots,N$,
    let 
    $\CG_i$ be encoded by
    $\big(d_i,m_i,(\BIv_{i,l})_{l=0:m_i},\allowbreak \widetilde{J}_i,\allowbreak\big(\varrho_{i,\tilde{j}},\allowbreak(\varkappa_{i,\tilde{j},\iota})_{\iota=0:\varrho_{i,\tilde{j}}},\varpi_{i,\tilde{j}}\big)_{\tilde{j}=1:\widetilde{J}_i},\overline{\eta}_i\big)$ as in \ref{enc:testfuncs},
    where $\overline{\eta}_i\ge \eta(\FS_i)$.
    Recall that $L_{c}>0$ is the Lipschitz constant of the cost functions $(c_i)_{i=1:N}$.
    Moreover, let $\big((\hat{y}_{i,0},\hat{\BIy}_i,\hat{\BIw}_i)_{i=1:N},\widehat{\CI}\big)$ be the output of 
    $\SolveLSIPDual\big(N,\allowbreak(\CG_i)_{i=0:N},(\bar{\BIg}_i,c_i)_{i=1:N},\epsilon_0,L_c\big)$, and
    subsequently
    let $(\hat{\theta}_i)_{i=1:N}$ be the output of 
    $\SolveLSIPPrimal\big(N,\allowbreak(\CG_i)_{i=0:N},(\bar{\BIg}_i,c_i)_{i=1:N},\widehat{\CI},\varsigma\big)$.
    Furthermore, 
    let us denote $m:=\max_{0\le i\le N}\{m_i\}$,
    $c_{\ssmin}:=\min_{1\le i\le N,\,0\le j\le m_i,\,0\le l\le m_0}\big\{c_i(\BIv_{i,j},\BIv_{0,l})\big\}$,
    $c_{\ssmax}:=\nobreak\max_{1\le i\le N,\,0\le j\le m_i,\,0\le l\le m_0}\big\{c_i(\BIv_{i,j},\BIv_{0,l})\big\}$,
    $\bar{g}_{\ssmin}:=\min_{1\le i\le N,\,0\le j\le m_i}\big\{\int_{\CX_i}g_{i,j}\DIFFX{\mu_i}\big\}$,
    and let $T_{\mathtt{Eval}}$ denote the number of arithmetic operations incurred when evaluating $c_i(\BIx,\BIz)$ at any $(\BIx,\BIz)\in\nobreak\CX_i\times\nobreak\CX_0$.
    Then, the following statements hold.\widowpenalty-5000
    \begin{enumerate}[label=(\roman*),beginpenalty=10000]
        \item\label{thms:proc-solve-LSIP-dual-complexity}
        It holds that $\bar{g}_{\ssmin}>0$, and that
        $\SolveLSIPDual\big(N,\allowbreak(\CG_i)_{i=0:N},(\bar{\BIg}_i,c_i)_{i=1:N},\epsilon_0,L_c\big)$ incurs
        $O\big((N^{\omega+1}m^{\omega+1}+N^2m^3T_{\mathtt{Eval}})\log\big(\frac{Nm(c_{\ssmax}-c_{\ssmin}+1)}{\epsilon_0\bar{g}_{\ssmin}}\big)\big)$
        arithmetic operations.
        Moreover, it holds that $|\widehat{\CI}|=O(Nm)$.

        \item\label{thms:proc-solve-LSIP-primal-complexity}
        $\SolveLSIPPrimal\big(N,(\CG_i)_{i=0:N},(\bar{\BIg}_i,c_i)_{i=1:N},\widehat{\CI},\varsigma\big)$
        incurs $O\big(N^{\omega}m^{\omega}\log(Nm)^2\log\big(\frac{Nm}{\varsigma}\big)+Nm^2T_{\mathtt{Eval}}\big)$
        arithmetic operations.

        \item\label{thms:proc-solve-LSIP-dual-feasibility}
        $(\hat{y}_{i,0},\hat{\BIy}_i,\hat{\BIw}_i)_{i=1:N}$ is feasible for \eqref{eqn:mt-tf-lsip}.

        \item\label{thms:proc-solve-LSIP-primal-feasibility}
        For $i=1,\ldots,N$,
        let $\hat{\mu}_i$ and $\hat{\nu}_i$ denote the marginals of 
        $\hat{\theta}_i$ on $\CX_i$ and $\CX_0$, respectively.
        Then, it holds for $i=1,\ldots,N$ that 
        $\support(\hat{\theta}_i)\subseteq V(\FS_i)\times V(\FS_0)$,
        and that 
        $\bar{\mu}_i\abovebelowset{\CG_i}{\varsigma}{\sim}\mu_i$, 
        $\bar{\nu}_i\abovebelowset{\CG_0}{\varsigma}{\sim}\bar{\nu}_1$.

        \item\label{thms:proc-solve-LSIP-gap}
        It holds that 
        \begin{align*}
            \Bigg(\sum_{i=1}^{N}\int_{\CX_i\times\CX_0}c_i\DIFFX{\hat{\theta}_i}\Bigg)
            -
            \Bigg(\sum_{i=1}^{N}\hat{y}_{i,0}+\langle\bar{\BIg}_i,\hat{\BIy}_i\rangle\Bigg)
            &\le \epsilon_{0} + (c_{\ssmax}-c_{\ssmin})\varsigma + L_{c}\Bigg(N\overline{\eta}_0+\sum_{i=1}^{N}\overline{\eta}_i\Bigg).
        \end{align*}

    \end{enumerate}
    Consequently, 
    $(\hat{\theta}_i)_{i=1:N}$,
    $(\hat{y}_{i,0},\hat{\BIy}_i,\hat{\BIw}_i)_{i=1:N}$
    is a pair of 
    $(\epsilon_{\sspar},\varsigma)$-optimizers of 
    \eqref{eqn:mt-tf-dual} and \eqref{eqn:mt-tf-lsip},
    where 
    $\epsilon_{\sspar}:=\epsilon_{0} + (c_{\ssmax}-c_{\ssmin})\varsigma + L_{c}\big(N\overline{\eta}_0+\sum_{i=1}^{N}\overline{\eta}_i\big)$.
\end{theorem}

\proof{Proof of Theorem~\ref{thm:proc-solve-LSIP}.}
    See Section~\ref{sapx:proof-complexity-lsip}.
\endproof

With the pair of 
$(\epsilon_{\sspar},\varsigma)$-optimizers 
$(\hat{\theta}_i)_{i=1:N}$, $(\hat{y}_{i,0},\hat{\BIy}_i,\hat{\BIw}_i)_{i=1:N}$
of \eqref{eqn:mt-tf-dual} and \eqref{eqn:mt-tf-lsip}
computed by 
$\SolveLSIPDual$ and $\SolveLSIPPrimal$,
the next subsection will introduce the procedures to construct an AME via Theorem~\ref{thm:mt-tf-approx}.

\subsection{Construction of an approximate matching equilibrium (AME)}\label{sapx:complexity-AME}

\begin{procedure}[t]
    \DontPrintSemicolon
    \caption{ConstructTransFuncs(){$\big(N,\CG_0,(\hat{\BIw}_i)_{i=1:N}\big)$} $\to$ $(\hat{\varphi}_i)_{i=1:N}$}
    \begin{footnotesize}
    \KwInAndOut{\!\!\!\!\!\begin{tabular}{lll}
        \textbf{Input} &
        $N\in\N\cap[2,\infty)$ & the number of agent populations \\ 
        & $\CG_0$ & the test functions on $\CX_0$ encoded via \ref{enc:testfuncs} \\
        & $\hat{\BIw}_i\in\R^{m_0}$ & the part of an approximate optimizer of \eqref{eqn:mt-tf-lsip} containing the coefficients with \\
        && respect to $\CG_0$ \\ 
        \textbf{Output} &
        $\hat{\varphi}_i:\CX_0\to\R$ & transfer function in an $\epsilon$-matching equilibrium encoded via \ref{enc:transfuncs}
    \end{tabular}\!\!\!\!\!\!\!\!\!\!\!\!\!\!}

    \nl%
    \For(\tcp*[f]{$O(Nm_0)$}){$i=1,\ldots,N$}{

        \nl%
        $(\hat{w}_{i,1},\ldots,\hat{w}_{i,m_0})^\TRANSP\leftarrow \hat{\BIw}_i$.\tcp*[f]{$O(m_0)$}

        \nl%
        Let $\hat{\varphi}_i$ be encoded by $\big(\CG_0,(\hat{w}_{i,l})_{l=1:m_0}\big)$.
    }

    \nl%
    \Return $(\hat{\varphi}_i)_{i=1:N}$.

    \end{footnotesize}
\end{procedure}

First, the construction of the approximately optimal transfer functions
$(\hat{\varphi}_i)_{i=1:N}$ from 
an approximate optimizer $(\hat{y}_{i,0},\hat{\BIy}_i,\hat{\BIw}_i)_{i=1:N}$
of \eqref{eqn:mt-tf-lsip} 
in Theorem~\ref{thm:mt-tf-approx} is straightforward,
as it amounts to interpreting the components of the vector $\hat{\BIw}_i$
as coefficients of the test functions $(g_{0,l})_{l=1:m_0}$ for $i=1,\ldots,N$,
which incurs no actual computation.
This is presented in the procedure $\ConstructTransFuncs$.

\begin{procedure}[t]
    \DontPrintSemicolon
    \caption{GetCoords(){$(C,\BIx)$} $\to$ $\big(\mathtt{Inside},(h_l)_{l=0:n}\big)$}
    \begin{footnotesize}
    \KwInAndOut{\!\!\!\!\!\begin{tabular}{lll}
        \textbf{Input} &
        $C$ & an $n$-simplex in $\R^{d}$ encoded by $\big(d,n,(\BIv_l)_{l=0:n}\big)$ as in \ref{enc:simplex} \\
        & $\BIx\in \R^d$ & any point to be examined \\ 
        \textbf{Output} &
        $\mathtt{Inside}$ & a boolean value indicating whether $\BIx\in C$ \\
        & $h_l\in [0,1]$ & if $\mathtt{Inside}=\mathtt{True}$, then $h_l$ is equal to the barycentric coordinates $\lambda_{\FF(C),\BIv_l}(\BIx)$ of $\BIx$ with \\
        &&respect to $\BIv_l$; if $\mathtt{Inside}=\mathtt{False}$, then $h_l=0$
    \end{tabular}\!\!\!\!\!\!\!\!\!\!\!\!\!\!}

    \nl%
    $\BV\leftarrow\left(\begin{smallmatrix}
        | & | & & | \\
        \BIv_0 & \BIv_1 & \cdots & \BIv_n \\
        | & | & & | \\
        1 & 1 & \cdots & 1
    \end{smallmatrix}\right)\in \R^{(d+1)\times (n+1)}$.\tcp*[f]{$O(d^2)$}

    \nl\label{alglin:get-coords-linear-system}%
    $\vphantom{\Big(\Big)}{\Blambda\leftarrow (\BV^\TRANSP\BV)}^{-1}\BV^\TRANSP\big(\begin{smallmatrix}
        \BIx \\ 1
    \end{smallmatrix}\big)$.\tcp*[f]{$O(d^\omega)$}

    \nl%
    \uIf(\tcp*[f]{$O(d^2)$}){$\BV\Blambda=\big(\begin{smallmatrix}
        \BIx \\ 1
    \end{smallmatrix}\big)$ and $\Blambda\in\R^{n+1}_+$}{
        
        \nl%
        $\mathtt{Inside}\leftarrow \mathtt{True}$, $(h_0,h_1,\ldots,h_n)^\TRANSP\leftarrow \Blambda$.\tcp*[f]{$O(d)$}
    }
    \nl%
    \Else{
        \nl%
        $\mathtt{Inside}\leftarrow \mathtt{False}$, $(h_0,h_1,\ldots,h_n)^\TRANSP\leftarrow \veczero_{n+1}$.\tcp*[f]{$O(d)$}
    }

    \nl%
    \Return{$\big(\mathtt{Inside},(h_l)_{l=0:n}\big)$.}
    \end{footnotesize}
\end{procedure}

\begin{procedure}[t]
    \DontPrintSemicolon
    \caption{EvalTransFunc(){$(\hat{\varphi}_i,\BIz)$} $\to$ $\varphi$}
    \begin{footnotesize}
    \KwInAndOut{\!\!\!\!\!\begin{tabular}{lll}
        \textbf{Input} &
        $\hat{\varphi}_i$ & the transfer function encoded by $\big(\CG_0,(w_{i,l})_{l=1:m_0}\big)$ as in \ref{enc:transfuncs},
        where $\CG_0$ is encoded by \\
        && $\big(d_0,m_0,(\BIv_{0,l})_{l=0:m_0},\widetilde{J}_0,\big(\varrho_{0,\tilde{j}},\allowbreak (\varkappa_{0,\tilde{j},\iota})_{\iota=0:\varrho_{0,\tilde{j}}},\varpi_{0,\tilde{j}}\big)_{\tilde{j}=1:\widetilde{J}_0},\overline{\eta}_0\big)$ as in \ref{enc:testfuncs} \\
        & $\BIz\in \CX_0$ & any point at which the transfer function is evaluated \\ 
        \textbf{Output} &
        $\varphi\in\R$ & the value of $\hat{\varphi}_i(\BIz)$
    \end{tabular}\!\!\!\!\!\!\!\!\!\!\!\!\!\!}

    \nl%
    $w_{i,0}\leftarrow 0$.

    \nl%
    \For(\tcp*[f]{$O(\widetilde{J}_0d_0^{\omega})$}){$\tilde{j}=1,\ldots,\widetilde{J}_0$}{

        \nl%
        Let $S_{0,\tilde{j}}$ be a simplex encoded by $(d_0,\varrho_{0,\tilde{j}},(\BIv_{0,\varkappa_{0,\tilde{j},\iota}})_{\iota=0:\varrho_{0,\tilde{j}}})$ as in \ref{enc:simplex}.\tcp*[f]{$O(d_0^2)$}

        \nl\label{alglin:eval-trans-func-get-coord}%
        $\big(\mathtt{Inside},(h_\iota)_{\iota=0:\varrho_{0,\tilde{j}}}\big)\leftarrow \GetCoords(S_{0,\tilde{j}},\BIz)$.\tcp*[f]{$O(d_0^\omega)$}

        \nl%
        \If{$\mathtt{Inside}=\mathtt{True}$}{

            \nl%
            $\varphi\leftarrow \sum_{\iota=0}^{\varrho_{0,\tilde{j}}}w_{i,\varkappa_{0,\tilde{j},\iota}}h_{\iota}$.\tcp*[f]{$O(d_0)$}

            \nl%
            \Return $\varphi$.
        }
    }

    \end{footnotesize}
\end{procedure}

After the transfer functions $(\hat{\varphi}_i)_{i=1:N}$ have been constructed,
each $\hat{\varphi}_i$ can be evaluated at any $\BIz\in\CX_0$ via the procedure $\EvalTransFunc$.
$\EvalTransFunc$ utilizes an auxiliary sub-procedure $\GetCoords$
for computing the barycentric coordinates 
(see Lemma~\ref{lem:proc-get-coords} in Section~\ref{sapx:proof-complexity-AME}
for the properties of $\GetCoords$).
We present in Lemma~\ref{lem:proc-construct-trans-funcs}
the properties of $\ConstructTransFuncs$ and $\EvalTransFunc$.

\begin{lemma}\label{lem:proc-construct-trans-funcs}
    Under Setting~\ref{sett:Euclidean}, let 
    $(\hat{y}_{i,0},\hat{\BIy}_i,\hat{\BIw}_i)_{i=1:N}$ be an arbitrary 
    feasible solution of \eqref{eqn:mt-tf-lsip}.
    Then, $\ConstructTransFuncs\big(\CG_0,\allowbreak(\hat{\BIw}_i)_{i=1:N}\big)$
    incurs $O(Nm_0)$ arithmetic operations and returns 
    $(\hat{\varphi}_i)_{i=1:N}$,
    where for $i=1,\ldots,N$, 
    $\hat{\varphi}_i(\BIz)=\langle\BIg_0(\BIz),\hat{\BIw}_i\rangle$ 
    $\forall \BIz\in\CX_0$.
    Subsequently,
    for any $i\in\nobreak\{1,\ldots,N\}$ and any $\BIz\in\CX_0$,
    $\EvalTransFunc(\hat{\varphi}_i,\BIz)$ incurs 
    $O\big(|\FM(\FS_0)|d_0^{\omega}\big)$ arithmetic operations and returns 
    $\hat{\varphi}_i(\BIz)$.
\end{lemma}

\proof{Proof of Lemma~\ref{lem:proc-construct-trans-funcs}.}
    See Section~\ref{sapx:proof-complexity-AME}.
\endproof

\begin{procedure}[t]
    \DontPrintSemicolon
    \caption{NormalizeRows(){$(n,n',\BP)$} $\to$ $\overrightarrow{\BP}$}
    \begin{footnotesize}
    \KwInAndOut{\!\!\!\!\!\begin{tabular}{lll}
        \textbf{Input} &
        $n,n'\in\N$ & the dimensions of the input matrix $\BP$ \\
        & $\BP\in\R_{+}^{n\times n'}$ & a matrix in which all entries sum to~1; its entries are denoted by $(p_{i,j})_{i=1:n,\,j=1:n'}$ \\
        \textbf{Output} &
        $\overrightarrow{\BP}\in\R_{+}^{n\times n'}$ & a matrix with each row summing to~1
    \end{tabular}\!\!\!\!\!\!\!\!\!\!\!\!\!\!}

    \nl%
    \For(\tcp*[f]{$O(nn')$}){$i=1,\ldots,n$}{

        \nl%
        \uIf(\tcp*[f]{$O(n')$}){$(p_{i,1},\ldots,p_{i,n'})^\TRANSP\ne\veczero_{n'}$}{

            \nl\label{alglin:normalize-rows-normalize}%
            $(q_{i,1},\ldots,q_{i,n'})^\TRANSP \leftarrow \big(\sum_{j=1}^{n'}p_{i,j}\big)^{-1}(p_{i,1},\ldots,p_{i,n'})^\TRANSP$.\tcp*[f]{$O(n')$}
        }
        \nl%
        \Else{
            \nl\label{alglin:normalize-rows-degenerate}%
            $(q_{i,1},\ldots,q_{i,n'})^\TRANSP \leftarrow \frac{1}{n'}\vecone_{n'}$.\tcp*[f]{$O(n')$}
        }
    }

    \nl%
    $\overrightarrow{\BP}\leftarrow(q_{i,j})_{i=1:n,\,j=1:n'}$.\tcp*[f]{$O(nn')$}

    \nl%
    \Return{$\overrightarrow{\BP}$.}
    \end{footnotesize}
\end{procedure}

\begin{procedure}
    \DontPrintSemicolon
    \caption{ConstructAllocs(){$\big(N,(\CX_i,\CG_i)_{i=0:N},(\mu_i,\mathtt{Aux}_i,\hat{\theta}_i)_{i=1:N}\big)$} $\to$ $\big((\hat{\gamma}_i)_{i=1:N},\hat{\nu}\big)$\!\!\!\!}
    \begin{footnotesize}
    \KwInAndOut{\!\!\!\!\!\begin{tabular}{lll}
        \textbf{Input} &
        $N\in\N\cap[2,\infty)$ & the number of agent populations \\
        & $\CX_i$ & the type space encoded via \ref{enc:spaces} \\
        & $\CG_i$ & the test functions on $\CX_i$ encoded by $\big(d_i,m_i,(\BIv_{i,l})_{l=0:m_i},\widetilde{J}_{i},\big(\varrho_{i,\tilde{j}},(\varkappa_{i,\tilde{j},\iota})_{\iota=0:\varrho_{i,\tilde{j}}},$ \\
        &&$\varpi_{i,\tilde{j}}\big)_{\tilde{j}=1:\widetilde{J}_{i}},\overline{\eta}_i\big)$ as in \ref{enc:testfuncs} \\
        & $\mu_i\in\CP(\CX_i)$ & the type measure on $\CX_i$; see \ref{enc:typemeasures} \\
        & $\mathtt{Aux}_i$ & the auxiliary data structure in the output of $\InitTypeMeasure(\CX_i,\mu_i,\CG_i)$ encoded \\
        && by $\big(\widetilde{R}_i,(F_{i,\tilde{r}},\tilde{\varpi}_{i,\tilde{r}})_{\tilde{r}=1:\widetilde{R}_i},\BP_i\big)$ as in \ref{enc:typemeasures-aux} \\
        & $\hat{\theta}_i$ & a discrete probability measure with support in $V(\FS_i)\times V(\FS_{0})$ encoded by \\
        && $\big(n_i,(j_{i,\iota},l_{i,\iota},\hat{\theta}_{i,\iota})_{\iota=1:n_i}\big)_{i=1:N}$ as in \ref{enc:par-primal} \\
        \textbf{Output} &
        $\hat{\gamma}_i\in\CP(\CX_i\times\CX_0)$ & coupling in team allocations encoded via \ref{enc:allocations} \\
        & $\hat{\nu}\in\CP(\CX_0)$ & quality measure in team allocations encoded via \ref{enc:qualitymeasure}
    \end{tabular}\!\!\!\!\!\!\!\!\!\!\!\!\!\!}

    \tcc{Let us denote $m:=\max_{0\le i\le N}\{m_i\}$, $\widetilde{R}:=\max_{1\le i\le N}\{\widetilde{R}_i\}$ in the following computational complexities.}

    \nl%
    \For(\tcp*[f]{$O\big(N(\widetilde{R}m^{\omega-1}+m^{\omega})\big)$}){$i=1,\ldots,N$}{

        \nl\label{alglin:construct-allocs-begin}%
        $\tilde{\theta}_{i,j,l}\leftarrow 0$ $\forall j\in\{0,1,\ldots,m_i\}$, 
        $\forall l\in\{0,1,\ldots,m_0\}$.\tcp*[f]{$O(m^2)$}

        \nl\label{alglin:construct-allocs-assignmatrix}%
        $\tilde{\theta}_{i,j_{i,\iota},l_{i,\iota}}\leftarrow \hat{\theta}_{i,\iota}$ $\forall \iota\in\{1,\ldots,n_i\}$.\tcp*[f]{$O(m^2)$}

        \nl%
        $\BTheta_i\leftarrow (\tilde{\theta}_{i,j,l})_{j=0:m_i,\,l=0:m_0}\in\R_{+}^{(m_i+1)\times(m_0+1)}$.\tcp*[f]{$O(m^2)$}

        \nl\label{alglin:construct-allocs-margvec}%
        $\Bmu_i\leftarrow \BP_i\vecone_{\widetilde{R}_i}$,
        $\hat{\Bmu}_i \leftarrow \BTheta_i\vecone_{m_0+1}$,
        $\hat{\Bnu}_i \leftarrow \BTheta_i^\TRANSP\vecone_{m_i+1}$.\tcp*[f]{$O(\widetilde{R}m+m^2)$}

        \nl\label{alglin:construct-allocs-coup-start}%
        \uIf(\tcp*[f]{$O(m)$}){$\hat{\Bmu}_i=\Bmu_i$}{

            \nl\label{alglin:construct-allocs-coup-same-first}%
            $\BA_i\leftarrow \diag(\hat{\Bmu}_i)$.\tcp*[f]{$O(m^2)$}
        }
        \nl%
        \Else{
            \nl\label{alglin:construct-allocs-coup-diff-first-1}%
            $\hat{\Bmu}_{i,\ssmin}\leftarrow \hat{\Bmu}_i \wedge \Bmu_i$ 
            (here $\wedge$ denotes component-wise minimum).\tcp*[f]{$O(m)$}

            \nl\label{alglin:construct-allocs-coup-diff-first-2}%
            $\BA_i\leftarrow \diag(\hat{\Bmu}_{i,\ssmin})+\big(1-\langle\hat{\Bmu}_{i,\ssmin},\vecone_{m_i+1}\rangle\big)^{-1}(\hat{\Bmu}_i-\hat{\Bmu}_{i,\ssmin})(\Bmu_{i}-\hat{\Bmu}_{i,\ssmin})^\TRANSP$.\tcp*[f]{$O(m^2)$}
        }

        \nl%
        \uIf(\tcp*[f]{$O(m)$}){$\hat{\Bnu}_1 = \hat{\Bnu}_i$}{

            \nl\label{alglin:construct-allocs-coup-same-second}%
            $\BB_i\leftarrow \diag(\hat{\Bnu}_1)$.\tcp*[f]{$O(m^2)$}

        }
        \nl%
        \Else{

            \nl\label{alglin:construct-allocs-coup-diff-second-1}%
            $\hat{\Bnu}_{i,\ssmin}\leftarrow \hat{\Bnu}_1 \wedge \hat{\Bnu}_i$ 
            (here $\wedge$ denotes component-wise minimum).\tcp*[f]{$O(m)$}

            \nl\label{alglin:construct-allocs-coup-end}%
            $\BB_i\leftarrow \diag(\hat{\Bnu}_{i,\ssmin})+\big(1-\langle\hat{\Bnu}_{i,\ssmin},\vecone_{m_0+1}\rangle\big)^{-1}(\hat{\Bnu}_1-\hat{\Bnu}_{i,\ssmin})(\hat{\Bnu}_{i}-\hat{\Bnu}_{i,\ssmin})^\TRANSP$.\tcp*[f]{$O(m^2)$}
        }

        \nl\label{alglin:construct-allocs-normalize-rows}%
        $\overleftarrow{\BTheta}_i\leftarrow\NormalizeRows(m_0+1,m_1+1,\BTheta_i^\TRANSP)$,
        $\overrightarrow{\BP}_i\leftarrow\NormalizeRows(m_i+1,\widetilde{R}_i,\BP_i)$,
        $\overrightarrow{\BA}_i\leftarrow\NormalizeRows(m_i+1,m_i+1,\BA_i)$,
        $\overrightarrow{\BB}_i\leftarrow\NormalizeRows(m_0+1,m_0+1,\BB_i)$.\tcp*[f]{$O(\widetilde{R}m+m^2)$}

        \nl\label{alglin:construct-allocs-matprod}%
        $\BT_i\leftarrow \overrightarrow{\BB}_{i}\overleftarrow{\BTheta}_i\overrightarrow{\BA}_i\overrightarrow{\BP}_i$.\tcp*[f]{$O(\widetilde{R}m^{\omega-1}+m^{\omega})$}

        \nl%
        \If{$i=1$}{

            \nl%
            $(\hat{\nu}_{1,0},\hat{\nu}_{1,1},\ldots,\hat{\nu}_{1,m_0})^\TRANSP\leftarrow \hat{\Bnu}_1$.\tcp*[f]{$O(m)$}

            \nl\label{alglin:construct-allocs-qualitymeasure}%
            Let $\hat{\nu}$ be encoded by 
            $\big(m_0+1,d_0,(\BIv_{0,l},\hat{\nu}_{1,l})_{l=0:m_0}\big)$ as in \ref{enc:qualitymeasure}.
        }

        \nl%
        Let $\hat{\gamma}_i$ be encoded by 
        $(\hat{\nu},\CX_i,\mu_i,\mathtt{Aux}_i,\BT_i)$ as in \ref{enc:allocations}.
    }

    \nl%
    \Return{$\big((\hat{\gamma}_i)_{i=1:N},\hat{\nu}\big)$.}
    \end{footnotesize}
\end{procedure}

Next, the procedure for the construction of approximately optimal 
team allocations $(\hat{\gamma}_i)_{i=1:N}$ and quality measure $\hat{\nu}$
is presented in $\ConstructAllocs$.
It utilizes an auxiliary sub-procedure $\NormalizeRows$
for normalizing each row of an input matrix to sum to~1, while 
handling the edge-case where an entire row of the input matrix consists of zeros (see Lemma~\ref{lem:proc-normalize-rows} in Section~\ref{sapx:proof-complexity-AME} for the properties of $\NormalizeRows$).

\begin{procedure}
    \DontPrintSemicolon
    \caption{SampQual(){$\big(\hat{\nu},(U_t)_{t\in\N}\big)$} $\to$ $(\BIz,\overline{t})$}
    \begin{footnotesize}
    \KwInAndOut{\!\!\!\!\!\begin{tabular}{lll}
        \textbf{Input} &
        $\hat{\nu}\in\CP(\CX_0)$ & quality measure in team allocations encoded by $\big(n_0,d_0,(\BIz_{0,l},\hat{\nu}_{0,l})_{l=0:n_0-1}\big)$ as in \ref{enc:qualitymeasure} \\
        & $(U_{t})_{t\in\N}\subset [0,1]$\! & an infinite sequence of independent realizations of $\mathrm{Uniform}[0,1]$ random variable \\
        \textbf{Output} &
        $\BIz\in\R^{d_0}$ & the generated random sample of quality \\
        & $\overline{t}\in\N_0$ & the number of realizations in $(U_t)_{t\in\N}$ consumed in the procedure
    \end{tabular}\!\!\!\!\!\!\!\!\!\!\!\!\!\!}

    \nl%
    $Q\leftarrow \min\big\{l\in\{0,1,\ldots,n_0-1\}:\sum_{s=0}^{l}\hat{\nu}_{0,s}\ge U_1\big\}$.\tcp*[f]{$O(n_0)$}

    \nl%
    \Return{$(\BIz_{0,Q},1)$.}
    \end{footnotesize}
\end{procedure}

\begin{procedure}[t]
    \DontPrintSemicolon
    \caption{SampAlloc(){$\big(\hat{\gamma}_i,(U_t)_{t\in\N}\big)$} $\to$ $(\BIx,\BIz,\overline{t})$}
    \begin{footnotesize}
    \KwInAndOut{\!\!\!\!\!\begin{tabular}{lll}
        \textbf{Input} &
        $\hat{\gamma}_i\in\CP(\CX_i\times\CX_0)$ & coupling in team allocations encoded by
        $(\hat{\nu},\CX_i,\mu_i,\mathtt{Aux}_i,\BT_i)$ as in \ref{enc:allocations}, where $\hat{\nu}$ is \\
        && encoded by $\big(n_0,d_0,(\BIz_{0,l},\hat{\nu}_{0,l})_{l=0:n_0-1}\big)$ as in \ref{enc:qualitymeasure},
        $\mathtt{Aux}_i$ is encoded by \\
        && $\big(\widetilde{R}_i,(F_{i,\tilde{r}},\tilde{\varpi}_{i,\tilde{r}})_{\tilde{r}=1:\widetilde{R}_i},\BP_i\big)$ as in \ref{enc:typemeasures-aux},
        $\mu_i$ is associated with the oracle procedures \\
        && $\mathtt{Int}$ and $\mathtt{Samp}$ in \ref{enc:typemeasures}, and the entries of $\BT_i$ are denoted by $(t_{l,\tilde{r}})_{l=0:n_0-1,\,\tilde{r}=1:\widetilde{R}_i}$ \\
        & $(U_{t})_{t\in\N}\subset [0,1]$ & an infinite sequence of independent realizations of $\mathrm{Uniform}[0,1]$ random variable \\
        \textbf{Output} &
        $\BIx\in\R^{d_i}$ & the generated random sample of type \\
        & $\BIz\in\R^{d_0}$ & the generated random sample of quality \\
        & $\overline{t}\in\N_0$ & the number of realizations in $(U_t)_{t\in\N}$ consumed in the procedure
    \end{tabular}\!\!\!\!\!\!\!\!\!\!\!\!\!\!}
    \tcc{Let $T_{\mathtt{Samp}}$ denote the number of arithmetic operations incurred by $\mathtt{Samp}$.}

    \nl\label{alglin:samp-alloc-qualityindex}%
    $Q\leftarrow \min\big\{l\in\{0,1,\ldots,n_0-1\}:\sum_{s=0}^{l}\hat{\nu}_{0,s}\ge U_1\big\}$.\tcp*[f]{$O(n_0)$}

    \nl\label{alglin:samp-alloc-faceindex}%
    $S\leftarrow \min\big\{\tilde{r}\in\{1,\ldots,\widetilde{R}_i\}:\sum_{\tilde{l}=1}^{\tilde{r}}t_{Q,\tilde{l}}\ge U_2\big\}$.\tcp*[f]{$O(\widetilde{R}_i)$}

    \nl\label{alglin:samp-alloc-condsamp}%
    $(\BIx,\overline{t}_0)\leftarrow \mathtt{Samp}\big(\CX_i,\mu_i,F_{i,S},\tilde{\varpi}_{i,S},(U_t)_{t\ge 3}\big)$.\tcp*[f]{$O(T_{\mathtt{Samp}})$}

    \nl%
    \Return{$(\BIx,\BIz_{0,Q},\overline{t}_0+2)$.}
    \end{footnotesize}
\end{procedure}

After $(\hat{\gamma}_i)_{i=1:N},\hat{\nu}$ have been constructed, 
one can generate any number of independent samples from 
$\hat{\nu}$ using the procedure $\SampQual$,
and for $i=1,\ldots,N$,
one can generate any number of independent samples from 
$\hat{\gamma}_i$ using the procedure $\SampAlloc$.
Specifically, 
assuming that one has access to an infinite sequence $(U_t)_{t\in\N}$ of independent realizations of $\mathrm{Uniform}[0,1]$ random variable,
$\SampQual\big(\hat{\nu},(U_t)_{t\in\N}\big)$
returns $(\BIz,\overline{t})$ 
where $\BIz\in\CX_0$ is a sample from $\hat{\nu}$ and $\overline{t}\in\N_{0}$ 
corresponds to the number of realizations in $(U_t)_{t\in\N}$ consumed in $\SampQual$,
that is, 
$\BIz$ is a Borel measurable function of $(U_t)_{t=1:\overline{t}}$.
Thus, the generation of $K\in\N$ independent samples from $\hat{\nu}$
can be carried out by 
initializing $t_0\leftarrow 0$,
calling 
$(\BIz_{k},\overline{t}_{k})\leftarrow \SampQual\big(\hat{\nu},(U_t)_{t\ge t_0+1}\big)$
followed by 
the update $t_0\leftarrow t_0+\overline{t}_{k}$
for $k=1,\ldots,K$,
and then gathering the samples 
$\{\BIz_{k}\}_{k=1:K}$.
Similarly, for each $i\in\{1,\ldots,N\}$,
$\SampAlloc\big(\hat{\gamma}_i,(U_t)_{t\in\N}\big)$
returns $(\BIx,\BIz,\overline{t})$
where $(\BIx,\BIz)\in\CX_i\times\CX_0$ is a sample from 
$\hat{\gamma}_i$ and $\overline{t}\in\N_{0}$
corresponds to the number of realizations in $(U_t)_{t\in\N}$ consumed in $\SampAlloc$,
that is,
$(\BIx,\BIz)$ is a Borel measurable function of $(U_t)_{t=1:\overline{t}}$.
Subsequently, the generation of $K\in\N$ independent samples from 
$\hat{\gamma}_i$ can be similarly carried out by iteratively calling 
$\SampAlloc$.

The following lemma presents the properties of 
$\ConstructAllocs$, $\SampQual$, and $\SampAlloc$.

\begin{lemma}\label{lem:proc-construct-allocs}
    Under Setting~\ref{sett:Euclidean},
    let us denote $m:=\max_{0\le i\le N}\{m_i\}$.
    For $i=1,\ldots,N$, let 
    $(\bar{\BIg}_i,\mathtt{Aux}_i)$ be
    the output of 
    $\InitTypeMeasure(\CX_i,\allowbreak\mu_i,\allowbreak\CG_i)$,
    where $\mathtt{Aux}_i$ is encoded by 
    $\big(\widetilde{R}_i,\allowbreak(F_{i,\tilde{r}},\tilde{\varpi}_{i,\tilde{r}})_{\tilde{r}=1:\widetilde{R}_i},\allowbreak\BP_i\big)$ as in \ref{enc:typemeasures-aux}.
    Let us denote $\widetilde{R}:=\max_{1\le i\le N}\{\widetilde{R}_i\}$.
    Moreover, for $i=1,\ldots,N$,
    let $\hat{\theta}_i\in\CP(\CX_i\times\CX_0)$ satisfy 
    $\support(\hat{\theta}_i)\subseteq V(\FS_i)\times V(\FS_0)$,
    let $\hat{\mu}_i$, $\hat{\nu}_i$ denote the marginals of $\hat{\theta}_i$ on $\CX_i$ and~$\CX_0$, respectively,
    and assume that 
    $\hat{\mu}_i\abovebelowset{\CG_i}{\varsigma}{\sim}\mu_i$, 
    $\hat{\nu}_i\abovebelowset{\CG_0}{\varsigma}{\sim}\hat{\nu}_1$
    $\forall 1\le i\le N$ hold with respect to 
    some $\varsigma>0$.
    Furthermore, 
    let 
    $\big((\hat{\gamma}_i)_{i=1:N},\hat{\nu}\big)$
    be the output of 
    $\ConstructAllocs\big(N,\allowbreak(\CX_i,\CG_i)_{i=0:N},\allowbreak(\mu_i,\mathtt{Aux}_i,\hat{\theta}_i)_{i=1:N}\big)$,
    let $(\Omega,\CF,\PROB)$ be a probability space and 
    let $(U_t)_{t\in\N}$ be an infinite sequence of independent $\mathrm{Uniform}[0,1]$ random variables on $\Omega$.
    Then, the following statements hold.\widowpenalty-1000
    \begin{enumerate}[label=(\roman*),beginpenalty=10000]
        \item\label{lems:proc-construct-allocs-complexity}
        $\ConstructAllocs\big(N,\allowbreak(\CX_i,\CG_i)_{i=0:N},\allowbreak(\mu_i,\mathtt{Aux}_i,\hat{\theta}_i)_{i=1:N}\big)$ incurs 
        $O\big(N(\widetilde{R}m^{\omega-1}+m^{\omega})\big)$ arithmetic operations.
        
        \item\label{lems:proc-construct-allocs-samp-qual}
        Let $(\widetilde{\BIZ},\overline{t})$ be the output of $\SampQual\big(\hat{\nu},(U_t)_{t\in\N}\big)$.
        Then, $\SampQual\big(\hat{\nu},(U_t)_{t\in\N}\big)$ incurs $O(m)$ arithmetic operations.
        Moreover, the law of the random variable $\widetilde{\BIZ}:\Omega\to\CX_0$ is equal to $\hat{\nu}$.
        
        \item\label{lems:proc-construct-allocs-samp-alloc}
        Let $T_{\mathtt{Samp}}$ denote the number of arithmetic operations incurred by the oracle procedure $\mathtt{Samp}$ in \ref{enc:typemeasures}.
        Let $i\in\nobreak\{1,\ldots,N\}$ be arbitrary,
        let 
        $(\widetilde{\BIX}_i,\widetilde{\BIZ},\overline{t})$ be the output of 
        $\SampAlloc\big(\hat{\gamma}_i,\allowbreak(U_t)_{t\in\N}\big)$,
        and let 
        $\widetilde{Q}:\nobreak\Omega\to\{0,1,\ldots,n_0-1\}$,
        $\widetilde{S}:\Omega\to\{1,\ldots,\widetilde{R}_i\}$
        be the random variables generated in 
        \mbox{Lines~\ref{alglin:samp-alloc-qualityindex}--\ref{alglin:samp-alloc-faceindex}} of
        $\SampAlloc$.
        Then, $\SampAlloc\big(\hat{\gamma}_i,(U_t)_{t\in\N}\big)$
        incurs $O(\widetilde{R}+m+T_{\mathtt{Samp}})$
        arithmetic operations.
        Moreover, 
        it holds that 
        $\PROB\big[\mu_i\big(\relint(F_{i,\widetilde{S}})\big)>\nobreak0\big]=\nobreak 1$, which guarantees that the call to $\mathtt{Samp}$
        in Line~\ref{alglin:samp-alloc-condsamp} in $\SampAlloc$ is well-defined.
        Furthermore,
        the law of the random variable $\widetilde{\BIX}_i:\nobreak\Omega\to\nobreak\CX_i$
        is equal to $\mu_i$,
        and the law of the random variable $\widetilde{\BIZ}:\Omega\to\CX_0$
        is equal to $\hat{\nu}$.
        Thus, $\hat{\gamma}_i\in\nobreak\CP(\CX_i\times\nobreak\CX_0)$ as characterized by the law of the random variable 
        $(\widetilde{\BIX}_i,\widetilde{\BIZ}):\Omega\to\CX_i\times\CX_0$ 
        satisfies $\hat{\gamma}_i\in\Gamma(\mu_i,\hat{\nu})$.

        \item\label{lems:proc-construct-allocs-error}
        One can construct random variables 
        $\BIZ:\Omega\to\CX_0$, 
        $(\BIZ_i:\Omega\to\CX_0)_{i=1:N}$,
        $(\BIX_i:\Omega\to\CX_i)_{i=1:N}$, and 
        $(\overbar{\BIX}_i:\Omega\to\CX_i)_{i=1:N}$ 
        (see the proof in Section~\ref{sapx:proof-complexity-AME} for the explicit constructions)
        such that 
        the law of $\BIZ$ is $\hat{\nu}$,
        and for $i=1,\ldots,N$,
        the law of $(\BIX_i,\BIZ_i):\Omega\to \CX_i\times\CX_0$
        is $\hat{\theta}_i$,
        the law of $\overbar{\BIX}_i$ is $\mu_i$,
        the law of $(\overbar{\BIX}_i,\BIZ):\Omega\to\CX_i\times\CX_0$
        is $\hat{\gamma}_i$,
        and it holds that 
        $\EXP\big[\|\BIZ-\BIZ_i\|_2\big]\le \overline{\epsilon}_0(\varsigma)$,
        $\EXP\big[\|\BIX_i-\overbar{\BIX}_i\|_2\big]\le \overline{\epsilon}_{i}(\varsigma)$.

    \end{enumerate}
    
\end{lemma}

\proof{Proof of Lemma~\ref{lem:proc-construct-allocs}.}
    See Section~\ref{sapx:proof-complexity-AME}.
\endproof

In particular, Lemma~\ref{lem:proc-construct-allocs}\ref{lems:proc-construct-allocs-error} shows that 
$\big((\hat{\gamma}_i)_{i=1:N},\hat{\nu}\big)$
returned by $\ConstructAllocs$ satisfy all the corresponding requirements 
in Theorem~\ref{thm:mt-tf-approx}.
Therefore, if we begin with 
a pair of $(\epsilon_{\sspar},\varsigma)$-optimizers of 
\eqref{eqn:mt-tf-dual} and \eqref{eqn:mt-tf-lsip}
constructed by $\SolveLSIPDual$ and $\SolveLSIPPrimal$
(recall Theorem~\ref{thm:proc-solve-LSIP}),
then $\ConstructTransFuncs$ and $\ConstructAllocs$ 
are concrete implementations of the construction of 
an AME
$(\hat{\varphi}_i,\hat{\gamma}_i)_{i=1:N},\hat{\nu}$ via Theorem~\ref{thm:mt-tf-approx}. 

We have now presented all the procedures used in Algorithm~\ref{alg:complexity-overall}.
One observes that Theorem~\ref{thm:complexity-explicit}
now directly follows by combining the results in 
Lemma~\ref{lem:proc-init-test-funcs}, 
Lemma~\ref{lem:proc-init-type-measure},
Theorem~\ref{thm:proc-solve-LSIP},
Lemma~\ref{lem:proc-construct-trans-funcs},
and 
Lemma~\ref{lem:proc-construct-allocs}.
The detailed proof of Theorem~\ref{thm:complexity-explicit} will be presented 
in Section~\ref{sapx:proof-numerics}.

\subsection{The particular class of $(\mu_i)_{i=1:N}$ in Setting~\ref{sett:Euclidean}{\normalfont\ref{settc:Euclidean-particular}}}\label{sapx:complexity-particular-case}

This subsection presents concrete implementations of the oracle procedures 
$\mathtt{Int}$ and $\mathtt{Samp}$ in \ref{enc:typemeasures} under the particular class of $(\mu_i)_{i=1:N}$ 
in Setting~\ref{sett:Euclidean}\ref{settc:Euclidean-particular}.
To that end, 
let us first present below the concrete encoding of 
the type measures $(\mu_i)_{i=1:N}$ under Setting~\ref{sett:Euclidean}\ref{settc:Euclidean-particular}.
\begin{enumerate}[label=(E\arabic*'),beginpenalty=10000,leftmargin=30pt]
    \setcounter{enumi}{3}
    \item\label{enc:typemeasures-particular}
    For $i=1,\ldots,N$,
    observe that the type measure
    $\mu_i\in\CP(\CX_i)$ is completely characterized by the affine functions 
    $(p_{i,K}:K\to\R_+)_{K\in\FK_i,\,K\ne\emptyset}$
    under Setting~\ref{sett:Euclidean}\ref{settc:Euclidean-particular},
    which are in turn completely characterized by 
    $\big(p_{i,K}(\BIu)\big)_{\BIu\in V(K),\, K\in\FK_i,\, K\ne\emptyset}$.
    Thus,
    given the encoding 
    $\big(d_i,q_i,\allowbreak(\BIu_{i,l})_{l=0:q_i},\allowbreak J_i,\allowbreak \big(\rho_{i,j},(\kappa_{i,j,\iota})_{\iota=0:\rho_{i,j}}\big)_{j=1:J_i}\big)$
    of $\CX_i$ in \ref{enc:spaces},
    we can efficiently encode $\mu_i$ by 
    an associative array $\mathtt{Map}_{\mu_i}$
    containing key-value pairs
    in which each key is a non-empty and increasing tuple
    $(\kappa_{i,j,\iota_r})_{r=0:n_{i,j}}$ with 
    $j\in\{1,\ldots,J_i\}$,
    $n_{i,j}\in\{0,1,\ldots, \rho_{i,j}\}$,
    and $\{\iota_r\}_{r=0:n_{i,j}}\subseteq \{0,1,\ldots,\rho_{i,j}\}$.
    The associated value 
    $\mathtt{Map}_{\mu_i}\big[(\kappa_{i,j,\iota_r})_{r=0:n_{i,j}}\big]$
    is a tuple 
    $\big(p_{i,K_{i,j,(\iota_r)}}(\BIu_{i,\kappa_{i,j,\iota_r}})\big)_{r=0:n_{i,j}}$
    where $K_{i,j,(\iota_r)}:=\conv\big(\{\BIu_{i,\kappa_{i,j,\iota_r}}\}_{r=0:n_{i,j}}\big)\in \FK_i$.
    Concretely, let 
    $\mathtt{Map}_{\mu_i}$ be a red-black tree 
    \citep*[Chapter~13]{ECcormen2009introduction} with integer-valued keys,
    where each key $(\kappa_{i,j,\iota_r})_{r=0:n_{i,j}}$ is 
    encoded by an integer 
    $\sum_{r=0}^{n_{i,j}}\big(J_i(d_i+\nobreak1)+\nobreak1\big)^{r}(\kappa_{i,j,\iota_r}+1)$.
    Since $q_i+1=V(\FK_i)\le J_i(d_i+1)$, 
    this integer encoding uniquely determines the key $(\kappa_{i,j,\iota_r})_{r=0:n_{i,j}}$.
    Moreover,
    it takes $O(d_i)$ arithmetic operations to encode any key
    into an integer.
    Subsequently,
    since $\mathtt{Map}_{\mu_i}$ contains at most $O(2^{d_i}J_i)$ entries,
    it follows from the properties of the red-black tree that 
    each search of $\mathtt{Map}_{\mu_i}$ for a given key incurs 
    $O\big(d_i+\log(J_i)\big)$ arithmetic operations.
\end{enumerate}

\begin{procedure}[t]
    \DontPrintSemicolon
    \caption{GetVol(){$(C)$} $\to$ $w$}
    \begin{footnotesize}
    \KwInAndOut{\!\!\!\!\!\begin{tabular}{lll}
        \textbf{Input} &
        $C$ & an $n$-simplex in $\R^{d}$ encoded by $\big(d,n,(\BIv_l)_{l=0:n}\big)$ as in \ref{enc:simplex} \\
        \textbf{Output} &
        $w\in\R_+$ & the $n$-dimensional volume $\volume_n(C)$
    \end{tabular}\!\!\!\!\!\!\!\!\!\!\!\!\!\!}

    \nl%
    \If{$n=0$}{
        \Return%
        $w\leftarrow 1$.
    }

    \nl%
    $\BV\leftarrow\left(\begin{smallmatrix}
        | & & | \\
        (\BIv_1-\BIv_0) & \cdots & (\BIv_n-\BIv_0) \\
        | & & | \\
    \end{smallmatrix}\right)\in \R^{d\times n}$.\tcp*[f]{$O(d^2)$}

    \nl\label{alglin:get-vol-QR}%
    $(\BQ,\BR)\leftarrow$ the QR decomposition of $\BV$,
    where $\BQ\in\R^{d\times d}$ is an orthonormal matrix, 
    $\BR=\left(\begin{smallmatrix}
    \BR_1 \\ \BO
    \end{smallmatrix}\right)$,
    $\BR_1\in\R^{n\times n}$ is an upper triangular matrix,
    and $\BO\in\R^{(d-n)\times n}$ is the all-zero $(d-n)$-by-$n$ matrix.\tcp*[f]{$O(d^3)$}

    \nl\label{alglin:get-vol-det}%
    $w\leftarrow \big|\prod_{i=1}^{n}r_{i,i}\big|/n!$, 
    where $r_{1,1},\ldots,r_{n,n}$ are the diagonal entries of $\BR_1$.\tcp*[f]{$O(d)$}

    \nl%
    \Return $w$.
    \end{footnotesize}
\end{procedure}

\begin{procedure}
    \DontPrintSemicolon
    \caption{IntSimplex(){$(\CX_i,\mu_i,C,\hat{j})$} $\to$ $(h_l)_{l=0:n}$}
    \begin{footnotesize}
    \KwInAndOut{\!\!\!\!\!\begin{tabular}{lll}
        \textbf{Input} & 
        $\CX_i$ & the type space encoded by $\big(d_i,q_i,\allowbreak(\BIu_{i,l})_{l=0:q_i},\allowbreak J_i,\allowbreak \big(\rho_{i,j},(\kappa_{i,j,\iota})_{\iota=0:\rho_{i,j}}\big)_{j=1:J_i}\big)$ as in \ref{enc:spaces} \\
        & $\mu_i\in\CP(\CX_i)$ & the type measure on $\CX_i$ encoded by an associative array $\mathtt{Map}_{\mu_i}$ as in \ref{enc:typemeasures-particular} \\
        & $C$ & an $n$-simplex encoded by $\big(d_i,n,(\BIv_l)_{l=0:n}\big)$ as in \ref{enc:simplex} \\
        & $\hat{j}\in\{1,\ldots,J_i\}$ & an index satisfying $\conv\big(\{\BIu_{i,\kappa_{i,\hat{j},\iota}}\}_{\iota=0:\rho_{i,\hat{j}}}\big)\supseteq C$ \\
        \textbf{Output} &
        $h_l\in\R_+$ & the integral $\int_{\CX_i}\lambda_{\FF(C),\BIv_l}\INDI_{\relint(C)}\DIFFX{\mu_i}$ where $\lambda_{\FF(C),\BIv_l}:\CX_i\to[0,1]$ is the barycentric \\
        && coordinate function with respect to $\BIv_l\in V(C)$
    \end{tabular}\!\!\!\!\!\!\!\!\!\!\!\!\!\!}

    \tcc{Let us denote $d:=\max_{0\le i\le N}\{d_i\}$ and $J:=\max_{0\le i\le N}\{J_i\}$ in the following computational complexities.}

    \nl\label{alglin:int-simplex-encode-outer-simplex}%
    Let $K_{i,\hat{j}}$ be a simplex encoded by $(d_i,\rho_{i,\hat{j}},(\BIu_{i,\kappa_{i,\hat{j},\iota}})_{\iota=0:\rho_{i,\hat{j}}})$ as in \ref{enc:simplex}.\tcp*[f]{$O(d^2)$}

    \nl\label{alglin:int-simplex-get-coords-forloop}%
    \For(\tcp*[f]{$O(d^{\omega+1})$}){$s=0,1,\ldots,n$}{

        \nl\label{alglin:int-simplex-get-coords}%
        $\big(\mathtt{Inside},(\zeta_{s,\iota})_{\iota=0:\rho_{i,\hat{j}}}\big)\leftarrow \GetCoords(K_{i,\hat{j}},\BIv_s)$.
        \tcp*[f]{$O(d^\omega)$}
    }

    \nl\label{alglin:int-simplex-search-face}%
    $\{\hat{\iota}_r\}_{r=0:\hat{n}}\leftarrow \big\{\iota\in\{0,1,\ldots,\rho_{i,\hat{j}}\}:\sum_{s=0}^{n}\zeta_{s,\iota}>0\big\}$ where $(\hat{\iota}_r)_{r=0:\hat{n}}$ is in ascending order.\tcp*[f]{$O(d^2)$}

    \nl%
    \If{$n<\hat{n}$}{
        \nl\label{alglin:int-simplex-degenerate}%
        \Return $(0)_{l=0:n}$.
    }

    \nl\label{alglin:int-simplex-get-dens}%
    $\big(p_{i,K_{i,\hat{j},(\hat{\iota}_r)}}(\BIu_{i,\kappa_{i,\hat{j},\hat{\iota}_r}})\big)_{r=0:\hat{n}}\leftarrow \mathtt{Map}_{\mu_i}\big[(\kappa_{i,\hat{j},\hat{\iota}_r})_{r=0:\hat{n}}\big]$.\tcp*[f]{$O\big(d+\log(J)\big)$}

    \nl\label{alglin:int-simplex-get-vol}%
    $w\leftarrow \GetVol(C)$.\tcp*[f]{$O(d^3)$}

    \nl%
    \For(\tcp*[f]{$O(d^3)$}){$l=0,1,\ldots,n$}{

        \nl\label{alglin:int-simplex-calc-int}%
        $h_l\leftarrow \frac{w}{(n+2)(n+1)}\Big[\sum_{r=0}^{\hat{n}}\big(\zeta_{l,\hat{\iota}_r}+\sum_{s=0}^{n}\zeta_{s,\hat{\iota}_r}\big)p_{i,K_{i,\hat{j},(\hat{\iota}_r)}}(\BIu_{i,\kappa_{i,\hat{j},\hat{\iota}_r}})\Big]$.\tcp*[f]{$O(d^2)$}
    }

    \nl \Return{$(h_l)_{l=0:n}$.}
    \end{footnotesize}
\end{procedure}

\begin{procedure}
    \DontPrintSemicolon
    \caption{SampSimplex(){$\big(\CX_i,\mu_i,C,\hat{j},(U_t)_{t\in\N}\big)$} $\to$ $(\BIx,\overline{t})$}
    \begin{footnotesize}
    \KwInAndOut{\!\!\!\!\!\begin{tabular}{lll}
        \textbf{Input} & 
        $\CX_i$ & the type space encoded by $\big(d_i,q_i,\allowbreak(\BIu_{i,l})_{l=0:q_i},\allowbreak J_i,\allowbreak \big(\rho_{i,j},(\kappa_{i,j,\iota})_{\iota=0:\rho_{i,j}}\big)_{j=1:J_i}\big)$ as in \ref{enc:spaces} \\
        & $\mu_i\in\CP(\CX_i)$ & the type measure on $\CX_i$ encoded by an associative array $\mathtt{Map}_{\mu_i}$ as in \ref{enc:typemeasures-particular} \\
        & $C$ & an $n$-simplex encoded by $\big(d_i,n,(\BIv_l)_{l=0:n}\big)$ as in \ref{enc:simplex}, where $\mu_i\big(\relint(C)\big)>0$ \\
        & $\hat{j}\in\{1,\ldots,J_i\}$ & an index satisfying $\conv\big(\{\BIu_{i,\kappa_{i,\hat{j},\iota}}\}_{\iota=0:\rho_{i,\hat{j}}}\big)\supseteq C$ \\
        & $(U_{t})_{t\in\N}\subset [0,1]$ & an infinite sequence of independent realizations of $\mathrm{Uniform}[0,1]$ random variable \\
        \textbf{Output} &
        $\BIx\in\R^{d_i}$ & the generated random sample from $\tilde{\mu}_{i,C}\in\CP(\CX_i)$ defined by $\frac{\DIFF\tilde{\mu}_{i,C}}{\DIFF\mu_i}:=\frac{\INDI_{\relint(C)}}{\mu_i(\relint(C))}$ \\
        & $\overline{t}\in\N_0$ & the number of realizations in $(U_t)_{t\in\N}$ consumed in the procedure
    \end{tabular}\!\!\!\!\!\!\!\!\!\!\!\!\!\!}

    \tcc{Let us denote $d:=\max_{0\le i\le N}\{d_i\}$ and $J:=\max_{0\le i\le N}\{J_i\}$ in the following computational complexities.}

    \nl\label{alglin:samp-simplex-order-init}%
    $U_{(0)}\leftarrow 0$.

    \nl\label{alglin:samp-simplex-order}%
    $(U_{(t)})_{t=1:n+1}\leftarrow$ $(U_t)_{t=1:n+1}$ sorted into ascending order, i.e., $(U_{(t)})_{t=1:n+1}$ satisfy
    $\{U_t\}_{t=1:n+1}=\{U_{(t)}\}_{t=1:n+1}$ and $U_{(t)}\le U_{(t+1)}$ 
    $\forall t\in\{1,\ldots,n\}$.\tcp*[f]{$O\big(d\log(d)\big)$}

    \nl\label{alglin:samp-simplex-encode-outer-simplex}%
    Let $K_{i,\hat{j}}$ be a simplex encoded by $(d_i,\rho_{i,\hat{j}},(\BIu_{i,\kappa_{i,\hat{j},\iota}})_{\iota=0:\rho_{i,\hat{j}}})$ as in \ref{enc:simplex}.\tcp*[f]{$O(d^2)$}

    \nl\label{alglin:samp-simplex-get-coords-forloop}\For(\tcp*[f]{$O(d^{\omega+1})$}){$s=0,1,\ldots,n$}{

        \nl\label{alglin:samp-simplex-get-coords}%
        $\big(\mathtt{Inside},(\zeta_{s,\iota})_{\iota=0:\rho_{i,\hat{j}}}\big)\leftarrow \GetCoords(K_{i,\hat{j}},\BIv_s)$.\tcp*[f]{$O(d^\omega)$}
    }

    \nl\label{alglin:samp-simplex-search-face}%
    $\{\hat{\iota}_r\}_{r=0:\hat{n}}\leftarrow \big\{\iota\in\{0,1,\ldots,\rho_{i,\hat{j}}\}:\sum_{s=0}^{n}\zeta_{s,\iota}>0\big\}$ where $(\hat{\iota}_r)_{r=0:\hat{n}}$ is in ascending order.\tcp*[f]{$O(d^2)$}

    \tcc{Under the assumption that $\mu_i\big(\relint(C)\big)>0$, it is guaranteed that $\hat{n}=n$; see Lemma~\ref{lem:proc-int-simplex}\ref{lems:proc-int-simplex-samp}.}

    \nl\label{alglin:samp-simplex-get-dens}%
    $\big(p_{i,K_{i,\hat{j},(\hat{\iota}_r)}}(\BIu_{i,\kappa_{i,\hat{j},\hat{\iota}_r}})\big)_{r=0:n}\leftarrow \mathtt{Map}_{\mu_i}\big[(\kappa_{i,\hat{j},\hat{\iota}_r})_{r=0:n}\big]$.\tcp*[f]{$O\big(d+\log(J)\big)$}

    \nl%
    \For(\tcp*[f]{$O(d^2)$}){$l=0,1,\ldots,n$}{

        \nl\label{alglin:samp-simplex-calc-int}%
        $q_l\leftarrow \sum_{r=0}^{\hat{n}}\zeta_{l,\hat{\iota}_r}p_{i,K_{i,\hat{j},(\hat{\iota}_r)}}(\BIu_{i,\kappa_{i,\hat{j},\hat{\iota}_r}})$.\tcp*[f]{$O(d)$}
    }

    \nl\label{alglin:samp-simplex-samp-component}%
    $S\leftarrow \min\big\{s\in\{0,1,\ldots,n\}:\sum_{l=0}^{s}q_{l}\ge \big(\sum_{l=0}^{n}q_{l}\big)U_{n+2}\big\}$.\tcp*[f]{$O(d)$}

    \nl\label{alglin:samp-simplex-samp-point}%
    $\BIx\leftarrow (1-U_{(n+1)})\BIv_{S}+\sum_{l=0}^{n}(U_{(l+1)}-U_{(l)})\BIv_{l}$.\tcp*[f]{$O(d^2)$}

    \nl%
    \Return{$(\BIx,n+2)$.}
    \end{footnotesize}
\end{procedure}

With the encoding of $(\mu_i)_{i=1:N}$ in \ref{enc:typemeasures-particular},
let us present the procedures 
$\IntSimplex$ and $\SampSimplex$,
which satisfy the requirements for the oracle procedures 
$\mathtt{Int}$ and $\mathtt{Samp}$ in \ref{enc:typemeasures}, respectively.
Specifically, 
$\IntSimplex$ utilizes the auxiliary sub-procedure 
$\GetCoords$ presented earlier in Section~\ref{sapx:complexity-AME},
as well as another auxiliary sub-procedure $\GetVol$
for computing the $n$-dimensional volume $\volume_n(C)$ of any $n$-simplex $C$ 
(see Lemma~\ref{lem:proc-get-vol} in Section~\ref{sapx:proof-complexity-particular-case} for the properties of $\GetVol$).
Moreover, $\SampSimplex$ also utilizes the auxiliary sub-procedure $\GetCoords$.
We present in the following lemma the properties of
$\IntSimplex$ and $\SampSimplex$.

\begin{lemma}\label{lem:proc-int-simplex}
    Under Setting~\ref{sett:Euclidean} and the particular class of $(\mu_i)_{i=1:N}$ in Setting~\ref{sett:Euclidean}\ref{settc:Euclidean-particular}, let us denote 
    $d:=\nobreak\max_{0\le i\le N}\{d_i\}$ and
    $J:=\max_{0\le i\le N}\big\{|\FM(\FK_i)|\big\}$.
    Let $i\in\{1,\ldots,N\}$ and $\hat{j}\in\{1,\ldots,J_i\}$ be arbitrary, 
    and let $C$ be an arbitrary $n$-simplex in $\R^{d_i}$ 
    where $n\in\{0,1,\ldots,d_i\}$,
    $V(C)=\{\BIv_l\}_{l=0:n}$,
    $(K_{i,j})_{j=1:J_i}$ is the enumeration of $\FM(\FK_i)$ in \ref{enc:spaces},
    and $C\subseteq K_{i,\hat{j}}$.
    Then, the following statements hold.
    \begin{enumerate}[label=(\roman*),beginpenalty=10000]
        \item\label{lems:proc-int-simplex-int}
        $\IntSimplex(\CX_i,\mu_i,C,\hat{j})$ incurs 
        $O\big(d^{\omega+1}+\log(J)\big)$ arithmetic operations and returns 
        $(h_l)_{l=0:n}$,
        where
        $h_l=\nobreak\int_{\CX_i}\lambda_{\FF(C),\BIv_l}\INDI_{\relint(C)}\DIFFX{\mu_i}$ 
        $\forall l\in\{0,1,\ldots,n\}$.
        In particular, $\IntSimplex$ satisfies the requirements for the 
        oracle procedure $\mathtt{Int}$ in \ref{enc:typemeasures}.

        \item\label{lems:proc-int-simplex-samp}
        Let $(\Omega,\CF,\PROB)$ be a probability space and 
        let $(U_t)_{t\in\N}$ be an infinite sequence of independent $\mathrm{Uniform}[0,1]$ random variables on $\Omega$.
        Let us assume in addition that $\mu_i\big(\relint(C)\big)>0$.
        Then, 
        $\SampSimplex\big(\CX_i,\mu_i,C,\hat{j},(U_t)_{t\in\N}\big)$
        incurs 
        $O\big(d^{\omega+1}+\log(J)\big)$ arithmetic operations and returns $(\BIX,\overline{t})$,
        where 
        $\hat{n}=n$ holds in Line~\ref{alglin:samp-simplex-search-face},
        and the law of 
        $\BIX:\Omega\to\CX_i$ is $\tilde{\mu}_{i,C}$
        which is defined by 
        $\frac{\DIFF\tilde{\mu}_{i,C}}{\DIFF\mu_i}:=\frac{\INDI_{\relint(C)}}{\mu_i(\relint(C))}$.
        In particular, $\SampSimplex$ satisfies the requirements for the oracle procedure $\mathtt{Samp}$ in \ref{enc:typemeasures}.
    \end{enumerate}
    
\end{lemma}

\proof{Proof of Lemma~\ref{lem:proc-int-simplex}.}
    See Section~\ref{sapx:proof-complexity-particular-case}.
\endproof

\section{Additional remarks about Algorithm~\ref{alg:cp-tf} and Algorithm~\ref{alg:mt-tf}}\label{apx:efficient-algo}

\vspace{6pt}

\subsection{Remarks about solving the LP problems in Algorithm~\ref{alg:cp-tf}}\label{sapx:efficient-algo-LP}

Line~\ref{alglin:cp-tf-lp} of Algorithm~\ref{alg:cp-tf} simultaneously solves the LP relaxation \eqref{eqn:mt-tf-algo-lp} and its dual LP problem \eqref{eqn:mt-tf-algo-lp-dual} via the dual simplex algorithm (see, e.g., \citep[Chapter~6.4]{ECvanderbei2020linear}) or the interior point algorithm (see, e.g., \citep[Chapter~18]{ECvanderbei2020linear}).
For any $r\in\N_{0}$,
let 
$d:=\max_{0\le i\le N}\{d_i\}$,
$q:=\nobreak\sum_{i=1}^N(1+\nobreak m_i+\nobreak m_0)$, 
$p:=\nobreak\sum_{i=1}^{N}|\CK^{(r)}_{i}|$,
and let us re-express \eqref{eqn:mt-tf-algo-lp} into a vectorized form as follows:
\begin{align*}
    \maximize_{\Bchi}\quad & \langle\bar{\BIg},\Bchi\rangle\\*
    \mathrm{subject~to}\quad & \BA\Bchi\le \BIc, \quad \BU\Bchi=\veczero_{m_0}, \quad  \Bchi\in\R^{q},
\end{align*}
where
$\bar{\BIg}:=(1,\bar{\BIg}_1^\TRANSP,\veczero_{m_0}^\TRANSP,\ldots,1,\bar{\BIg}_N^\TRANSP,\veczero_{m_0}^\TRANSP)^\TRANSP\in\R^{q}$,
$\BA\in\R^{p\times q}$,
and 
$\BU\in\R^{m_0\in q}$.
Observe that 
in the LP problem above,
each row of
the matrix $\BA$ has 
no more than $2d+1$ non-zero entries,
and hence the matrix $\BA$ has at most $p(2d+1)$ non-zero entries, 
whereas
the matrix $\BU$ has 
$Nm_0$ non-zero entries.
Therefore, the constraint matrix of \eqref{eqn:mt-tf-algo-lp} is sparse,
which can be exploited by modern LP solvers for efficiency gain;
see, e.g., \citep[Chapter~8.2]{ECvanderbei2020linear} 
and \citep[Chapter~11.8]{ECboyd2004convex}
for related discussions about 
how sparsity of the constraint matrix can be exploited 
in the dual simplex algorithm and the interior point algorithm.
Moreover, in the case where the dual simplex algorithm is used for solving 
$\CPDLPDUAL{r}$ and $\CPDLPPRIMAL{r}$,
observe that the constraint matrix of $\CPDLPDUAL{r+1}$
is formed by adding rows to the constraint matrix of $\CPDLPDUAL{r}$.
Hence, the optimizers of $\CPDLPDUAL{r}$ and $\CPDLPPRIMAL{r}$
computed by the dual simplex algorithm can be used to warm-start
the dual simplex algorithm when solving 
$\CPDLPDUAL{r+1}$ and $\CPDLPPRIMAL{r+1}$;
see, e.g., the reference manual of \citet{gurobi}.

\subsection{Remarks about solving the global minimization problems in Algorithm~\ref{alg:cp-tf}}\label{sapx:efficient-algo-global}

In Line~\ref{alglin:cp-tf-global-tol} of Algorithm~\ref{alg:cp-tf},
a positive sequence $(\tau^{(r)})_{r\in\N_0}$ is chosen,
where each $\tau^{(r)}$ corresponds to the error tolerance when solving the global minimization problem in Line~\ref{alglin:cp-tf-global} of Algorithm~\ref{alg:cp-tf} in iteration~$r\in\N_{0}$.
It is required that $\limsup_{r\to\infty}\tau^{(r)}<\frac{\epsilon_{\sspar}}{N}$
in order to guarantee the finite termination of Algorithm~\ref{alg:cp-tf}.
In practice, one typically starts with a large tolerance $\tau^{(0)}$
and then (possibly adaptively) decrease $\tau^{(r)}$ over the iterations.
This can often improve the computational efficiency compared to using a constant tolerance.

Line~\ref{alglin:cp-tf-global} of Algorithm~\ref{alg:cp-tf} computes a $\tau^{(r)}$-optimizer 
$(\tilde{\BIx}_{i}^{(r)},\tilde{\BIz}_{i}^{(r)})$
as well as a lower bound $\underline{\beta}_i^{(r)}$
of the following global minimization problem:
\begin{align}
    \begin{split}
        \minimize_{\BIx,\,\BIz}\quad & c_i(\BIx,\BIz)
            -\langle\BIg_i(\BIx),\BIy_i^{(r)}\rangle
            -\langle\BIg_0(\BIz),\BIw_i^{(r)}\rangle \\
        \mathrm{subject~to}\quad & \BIx\in\CX_i,\quad \BIz\in\CX_0,
    \end{split}
    \label{eqn:oracle-globalmin}
\end{align}
where $(y_{i,0}^{(r)},\BIy_{i}^{(r)},\BIw_{i}^{(r)})_{i=1:N}$ 
are computed in Line~\ref{alglin:cp-tf-lp}.
In the following, 
we present two methods to solve (\ref{eqn:oracle-globalmin})
in practice.

The first method is based on grid discretization.
Let $\FD_i$ be a subdivision of $\FS_i$ such that $\eta(\FD_i)\le\nobreak \frac{\tau^{(r)}}{2L_{c}}$, 
and let $\FD_0$ be a subdivision of $\FS_0$ such that $\eta(\FD_0)\le \frac{\tau^{(r)}}{2L_{c}}$.
Then, one checks that the function $\hat{c}_i:\CX_i\times\CX_0\to\R$ defined by
\begin{align*}
    \hat{c}_i(\BIx,\BIz):=\sum_{\BIu\in V(\FD_i)}\sum_{\BIv\in V(\FD_0)}
        \lambda_{\FD_i,\BIu}(\BIx)\lambda_{\FD_0,\BIv}(\BIz)c_i(\BIu,\BIv)
        \qquad \forall (\BIx,\BIz)\in\CX_i\times\CX_0
\end{align*}
satisfies 
$\hat{c}_i(\BIx,\BIz)=c_i(\BIx,\BIz)$
$\forall (\BIx,\BIz)\in V(\FD_i)\times V(\FD_0)$.
Moreover, it holds that  
\begin{align}
    \max_{(\BIx,\BIz)\in\CX_i\times\CX_0}\big\{\big|c_i(\BIx,\BIz)-\hat{c}_i(\BIx,\BIz)\big|\big\}\le\nobreak \tau^{(r)}.
    \label{eqn:global-grid-discretization}
\end{align}
To see this,
let us fix arbitrary $K_i\in \FD_i$, $K_0\in \FD_0$,
fix an arbitrary $(\BIx,\BIz)\in K_i\times K_0$, 
and observe that 
\begin{align*}
    c_i(\BIx,\BIz)&\le c_i(\BIu,\BIv)+L_{c}\big(\eta(\FD_i)+\eta(\FD_0)\big)
    \le c_i(\BIu,\BIv)+\tau^{(r)} 
    & \forall (\BIu,\BIv)\in V(K_i)\times V(K_0),\\
    c_i(\BIx,\BIz)&\ge c_i(\BIu,\BIv)-L_{c}\big(\eta(\FD_i)+\eta(\FD_0)\big)
    \ge c_i(\BIu,\BIv)-\tau^{(r)} 
    & \forall (\BIu,\BIv)\in V(K_i)\times V(K_0).
\end{align*}
Multiplying the above inequalities by 
$\lambda_{\FD_i,\BIu}(\BIx)\lambda_{\FD_0,\BIv}(\BIz)$
and subsequently summing over 
all $(\BIu,\BIv)\in V(\FD_i)\times V(\FD_0)$ leads to 
$\big|c_i(\BIx,\BIz)-\hat{c}_i(\BIx,\BIz)\big|\le \tau^{(r)}$.
Since $\CX_i=\bigcup_{K\in \FD_i}K$, $\CX_0=\bigcup_{K\in\FD_0}K$,
this shows that (\ref{eqn:global-grid-discretization}) indeed holds.
Consequently, we can let 
\begin{align*}
    (\tilde{\BIx}_{i}^{(r)},\tilde{\BIz}_{i}^{(r)})\in
    \argmin_{(\BIx,\BIz)\in V(\FD_i)\times V(\FD_0)}\big\{c_i(\BIx,\BIz)-\langle\BIg_i(\BIx),\BIy_{i}^{(r)}\rangle-\langle\BIg_0(\BIz),\BIw_{i}^{(r)}\rangle\big\},
\end{align*}
and let 
$\underline{\beta}_{i}^{(r)}:=c_i(\tilde{\BIx}_{i}^{(r)},\tilde{\BIz}_{i}^{(r)})-\langle\BIg_i(\tilde{\BIx}_{i}^{(r)}),\BIy_{i}^{(r)}\rangle-\langle\BIg_0(\tilde{\BIz}_{i}^{(r)}),\BIw_{i}^{(r)}\rangle-\tau^{(r)}$.
Since
$\CX_i\ni \BIx\mapsto \hat{c}_{i}(\BIx,\BIz)\in\nobreak\R$ 
is continuous and piece-wise affine on each $K\in\FD_i$ for every $\BIz\in\CX_0$,
$\CX_0\ni \BIz\mapsto \hat{c}_{i}(\BIx,\BIz)\in\R$
is continuous and piece-wise affine on each $K\in\FD_0$ for every 
$\BIx\in\CX_i$,
$\CX_i\ni \BIx\mapsto \langle\BIg_i(\BIx),\BIy_{i}^{(r)}\rangle\in\nobreak\R$
is continuous and piece-wise affine on each $K\in\FD_i$,
and 
$\CX_0\ni \BIz\mapsto \langle\BIg_0(\BIz),\BIw_{i}^{(r)}\rangle\in\nobreak\R$
is continuous and piece-wise affine on each $K\in\FD_0$,
one checks that 
$(\tilde{\BIx}_{i}^{(r)},\tilde{\BIz}_{i}^{(r)})$
is an optimizer of the problem
$\min_{(\BIx,\BIz)\in\CX_i\times\CX_0}\big\{\hat{c}_i(\BIx,\BIz)-\langle\BIg_i(\BIx),\BIy_{i}^{(r)}\rangle-\langle\BIg_0(\BIz),\BIw_{i}^{(r)}\rangle\big\}$
and therefore 
the condition 
$\underline{\beta}_i^{(r)}
\le \min_{\BIx\in\CX_i,\,\BIz\in\CX_0}\big\{c_i(\BIx,\BIz)-\langle\BIg_i(\BIx),\BIy_i^{(r)}\rangle-\langle\BIg_0(\BIz),\BIw_i^{(r)}\rangle\big\}
\le \allowbreak c_i(\tilde{\BIx}_i^{(r)},\tilde{\BIz}_i^{(r)})\allowbreak-\langle\BIg_i(\tilde{\BIx}_i^{(r)}),\BIy_i^{(r)}\rangle-\langle\BIg_0(\tilde{\BIz}_i^{(r)}),\BIw_i^{(r)}\rangle
\le \underline{\beta}_i^{(r)} + \tau^{(r)}$
in Line~\ref{alglin:cp-tf-global} is satisfied.
Overall, this method requires the construction of 
two simplicial complexes $\FD_i$ and $\FD_0$ which are subdivisions of 
$\FS_i$ and $\FS_0$, respectively,
and the evaluation of 
$c_i(\BIx,\BIz)-\langle\BIg_i(\BIx),\BIy_i^{(r)}\rangle-\langle\BIg_0(\BIz),\BIw_i^{(r)}\rangle$
at the finitely many points in
$V(\FD_i)\times V(\FD_0)$.

The second method
requires the cost function $c_i$ to possess certain structures, e.g.,
a piece-wise convex quadratic structure, 
and is based on 
the reformulation of (\ref{eqn:oracle-globalmin}) into 
mixed-integer programming (MIP) problems.
To simplify the presentation, 
we can turn $c_i(\BIx,\BIz)$, $\langle\BIg_i(\BIx),\BIy_{i}^{(r)}\rangle$, and $\langle\BIg_0(\BIz),\BIw_{i}^{(r)}\rangle$ 
in the objective of (\ref{eqn:oracle-globalmin})
into constraints via their epigraphs and hypographs in order to make the objective of (\ref{eqn:oracle-globalmin}) linear:
\begin{align}
    \begin{split}
        \minimize_{\scalebox{0.62}{$\BIx,\,\BIz,\,\beta_0,\,\beta_1,\,\beta_2$}} \quad & \beta_0-\beta_1-\beta_2 \\
        \mathrm{subject~to} \quad & (\beta_{0},\BIx,\BIz) \in \epigraph(c_i), \\
        & (\beta_{1},\BIx)\in \hypograph\big(\langle\BIg_i(\,\cdot\,),\BIy_{i}^{(r)}\rangle\big),\\
        & (\beta_{2},\BIz)\in \hypograph\big(\langle\BIg_0(\,\cdot\,),\BIw_{i}^{(r)}\rangle\big),\\
        & \BIx\in\R^{d_i},\; \BIz\in\R^{d_0},\; \beta_0,\beta_1,\beta_2\in\R.
    \end{split}
    \label{eqn:oracle-globalmin-epigraph}
\end{align}
In the problem above, 
$\epigraph(c_i):=\big\{(\beta,\BIx,\BIz)\in \R\times\CX_i\times\CX_0:
    c_i(\BIx,\BIz)\le \beta\big\}$
denotes the epigraph of the function $c_i$,
and 
$\hypograph\big(\langle\BIg_i(\,\cdot\,),\BIy_{i}^{(r)}\rangle\big)
:=\big\{(\beta,\BIx)\in\R\times\CX_i:
    \langle\BIg_i(\BIx),\BIy_{i}^{(r)}\rangle\ge \beta\big\}$,
$\hypograph\big(\langle\BIg_0(\,\cdot\,),\BIw_{i}^{(r)}\rangle\big)
:=\big\{(\beta,\BIz)\in\R\times\CX_0:
    \langle\BIg_0(\BIz),\BIw_{i}^{(r)}\rangle\ge \beta\big\}$
denote the hypographs of the functions 
$\CX_i\ni\BIx\mapsto \langle\BIg_i(\BIx),\BIy_{i}^{(r)}\rangle\in\nobreak\R$,
$\CX_0\ni\BIz\mapsto \langle\BIg_0(\BIz),\BIw_{i}^{(r)}\rangle\in\nobreak\R$,
respectively.
Subsequently, it suffices to tractably formulate each of the constraints 
$(\beta_{0},\BIx,\BIz) \in \epigraph(c_i)$, 
$(\beta_{1},\BIx)\in \hypograph\big(\langle\BIg_i(\,\cdot\,),\BIy_{i}^{(r)}\rangle\big)$, 
$(\beta_{2},\BIz)\in \hypograph\big(\langle\BIg_0(\,\cdot\,),\BIw_{i}^{(r)}\rangle\big)$ 
in (\ref{eqn:oracle-globalmin-epigraph}), e.g., into linear or quadratic constraints possibly involving integer-valued auxiliary variables.
In the following, we will discuss in detail the formulation of 
the constraint 
$(\beta_{2},\BIz)\in \hypograph\big(\langle\BIg_0(\,\cdot\,),\BIw_{i}^{(r)}\rangle\big)$
into linear constraints with binary-valued auxiliary variables.
The formulation of the constraint $(\beta_{1},\BIx)\in \hypograph\big(\langle\BIg_i(\,\cdot\,),\BIy_{i}^{(r)}\rangle\big)$ 
is completely analogous.
Needless to say, the formulation of the constraint $(\beta_{0},\BIx,\BIz) \in \epigraph(c_i)$ depends on the specific form of the cost function~$c_i$.
However, we would like to remark that when $c_i$ has a continuous piece-wise affine structure (which is the case in our Experiment~1 and Experiment~3),
the formulation of the constraint $(\beta_{0},\BIx,\BIz) \in \epigraph(c_i)$
into linear constraints with binary-valued auxiliary variables 
can be carried out via similar techniques.
The formulations of the constraint $(\beta_{0},\BIx,\BIz) \in \epigraph(c_i)$
under the settings of Experiment~1 and Experiment~3 
are presented in 
Section~\ref{sapx:experiments-biz} and Section~\ref{sapx:experiments-1d},
respectively.

Let us 
denote 
$\BIw_{i}^{(r)}=(w_{i,1}^{(r)},\ldots,w_{i,m_0}^{(r)})^\TRANSP$
and denote 
$w_{i,\BIv}^{(r)}:=\sum_{j=1}^{m_0}w_{i,j}^{(r)}\INDI_{\{\BIv=\BIv_{0,j}\}}$
$\forall \BIv\in V(\FS_0)$.
One can adopt any formulation by \citet*{ECvielma2010mixed} (see \citep[Sections~3.1.1, 3.1.2, 3.2.1, 3.2.2, 3.3, \& 3.4]{ECvielma2010mixed})
to represent 
$\hypograph\big(\langle\BIg_0(\,\cdot\,),\BIw_{i}^{(r)}\rangle\big)$.
Let us demonstrate one formulation via the so-called logarithmic disaggregated convex combination (DLog) model \citep[Section~3.1.2]{ECvielma2010mixed}, 
which is the one that we have implemented in our Experiment~1 and Experiment~3. 
Let $B:=\big\lceil\log_2(|\FM(\FS_0)|)\big\rceil$,\footnote{We use 
$\lceil\,\cdot\,\rceil$ to denote the ceiling functions, that is, 
$\lceil a\rceil$ is the smallest integer greater than or equal to $a\in\R$.} 
let $\BIb:\FM(\FS_0)\to\{0,1\}^{B}$ be an arbitrary injective function,
let $b_l(C)\in\{0,1\}$ denote the $l$-th component of $\BIb(C)$ 
for any $C\in\nobreak\FM(\FS_0)$, for $l=\nobreak1,\ldots,B$,
and define 
$\FC_{l,0}:=\nobreak\big\{C\in\FM(\FS_0):b_l(C)=0\big\}$,
$\FC_{l,1}:=\nobreak\big\{C\in\FM(\FS_0):b_l(C)=1\big\}$ 
for $l=1,\ldots,B$.
By introducing at most $(d_0+1)|\FM(\FS_0)|$ continuous auxiliary variables 
$(\lambda_{C,\BIv})_{\BIv\in V(C),\,C\in\FM(\FS_0)}$ and $B$ binary-valued auxiliary variables $(\iota_l)_{l=1:B}$, 
$\hypograph\big(\langle\BIg_0(\,\cdot\,),\BIw_{i}^{(r)}\rangle\big)$ can be formulated as follows:
\begin{align}
    (\beta_{2},\BIz)\in \hypograph\big(\langle\BIg_0(\,\cdot\,),\BIw_{i}^{(r)}\rangle\big) 
    \quad \Leftrightarrow \quad 
    & \exists (\lambda_{C,\BIv})_{\BIv\in V(C),\,C\in\FM(\FS_0)},\; 
    \exists (\iota_l)_{l=1:B}:\label{eqn:DLog-formulation} \\
    & \begin{cases}
        \lambda_{C,\BIv}\ge 0 & \forall \BIv\in V(C),\; \forall C\in\FM(\FS_0),\\
        \iota_l\in\{0,1\} &  \hspace{67.5pt}\forall 1\le l\le B, \\
        \sum_{C\in \FM(\FS_0)}\sum_{\BIv\in V(C)}\lambda_{C,\BIv}=1, \\
        \sum_{C\in\FC_{l,1}}\sum_{\BIv\in V(C)}\lambda_{C,\BIv}\le \iota_l & \hspace{67.5pt}\forall 1\le l\le B, \\
        \sum_{C\in\FC_{l,0}}\sum_{\BIv\in V(C)}\lambda_{C,\BIv}\le 1-\iota_l & \hspace{67.5pt}\forall 1\le l\le B,\\
        \sum_{C\in\FM(\FS_0)}\sum_{\BIv\in V(C)}\lambda_{C,\BIv}\BIv=\BIz, \\
        \sum_{C\in\FM(\FS_0)}\sum_{\BIv\in V(C)}\lambda_{C,\BIv}w_{i,\BIv}^{(r)}\ge \beta_{2}.
    \end{cases}\nonumber
\end{align}

In the special case where $\CX_0$ is a compact one-dimensional interval $[\underline{\kappa},\overline{\kappa}]\subset\R$ and 
$\FM(\FS_0)=\big\{[\kappa_0,\kappa_1],\ldots,[\kappa_{m_0-1},\kappa_{m_0}]\big\}$,
where $\underline{\kappa}=\kappa_0<\kappa_1<\cdots<\kappa_{m_0-1}<\kappa_{m_0}=\overline{\kappa}$, 
one may instead adopt the so-called logarithmic convex combination (Log) model \citep[Section~3.2.2]{ECvielma2010mixed}, which results in fewer auxiliary variables.
Let us denote $\BIw_{i}^{(r)}=(w_{i,1}^{(r)},\ldots,w_{i,m_0}^{(r)})^\TRANSP$, 
let $B:=\big\lceil\log_2(m_0)\big\rceil$, 
let $(\BIb_j)_{j=1:m_0}\subseteq\{0,1\}^B$ be a sequence of distinct binary-valued vectors such that 
${\|\BIb_{j+1}-\BIb_{j}\|_1=1}$ $\forall 1\le\nobreak j\le\nobreak m_0-\nobreak1$,
and denote 
$\BIb_j=(b_{j,1},\ldots,b_{j,B})^\TRANSP$
$\forall 1\le j\le m_0$.
Moreover, 
for $l=1,\ldots,B$, 
let us define 
$\FC_{l,0}\subset \{0,1,\ldots,m_0\}$, 
$\FC_{l,1}\subset \{0,1,\ldots,m_0\}$ as follows:
\begin{align*}
    j\in \FC_{l,\iota} \quad \Leftrightarrow \quad &(j=0 \text{ and } b_{1,l}=\iota) \text{ or } (1\le j\le m_0-1 \text{ and } b_{j,l}=b_{j+1,l}=\iota) \text{ or } (j=m_0 \text{ and } b_{m_0,l}=\iota)\\
    & \hspace{230pt} \forall 0\le j\le m_0,\; \forall \iota\in\{0,1\},\; \forall 1\le l\le B.
\end{align*}
Subsequently, by introducing $m_0+1$ continuous auxiliary variables $(\lambda_j)_{j=0:m_0}$ and 
$B$ binary-valued auxiliary variables $(\iota_l)_{l=1:B}$,
$\hypograph\big(\langle\BIg_0(\,\cdot\,),\BIw_{i}^{(r)}\rangle\big)$ can be formulated as follows:
\begin{align}
    (\beta_{2},z)\in \hypograph\big(\langle\BIg_0(\,\cdot\,),\BIw_{i}^{(r)}\rangle\big) 
    \quad \Leftrightarrow \quad 
    & \exists (\lambda_j)_{j=0:m_0},\; 
    \exists (\iota_l)_{l=1:B}:\label{eqn:Log-formulation} \\
    & \begin{cases}
        \lambda_j\ge 0 & \forall 0\le j\le m_0,\\
        \iota_l\in\{0,1\} & \hspace{6.8pt} \forall 1\le l\le B, \\
        \sum_{j=0}^{m_0}\lambda_{j}=1, \\
        \sum_{j\in\FC_{l,0}}\lambda_{j}\le \iota_l & \hspace{6.8pt} \forall 1\le l\le B, \\
        \sum_{j\in\FC_{l,1}}\lambda_{j}\le 1-\iota_l & \hspace{6.8pt} \forall 1\le l\le B,\\
        \sum_{j=0}^{m_0}\lambda_j\kappa_{j}=z, \\
        \sum_{j=0}^{m_0}\lambda_{j}w_{i,j}^{(r)}\ge \beta_{2}.
    \end{cases}\nonumber
\end{align}

Once (\ref{eqn:oracle-globalmin-epigraph}) has been formulated into 
an MIP problem via the techniques presented above,
it can be efficiently solved by
state-of-the-art MIP solvers such as the Gurobi optimizer \citep{ECgurobi},
which utilize highly optimized branch-and-cut algorithms.
When the tolerance of the so-called absolute MIP gap is set to $\tau^{(r)}$,
branch-and-cut algorithms also naturally compute a lower bound 
$\underline{\beta}_i^{(r)}$, 
along with the computed approximate optimizer 
$(\tilde{\BIx}_{i}^{(r)},\tilde{\BIz}_{i}^{(r)})$, 
that satisfies the condition
$\underline{\beta}_i^{(r)}
\le\nobreak \min_{\BIx\in\CX_i,\,\BIz\in\CX_0}\big\{c_i(\BIx,\BIz)-\langle\BIg_i(\BIx),\BIy_i^{(r)}\rangle-\langle\BIg_0(\BIz),\BIw_i^{(r)}\rangle\big\}
\le \allowbreak c_i(\tilde{\BIx}_i^{(r)},\tilde{\BIz}_i^{(r)})\allowbreak-\langle\BIg_i(\tilde{\BIx}_i^{(r)}),\BIy_i^{(r)}\rangle-\langle\BIg_0(\tilde{\BIz}_i^{(r)}),\BIw_i^{(r)}\rangle
\le \underline{\beta}_i^{(r)} + \tau^{(r)}$
in Line~\ref{alglin:cp-tf-global} of Algorithm~\ref{alg:cp-tf}.

Additionally, we would like to mention that 
certain structures in the cost function $c_i$ can be exploited to 
simplify the formulation of (\ref{eqn:oracle-globalmin}).
For example, suppose that the function 
$\CX_i\ni\nobreak\BIx\mapsto c_i(\BIx,\BIz)\in\R$ is concave for every $\BIz\in\CX_0$.
Then, for every $\BIz\in\CX_0$, the function 
$\CX_i\ni\BIx\mapsto c_i(\BIx,\BIz) - \langle\BIg_i(\BIx),\BIy_{i}^{(r)}\rangle - \langle\BIg_0(\BIz),\BIw_{i}^{(r)}\rangle\in\R$ is continuous and piece-wise concave on each $C\in\FM(\FS_i)$,
and thus attains a minimizer in $V(\FS_i)$.
Consequently, we can simplify (\ref{eqn:oracle-globalmin}) by first solving the sub-problem: 
$\min_{\BIz\in\CX_0}\big\{c_i(\BIv_{i,j},\BIz)-\langle\BIg_0(\BIz),\BIw_{i}^{(r)}\rangle\big\}$
for each $j\in\{0,1,\ldots,m_i\}$
via, for example, the grid discretization method or the MIP-based method presented earlier,
storing the computed minimum values $(\tilde{\beta}_{j})_{j=0:m_i}$, 
and then computing 
$\min_{0\le j\le m_i}\big\{\tilde{\beta}_{j}- \langle\BIg_i(\BIv_{i,j}),\BIy_{i}^{(r)}\rangle\big\}$.
This simplification analogously applies to the case where the function
$\CX_0\ni\BIz\mapsto c_i(\BIx,\BIz)\in\R$ is concave for every $\BIx\in\CX_i$.
This technique is applied in our Experiment~2 
to solve (\ref{eqn:oracle-globalmin}) exactly;
see Section~\ref{sapx:experiments-WB}
for concrete details.

Lastly, the set $\widetilde{\CK}^{(r)}_i$ 
in Line~\ref{alglin:cp-tf-aggregate} of Algorithm~\ref{alg:cp-tf}
is typically chosen to be a set of sub-optimal solutions of the global minimization problem (\ref{eqn:oracle-globalmin}).
For example, in the grid discretization method and in typical branch-and-cut algorithms for MIP problems discussed above,
such a set of sub-optimal solutions is naturally obtained.

\subsection{Remarks about Algorithm~\ref{alg:mt-tf}}\label{sapx:efficient-algo-coupling}

Let us first recall in the following proposition the numerical constructions of $W_1$-optimal couplings 
via classical results about discrete optimal transport 
(see, e.g., \citep[Section~2.3]{ECpeyre2019computational} and \citep[Section~1.3]{ECbenamou2021optimal}), 
semi-discrete optimal transport 
(see, e.g., \citep{EClevy2015numerical} and \citep[Section~5.2]{ECpeyre2019computational}), 
and one-dimensional optimal transport 
(see, e.g., \citep[Section~3.1]{ECrachev1998mass}).

\begin{proposition}[Construction of $W_1$-optimal coupling]\label{prop:OT-construction}%
    Under Setting~\ref{sett:Euclidean},
    let $i\in\nobreak\{0,1,\ldots,N\}$ be arbitrary,
    and let $\mu=\sum_{j=1}^{n}p_j\delta_{\BIx_j}$
    where $n\in\N$, 
    $(p_j)_{j=1:n}\subset(0,1]$ satisfy 
    $\sum_{j=1}^{n}p_j=1$,
    and $(\BIx_j)_{j=1:n}\subseteq \CX_i$ are distinct points.
    Let $(\Omega,\CF,\PROB)$ be a probability space
    and let 
    $\BIX_i:\Omega\to\CX_i$ be a random variable that satisfies 
    $\PROB[\BIX_i=\BIx_j]=\nobreak p_{j}$ $\forall 1\le\nobreak j\le\nobreak n$.
    Moreover, let $\nu\in\CP(\CX_i)$,
    and suppose that any one of the following assumptions hold:
    \begin{enumerate}[label=\normalfont{(A\arabic*)},leftmargin=28pt,beginpenalty=10000]

    \item\label{propc:OT-construction-discrete}%
    \underline{Discrete-to-discrete.}
    $\nu=\sum_{l=1}^{k}q_l\delta_{\BIy_l}$ for $k\in\N$, $(q_l)_{l=1:k}\subset(0,1]$ with 
    $\sum_{l=1}^{k}q_l=\nobreak 1$, and 
    distinct points $(\BIy_l)_{l=1:k}\subseteq \CX_i$.

    \item\label{propc:OT-construction-semidiscrete}%
    \underline{Discrete-to-continuous.}
    $\nu$ is absolutely continuous with respect to the Lebesgue measure on~$\CX_i$. 

    \item\label{propc:OT-construction-1d}%
    \underline{One-dimensional.}
    $\CX_i\subset\R$, i.e., $d_i=1$.

    \end{enumerate}

    Let the random variable $\overbar{\BIX}_i:\Omega\to\CX_i$ be constructed by the procedures below.
    \begin{itemize}[leftmargin=10pt, beginpenalty=10000]
    \item \underline{Discrete-to-discrete.}
    If \ref{propc:OT-construction-discrete} holds,
    let $(\gamma^\star_{j,l})_{j=1:n,\,l=1:k}$ be an optimizer of the following LP problem:
    \begin{align}
        \begin{split}
            \minimize_{(\gamma_{j,l})}\quad & \sum_{j=1}^{n}\sum_{l=1}^{k}\|\BIx_j-\BIy_l\|_2\gamma_{j,l}\\
            \mathrm{subject~to}\quad & \sum_{l=1}^{k} \gamma_{j,l}=p_j \quad\forall 1\le j\le n,\; \qquad 
            \sum_{j=1}^{n} \gamma_{j,l}=q_l \quad\forall 1\le l\le k,\\
            & \gamma_{j,l}\in\R_{+} \quad \forall 1\le j\le n,\;\forall 1\le l\le k.
        \end{split}
        \label{eqn:discrete-OT}
    \end{align}
    Subsequently, let $\overbar{\BIX}_i:\Omega\to\CX_i$ be such that 
    $\PROB[\BIX_i=\BIx_j,\;\overbar{\BIX}_i=\BIy_l]=\gamma^\star_{j,l}$
    $\forall 1\le j\le n$, $\forall 1\le l\le k$.

    \item \underline{Discrete-to-continuous.}
    If \ref{propc:OT-construction-semidiscrete} holds,
    let $\big(\phi^\star_{j}\big)_{j=1:n}\subset\R$ be an optimizer of the following concave maximization problem (which always exists; see, e.g., \citep[Proposition~3.4]{ECneufeld2022v7numerical}):
    $\max_{\phi_{1},\ldots,\phi_{n}\in\R}
    \Big\{\big(\sum_{j=1}^{n}p_{j}\phi_{j}\big)-\int_{\CX_i}\max_{1\le j\le n}\big\{\phi_{j}-\|\BIx_j-\BIx\|_2\big\}\DIFFM{\nu}{\DIFF \BIx}\Big\}$.
    Subsequently, let
    $V_{j}:=\Big\{\BIx\in\nobreak\CX_i:\phi_{j}^\star-\|\BIx_j-\BIx\|_2=\max_{1\le l\le n}\big\{\phi_{l}^\star-\|\BIx_l-\BIx\|_2\big\}\Big\}$
    $\forall 1\le j\le n$, and let 
    $\overbar{\BIX}_i:\Omega\to\CX_i$ be such that 
    $\PROB[\overbar{\BIX}_i\in E|\BIX_i=\BIx_j]=\frac{\nu(E\,\cap\, V_{j})}{\nu(V_{j})}$ $\forall E\in\CB(\CX_i)$, $\forall 1\le j\le n$;
    note that 
    \citep[Proposition~3.4(ii)]{ECneufeld2022v7numerical}
    guarantees $\nu(V_j)=p_j>0$ $\forall 1\le j\le n$.

    \item \underline{One-dimensional.}
    If \ref{propc:OT-construction-1d} holds,
    let 
    $\overbar{\BIX}_i:\Omega\to\CX_i$ be constructed via the following steps.
    \begin{itemize}[beginpenalty=1000]
        \item Define $J:\Omega\to\{1,\ldots,n\}$ by 
        $J:=\sum_{j=1}^{n}j\INDI_{\{\BIX_i=\BIx_j\}}$.

        \item Sort the sequence 
        $(\BIx_{1},\ldots,\BIx_{n})$
        of real numbers into ascending order 
        $\BIx_{\sigma(1)}< \BIx_{\sigma(2)}< \cdots<\nobreak \BIx_{\sigma(n)}$,
        i.e.,
        $\sigma:\{1,\ldots,n\}\to\{1,\ldots,n\}$ is the unique bijection satisfying
        $\BIx_{\sigma(1)}< \BIx_{\sigma(2)}< \cdots< \BIx_{\sigma(n)}$.
        
        \item Define 
        $F_j:=\sum_{l=1}^{\sigma^{-1}(j)-1}p_{\sigma(l)}$ $\forall 1\le j\le n$.
        
        \item Let $U:\Omega\to[0,1]$ be a $\mathrm{Uniform}[0,1]$ random variable that is independent of $\BIX_i$. 
        
        \item Define $\overbar{\BIX}_i:=F^{-1}_{\nu}(F_{J}+p_{J}U)$,
        where
        $F_{\nu}^{-1}(u):=\inf\big\{x\in\CX_i:{\nu\big(\CX_i\cap(-\infty,x]\big)\ge u}\big\}$ 
        $\forall u\in\nobreak[0,1]$.
    \end{itemize}

    \end{itemize}
    Then, in all three cases, 
    the law of $\overbar{\BIX}_i$ is equal to $\nu$, 
    and it holds that
    $\EXP\big[\|\BIX_i-\overbar{\BIX}_i\|_2\big]=W_1(\mu,\nu)$.
\end{proposition}

\proof{Proof of Proposition~\ref{prop:OT-construction}.}
    See Section~\ref{sapx:proof-efficient-algo}.
\endproof

Subsequently, the construction of the random variable $\BIZ_i:\Omega\to\CX_0$ in
Line~\ref{alglin:mt-tf-binding1} of Algorithm~\ref{alg:mt-tf}
can be carried out via Proposition~\ref{prop:OT-construction} since both $\hat{\nu}_{1}$ and $\hat{\nu}_i$ have finite supports.
The construction of the random variable $\overbar{\BIX}_i:\Omega\to\CX_i$ in Line~\ref{alglin:mt-tf-reassembly} can also be 
carried out via Proposition~\ref{prop:OT-construction} due to the assumptions of Setting~\ref{sett:typemeasure}.
Moreover, the random variable $\BIX_i:\Omega\to\CX_i$
in Line~\ref{alglin:mt-tf-binding2} is well-defined due to the finite support of $\hat{\nu}_i$.

In Line~\ref{alglin:mt-tf-primal2},
the existence of the random variable 
$\overbar{\BIZ}:\Omega\to\CX_0$ follows from 
\citep[Proposition~7.33]{ECbertsekas1978stochastic}.
Note that in order to generate a random sample from $\tilde{\nu}$,
one would need to first generate realizations $\BIx_1,\ldots,\BIx_N$
of the random variables $\BIX_1,\ldots,\BIX_N$,
and subsequently solve the problem: 
$\min_{\BIz\in\CX_0}\big\{\sum_{i=1}^{N}c_i(\BIx_i,\BIz)\big\}$.
This is a global minimization problem, 
and the techniques mentioned in Section~\ref{sapx:efficient-algo-global}
can be applied to efficiently solve it.

To numerically implement Algorithm~\ref{alg:mt-tf}, 
the probability space $(\Omega,\CF,\PROB)$ and the constructed random variables 
$\BIZ,\BIX_1,\ldots,\BIX_N,\overbar{\BIZ},\overbar{\BIX}_1,\ldots,\overbar{\BIX}_N$ 
allow us to generate independent random samples from the probability measures $\hat{\nu}$, $(\hat{\gamma}_i)_{i=1:N}$, $\tilde{\nu}$, and $(\tilde{\gamma}_i)_{i=1:N}$.
Hence, one can compute highly accurate approximations of 
$\hat{\alpha}_{\MT}^{\ssUB}$, 
$\tilde{\alpha}_{\MT}^{\ssUB}$, 
$\hat{\epsilon}_{\sssub}$, 
$\tilde{\epsilon}_{\sssub}$ 
via Monte Carlo integration.

\section{Discussions about the connection between matching for teams and multi-marginal optimal transport}\label{apx:mmot}

In this section, we discuss how our results in the paper contribute to the multi-marginal optimal transport (MMOT) literature.
\citet[Proposition~3]{ECcarlier2010matching} have shown that 
\eqref{eqn:mt-primalopt}
is related to the following MMOT problem:
\begin{align}
    \inf_{\mu\in\Gamma(\mu_1,\ldots,\mu_N)}\bigg\{\int_{\CX_1\times\cdots\times\CX_N}\bar{c}\DIFFX{\mu}\bigg\}
    \label{eqn:mt-mmot},
\end{align}
where $\Gamma(\mu_1,\ldots,\mu_N):=\big\{\mu\in\CP(\CX_1\times\cdots\times\CX_N):$ the marginal of $\mu$ on $\CX_i$ is $\mu_i$ $\forall i\in\{1,\ldots,N\}\big\}$,
and 
$\bar{c}(x_1,\ldots,x_N):=\min_{z\in\CX_0}\big\{\sum_{i=1}^{N}c_i(x_i,z)\big\}$
$\forall x_1\in\CX_1,\ldots,\forall x_N\in\CX_N$.
Specifically,
\citep[Proposition~3]{ECcarlier2010matching}
states that the optimal value of (\ref{eqn:mt-mmot}) coincides with 
the optimal value of~\eqref{eqn:mt-primalopt}.
We show in the following corollary that the law of the random variable
$(\bar{X}_1,\ldots,\bar{X}_N):\Omega\to\CX_1\times\cdots\times\CX_N$
in Theorem~\ref{thm:mt-tf-approx}
is an approximate optimizer of (\ref{eqn:mt-mmot}).

\begin{corollary}\label{cor:mmot}
    Under the assumptions of Theorem~\ref{thm:mt-tf-approx},
    the law $\mu\in\CP(\CX_1\times\cdots\times\CX_N)$ of the random variable
    $(\bar{X}_1,\ldots,\bar{X}_N):\Omega\to\CX_1\times\cdots\times\CX_N$ 
    is an $\epsilon_{\ssapx}$-optimizer of (\ref{eqn:mt-mmot}).
\end{corollary}

\proof{Proof of Corollary~\ref{cor:mmot}}
    See Section~\ref{sapx:proof-mmot}.
\endproof

\begin{procedure}[t]
    \DontPrintSemicolon
    \caption{SampMMOT(){$\big((\hat{\gamma}_i)_{i=1:N},(U_t)_{t\in\N}\big)$} $\to$ $\big((\BIx_i)_{i=1:N},\overline{t}\big)$}
    \begin{footnotesize}
    \KwInAndOut{\!\!\!\!\!\begin{tabular}{lll}
        \textbf{Input} &
        $\hat{\gamma}_i\in\CP(\CX_i\times\CX_0)$ & coupling in team allocations encoded by
        $(\hat{\nu},\CX_i,\mu_i,\mathtt{Aux}_i,\BT_i)$ as in \ref{enc:allocations}, where $\hat{\nu}$ is \\
        && encoded by $\big(n_0,d_0,(\BIz_{0,l},\hat{\nu}_{0,l})_{l=0:n_0-1}\big)$ as in \ref{enc:qualitymeasure},
        $\mathtt{Aux}_i$ is encoded by \\
        && $\big(\widetilde{R}_i,(F_{i,\tilde{r}},\tilde{\varpi}_{i,\tilde{r}})_{\tilde{r}=1:\widetilde{R}_i},\BP_i\big)$ as in \ref{enc:typemeasures-aux},
        $\mu_i$ is associated with the oracle procedures \\
        && $\mathtt{Int}$ and $\mathtt{Samp}$ in \ref{enc:typemeasures}, and the entries of $\BT_i$ are denoted by $(t_{l,\tilde{r}})_{l=0:n_0-1,\,\tilde{r}=1:\widetilde{R}_i}$ \\
        & $(U_{t})_{t\in\N}\subset [0,1]$ & an infinite sequence of independent realizations of $\mathrm{Uniform}[0,1]$ random variable \\
        \textbf{Output} &
        $\BIx_i\in\R^{d_i}$ & a component of the generated random sample \\
        & $\overline{t}\in\N_0$ & the number of realizations in $(U_t)_{t\in\N}$ consumed in the procedure
    \end{tabular}\!\!\!\!\!\!\!\!\!\!\!\!\!\!}
    \tcc{Let $T_{\mathtt{Samp}}$ denote the number of arithmetic operations incurred by $\mathtt{Samp}$.}

    \nl%
    $Q\leftarrow \min\big\{l\in\{0,1,\ldots,n_0-1\}:\sum_{s=0}^{l}\hat{\nu}_{0,s}\ge U_1\big\}$, 
    $\overline{t}_0\leftarrow 1$.\tcp*[f]{$O(n_0)$}

    \nl%
    \For(\tcp*[f]{$O\big(N(\widetilde{R}_i+T_{\mathtt{Samp}})\big)$}){$i=1,\ldots,N$}{

        \nl%
        $S_i\leftarrow \min\big\{\tilde{r}\in\{1,\ldots,\widetilde{R}_i\}:\sum_{\tilde{l}=1}^{\tilde{r}}t_{Q,\tilde{l}}\ge U_{\overline{t}_0+1}\big\}$.\tcp*[f]{$O(\widetilde{R}_i)$}

        \nl%
        $(\BIx_i,\overline{t})\leftarrow \mathtt{Samp}\big(\CX_i,\mu_i,F_{i,S_i},\tilde{\varpi}_{i,S_i},(U_t)_{t\ge \overline{t}_0+2}\big)$.\tcp*[f]{$O(T_{\mathtt{Samp}})$}

        \nl%
        $\overline{t}_0\leftarrow \overline{t}_0+\overline{t}+1$.
    }

    \nl%
    \Return{$\big((\BIx_i)_{i=1:N},\overline{t}_0\big)$.}
    \end{footnotesize}
\end{procedure}

Under Setting~\ref{sett:Euclidean},
since Algorithm~\ref{alg:complexity-overall} and Algorithms~\ref{alg:cp-tf}--\ref{alg:mt-tf} use Theorem~\ref{thm:mt-tf-approx} as a template,
we can extend the results in 
Theorem~\ref{thm:complexity-explicit} 
and Theorem~\ref{thm:mt-tf} to approximately solve the MMOT problem (\ref{eqn:mt-mmot}).
Specifically,
$(\hat{\gamma}_i)_{i=1:N}$ returned by Algorithm~\ref{alg:complexity-overall}
contain all the information to encode an $\epsilon$-optimizer 
$\hat{\mu}$ of (\ref{eqn:mt-mmot}),
and a random sample from $\hat{\mu}$ can be generated 
by the procedure $\SampMMOT$ that we present here.
Therefore,
we are able to extend the computational complexity results in Theorem~\ref{thm:complexity-explicit}
to the computation of an $\epsilon$-optimizer of (\ref{eqn:mt-mmot}),
which we present in the following corollary.

\begin{corollary}\label{cor:complexity-mmot}
    Under the assumptions of Theorem~\ref{thm:complexity-explicit}, 
    the following statements hold.
    \begin{enumerate}[label=(\roman*),beginpenalty=10000]
        \item\label{cors:complexity-mmot-general}
        For any $\epsilon\in(0,1)$,
        Algorithm~\ref{alg:complexity-overall} incurs
        $O\big(\big(\frac{64NDd\overline{L}_{c}}{\epsilon}\big)^{(\omega+1)d}N^{\omega+1}J^{\omega+1}d^{\omega+2}\log\big(\frac{NJDd\overline{L}_{c}}{I_{\ssmin}\epsilon}\big)
        +\big(\frac{64NDd\overline{L}_{c}}{\epsilon}\big)^{3d}N^{3}J^{3}d^{4}\log\big(\frac{NJDd\overline{L}_{c}}{I_{\ssmin}\epsilon}\big)T_{\mathtt{Eval}}
        +\big(\frac{128NDd\overline{L}_{c}}{\epsilon}\big)^{d}NJT_{\mathtt{Int}}\big)$
        arithmetic operations
        to compute 
        an $\epsilon$-optimal solution $\hat{\mu}$ 
        of (\ref{eqn:mt-mmot}).
        Subsequently, 
        a random sample from $\hat{\mu}$ can be generated by the procedure $\SampMMOT$, which incurs
        $O\big(\big(\frac{128NDd\overline{L}_{c}}{\epsilon}\big)^{d}NJ+NT_{\mathtt{Samp}}\big)$
        arithmetic operations.

        \item\label{cors:complexity-mmot-particular-case}
        Under the particular class of $(\mu_i)_{i=1:N}$ in Setting~\ref{sett:Euclidean}\ref{settc:Euclidean-particular},
        statement~\ref{cors:complexity-mmot-general}
        holds with respect to 
        $T_{\mathtt{Int}}=O\big(d^{\omega+1}+\log(J)\big)$ and 
        $T_{\mathtt{Samp}}=O\big(d^{\omega+1}+\log(J)\big)$.
    \end{enumerate}
\end{corollary}

\proof{Proof of Corollary~\ref{cor:complexity-mmot}}
    See Section~\ref{sapx:proof-mmot}.
\endproof

Similarly, by letting $\hat{\mu}\leftarrow $ the law of $(\overbar{\BIX}_1,\ldots,\overbar{\BIX}_N)$ in Algorithm~\ref{alg:mt-tf},
we can extend all the results in Theorem~\ref{thm:mt-tf} to the computation of an approximate optimizer of (\ref{eqn:mt-mmot}).
Since this is a direct consequence of Theorem~\ref{thm:mt-tf} and Corollary~\ref{cor:mmot}, we omit the details here.

In summary, for MMOT problems in which the cost function possesses a minimum-of-sum structure:
$\bar{c}(x_1,\ldots,x_N):=\min_{z\in\CX_0}\big\{\sum_{i=1}^{N}c_i(x_i,z)\big\}$
$\forall x_1\in\CX_1,\ldots,\forall x_N\in\CX_N$,
and
the spaces $(\CX_i)_{i=0:N}$,
the marginals $(\mu_i)_{i=1:N}$,
and the functions $(c_i)_{i=1:N}$
satisfy Setting~\ref{sett:Euclidean},
our results in the paper 
provide an explicit computational complexity bound (i.e., Corollary~\ref{cor:complexity-mmot})
and a practically efficient algorithm (i.e., Algorithm~\ref{alg:mt-tf})
for computing approximate optimizers.

\section{Additional details of the numerical experiments}\label{apx:experiments}

\vspace{6pt}

\subsection{Experiment~1: retail business problem}\label{sapx:experiments-biz}

The list below shows the concrete setting of Experiment~1.
\begin{itemize}[leftmargin=10pt, beginpenalty=1000]
\item $N=5$ and
$\CX_0=[-2,2]^2\subset\R^2$.
For $i=1,\ldots,4$, $\CX_i=[-2,2]^2\setminus (-1,1)^2\subset \R^2$, and $\CX_5=[-2,2]\times[-3,-2]\subset\R^2$. Moreover, $(\CX_i)_{i=0:5}$ are equipped with the norm $\|\cdot\|_2$. 

\item For $i=1,\ldots,5$, $\mu_i\in\CP(\CX_i)$ is absolutely continuous with respect to the Lebesgue measure on $\CX_i$ and $\support(\mu_i)=\CX_i$. The probability density function of $\mu_i$ is a continuous piece-wise affine function on $\CX_i$. 

\item For $i=1,\ldots,4$, $c_i:\CX_i\times\CX_0\to\R$ is defined by
\begin{align*}
    c_i(\BIx,\BIz)&:=\min_{1\le j\le 5,\,1\le j'\le 5}\big\{c_{\mathrm{walk}}\|\BIx-\BIu_j\|_1+c_{\mathrm{walk}}\|\BIz-\BIu_{j'}\|_1+c_{\mathrm{train}}|j-j'|\big\} \wedge c_{\mathrm{walk}}\|\BIx-\BIz\|_1 \\
    &\hspace{295pt} \forall (\BIx,\BIz)\in\CX_i\times\CX_0,
\end{align*}
where $c_{\mathrm{walk}}=0.15$, $c_{\mathrm{train}}=0.015$.
Moreover, 
$c_5:\CX_5\times\CX_0\to \R$ is defined by $c_5(\BIx,\BIz):=c_{\mathrm{restock}}\|\BIx-\nobreak\BIz\|_1$ $\forall(\BIx,\BIz)\in\CX_5\times\CX_0$,
where $c_{\mathrm{restock}}=0.4$.
\end{itemize}
Observe that 
\begin{align*}
    \big|c_i(\BIx,\BIz)-c_i(\BIx',\BIz')\big|
    &\le c_{\mathrm{walk}}\big(\|\BIx-\nobreak\BIx'\|_1+\|\BIz-\BIz'\|_1\big)\\
    &\le \sqrt{2}c_{\mathrm{walk}}\big(\|\BIx-\nobreak\BIx'\|_2+\|\BIz-\BIz'\|_2\big)
    & \forall (\BIx,\BIz),(\BIx',\BIz')\in\CX_i\times\CX_0,\; 
    \forall 1\le i\le 4,
\end{align*}
and that 
\begin{align*}
    \big|c_5(\BIx,\BIz)-c_5(\BIx',\BIz')\big|
    &\le \sqrt{2}c_{\mathrm{restock}}\big(\|\BIx-\nobreak\BIx'\|_2+\|\BIz-\BIz'\|_2\big)
    & \forall (\BIx,\BIz),(\BIx',\BIz')\in\CX_5\times\CX_0.
\end{align*}
Therefore, one checks that 
Setting~\ref{sett:Euclidean} holds with respect to 
$L_{c}\leftarrow \sqrt{2}(c_{\mathrm{walk}}\vee c_{\mathrm{restock}})$.

In the following, we will present the mixed-integer programming (MIP) formulation of the global minimization problem 
$\min_{\BIx\in\CX_i,\,\BIz\in\CX_0}\big\{c_i(\BIx,\BIz)-\langle\BIg_i(\BIx),\BIy_i^{(r)}\rangle-\langle\BIg_0(\BIz),\BIw_i^{(r)}\rangle\big\}$
in Line~\ref{alglin:cp-tf-global} of Algorithm~\ref{alg:cp-tf}.
As mentioned in Section~\ref{sapx:efficient-algo-global},
we reformulate it via the epigraph 
$\epigraph(c_i)$ and the hypographs 
$\hypograph\big(\langle\BIg_i(\,\cdot\,),\BIy_i^{(r)}\rangle\big)$,
$\hypograph\big(\langle\BIg_0(\,\cdot\,),\BIw_i^{(r)}\rangle\big)$
as in (\ref{eqn:oracle-globalmin-epigraph}),
where the constraints 
$(\beta_{1},\BIx)\in \hypograph\big(\langle\BIg_i(\,\cdot\,),\BIy_{i}^{(r)}\rangle\big)$ and 
$(\beta_{2},\BIz)\in \hypograph\big(\langle\BIg_0(\,\cdot\,),\BIw_{i}^{(r)}\rangle\big)$
are formulated into linear constraints with binary-valued auxiliary variables
via the so-called logarithmic disaggregated convex combination 
(DLog) model
\citep[Section~3.1.2]{ECvielma2010mixed}
as in (\ref{eqn:DLog-formulation}).

For $i\in\{1,\ldots,4\}$,
observe that 
$c_i$ is the minimum of 21~convex functions,
and thus the constraint
$(\beta_0,\BIx,\BIz)\in\epigraph(c_i)$
can be formulated 
by introducing 21~continuous auxiliary variables 
$s_{0,0},(s_{j,j'})_{j,j'\in\{1,\ldots,5\},\,j\ne j'}$
and 21~binary-valued auxiliary variables
$\iota_{0,0},(\iota_{j,j'})_{j,j'\in\{1,\ldots,5\},\,j\ne j'}$
as follows:
\begin{align*}
    &(\beta_{0},\BIx,\BIz) \in \epigraph(c_i) \\
     \Leftrightarrow \quad & \exists s_{0,0},\; \exists\iota_{0,0},\;\exists (s_{j,j'}, \iota_{j,j'})_{j,j'\in\{1,\ldots,5\},\,j\ne j'}:\\
    &\begin{cases}
        s_{0,0} \ge 0,\; \iota_{0,0} \in \{0,1\}, \\
        s_{j,j'}\ge 0,\; \iota_{j,j'}\in\{0,1\} 
        & \forall j,j'\in\{1,\ldots,5\},\;j\ne j',\\
        \beta_{0} \ge c_{\mathrm{walk}}\|\BIx-\BIz\|_1 - s_{0,0},\\
        \beta_{0} \ge c_{\mathrm{walk}}\|\BIx-\BIu_j\|_1+c_{\mathrm{walk}}\|\BIz-\BIu_{j'}\|_1 + c_{\mathrm{train}}|j-j'| - s_{j,j'}
        & \forall j,j'\in\{1,\ldots,5\},\;j\ne j',\\
        s_{0,0} \le M(1-\iota_{0,0}), \\
        s_{j,j'} \le M(1 - \iota_{j,j'}) 
        & \forall j,j'\in\{1,\ldots,5\},\;j\ne j',\\
        \iota_{0,0} + \sum_{j,j'\in\{1,\ldots,5\},\, j\ne j'}\iota_{j,j'} = 1,
    \end{cases}
\end{align*}
where $M:=c_{\mathrm{walk}}\Big[2\max_{j,j'\in\{1,\ldots,5\},\, j\ne j'}\big\{\|\BIu_j-\BIu_{j'}\|_1\big\}\vee 2\Big(\max_{\BIz\in\CX_0}\big\{\|\BIz\|_1\big\}+\max_{1\le j\le 5}\big\{\|\BIu_j\|_1\big\}\Big)\Big]+4c_{\mathrm{train}}$.
Note that the non-linear constraints 
$\beta_{0} \ge c_{\mathrm{walk}}\|\BIx-\BIz\|_1-s_{0,0}$ and
$\beta_{0} \ge c_{\mathrm{walk}}\|\BIx-\nobreak\BIu_j\|_1+c_{\mathrm{walk}}\|\BIz-\nobreak\BIu_{j'}\|_1 + c_{\mathrm{train}}|j-j'|-s_{j,j'}$ 
$\forall j,j'\in\{1,\ldots,5\}$, $j\ne j'$
can be linearized through introducing 22~additional continuous auxiliary variables.
In total,
our MIP formulation 
of $\min_{\BIx\in\CX_i,\,\BIz\in\CX_0}\big\{c_i(\BIx,\BIz)-\langle\BIg_i(\BIx),\BIy_i^{(r)}\rangle-\langle\BIg_0(\BIz),\BIw_i^{(r)}\rangle\big\}$
involves 
$\big(3\big(|\FM(\FS_i)|+|\FM(\FS_0)|\big)+48\big)$ continuous decision variables and
$\big(\big\lceil\log_2(|\FM(\FS_i)|)\big\rceil+\big\lceil\log_2(|\FM(\FS_0)|)\big\rceil+21\big)$ binary-valued decision variables 
for $i\in\nobreak\{1,\ldots,4\}$.
In the $i=5$ case, 
it is straightforward to formulate 
the constraint $(\beta_0,\BIx,\BIz)\in\epigraph(c_i)$
using 2~continuous auxiliary variables and no binary-valued auxiliary variable,
and thus 
our MIP formulation 
of $\min_{\BIx\in\CX_5,\,\BIz\in\CX_0}\big\{c_5(\BIx,\BIz)-\langle\BIg_5(\BIx),\BIy_5^{(r)}\rangle-\langle\BIg_0(\BIz),\BIw_5^{(r)}\rangle\big\}$
involves 
$\big(3\big(|\FM(\FS_i)|+|\FM(\FS_0)|\big)+7\big)$ continuous decision variables and
$\big(\big\lceil\log_2(|\FM(\FS_i)|)\big\rceil+\big\lceil\log_2(|\FM(\FS_0)|)\big\rceil\big)$ binary-valued decision variables.
Overall,
the largest instance of MIP problem solved in Experiment~1
contains 43440 continuous decision variables and 47 binary-valued decision variables.

The resulting linear MIP problem is subsequently solved 
by the MIP solver of the Gurobi optimizer \citep{ECgurobi}, 
where we set the absolute MIP gap tolerance to
$\tau^{(r)}\leftarrow 10^{-10}$ $\forall r\in\N_{0}$.

\subsection{Experiment~2: 2-Wasserstein barycenter}\label{sapx:experiments-WB}

The list below shows the concrete setting of Experiment~2.
\begin{itemize}[leftmargin=10pt, beginpenalty=1000]
    
\item $N=20$. 
For $i=1,\ldots,20$, $\CX_i\subset\R^2$ is the union of two triangles.
Moreover, $\CX_0\subset\R^2$ is the Minkowski sum $\frac{1}{20}\sum_{i=1}^{20}\conv(\CX_i)$.
Furthermore, $(\CX_i)_{i=0:20}$ are equipped with the norm $\|\cdot\|_2$.

\item For $i=1,\ldots,20$, $\mu_i\in\CP(\CX_i)$ is absolutely continuous with respect to the Lebesgue measure on $\CX_i$. The probability density function of $\mu_i$ is a continuous piece-wise affine function on $\CX_i$.

\item For $i=1,\ldots,20$, $c_i:\CX_i\times\CX_0\to\R$ is given by $c_i(\BIx,\BIz):=\frac{1}{N}\big(\|\BIz\|^2_2-2\langle\BIx,\BIz\rangle\big)=\frac{1}{N}\|\BIx-\BIz\|_2^2-\frac{1}{N}\|\BIx\|_2^2$ $\forall\BIx\in\CX_i$, $\forall\BIz\in\CX_0$.

\end{itemize}
Note that we have subtracted the term $\frac{1}{N}\|\BIx\|_2^2$
from the quadratic cost function $\frac{1}{N}\|\BIx-\BIz\|_2^2$ in $c_i(\BIx,\BIz)$ for $i=1,\ldots,20$,
since these terms will only shift the common optimal value of \eqref{eqn:mt-primalopt} and \eqref{eqn:mt-dualopt} by an additive constant 
$C:=\frac{1}{N}\sum_{i=1}^N\int_{\CX_i}\|\BIx\|_2^2\DIFFM{\mu_i}{\DIFF\BIx}$
and have no impact on the matching equilibria.
Because of this, 
we add the constant~$C$ to the lower and upper bounds 
$\hat{\alpha}_{\MT}^{\ssLB}$, 
$\hat{\alpha}_{\MT}^{\ssUB}$,
$\tilde{\alpha}_{\MT}^{\ssUB}$
computed by Algorithm~\ref{alg:mt-tf}
so that they correspond to lower and upper bounds for  
the optimal value of the \mbox{2-Wasserstein} barycenter problem, i.e., 
$\inf_{\nu\in\CP(\CX_0)}\big\{\frac{1}{N}\sum_{i=1}^N W_2(\mu_i,\nu)^2\big\}$.
Moreover, observe that 
\begin{align*}
    \big|c_i(\BIx,\BIz)-c_i(\BIx',\BIz')\big|
    &\le \frac{2}{N}\Big(\max_{\bar{\BIx}\in\CX_i}\big\{\|\bar{\BIx}\|_2\big\}
    +\max_{\bar{\BIz}\in\CX_0}\big\{\|\bar{\BIz}\|_2\big\}\Big) \|\BIz-\BIz'\|_2
    + \frac{2}{N}\max_{\bar{\BIz}\in\CX_0}\big\{\|\bar{\BIz}\|_2\big\} \|\BIx-\BIx'\|_2\\
    & \hspace{180pt} \forall (\BIx,\BIz),(\BIx',\BIz')\in\CX_i\times\CX_0,\; \forall 1\le i\le 20,
\end{align*}
and hence Setting~\ref{sett:Euclidean} holds with respect 
to $L_{c}\leftarrow \frac{4}{N}\max_{0\le i\le N}\big\{\max_{\BIx\in\CX_i}\big\{\|\BIx\|_2\big\}\big\}$.

Next, let us present the detailed implementation of 
the global minimization problem in Line~\ref{alglin:cp-tf-global} of Algorithm~\ref{alg:cp-tf}.
Notice that since the function 
$\CX_i\ni\BIx\mapsto c_i(\BIx,\BIz)\in\R$ is concave for every $\BIz\in\CX_0$, 
we can simplify $\min_{\BIx\in\CX_i,\,\BIz\in\CX_0}\big\{c_i(\BIx,\BIz)-\langle\BIg_i(\BIx),\BIy_i^{(r)}\rangle-\langle\BIg_0(\BIz),\BIw_i^{(r)}\rangle\big\}$ into 
$\min_{\BIx\in V(\FS_i)}\big\{\min_{\BIz\in\CX_0}\big\{c_i(\BIx,\BIz)-\langle\BIg_0(\BIz),\BIw_i^{(r)}\rangle\big\} - \langle\BIg_i(\BIx),\BIy_i^{(r)}\rangle\big\}$, as discussed in Section~\ref{sapx:efficient-algo-global}.
Moreover, since the function 
$\CX_0\ni\nobreak\BIz\mapsto \langle\BIg_0(\BIz),\BIw_i^{(r)}\rangle\in\R$ is continuous and piece-wise affine on $\relint(K)$ for each non-empty $K\in\FS_0$, 
and since $\bigcup_{K\in\FS_0,\,K\ne\emptyset}\relint(K)=\CX_0$,
we can further simplify 
$\min_{\BIz\in\CX_0}\big\{c_i(\BIx,\BIz)-\langle\BIg_0(\BIz),\BIw_i^{(r)}\rangle\big\}$
into 
$\min_{K\in\FS_0,\,K\ne\emptyset}\big\{\min_{\BIz\in \relint(K)}\big\{c_i(\BIx,\BIz)-\langle\BIg_0(\BIz),\BIw_i^{(r)}\rangle\big\}\big\}$.
Since each non-empty $K\in\FS_0$
is a single point if $\dim(K)=0$, 
a line segment if $\dim(K)=1$, 
or a triangle if $\dim(K)=2$,
the innermost minimization problem $\min_{\BIz\in \relint(K)}\big\{c_i(\BIx,\BIz)-\langle\BIg_0(\BIz),\BIw_i^{(r)}\rangle\big\}$
can be solved exactly in each of these cases.

\renewcommand{\thealgocf}{EC.\arabic{algocf}}
\setcounter{algocf}{0}
\begin{algorithm}[t]
\caption{\textbf{Implementation of Line~\ref{alglin:cp-tf-global}
of Algorithm~\ref{alg:cp-tf} in Experiment~2}}\label{alg:WB-global}
\begin{footnotesize}

\KwIn{$N$, $\CX_i$, $\FS_i$, $\CG_i$, $\CX_0$, $\FS_0$, $\CG_0$ (see Setting~\ref{sett:Euclidean}), $\BIy_i^{(r)}\in\R^{m_i}$, $\BIw_i^{(r)}\in\R^{m_0}$.}

\nl%
\For{$j=0,1,\ldots,m_i$}{

    \nl%
    \For{each $K\in \FS_0$}{

        \nl\uIf{$K=\emptyset$}{
            \nl $\beta_{j,K}\leftarrow \infty$, $\BIz_{j,K}\leftarrow \veczero$.
        }
        \nl\uElseIf{$\dim(K)=0$}{
            \nl%
            $\BIv_0\leftarrow$ the unique point in $K$.

            \nl%
            $\beta_{j,K}\leftarrow c_i(\BIv_{i,j},\BIv_0)-\langle\BIg_0(\BIv_0),\BIw_i^{(r)}\rangle$,
            $\BIz_{j,K}\leftarrow \BIv_0$.
        }
        \nl\uElseIf{$\dim(K)=1$}{
            \nl%
            $(\BIv_0,\BIv_1)\leftarrow$ any enumeration of $V(K)$.

            \nl\label{alglin:WB-global-linemin}%
            $\lambda_0^{\star}\leftarrow \argmin_{\lambda_0\in\R}\Big\{
                \frac{1}{N}\big\|\lambda_0\BIv_0+(1-\lambda_0)\BIv_1\big\|_2^2
                -\frac{2}{N}\langle\BIv_{i,j},\lambda_0\BIv_0+(1-\lambda_0)\BIv_1\rangle
                -\lambda_0\langle\BIg_0(\BIv_0),\BIw_i^{(r)}\rangle
                -\allowbreak(1-\lambda_0)\langle\BIg_0(\BIv_1),\BIw_i^{(r)}\rangle
            \Big\}$.

            \nl\uIf{$0<\lambda_0^\star<1$}{
                \nl%
                $\BIz_{j,K}\leftarrow \lambda_0^\star\BIv_0+(1-\lambda_0^\star)\BIv_1$.

                \nl%
                $\beta_{j,K}\leftarrow 
                \frac{1}{N}\|\BIz_{j,K}\|_2^2
                -\frac{2}{N}\langle\BIv_{i,j},\BIz_{j,K}\rangle
                -\langle\BIg_0(\BIz_{j,K}),\BIw_i^{(r)}\rangle$.
            }
            \nl\Else{
                \nl%
                $\beta_{j,K}\leftarrow \infty$,
                $\BIz_{j,K}\leftarrow \veczero$.
            }
        }
        \nl\Else{
            \nl%
            $(\BIv_0,\BIv_1,\BIv_2)\leftarrow$ any enumeration of $V(K)$.

            \nl\label{alglin:WB-global-trianglemin}%
            $(\lambda_0^\star,\lambda_1^\star)\leftarrow \argmin_{\lambda_0\in\R,\,\lambda_1\in\R}\Big\{
                \frac{1}{N}\big\|\lambda_0\BIv_0+\lambda_1\BIv_1+(1-\lambda_0-\lambda_1)\BIv_2\big\|_2^2
                -\frac{2}{N}\langle\BIv_{i,j},\lambda_0\BIv_0+\lambda_1\BIv_1+\allowbreak(1-\lambda_0-\lambda_1)\BIv_2\rangle
                -\lambda_0\langle\BIg_0(\BIv_0),\BIw_i^{(r)}\rangle-\lambda_1\langle\BIg_0(\BIv_1),\BIw_i^{(r)}\rangle
                -\allowbreak(1-\lambda_0-\lambda_1)\langle\BIg_0(\BIv_2),\BIw_i^{(r)}\rangle
            \Big\}$.

            \nl\uIf{$\lambda_0^\star>0$, $\lambda_1^\star>0$, $1-\lambda_0^\star-\lambda_1^\star>0$}{
                \nl%
                $\BIz_{j,K}\leftarrow \lambda_0^\star\BIv_0+\lambda_1^\star\BIv_1+(1-\lambda_0^\star-\lambda_1^\star)\BIv_2$.

                \nl%
                $\beta_{j,K}\leftarrow 
                \frac{1}{N}\|\BIz_{j,K}\|_2^2
                -\frac{2}{N}\langle\BIv_{i,j},\BIz_{j,K}\rangle
                -\langle\BIg_0(\BIz_{j,K}),\BIw_i^{(r)}\rangle$.
            }
            \nl\Else{
                \nl%
                $\beta_{j,K}\leftarrow \infty$,
                $\BIz_{j,K}\leftarrow \veczero$.
            }

        }
    }
}

\nl\label{alglin:WB-global-overallmin}%
$\underline{\beta}_i^{(r)}\leftarrow \min_{j\in\{0,1,\ldots,m_i\},\, K\in \FS_0}\big\{\beta_{j,K}-\langle\BIg_i(\BIv_{i,j}),\BIy_i^{(r)}\rangle\big\}$,
let $(j^\star,K^\star)\in \argmin_{j\in\{0,1,\ldots,m_i\},\, K\in \FS_0}\big\{\beta_{j,K}-\allowbreak\langle\BIg_i(\BIv_{i,j}),\BIy_i^{(r)}\rangle\big\}$.

\nl%
$(\tilde{\BIx}_{i}^{(r)},\tilde{\BIz}_{i}^{(r)})\leftarrow (\BIv_{i,j^\star},\BIz_{j^\star,K^\star})$.

\KwOut{$(\tilde{\BIx}_i^{(r)},\tilde{\BIz}_i^{(r)})$,
$\underline{\beta}_i^{(r)}$.} 
\end{footnotesize}
\end{algorithm}

We present the detailed implementation of Line~\ref{alglin:cp-tf-global}
of Algorithm~\ref{alg:cp-tf} in 
Algorithm~\ref{alg:WB-global}.
Note that the minimization problems in Line~\ref{alglin:WB-global-linemin} and Line~\ref{alglin:WB-global-trianglemin}
are unconstrained minimizations of convex quadratic functions and can be solved in closed-form.
In addition, a finite set $\widetilde{\CK}^{(r)}_i\subseteq\CX_i\times\CX_0$ of sub-optimal solutions of $\min_{\BIx\in\CX_i,\,\BIz\in\CX_0}\big\{c_i(\BIx,\BIz)-\langle\BIg_i(\BIx),\BIy_i^{(r)}\rangle-\langle\BIg_0(\BIz),\BIw_i^{(r)}\rangle\big\}$ in Line~\ref{alglin:cp-tf-aggregate} of Algorithm~\ref{alg:cp-tf}
can be naturally generated by Algorithm~\ref{alg:WB-global}
by first gathering sub-optimal solutions $\widetilde{\CI}\subseteq \{0,1,\ldots,m_i\}\times\nobreak\FS_0$ of 
$\min_{j\in\{0,1,\ldots,m_i\},\, K\in \FS_0}\big\{\beta_{j,K}-\langle\BIg_i(\BIv_{i,j}),\BIy_i^{(r)}\rangle\big\}$
in Line~\ref{alglin:WB-global-overallmin},
and then setting 
$\widetilde{\CK}_i^{(r)}\leftarrow \big\{(\BIv_{i,j},\BIz_{j,K}):(j,K)\in \widetilde{\CI}\big\}$.

In the following, let us present the detailed procedures that we carried out in order to compare Algorithm~\ref{alg:mt-tf} with five state-of-the-art \mbox{2-Wasserstein} barycenter algorithms.
The first four algorithms by \citet{ECstaib2017parallel}, \citet{ECfan2020scalable}, \citet{ECkorotin2021continuous}, and \citet{ECkorotin2022wasserstein} only take independent random samples from $(\mu_i)_{i=1:N}$ as inputs in the process of approximating the \mbox{2-Wasserstein} barycenter.
Hence, we pre-generated a large collection of independent random samples from $(\mu_i)_{i=1:N}$, and we used this collection in the computation of approximate \mbox{2-Wasserstein} barycenters via these four algorithms.
The MMOT formulation based algorithm of \citet{ECneufeld2022v7numerical} takes the integrals of continuous test functions with respect to $(\mu_i)_{i=1:N}$ as inputs, which is the same as Algorithm~\ref{alg:mt-tf}.

The code for the algorithm of \citet{ECstaib2017parallel} is at: \url{https://github.com/mstaib/stochastic-barycenter-code}.
We fix the support of the approximate \mbox{2-Wasserstein} barycenter to the finite set $V(\FS_0)$, where $\FS_0$ is the subdivision of $\FK_0$ with the smallest mesh size $\eta(\FS_0)$ under Setting~\ref{sett:Euclidean}. 
The resulting size of the fixed support is $|V(\FS_0)|=9020$.
The algorithm was run with one master thread and 20 worker threads each representing an input measure $\mu_i$, where each worker thread was run for $10^6$ iterations, resulting in $2\times10^{7}$ total iterations. 
We took the approximate \mbox{2-Wasserstein} barycenter computed after $10^7$ total iterations instead of waiting until all $2\times 10^7$ total iterations are completed.
This was to prevent an issue where some worker threads would stop gradient-based updates before all $2\times 10^7$ total iterations were completed due to being faster than the other worker threads.

The code for the algorithm of \citet{ECfan2020scalable} is at: \url{https://github.com/sbyebss/Scalable-Wasserstein-Barycenter}.
We use the same neural network configurations, the same training configurations, and the same sample sizes for training the neural networks as used by the file \texttt{G2G\_sameW\_3loop.py}.
We ran Algorithm~1 of \citep{ECfan2020scalable} for 500 iterations (the outermost for-loop). 
Subsequently, we randomly generated 10000 independent samples from the trained generative neural network (GNN) and evaluated the resulting empirical measure as the approximate \mbox{2-Wasserstein} barycenter.

The code for the algorithm of \citet{ECkorotin2021continuous} is at: \url{https://github.com/iamalexkorotin/Wasserstein2Barycenters}.
We adopt the same setting as the file \texttt{CW2B\_toy\_experiments.ipynb}, where we use the same configurations of the input convex neural networks (ICNNs), the same training configurations, and the same sample sizes for training. 
We ran Algorithm~1 of \citep{ECkorotin2021continuous} for 3000 iterations (the outermost for-loop). 
Next, we randomly generated 10000 independent samples from the pushforward of $\mu_1$ by the gradient of the convex conjugate of the ICNN associated with $\mu_1$, i.e., $\nabla\psi_1^{\dagger}\sharp\mu_1$ in the notation of \citep{ECkorotin2021continuous},
and we evaluated the resulting empirical measure as the approximate \mbox{2-Wasserstein} barycenter.

The code for the algorithm of \citet{ECkorotin2022wasserstein} is at: \url{https://github.com/iamalexkorotin/WassersteinIterativeNetworks}.
The neural network configurations, the training configurations, and the sample sizes for training are taken to be identical to those specified in the file \texttt{WIN\_location\_scatter.ipynb}.
We ran Algorithm~1 of \citep{ECkorotin2022wasserstein} for 300 iterations (the outermost loop). 
Subsequently, we randomly generated 10000 independent samples from the trained GNN and evaluated the resulting empirical measure as the approximate \mbox{2-Wasserstein} barycenter.

The code for the algorithm of \citet{ECneufeld2022v7numerical} is at:
\url{https://github.com/qikunxiang/MultiMarginalOptimalTransport}.
Since Algorithm~2 of \citet{ECneufeld2022v7numerical} takes integrals of finite sets of continuous test functions on $(\CX_i)_{i=1:N}$ with respect to $(\mu_i)_{i=1:N}$ as inputs, 
we took the combinations of test functions $(\CG_i)_{i=1:N}$ used in our algorithm, and used them in Algorithm~2 of \citet{ECneufeld2022v7numerical}.
However, the largest instance of linear semi-infinite programming (LSIP) formulation that we managed to solve with their algorithm contains $|\CG_i|=288$ test functions for each~$\mu_i$.
For LSIP instances larger than this, the convergence of their algorithm became extremely slow. 
Therefore, the results shown in Section~\ref{ssec:experiment-WB} are based on this largest LSIP formulation.

\subsection{Experiment~3: one-dimensional type spaces}\label{sapx:experiments-1d}

In this subsection, 
we will present the MIP formulation of the global minimization problem 
$\min_{x\in\CX_i,\,\BIz\in\CX_0}\big\{c_i(x,\BIz)-\langle\BIg_i(x),\BIy_i^{(r)}\rangle-\langle\BIg_0(\BIz),\BIw_i^{(r)}\rangle\big\}$
in Line~\ref{alglin:cp-tf-global} of Algorithm~\ref{alg:cp-tf}.
As mentioned in Section~\ref{sapx:efficient-algo-global},
we reformulate it via the epigraph 
$\epigraph(c_i)$ and the hypographs 
$\hypograph\big(\langle\BIg_i(\,\cdot\,),\BIy_i^{(r)}\rangle\big)$,
$\hypograph\big(\langle\BIg_0(\,\cdot\,),\BIw_i^{(r)}\rangle\big)$
as in (\ref{eqn:oracle-globalmin-epigraph}),
where the constraint
$(\beta_{1},x)\in \hypograph\big(\langle\BIg_i(\,\cdot\,),\BIy_{i}^{(r)}\rangle\big)$ 
is formulated into linear constraints with binary-valued auxiliary variables 
via the so-called logarithmic convex combination (Log) model 
\citep[Section~3.2.2]{ECvielma2010mixed}
as in (\ref{eqn:Log-formulation}) because 
$\CX_i$ is a compact one-dimensional interval,
and the constraint
$(\beta_{2},\BIz)\in \hypograph\big(\langle\BIg_0(\,\cdot\,),\BIw_{i}^{(r)}\rangle\big)$
is formulated into linear constraints with binary-valued auxiliary variables
via the so-called logarithmic disaggregated convex combination 
(DLog) model
\citep[Section~3.1.2]{ECvielma2010mixed}
as in (\ref{eqn:DLog-formulation}).
The cost function $c_i(x,\BIz):=\frac{1}{N}\big(\big(|x-\langle\BIs_i,\BIz\rangle|\wedge \kappa_{i,2}\big)-\kappa_{i,1}\big)^+$ $\forall x\in\CX_i$, $\forall\BIz\in\CX_0$ in Experiment~3
can be represented as $c_i(x,\BIz)=l_i(x-\langle\BIs_i,\BIz\rangle)$ $\forall x\in\CX_i$, $\forall \BIz\in\CX_0$, where 
$\R\ni f\mapsto l_i(f):=\frac{1}{N}\big(\big(|f|\wedge \kappa_{i,2}\big)-\kappa_{i,1}\big)^+\in\R$ is a continuous piece-wise affine function.
Thus, we can formulate the constraint $(\beta_0,x,\BIz)\in\epigraph(c_i)$ as follows:
\begin{align*}
    (\beta_{0},x,\BIz)\in\epigraph(c_i) \quad \Leftrightarrow \quad & \exists f:\; \begin{cases}
        f=x-\langle\BIs_i,\BIz\rangle,\\
        (\beta_{0},f) \in \epigraph(l_i):=\big\{(\beta,f'):l_i(f')\le \beta\big\}.
    \end{cases}
\end{align*}
Since 
$\CX_i=[0,1]$, $\CX_0=\big\{(z_1,z_2)^\TRANSP:z_1+z_2 \le 1,\; z_1,z_2\in\R_{+}\big\}$, and $\|\BIs_i\|_2=1$,
we can restrict the domain of $l_i$ to $[-1,2]$ and formulate 
the constraint $(\beta_0,f)\in\epigraph(l_i)$ 
via the Log formulation \citep[Section~3.2.2]{ECvielma2010mixed} analogous to (\ref{eqn:Log-formulation}).
Overall, our MIP formulation of 
$\min_{x\in\CX_i,\,\BIz\in\CX_0}\big\{c_i(x,\BIz)-\langle\BIg_i(x),\BIy_i^{(r)}\rangle-\langle\BIg_0(\BIz),\BIw_i^{(r)}\rangle\big\}$
involves 
$\big(m_i+3|\FM(\FS_0)|+10\big)=3131$ continuous decision variables and 
$\big(\big\lceil\log_2(m_i)\big\rceil+\big\lceil\log_2(|\FM(\FS_0)|)\big\rceil+3\big)=19$ binary-valued decision variables.

Same as in Experiment~1, the resulting linear MIP problem is subsequently solved via the Gurobi optimizer \citep{ECgurobi}, where we set the absolute MIP gap tolerance to $\tau^{(r)}\leftarrow 10^{-10}$ $\forall r\in\N_0$.

\section{Proofs of theoretical results}\label{apx:proof}
\vspace{6pt}

\subsection{Proofs of results in Section~\ref{sec:model}}\label{sapx:proof-discussion-assumption}

\proof{Proof of Theorem~\ref{thm:mt-tf-equilibrium}.}
To begin, 
it holds by definition that
\begin{align}
    \tilde{\varphi}_i^{(l)}(z)&\ge \varphi_i^{(l)}(z) & \forall z\in\CX_0,\; \forall 1\le i\le N-1,\; \forall l\in\N,\label{eqn:mt-tf-equilibrium-proof-doublectrans}\\
    \tilde{\varphi}{_i^{(l)}}^{c_i}(x)&={\varphi_i^{(l)}}^{c_i}(x) & \forall x\in\CX_i,\; \forall 1\le i\le N-1,\; \forall l\in\N.\label{eqn:mt-tf-equilibrium-proof-triplectrans}
\end{align}
Moreover, for any $l\in\N$ and for $i=1,\ldots,N-1$, it follows from the $L_c$-Lipschitz continuity of $c_i$ that 
\begin{align*}
    \tilde{\varphi}_i^{(l)}(z)\le c_i(x,z) - {\varphi_i^{(l)}}^{c_i}(x)
    \le c_i(x,z') + L_cd_{\CX_0}(z,z') - {\varphi_i^{(l)}}^{c_i}(x) 
    \qquad \forall x\in\CX_i,\; \forall z,z'\in\CX_0.
\end{align*}
Thus, we get 
\begin{align*}
    \tilde{\varphi}_i^{(l)}(z)-L_cd_{\CX_0}(z,z') 
    \le \inf_{x\in\CX_i}\big\{c_i(x,z') - {\varphi_i^{(l)}}^{c_i}(x)\big\}
    =\tilde{\varphi}_i^{(l)}(z') \qquad \forall z,z'\in\CX_0,
\end{align*}
which shows that $\tilde{\varphi}_i^{(l)}$ is $L_c$-Lipschitz continuous
for $i=1,\ldots,N-1$ and for all $l\in\N$.
Consequently, it holds that $\tilde{\varphi}_N^{(l)}$
is $(N-1)L_c$-Lipschitz continuous for all $l\in\N$,
and that (\ref{eqn:mt-tf-equilibrium-proof-doublectrans}) implies
$\tilde{\varphi}_N^{(l)}(z)={-\sum_{i=1}^{N-1}}\tilde{\varphi}_i^{(l)}(z)\le {-\sum_{i=1}^{N-1}}\varphi_i^{(l)}(z)=\varphi_N^{(l)}(z)$
$\forall z\in\CX_0$, $\forall l\in\N$.
We hence get
\begin{align}
    \tilde{\varphi}{_N^{(l)}}^{c_N}(x)\ge {\varphi_N^{(l)}}^{c_N}(x) \qquad \forall x\in\CX_N,\; \forall l\in\N.\label{eqn:mt-tf-equilibrium-proof-triplectrans-N}
\end{align}

To prove statement~\ref{thms:mt-tf-equilibrium-invariance},
let us define 
$\kappa^{(l)}_i:=-\min_{z\in\CX_0}\big\{\tilde{\varphi}_i^{(l)}(z)\big\}$
$\forall 1\le i\le N-1$, $\forall l\in\N$
and define
$\kappa^{(l)}_{N}:=-\sum_{i=1}^{N-1}\kappa^{(l)}_i$ 
$\forall l\in\N$.
It holds that $\sum_{i=1}^{N}(\tilde{\varphi}_i^{(l)} + \kappa_i^{(l)})=\nobreak0$ $\forall l\in\N$,
and it follows from (\ref{eqn:mt-tf-equilibrium-proof-triplectrans}) and (\ref{eqn:mt-tf-equilibrium-proof-triplectrans-N}) that 
\begin{align*}
    \Bigg(\sum_{i=1}^{N}\int_{\CX_i\times\CX_0}c_i\DIFFX{\gamma_i^{(l)}}\Bigg)
    - \Bigg(\sum_{i=1}^{N}\int_{\CX_i}(\tilde{\varphi}_i^{(l)} + \kappa_i^{(l)})\big.^{c_i}\DIFFX{\mu_i}\Bigg)
    &=\Bigg(\sum_{i=1}^{N}\int_{\CX_i\times\CX_0}c_i\DIFFX{\gamma_i^{(l)}}\Bigg)
    - \Bigg(\sum_{i=1}^{N}\int_{\CX_i}\tilde{\varphi}{_i^{(l)}}^{c_i}\DIFFX{\mu_i}\Bigg)\\
    &\le \Bigg(\sum_{i=1}^{N}\int_{\CX_i\times\CX_0}c_i\DIFFX{\gamma_i^{(l)}}\Bigg)
    - \Bigg(\sum_{i=1}^{N}\int_{\CX_i}\varphi{_i^{(l)}}^{c_i}\DIFFX{\mu_i}\Bigg)\\
    &\le \epsilon^{(l)} \hspace{163pt} \forall l\in\N,
\end{align*}
which proves that $\big(\tilde{\varphi}^{(l)}_i+\kappa^{(l)}_i,{\gamma}_i^{(l)})_{i=1:N},\allowbreak{\nu}^{(l)}$ is an $\epsilon^{(l)}$-matching equilibrium for each $l\in\N$.
Moreover, for $i=1,\ldots,N-1$, it holds that 
$\min_{z\in\CX_0}\big\{\tilde{\varphi}_i^{(l)}(z)+\kappa_i^{(l)}\big\}=\nobreak 0$,
and hence the functions 
$(\tilde{\varphi}_i^{(l)}+\kappa_i^{(l)})_{l\in\N}$ are uniformly bounded by
$L_c\max_{z,z'\in\CX_0}\big\{d_{\CX_0}(z,z')\big\}<\nobreak\infty$.
Consequently, 
the functions $(\tilde{\varphi}_i^{(l)}+\kappa_i^{(l)})_{i=1:N,\,l\in\N}$
are uniformly bounded by 
$(N-1)L_c\max_{z,z'\in\CX_0}\big\{d_{\CX_0}(z,z')\big\}<\nobreak\infty$.
The proof of statement~\ref{thms:mt-tf-equilibrium-invariance} is now complete.

Now, for each $i\in\{1,\ldots,N\}$,
since the set of functions 
$(\tilde{\varphi}_i^{(l)}+\kappa_i^{(l)})_{l\in\N}$
on the compact metric space $\CX_0$ 
are uniformly bounded and Lipschitz continuous with respect to the same Lipschitz constant,
it follows from the Arzel{\`a}--Ascoli theorem that 
$(\tilde{\varphi}_i^{(l)}+\kappa_i^{(l)})_{l\in\N}$ is precompact in the topology induced by the uniform norm.
Moreover, since $(\CP(\CX_0),W_1)$ and 
$(\CP(\CX_1\times\nobreak\CX_0),W_1),\ldots,(\CP(\CX_N\times\nobreak\CX_0),W_1)$ are compact metric spaces 
(see, e.g., \citep[Remark~6.19]{ECvillani2008optimal}),
it holds that 
$(\nu^{(l)})_{l\in\N}$,
$(\gamma_1^{(l)})_{l\in\N},\ldots,(\gamma_N^{(l)})_{l\in\N}$
are precompact in their respective topologies.
Therefore, it holds that every subsequence of 
$\Big(\big(\tilde{\varphi}^{(l)}_i+\kappa^{(l)}_i,{\gamma}_i^{(l)})_{i=1:N},\allowbreak{\nu}^{(l)}\Big)_{l\in\N}$
admits a further subsequence
$\Big(\big(\tilde{\varphi}^{(l_t)}_i+\kappa^{(l_t)}_i,\allowbreak{\gamma}_i^{(l_t)})_{i=1:N},\allowbreak{\nu}^{(l_t)}\Big)_{t\in\N}$
where
$({\nu}^{(l_t)})_{t\in\N}$ converges in $(\CP(\CX_0),W_1)$ to some~${\nu}^{(\infty)}$, 
and for $i=1,\ldots,N$, 
$\big(\tilde{\varphi}_i^{(l_t)}+\kappa_i^{(l_t)}\big)_{t\in\N}$ converges uniformly to some $\tilde{\varphi}_i^{(\infty)}\in\CC(\CX_0)$, 
$({\gamma}_i^{(l_t)})_{t\in\N}$ converges in $(\CP({\CX_i\times\CX_0}),W_1)$ to some ${\gamma}^{(\infty)}_i$.
One checks that 
$\sum_{i=1}^{N}\tilde{\varphi}_i^{(\infty)}=\nobreak0$,
and that 
$\gamma_i^{(\infty)}\in\Gamma(\mu_i,\nu^{(\infty)})$
for $i=1,\ldots,N$.
Furthermore, 
since \citep[Corollary~1]{ECcarlier2010matching} implies that the mapping
$\CC(\CX_0)\ni\varphi\mapsto\int_{\CX_i}\varphi^{c_i}\DIFFX{\mu_i}\in\R$
is continuous for $i=1,\ldots,N$,
it holds that 
\begin{align*}
    \hspace{20pt}&\hspace{-20pt}
    \Bigg(\sum_{i=1}^{N}\int_{\CX_i\times\CX_0}c_i\DIFFX{\gamma_i^{(\infty)}}\Bigg)
    - \Bigg(\sum_{i=1}^{N}\int_{\CX_i}\tilde{\varphi}{_i^{(\infty)}}^{c_i}\DIFFX{\mu_i}\Bigg)\\
    &=\Bigg(\sum_{i=1}^{N}\lim_{t\to\infty}\int_{\CX_i\times\CX_0}c_i\DIFFX{\gamma_i^{(l_t)}}\Bigg)
    - \Bigg(\sum_{i=1}^{N}\lim_{t\to\infty}\int_{\CX_i}(\tilde{\varphi}_i^{(l_t)}+\kappa_i^{(l_t)})\big.^{c_i}\DIFFX{\mu_i}\Bigg)\\
    &\le \limsup_{t\to\infty} \Bigg(\sum_{i=1}^{N}\int_{\CX_i\times\CX_0}c_i\DIFFX{\gamma_i^{(l_t)}}\Bigg)
    - \Bigg(\sum_{i=1}^{N}\int_{\CX_i}(\tilde{\varphi}_i^{(l_t)}+\kappa_i^{(l_t)})\big.^{c_i}\DIFFX{\mu_i}\Bigg)\le \limsup_{t\to\infty} \epsilon^{(l)}=0.
\end{align*}
It thus follows that 
$(\tilde{\varphi}^{(\infty)}_i,{\gamma}_i^{(\infty)})_{i=1:N},{\nu}^{(\infty)}$ is a matching equilibrium.
The proof is now complete.
\endproof

\proof{Proof of Proposition~\ref{prop:perturbation}.}
Since 
$\alpha_{\MT}^\star$ is the optimal value of \eqref{eqn:mt-primalopt} with respect to 
$c_i\leftarrow c_i$ $\forall 1\le\nobreak i\le\nobreak N$,
and 
$\tilde{\alpha}_{\MT}^\star$ is the optimal value of \eqref{eqn:mt-primalopt} with respect to 
$c_i\leftarrow \tilde{c}_i$ $\forall 1\le i\le N$, we have 
$\big(\sum_{i=1}^{N}\int_{\CX_i\times\CX_0}c_i\DIFFX{\gamma_i}\big)-\epsilon 
\le \alpha_{\MT}^\star
\le \big(\sum_{i=1}^{N}\int_{\CX_i}\varphi_i^{c_i}\DIFFX{\mu_i}\big)+\epsilon$
and 
$\sum_{i=1}^{N}\int_{\CX_i}\varphi_i^{\tilde{c}_i}\DIFFX{\mu_i}
\le\tilde{\alpha}_{\MT}^\star
\le\sum_{i=1}^{N}\int_{\CX_i\times\CX_0}\tilde{c}_i\DIFFX{\gamma_i}$.
Combining these two inequalities completes the proof.
\endproof

\subsection{Proofs of results in Section~\ref{sec:parametric}}\label{sapx:proof-parametric}

\proof{Proof of Theorem~\ref{thm:mt-tf-duality}.}
Let us first introduce some notations in order to express the LSIP problem 
\eqref{eqn:mt-tf-lsip} into a concise vectorized form.
Let us denote $q:=\sum_{i=1}^N(1+m_i+m_0)$,
arrange the $q$ decision variables of \eqref{eqn:mt-tf-lsip} into a vector $\Bchi:=\big(y_{1,0},\BIy_1^\TRANSP,\BIw_1^\TRANSP,\ldots,y_{N,0},\BIy_N^\TRANSP,\BIw_N^\TRANSP\big)^\TRANSP\in\R^{q}$, 
and denote the objective vector of \eqref{eqn:mt-tf-lsip} by 
$\bar{\BIg}:=\big(1,\bar{\BIg}_1^\TRANSP,\veczero_{m_0}^\TRANSP,\ldots,1,\bar{\BIg}_N^\TRANSP,\veczero_{m_0}^\TRANSP\big)^\TRANSP\in\R^{q}$.
For $i=1,\ldots,N$, let us define $\veczero^{(i)}:=\veczero_{1+m_i+m_0}$ for notational simplicity, and let $\BIa_i:\CX_i\times\CX_0\to\R^{q}$ be defined as follows:
\begin{align*}
\BIa_1(x,z)&:=\big(1, \BIg_1(x)^\TRANSP, \BIg_0(z)^\TRANSP, \veczero^{(2)\TRANSP},\ldots, \veczero^{(N)\TRANSP}\big)^\TRANSP & \forall (x,z)\in\CX_1\times\CX_0,\\
\BIa_2(x,z)&:=\big(\veczero^{(1)\TRANSP},1, \BIg_2(x)^\TRANSP, \BIg_0(z)^\TRANSP, \veczero^{(3)\TRANSP},\ldots, \veczero^{(N)\TRANSP}\big)^\TRANSP &  \forall (x,z)\in\CX_2\times\CX_0,\\
&\;\;\;\vdots\\
\BIa_N(x,z)&:=\big(\veczero^{(1)\TRANSP},\ldots,\veczero^{(N-1)\TRANSP},1, \BIg_N(x)^\TRANSP, \BIg_0(z)^\TRANSP\big)^\TRANSP & \forall (x,z)\in\CX_N\times\CX_0.
\end{align*}
Let $\BIe_l$ denote the $l$-th standard basis vector of $\R^{m_0}$
for $l=1,\ldots,m_0$, 
and let us define $\BIu^{(l)}:=\big(0,\veczero_{m_1}^\TRANSP,\BIe_l^\TRANSP,0,\veczero_{m_2}^\TRANSP,\BIe_l^\TRANSP,\ldots,0,\veczero_{m_N}^\TRANSP,\BIe_l^\TRANSP\big)^\TRANSP\in\R^{q}$ $\forall 1\le l\le m_0$.
With these notations, we can now equivalently express \eqref{eqn:mt-tf-lsip} concisely as follows: 
\begin{align}
\begin{split}
\maximize_{\Bchi}\quad & \langle\bar{\BIg},\Bchi\rangle\\
\mathrm{subject~to}\quad & \langle\BIa_i(x,z),\Bchi\rangle  \le c_i(x,z) \quad \forall (x,z)\in\CX_i\times\CX_0,\; \forall 1\le i\le N,\\
& \langle\BIu^{(l)},\Bchi\rangle=0 \quad\forall 1\le l\le m_0, \qquad \Bchi\in\R^{q}.
\end{split}
\label{eqn:mt-tf-lsip-concise}
\end{align}
Hence, it holds that $\eqref{eqn:mt-tf-lsip-concise}=\eqref{eqn:mt-tf-lsip}$.
The so-called \textit{Haar's dual optimization problem} of (\ref{eqn:mt-tf-lsip-concise}) (see, e.g., \citep[p.49]{ECgoberna1998linear}) is given by:
\begin{align}
\begin{split}
\minimize_{\scalebox{0.5}{$(\theta_{i,j},x_{i,j},z_{i,j}),(\xi_l)$}}\quad & \sum_{i=1}^N\sum_{j=1}^{n_i}\theta_{i,j}c_i(x_{i,j},z_{i,j})\\
\mathrm{subject~to}\quad & \left(\sum_{i=1}^N\sum_{j=1}^{n_i} \theta_{i,j}\BIa_i(x_{i,j},z_{i,j})\right)-\left(\sum_{l=1}^{m_0}\xi_l\BIu^{(l)}\right)=\bar{\BIg},\\
& n_i\in\N,\; 
(\theta_{i,j})_{j=1:n_i}\subset \R_+,\;
(x_{i,j})_{j=1:n_i}\subseteq\CX_i, \;
(z_{i,j})_{j=1:n_i}\subseteq\CX_0,
\quad \forall 1\le i\le N, \\
& (\xi_l)_{l=1:m_0}\subset\R.
\end{split}
\label{eqn:mt-tf-lsip-concise-dual}
\end{align}
In order to prove the strong duality between (\ref{eqn:mt-tf-lsip-concise}) and (\ref{eqn:mt-tf-lsip-concise-dual}), 
\citet[Theorem~4.5 \& Theorem~8.2]{ECgoberna1998linear} showed that it is sufficient to prove the closedness of 
the so-called \textit{second-moment cone} $C$ defined as follows (see \citep[p.81]{ECgoberna1998linear}):
\begin{align*}
{K}_{1,i}&:=\conv\Big(\Big\{\big(\BIa_i(x,z)^\TRANSP, c_i(x,z)\big)^\TRANSP:(x,z)\in\CX_i\times\CX_0\Big\}\Big)\subset\R^{q+1} 
&& \forall 1\le i\le N,\\
{C}_{1,i}&:=\cone({K}_{1,i})\subseteq\R^{q+1} 
&& \forall 1\le i\le N,\\
{C}_2&:=\cone\Big(\Big\{\big(\iota\BIu^{(l)\TRANSP},0\big)^\TRANSP:\iota\in\{-1,1\},\;1\le l\le m_0\Big\}\Big)\subset\R^{q+1},\\
{C}&:={C}_{1,1}+\cdots+{C}_{1,N}+{C}_{2}\subseteq\R^{q+1}.
\end{align*}

In order to get the strong duality $\eqref{eqn:mt-tf-dual}= \eqref{eqn:mt-tf-lsip}$,
we will divide the rest of the proof into the following 5 steps.
\begin{itemize}
    \item Step~1: establishing the weak duality $\eqref{eqn:mt-tf-dual}\ge \eqref{eqn:mt-tf-lsip}$.
    
    \item Step~2: showing that $C_{1,i}$ is closed for $i=1,\ldots,N$.
    
    \item Step~3: showing that $C$ is closed.
    
    \item Step~4: establishing the strong duality $\eqref{eqn:mt-tf-lsip-concise}=\eqref{eqn:mt-tf-lsip-concise-dual}$.
    
    \item Step~5: showing that $\eqref{eqn:mt-tf-lsip-concise-dual}\ge\eqref{eqn:mt-tf-dual}$.
\end{itemize}

\underline{Step~1.}
It follows from the compactness of $(\CX_i)_{i=0:N}$ and the continuity of $(c_i)_{i=1:N}$ that \eqref{eqn:mt-tf-lsip} is feasible. 
Observe that \eqref{eqn:mt-tf-dual} is also feasible. 
Let us fix an arbitrary feasible solution $(y_{i,0},\BIy_i,\BIw_i)_{i=1:N}$ of \eqref{eqn:mt-tf-lsip} as well as an arbitrary feasible solution $(\theta_i)_{i=1:N}$ of \eqref{eqn:mt-tf-dual},
and let us denote 
$\BIy_i=(y_{i,1},\ldots,y_{i,m_i})^\TRANSP$, $\BIw_i=(w_{i,1},\ldots,w_{i,m_0})^\TRANSP$ for $i=1,\ldots,N$.
By the constraints of \eqref{eqn:mt-tf-lsip}, it holds that 
$\sum_{i=1}^N\BIw_i=\veczero_{m_0}$, and 
$y_{i,0}+\langle\BIg_i(x),\BIy_i\rangle+\langle\BIg_0(z),\BIw_i\rangle  \le c_i(x,z)$ $\forall (x,z)\in\CX_i\times\CX_0$, $\forall 1\le i\le N$. 
Moreover, by the constraints of \eqref{eqn:mt-tf-dual}, it holds that $\theta_i\in\Gamma(\bar{\mu}_i,\bar{\nu}_i)$ for some $(\bar{\mu}_i,\bar{\nu}_i)_{i=1:N}$ satisfying $\bar{\mu}_i\overset{\CG_i}{\sim}\mu_i$, $\bar{\nu}_i\overset{\CG_0}{\sim}\bar{\nu}_1$ $\forall 1\le i\le N$. 
In particular, this implies that
\begin{align*}
    \int_{\CX_i\times\CX_0}\langle\BIg_i(x),\BIy_i\rangle\DIFFM{\theta_i}{\DIFF x,\DIFF z}
    &=\int_{\CX_i}\langle\BIg_i(x),\BIy_i\rangle\DIFFM{\bar{\mu}_i}{\DIFF x}
    =\sum_{j=1}^{m_i}y_{i,j}\int_{\CX_i}g_{i,j}\DIFFX{\mu_i}=\langle\bar{\BIg}_i,\BIy_i\rangle \qquad \forall 1\le i\le N,
\end{align*}
and that
\begin{align*}
    \int_{\CX_i\times\CX_0}\langle\BIg_0(z),\BIw_i\rangle\DIFFM{\theta_i}{\DIFF x,\DIFF z}
    &=\int_{\CX_0}\langle\BIg_0(z),\BIw_i\rangle\DIFFM{\bar{\nu}_i}{\DIFF z}
    =\sum_{l=1}^{m_0}w_{i,l}\int_{\CX_0}g_{0,l}\DIFFX{\bar{\nu}_1}\\
    &=\int_{\CX_0}\langle\BIg_0(z),\BIw_i\rangle\DIFFM{\bar{\nu}_1}{\DIFF z} \qquad\qquad\quad \forall 1\le i\le N.
\end{align*}
Consequently, we obtain
\begin{align*}
\sum_{i=1}^N\int_{\CX_i\times\CX_0}c_i(x,z)\DIFFM{\theta_i}{\DIFF x,\DIFF z}&
\ge \sum_{i=1}^N\int_{\CX_i\times\CX_0}y_{i,0}+\langle\BIg_i(x),\BIy_i\rangle+\langle\BIg_0(z),\BIw_i\rangle\DIFFM{\theta_i}{\DIFF x,\DIFF z}\\
&=\Bigg(\sum_{i=1}^N y_{i,0}+\langle\BIy_i,\bar{\BIg}_i\rangle\Bigg)+\int_{\CX_0}\big\langle\BIg_0(z),\textstyle\sum_{i=1}^N\BIw_i\big\rangle\DIFFM{\bar{\nu}_1}{\DIFF z}\\
&=\sum_{i=1}^N y_{i,0}+\langle\BIy_i,\bar{\BIg}_i\rangle.
\end{align*}
Taking the supremum over all $(y_{i,0},\BIy_i,\BIw_i)_{i=1:N}$ feasible for \eqref{eqn:mt-tf-lsip} and taking the infimum over all $(\theta_i)_{i=1:N}$ feasible for \eqref{eqn:mt-tf-dual} in the inequality above proves  $\eqref{eqn:mt-tf-dual}\ge \eqref{eqn:mt-tf-lsip}$.

\underline{Step~2.}
Let us fix an arbitrary $i\in\{1,\ldots,N\}$.
Notice that the continuity of 
$\BIg_i$, $\BIg_0$, and $c_i$, 
the compactness of $\CX_i\times\CX_0$, 
and \citep[Theorem~17.2]{ECrockafellar1970convex} imply that ${K}_{1,i}$ is compact. 
In addition, observe that ${K}_{1,i}$ does not contain the origin. 
It thus follows from \citep[Corollary~9.6.1]{ECrockafellar1970convex} that ${C}_{1,i}$ is closed. 

\underline{Step~3.}
${C}_2$ is closed since it is a subspace of $\R^{q+1}$.
To prove the closedness of ${C}$, we need to verify an additional condition in \citep[Corollary~9.1.3]{ECrockafellar1970convex} that
whenever ${\BIc}_{1,1}\in\nobreak {C}_{1,1},\ldots,{\BIc}_{1,N}\in\nobreak {C}_{1,N},\allowbreak{\BIc}_2\in {C}_2$
satisfy
${\BIc}_{1,1}+\cdots+{\BIc}_{1,N}+{\BIc}_2=\veczero$, 
it holds that ${\BIc}_{1,1}=\cdots={\BIc}_{1,N}={\BIc}_2=\nobreak\veczero$.
To that end, 
let us fix arbitrary 
${\BIc}_{1,1}\in {C}_{1,1},\ldots,{\BIc}_{1,N}\in {C}_{1,N},{\BIc}_2\in {C}_2$
that satisfy
${\BIc}_{1,1}+\cdots+{\BIc}_{1,N}+{\BIc}_2=\veczero$,
and let $s_1:=1$, $s_i:=s_{i-1}+1+m_{i-1}+m_0$ $\forall 2\le i\le N$. 
Then, for $i=1,\ldots,N$, 
there exist 
${\BIk}_{1,i}\in{K}_{1,i}$ and 
$\beta_{1,i}\ge 0$ such that 
${\BIc}_{1,i}=\beta_{1,i}{\BIk}_{1,i}$.
It follows from the definition of $\BIa_i(\,\cdot\,,\cdot\,)$ that the \mbox{$s_i$-th} component of ${\BIk}_{1,i}$ is equal to~$1$ and that the \mbox{$s_{i'}$-th} component of ${\BIk}_{1,i}$ is equal to~$0$ for $i'\ne i$. 
Moreover, it follows from the definition of $\big(\BIu^{(l)}\big)_{l=1:m_0}$ that 
the \mbox{$s_i$-th} component of ${\BIc}_2$ is equal to~$0$ for $i=1,\ldots,N$. 
Consequently, the $s_i$-th component of 
${\BIc}_{1,1}+\cdots+{\BIc}_{1,N}+{\BIc}_2$ is equal to $\beta_{1,i}$ for $i=1,\ldots,N$, 
implying that 
$\beta_{1,1}=\cdots=\beta_{1,N}=0$, 
${\BIc}_{1,1}=\cdots={\BIc}_{1,N}=\veczero$, 
and hence ${\BIc}_2=\veczero$ as well.
Now, 
since $C_{1,1},\ldots,C_{1,N},C_{2}$ are all closed,
we can apply \citep[Corollary~9.1.3]{ECrockafellar1970convex} to prove that ${C}$ is closed. 

\underline{Step~4.}
Since $C$ is closed, and since (\ref{eqn:mt-tf-lsip-concise}) is feasible,
it follows from \citep[Theorem~4.5]{ECgoberna1998linear} (with respect to $N\leftarrow {C}$, $K\leftarrow \cone\big({C}\cup (\veczero_{q}^\TRANSP,1)^\TRANSP\big)$ in the notation of \citep{ECgoberna1998linear}) that $\cone\big({C}\cup (\veczero_{q}^\TRANSP,1)^\TRANSP\big)$ is closed,
and it then follows from \citep[Theorem~8.2]{ECgoberna1998linear} that the optimal values of (\ref{eqn:mt-tf-lsip-concise}) and (\ref{eqn:mt-tf-lsip-concise-dual}) are identical (see the last three cases in \citep[Table~8.1]{ECgoberna1998linear}).

\underline{Step~5.}
Since we have shown that 
$\eqref{eqn:mt-tf-lsip}=\eqref{eqn:mt-tf-lsip-concise-dual}$, 
it holds that (\ref{eqn:mt-tf-lsip-concise-dual}) is feasible.
Thus, let us fix an arbitrary feasible solution $(\theta_{i,j},x_{i,j},z_{i,j})_{j=1:n_i,\,i=1:N}$, $(\xi_l)_{l=1:m_0}$ of (\ref{eqn:mt-tf-lsip-concise-dual}) and characterize its properties. 
The constraints of (\ref{eqn:mt-tf-lsip-concise-dual}) implies that
\begin{align}
\Bigg(\sum_{i=1}^N\sum_{j=1}^{n_i} \theta_{i,j}\BIa_i(x_{i,j},z_{i,j})\Bigg)
-\Bigg(\sum_{l=1}^{m_0}\xi_l\BIu^{(l)}\Bigg)=\bar{\BIg}.
\label{eqn:mt-tf-duality-proof-dualconstraint}
\end{align}
Consequently, it follows from the definitions of $\bar{\BIg}$, $(\BIa(\,\cdot\,,\cdot\,))_{i=1:N}$, $(\BIu^{(l)})_{l=1:m_0}$ and a component-wise expansion of (\ref{eqn:mt-tf-duality-proof-dualconstraint}) that
\begin{align}
\sum_{j=1}^{n_i}\theta_{i,j}&=1 
&\forall 1\le i\le N,
\label{eqn:mt-tf-duality-proof-dualconstraint1}\\
\sum_{j=1}^{n_i}\theta_{i,j}\BIg_i(x_{i,j})&=\bar{\BIg}_i 
&\forall 1\le i\le N,
\label{eqn:mt-tf-duality-proof-dualconstraint2}\\
\sum_{j=1}^{n_i}\theta_{i,j}\BIg_0(z_{i,j})&=\sum_{l=1}^{m_0} \xi_l\BIe_l &\forall 1\le i\le N.
\label{eqn:mt-tf-duality-proof-dualconstraint3}
\end{align}
Subsequently, let us define $\theta_{i}:=\sum_{j=1}^{n_i}\theta_{i,j}\delta_{(x_{i,j},z_{i,j})}$ $\forall 1\le i\le N$. 
By (\ref{eqn:mt-tf-duality-proof-dualconstraint1}) and $(\theta_{i,j})_{j=1:n_i}\subset\nobreak\R_+$, it holds that $\theta_i\in\CP(\CX_i\times\CX_0)$. 
Let $\bar{\mu}_i$ and $\bar{\nu}_i$ denote the marginals of $\theta_i$ on $\CX_i$ and $\CX_0$, respectively. 
Then, 
(\ref{eqn:mt-tf-duality-proof-dualconstraint2}) implies that 
$\int_{\CX_i}g_{i,j}\DIFFX{\bar{\mu}_i}=\int_{\CX_i\times\CX_0}g_{i,j}(x)\DIFFM{\theta_i}{\DIFF x,\DIFF z}=\sum_{k=1}^{n_i}\theta_{i,k}g_{i,j}(x_{i,k})=\int_{\CX_i}g_{i,j}\DIFFX{\mu_i}$
$\forall 1\le j\le m_i$, $\forall 1\le i\le N$,
and it hence holds that 
$\bar{\mu}_i\overset{\CG_i}{\sim}\mu_i$ $\forall 1\le i\le N$. 
Moreover, (\ref{eqn:mt-tf-duality-proof-dualconstraint3}) implies that 
$\int_{\CX_0}g_{0,l}\DIFFX{\bar{\nu}_i}
=\int_{\CX_i\times\CX_i}g_{0,l}(z)\DIFFM{\theta_i}{\DIFF x,\DIFF z}
=\sum_{k=1}^{n_i}\theta_{i,k}g_{0,l}(z_{i,k})=\xi_l$
$\forall 1\le l\le m_0$, $\forall 1\le\nobreak i\le\nobreak N$.
This shows that $\int_{\CX_0}g_{0,l}\DIFFX{\bar{\nu}_1}=\cdots=\int_{\CX_0}g_{0,l}\DIFFX{\bar{\nu}_N}=\xi_l$ 
$\forall 1\le l\le m_0$
and hence $\bar{\nu}_1\overset{\CG_0}{\sim}\cdots\overset{\CG_0}{\sim}\bar{\nu}_N$. 
The above analysis shows that $(\theta_i)_{i=1:N}$ is a feasible solution of \eqref{eqn:mt-tf-dual}. 
Furthermore, it holds that
$\sum_{i=1}^N\int_{\CX_i\times\CX_0}c_i(x,z)\DIFFM{\theta_i}{\DIFF x,\DIFF z}= \sum_{i=1}^N\sum_{k=1}^{n_i}\theta_{i,k}c_i(x_{i,k},z_{i,k})$.
Therefore, 
taking the infimum over all $(\theta_{i,j},x_{i,j},z_{i,j})_{j=1:n_i,\,i=1:N}$, $(\xi_l)_{l=1:m_0}$ feasible for (\ref{eqn:mt-tf-lsip-concise-dual}) shows that
$\eqref{eqn:mt-tf-lsip-concise-dual}\ge\eqref{eqn:mt-tf-dual}$. 
The proof is now complete. 
\endproof

\proof{Proof of Theorem~\ref{thm:mt-tf-approx}.}
Observe that
the continuity of $\BIg_0:\CX_0\to\R$ and the constraints of \eqref{eqn:mt-tf-lsip}
guarantee that 
$(\hat{\varphi}_i)_{i=1:N}$ are continuous
and $\sum_{i=1}^{N}\hat{\varphi}_i=\nobreak0$.
Moreover, the construction of $\hat{\nu}$, $\tilde{\nu}$, $(\hat{\gamma}_i,\tilde{\gamma}_i)_{i=1:N}$
guarantee that 
$\hat{\gamma}_i\in\Gamma(\mu_i,\hat{\nu})$
and 
$\tilde{\gamma}_i\in\Gamma(\mu_i,\tilde{\nu})$
for $i=1,\ldots, N$.
Therefore, 
in view of the assumption
$\big(\sum_{i=1}^{N}\int_{\CX_i\times\CX_0}c_i\DIFFX{\hat{\theta}_i}\big)-\big(\sum_{i=1}^{N}\hat{y}_{i,0}+\langle\bar{\BIg}_i,\hat{\BIy}_i\rangle\big)\le \epsilon_{\sspar}$,
we will prove that 
$(\hat{\varphi}_i,\hat{\gamma}_i)_{i=1:N},\hat{\nu}$ and $(\hat{\varphi}_i,\tilde{\gamma}_i)_{i=1:N},\tilde{\nu}$ are $\epsilon_{\ssapx}$-matching equilibria
by establishing the following inequalities:
\begin{align}
    \sum_{i=1}^{N}\int_{\CX_i\times\CX_0}c_i\DIFFX{\tilde{\gamma}_i}
    \le \sum_{i=1}^{N}\int_{\CX_i\times\CX_0}c_i\DIFFX{\hat{\gamma}_i}
    &\le \sum_{i=1}^{N}\bigg(\int_{\CX_i\times\CX_0}c_i\DIFFX{\hat{\theta}_i}+ L_c\big(\overline{\epsilon}_i(\varsigma)+\overline{\epsilon}_0(\varsigma)\big)\bigg),\label{eqn:mt-tf-approx-proof-ineq1}\\
    \sum_{i=1}^{N}\int_{\CX_i}\hat{\varphi}_i^{c_i}\DIFFX{\mu_i}
    &\ge \sum_{i=1}^{N}\hat{y}_{i,0}+\langle\bar{\BIg}_i,\hat{\BIy}_i\rangle.\label{eqn:mt-tf-approx-proof-ineq2}
\end{align}

First, the definitions of $(\hat{\gamma}_i)_{i=1:N}$, 
the properties of $Z$, $(Z_i,X_i,\bar{X}_i)_{i=1:N}$,
and the $L_c$-Lipschitz continuity of $(c_i)_{i=1:N}$ imply that
\begin{align*}
    \int_{\CX_i\times\CX_0}c_i\DIFFX{\hat{\gamma}_i}
    = \EXP\big[c_i(\bar{X}_i,Z)\big]
    &\le \EXP\big[c_i(X_i,Z_i)\big] 
    + L_c\EXP\big[d_{\CX_i}(X_i,\bar{X}_i)\big] 
    + L_c\EXP\big[d_{\CX_0}(Z,Z_i)\big]\\
    &\le \int_{\CX_i\times\CX_0}c_i\DIFFX{\hat{\theta}_i}
    + L_c\big(\overline{\epsilon}_i(\varsigma)+\overline{\epsilon}_0(\varsigma)\big)
    \qquad\qquad\qquad \forall 1\le i\le N.
\end{align*}
Moreover, it follows from the property of the random variable $\bar{Z}$ that 
\begin{align*}
    \sum_{i=1}^{N}\int_{\CX_i\times\CX_0}c_i\DIFFX{\tilde{\gamma}_i}
    &=\EXP\Bigg[\sum_{i=1}^{N}c_i(\bar{X}_i,\bar{Z})\Bigg]
    \le \EXP\Bigg[\sum_{i=1}^{N}c_i(\bar{X}_i,Z)\Bigg]
    =\sum_{i=1}^{N}\int_{\CX_i\times\CX_0}c_i\DIFFX{\hat{\gamma}_i}.
\end{align*}
Thus, (\ref{eqn:mt-tf-approx-proof-ineq1}) holds.

Lastly, fixing an arbitrary $i\in\{1,\ldots,N\}$,
we obtain from the constraints of \eqref{eqn:mt-tf-lsip} that 
$\hat{y}_{i,0}+\langle\BIg_i(x),\hat{\BIy}_i\rangle\le \inf_{z\in\CX_0}\big\{c_i(x,z)-\hat{\varphi}_i(z)\big\}=\hat{\varphi}_i^{c_i}(x)$ 
$\forall x\in\CX_i$,
and integrating both sides with respect to~$\mu_i$
then yields 
$\hat{y}_{i,0}+\langle\bar{\BIg}_i,\hat{\BIy}_i\rangle\le \int_{\CX_i}\hat{\varphi}_i^{c_i}\DIFFX{\mu_i}$.
This establishes (\ref{eqn:mt-tf-approx-proof-ineq2}) and completes the proof.
\endproof

\subsection{Proofs of results in Section~\ref{sec:numerics}}\label{sapx:proof-numerics}

\proof{Proof of Theorem~\ref{thm:error-control}.}
    Let us first prove statement~\ref{thms:error-control-W1}.
    To that end, let us fix an arbitrary $i\in\nobreak\{0,1,\ldots,N\}$ and an arbitrary $\varsigma\ge 0$,
    let $\FS_i$ be an arbitrary subdivision of $\FK_i$, and let 
    $\CG_i$, $\eta(\FS_i)$, $\overline{\epsilon}_i(\varsigma)$ be defined by Setting~\ref{sett:Euclidean}.
    Subsequently,
    it follows from \citep[Proposition~6.3]{ECneufeld2022v7numerical} that $(g_{i,j})_{j=0:m_i}$ is a so-called \textit{vertex interpolation function set} with respect to~$\FS_i$. 
    Let us fix two arbitrary probability measures $\nu,\nu'\in\CP(\CX_i)$ such that $\nu\abovebelowset{\CG_i}{\varsigma}{\sim}\nu'$.
    Subsequently, it follows from the proof of \citep[Proposition~6.5]{ECneufeld2022v7numerical} that 
    there exist $\hat{\nu},\hat{\nu}'\in\CP(\CX_i)$ 
    with $\support(\hat{\nu})\subseteq V(\FS_i)$, $\support(\hat{\nu}')\subseteq V(\FS_i)$
    satisfying 
    \begin{align}
        W_1(\nu,\hat{\nu})&\le \eta(\FS_i),\label{eqn:error-control-proof-ineq1}\\
        W_1(\nu',\hat{\nu}')&\le \eta(\FS_i),\label{eqn:error-control-proof-ineq2}\\
        \hat{\nu}\big(\{\BIv_{i,j}\}\big)&=\int_{\CX_i}g_{i,j}\DIFFX{\nu} \hspace{2.7pt}\qquad \forall 0\le j\le m_i,\label{eqn:error-control-proof-moment1}\\
        \hat{\nu}'\big(\{\BIv_{i,j}\}\big)&=\int_{\CX_i}g_{i,j}\DIFFX{\nu'} \qquad \forall 0\le j\le m_i.\label{eqn:error-control-proof-moment2}
    \end{align}
    Recall that $\nu\abovebelowset{\CG_i}{\varsigma}{\sim}\nu'$ implies 
    \begin{align}
        \sum_{j=1}^{m_i}\bigg|\int_{\CX_i}g_{i,j}\DIFFX{\nu}-\int_{\CX_i}g_{i,j}\DIFFX{\nu'}\bigg|
        &=\sum_{g\in\CG_i}\bigg|\int_{\CX_i}g\DIFFX{\nu}-\int_{\CX_i}g\DIFFX{\nu'}\bigg|
        \le \varsigma.
        \label{eqn:error-control-proof-approx1}
    \end{align}
    Since $g_{i,0}=1-\big(\sum_{j=1}^{m_i}g_{i,j}\big)$ holds by definition, we get
    \begin{align}
        \begin{split}
            \bigg|\int_{\CX_i}g_{i,0}\DIFFX{\nu}-\int_{\CX_i}g_{i,0}\DIFFX{\nu'}\bigg|
            &\le \sum_{j=1}^{m_i}\bigg|\int_{\CX_i}g_{i,j}\DIFFX{\nu}-\int_{\CX_i}g_{i,j}\DIFFX{\nu'}\bigg|
            \le \varsigma.
        \end{split}
        \label{eqn:error-control-proof-approx2}
    \end{align}
    Consequently, combining (\ref{eqn:error-control-proof-moment1})--(\ref{eqn:error-control-proof-approx2})
    yields
    \begin{align*}
        \sum_{\BIv\in V(\FS_i)}\big|\hat{\nu}\big(\{\BIv\}\big)-\hat{\nu}'\big(\{\BIv\}\big)\big|=\sum_{j=0}^{m_i}\bigg|\int_{\CX_i}g_{i,j}\DIFFX{\nu}-\int_{\CX_i}g_{i,j}\DIFFX{\nu'}\bigg|\le 2\varsigma,
    \end{align*}
    which implies that $\|\hat{\nu}-\hat{\nu}'\|_{\mathrm{TV}}\le \varsigma$, 
    where $\|\hat{\nu}-\nobreak\hat{\nu}'\|_{\mathrm{TV}}$ denotes the total variation distance between $\hat{\nu}$ and~$\hat{\nu}'$.
    Therefore,
    it follows from 
    (\ref{eqn:error-control-proof-ineq1}), (\ref{eqn:error-control-proof-ineq2}), and
    the upper bound for the $W_1$-distance through the total variation distance (see, e.g., \citep[Particular Case~6.16]{ECvillani2008optimal}) that
    $W_1(\nu,\nu')\le W_1(\nu,\hat{\nu})+W_1(\nu',\hat{\nu}')+W_1(\hat{\nu},\hat{\nu}')\le \overline{\epsilon}_{i}(\varsigma)$. 
    The proof of statement~\ref{thms:error-control-W1} is now complete.

    In statement~\ref{thms:error-control-arbitrary},
    for $i=0,1,\ldots,N$,
    one uses the procedure 
    $\mathtt{InitTestFuncs}\big(\CX_i,\frac{\epsilon}{16N(L_c\vee 1)}\big)$ and 
    subsequently construct $\FS_i$ by Lemma~\ref{lem:proc-init-test-funcs}.
    It thus holds by Lemma~\ref{lem:proc-init-test-funcs}\ref{lems:proc-init-test-funcs-subdivision}
    that $\FS_i$ is a simplicial complex that is a subdivision of $\FK_i$ 
    and that $\eta(\FS_i)\le \frac{\epsilon}{16NL_c}$,
    $|V(\FS_i)|-1\le \big(\frac{64ND_id_i(L_c\vee 1)}{\epsilon}\big)^{d_i}|\FM(\FK_i)|(d_i+1)$
    for $i=0,1,\ldots,N$.
    Then, the definition of $(\overline{\epsilon}_i(\,\cdot\,))_{i=0:N}$ in Setting~\ref{sett:Euclidean}
    along with
    $\varsigma\leftarrow\nobreak \frac{\epsilon}{8NL_c}\Big(\max_{0\le i\le N}\big\{\max_{\BIx,\BIx'\in\CX_i}\big\{\|\BIx-\BIx'\|_2\big\}\big\}\Big)^{-1}$
    yields $\overline{\epsilon}_i(\varsigma)\le \frac{\epsilon}{4NL_c}$
    $\forall i\in\{0,1,\ldots,N\}$.
    Moreover, statement~\ref{thms:error-control-W1}
    shows that $(\overline{\epsilon}_i(\varsigma))_{i=0:N}$
    satisfy the conditions of Theorem~\ref{thm:mt-tf-approx}.
    Thus, with the test functions $(\CG_i)_{i=0:N}$ defined in Setting~\ref{sett:Euclidean} 
    with respect to $(\FS_i)_{i=0:N}$
    and with $\epsilon_{\sspar}\leftarrow \frac{\epsilon}{2}$,
    Theorem~\ref{thm:mt-tf-approx} yields 
    two $\epsilon_{\ssapx}$-matching equilibria 
    $(\hat{\varphi}_i,\hat{\gamma}_i)_{i=1:N},\hat{\nu}$ and $(\hat{\varphi}_i,\tilde{\gamma}_i)_{i=1:N},\tilde{\nu}$,
    where $\epsilon_{\ssapx}:=\epsilon_{\sspar}
    +L_c\big(N\overline{\epsilon}_0(\varsigma)+\sum_{i=1}^N\overline{\epsilon}_i(\varsigma)\big)\le \epsilon$.
    The proof is now complete.
\endproof

\proof{Proof of Theorem~\ref{thm:complexity-explicit}.}
    Throughout this proof, let us denote 
    $\Lambda:=\big(\frac{64NDd\overline{L}_{c}}{\epsilon}\big)^{d}J$ 
    to simplify the computational complexities.
    To begin, 
    substituting $\eta_{\ssUB}\leftarrow \frac{\epsilon}{16N(L_{c}\vee 1)}$
    into Lemma~\ref{lem:proc-init-test-funcs}\ref{lems:proc-init-test-funcs-complexity}
    shows that 
    Line~\ref{alglin:complexity-overall-init-test-funcs-forloop} incurs 
    $O\Big(N\big(\frac{64NDd\overline{L}_{c}}{\epsilon}\big)^{2d}J^2d^3\Big)=O(\Lambda^2Nd^3)$ 
    arithmetic operations.
    Moreover, 
    for $i=\nobreak1,\ldots,N$,
    let 
    $(\CG_i,\overline{D}_i)$ be the output 
    of $\mathtt{InitTestFuncs}$
    in Line~\ref{alglin:complexity-overall-init-test-funcs-forloop} 
    where $\CG_i$ is encoded by
    $\big(d_i,m_i,(\BIv_{i,l})_{l=0:m_i},\allowbreak \widetilde{J}_i,\allowbreak\big(\varrho_{i,\tilde{j}},\allowbreak(\varkappa_{i,\tilde{j},\iota})_{\iota=0:\varrho_{i,\tilde{j}}},\varpi_{i,\tilde{j}}\big)_{\tilde{j}=1:\widetilde{J}_i},\overline{\eta}_i\big)$ as in \ref{enc:testfuncs}.
    Then,
    Lemma~\ref{lem:proc-init-test-funcs}\ref{lems:proc-init-test-funcs-subdivision}
    shows that 
    $\max_{\BIx,\BIx'\in\CX_i}\big\{\|\BIx-\nobreak\BIx'\|_2\big\}\le\overline{D}_i$ $\forall i\in\nobreak\{0,1,\ldots,N\}$,
    $\max_{0\le i\le N}\{\overline{D}_i\}=O(D)$,
    $\max_{0\le i\le N}\{\overline{\eta}_i\}\le \frac{\epsilon}{16NL_c}$,
    $m:=\max_{0\le i\le N}\{m_i\}=O(\Lambda d)$,
    $\widetilde{J}:=\max_{0\le i\le N}\big\{|\FM(\FS_i)|\big\}=O(\Lambda)$.
    Furthermore,
    Lemma~\ref{lem:proc-init-test-funcs}\ref{lems:proc-init-test-funcs-LB}
    implies that 
    \begin{align}
        \bar{g}_{\ssmin}&:=\min_{1\le i\le N,\,0\le j\le m_i}\bigg\{\int_{\CX_i}g_{i,j}\DIFFX{\mu_i}\bigg\}
        \ge \frac{I_{\ssmin}}{d+2}\bigg(\frac{\epsilon}{64NDd\overline{L}_{c}}\bigg)^{d+1}> \Lambda^{-3}I_{\ssmin}.
        \label{eqn:complexity-explicit-proof-gmin}
    \end{align}
    Subsequently, 
    Lemma~\ref{lem:proc-init-type-measure}\ref{lems:proc-init-type-measure-complexity} shows that 
    Line~\ref{alglin:complexity-overall-init-type-measure-forloop} incurs 
    $O\big(N2^d\widetilde{J}(T_{\mathtt{Int}}+\widetilde{J}d)\big)=O(\Lambda^2 N2^d d + \Lambda N2^d T_{\mathtt{Int}})$
    arithmetic operations.
    Moreover, letting 
    $\big(\widetilde{R}_i,(F_{i,\tilde{r}},\tilde{\varpi}_{i,\tilde{r}})_{\tilde{r}=1:\widetilde{R}_i},\BP_i\big)$ be the encoding of $\mathtt{Aux}_i$ as in \ref{enc:typemeasures-aux} for $i=1,\ldots,N$,
    we get 
    $\widetilde{R}:=\max_{1\le i\le N}\{\widetilde{R}_i\}=O(\Lambda 2^d)$.

    Next, let us apply Theorem~\ref{thm:proc-solve-LSIP} to derive the number of arithmetic operations incurred by Lines~\ref{alglin:complexity-overall-solve-LSIP-dual}--\ref{alglin:complexity-overall-solve-LSIP-primal}.
    To that end,
    let us first denote 
    $c_{\ssmin}:=\min_{1\le i\le N,\,0\le j\le m_i,\,0\le l\le m_0}\big\{c_i(\BIv_{i,j},\BIv_{0,l})\big\}$,
    $c_{\ssmax}:=\nobreak\max_{1\le i\le N,\,0\le j\le m_i,\,0\le l\le m_0}\big\{c_i(\BIv_{i,j},\BIv_{0,l})\big\}$,
    and observe that the $L_c$-Lipschitz continuity of $(c_i)_{i=1:N}$ implies that
    \begin{align}
        \begin{split}
        c_{\ssmax}-c_{\ssmin}
        &\le \max_{1\le i\le N}\!\Bigg\{L_c\Bigg(\max_{\BIx,\BIx'\in\CX_i}\big\{\|\BIx-\BIx'\|_2\big\}+\max_{\BIz,\BIz'\in\CX_0}\big\{\|\BIz-\BIz'\|_2\big\}\Bigg)\Bigg\}\le 2L_c\max_{0\le i\le N}\{\overline{D}_i\}.
        \end{split}
        \label{eqn:complexity-explicit-proof-cmaxcmin}
    \end{align}
    Now, combining (\ref{eqn:complexity-explicit-proof-gmin}), 
    (\ref{eqn:complexity-explicit-proof-cmaxcmin}), 
    and 
    $\epsilon_0:=\frac{\epsilon}{4}$, $\varsigma:=\frac{\epsilon}{8N(L_c\vee 1)\max_{0\le i\le N}\{\overline{D}_i\}}$
    shows that 
    \begin{align*}
        \log\bigg(\frac{Nm(c_{\ssmax}-c_{\ssmin}+1)}{\epsilon_0\bar{g}_{\ssmin}}\bigg)
        &=O\bigg(\log\bigg(\frac{\Lambda}{I_{\ssmin}}\bigg)\bigg),\\ 
        \log(Nm)^2\log\bigg(\frac{Nm}{\varsigma}\bigg)&=O\big(\log(\Lambda)\big)^3.
    \end{align*}
    Because of these rates, 
    Theorem~\ref{thm:proc-solve-LSIP}\ref{thms:proc-solve-LSIP-dual-complexity} 
    with respect to $\epsilon_0\leftarrow \frac{\epsilon}{4}$
    implies that
    Line~\ref{alglin:complexity-overall-solve-LSIP-dual} incurs
    $O\Big(\Lambda^{\omega+1}N^{\omega+1}d^{\omega+1}\log\big(\frac{\Lambda}{I_{\ssmin}}\big)+\Lambda^{3}N^2d^3\log\big(\frac{\Lambda}{I_{\ssmin}}\big)T_{\mathtt{Eval}}\Big)$
    arithmetic operations,
    and
    Theorem~\ref{thm:proc-solve-LSIP}\ref{thms:proc-solve-LSIP-primal-complexity}
    with respect to $\varsigma\leftarrow \frac{\epsilon}{8N(L_c\vee 1)\max_{0\le i\le N}\{\overline{D}_i\}}$
    implies that 
    Line~\ref{alglin:complexity-overall-solve-LSIP-primal}
    incurs 
    $O\big(\Lambda^{\omega}N^{\omega}d^{\omega}\log(\Lambda)^3+\Lambda^2Nd^2T_{\mathtt{Eval}}\big)$
    arithmetic operations.
    Moreover, 
    (\ref{eqn:complexity-explicit-proof-cmaxcmin}) and
    Theorem~\ref{thm:proc-solve-LSIP}\ref{thms:proc-solve-LSIP-gap}
    guarantee that 
    $(\hat{\theta}_i)_{i=1:N}$ and
    $(\hat{y}_{i,0},\hat{\BIy}_i,\hat{\BIw}_i)_{i=1:N}$
    computed by 
    Lines~\ref{alglin:complexity-overall-solve-LSIP-dual}--\ref{alglin:complexity-overall-solve-LSIP-primal}
    constitute a pair of 
    $(\epsilon_{\sspar},\varsigma)$-optimizers of 
    \eqref{eqn:mt-tf-dual} and \eqref{eqn:mt-tf-lsip},
    where 
    \begin{align}
        \epsilon_{\sspar}
        &:=\epsilon_{0} + (c_{\ssmax}-c_{\ssmin})\varsigma + L_{c}\Bigg(N\overline{\eta}_0+\sum_{i=1}^{N}\overline{\eta}_i\Bigg)
        \le \frac{\epsilon}{4}+\frac{\epsilon}{4N}+\frac{\epsilon}{8}
        \le\frac{\epsilon}{2}.
        \label{eqn:complexity-explicit-proof-epspar}
    \end{align}

    Furthermore, 
    Lemma~\ref{lem:proc-construct-trans-funcs}
    shows that Line~\ref{alglin:complexity-overall-construct-trans-funcs}
    incurs 
    $O(Nm_0)=O(\Lambda Nd)$
    arithmetic operations,
    whereas 
    Lemma~\ref{lem:proc-construct-allocs}\ref{lems:proc-construct-allocs-complexity}
    shows that 
    Line~\ref{alglin:complexity-overall-construct-allocs} 
    incurs
    $O\big(N(\widetilde{R}m^{\omega-1}+m^{\omega})\big)=O(\Lambda^{\omega}N2^dd^{\omega-1})$
    arithmetic operations.
    Since 
    it holds for $i=0,1,\ldots,N$ that
    $\overline{\epsilon}_{i}(\varsigma)
    =2\eta(\FS_i)+\varsigma\max_{\BIx,\BIx'\in\CX_i}\big\{\|\BIx-\BIx'\|_2\big\}
    \le\nobreak \frac{\epsilon}{4NL_c}$,
    we can conclude by 
    Theorem~\ref{thm:proc-solve-LSIP}\ref{thms:proc-solve-LSIP-gap},
    (\ref{eqn:complexity-explicit-proof-epspar}),
    Lemma~\ref{lem:proc-construct-trans-funcs},
    Lemma~\ref{lem:proc-construct-allocs}\ref{lems:proc-construct-allocs-error},
    and Theorem~\ref{thm:mt-tf-approx}
    that 
    $(\hat{\varphi}_i,\hat{\gamma}_i)_{i=1:N},\hat{\nu}$
    computed by Algorithm~\ref{alg:complexity-overall}
    is an $\epsilon_{\ssapx}$-matching equilibrium
    with 
    $\epsilon_{\ssapx}=\nobreak\epsilon_{\sspar}+L_{c}\big(N\overline{\epsilon}_{0}(\varsigma)+\sum_{i=1}^{N}\overline{\epsilon}_{i}(\varsigma)\big)
    \le \frac{\epsilon}{2}+\frac{\epsilon}{2}=\epsilon$.
    This shows that 
    $(\hat{\varphi}_i,\hat{\gamma}_i)_{i=1:N},\hat{\nu}$
    is indeed an \mbox{$\epsilon$-matching} equilibrium.

    Now, let us summarize the number of arithmetic operations incurred by each step of Algorithm~\ref{alg:complexity-overall} in the following list.
    \begin{itemize}
        \item \textbf{Line~\ref{alglin:complexity-overall-init-test-funcs-forloop}:}
        $O(\Lambda^2Nd^3)$.
        
        \item \textbf{Line~\ref{alglin:complexity-overall-init-type-measure-forloop}:}
        $O(\Lambda^2N2^dd+\Lambda N2^dT_{\mathtt{Int}})$.
        
        \item \textbf{Line~\ref{alglin:complexity-overall-solve-LSIP-dual}:}
        $O\Big(\Lambda^{\omega+1}N^{\omega+1}d^{\omega+1}\log\big(\frac{\Lambda}{I_{\ssmin}}\big)+\Lambda^{3}N^2d^3\log\big(\frac{\Lambda}{I_{\ssmin}}\big)T_{\mathtt{Eval}}\Big)$.
        
        \item \textbf{Line~\ref{alglin:complexity-overall-solve-LSIP-primal}:}
        $O\big(\Lambda^{\omega}N^{\omega}d^{\omega}\log(\Lambda)^3+\Lambda^2Nd^2T_{\mathtt{Eval}}\big)$.
        
        \item \textbf{Line~\ref{alglin:complexity-overall-construct-trans-funcs}:}
        $O(\Lambda Nd)$.
        
        \item \textbf{Line~\ref{alglin:complexity-overall-construct-allocs}:}
        $O(\Lambda^{\omega}N2^dd^{\omega-1})$.
    \end{itemize}
    Therefore, we conclude that Algorithm~\ref{alg:complexity-overall}
    computes an $\epsilon$-matching equilibrium 
    $(\hat{\varphi}_i,\hat{\gamma}_i)_{i=1:N},\hat{\nu}$
    while incurring 
    $O\Big(\Lambda^{\omega+1}N^{\omega+1}d^{\omega+1}\log\big(\frac{\Lambda}{I_{\ssmin}}\big)+\Lambda^{3}N^2d^3\log\big(\frac{\Lambda}{I_{\ssmin}}\big)T_{\mathtt{Eval}}+\Lambda N 2^dT_{\mathtt{Int}}\Big)$
    arithmetic operations.
    Lastly, 
    Lemma~\ref{lem:proc-construct-allocs}\ref{lems:proc-construct-allocs-samp-qual} 
    implies that 
    generating a random sample from $\hat{\nu}$ incurs 
    $O(m)=O(\Lambda d)$ 
    arithmetic operations,
    and for $i=1,\ldots,N$,
    Lemma~\ref{lem:proc-construct-trans-funcs} 
    implies that evaluating $\hat{\varphi}_i(\BIz)$ at any $\BIz\in\nobreak\CX_0$
    incurs 
    $O\big(|\FM(\FS_0)|d_0^{\omega}\big)=O(\Lambda d^{\omega})$
    arithmetic operations,
    whereas 
    Lemma~\ref{lem:proc-construct-allocs}\ref{lems:proc-construct-allocs-samp-alloc} implies that
    generating a random sample from $\hat{\gamma}_i$ incurs
    $O(\widetilde{R}+m+T_{\mathtt{Samp}})=O(\Lambda 2^d+T_{\mathtt{Samp}})$
    arithmetic operations.
    The proof of statement~\ref{thms:complexity-explicit-general}
    is now complete.

    Finally, 
    under the particular class of $(\mu_i)_{i=1:N}$ in Setting~\ref{sett:Euclidean}\ref{settc:Euclidean-particular},
    Lemma~\ref{lem:proc-int-simplex}\ref{lems:proc-int-simplex-int}
    implies that  
    $T_{\mathtt{Int}}=O\big(d^{\omega+1}+\log(J)\big)$,
    and 
    Lemma~\ref{lem:proc-int-simplex}\ref{lems:proc-int-simplex-samp}
    implies that  
    $T_{\mathtt{Samp}}=O\big(d^{\omega+1}+\log(J)\big)$.
    The proof is now complete.
\endproof

\proof{Proof of Proposition~\ref{prop:cp-tf-properties}.}
To begin,
statement~\ref{props:cp-tf-boundedness}
is a direct consequence of the results in Section~\ref{sapx:proof-complexity-lsip}.
Let $(y_{i,0}^\star,\BIy_{i}^\star,\BIw_{i}^\star)_{i=1:N}$ be an arbitrary
optimizer of $\CPDLPDUAL{0}$
with respect to 
$\CK_i^{(0)}\leftarrow\nobreak\big\{(\BIv_{i,j},\BIv_{0,0}):j\in\{0,1,\ldots,m_i\}\big\}
\cup \big\{(\BIv_{i,0},\BIv_{0,l}):l\in\{1,\ldots,m_0\}\big\}$ 
$\forall i\in\nobreak\{1,\ldots,N\}$.
Let us define $c_{\ssmax}:=\max_{1\le i\le N,\, 0\le j\le m_i,\, 0\le l\le m_0}\big\{c_i(\BIv_{i,j},\BIv_{0,l})\big\}$.
One observes that 
$(y_{i,0}^\star-c_{\ssmax},\BIy_{i}^\star,\BIw_{i}^\star)_{i=1:N}$
is feasible for the LP problem
\LPTAG{$\CV$} in Section~\ref{sapx:proof-complexity-lsip}.
Consequently, 
Lemma~\ref{lem:lp-any-index}\ref{lems:lp-any-index-Linfinity-bound}
in Section~\ref{sapx:proof-complexity-lsip}
provides an explicit upper bound for 
$\big\|(y_{1,0}^{\star}-c_{\ssmax},\BIy_{1}^{\star\TRANSP},\BIw_{1}^{\star\TRANSP},\ldots,y_{N,0}^{\star}-c_{\ssmax},\BIy_{N}^{\star\TRANSP},\BIw_{N}^{\star\TRANSP})^\TRANSP\big\|_{\infty}$
that does not depend on 
$(y_{i,0}^\star,\BIy_{i}^\star,\BIw_{i}^\star)_{i=1:N}$,
which proves statement~\ref{props:cp-tf-boundedness}.

Next, let us prove statement~\ref{props:cp-tf-termination}.
Let $q:=\nobreak\sum_{i=1}^N(1+\nobreak m_i+\nobreak m_0)$, 
and let us adopt the concise notations $\bar{\BIg}$, $(\BIa_i(\,\cdot\,,\cdot\,))_{i=1:N}$, $(\BIu^{(j)})_{j=1:m_0}$ defined in the proof of Theorem~\ref{thm:mt-tf-duality} 
in Section~\ref{sapx:proof-parametric} 
to re-express  \eqref{eqn:mt-tf-algo-lp} as follows:
\begin{align}
\begin{split}
\maximize_{\Bchi}\quad & \langle\bar{\BIg},\Bchi\rangle\\
\mathrm{subject~to}\quad & \langle\BIa_i(\BIx,\BIz),\Bchi\rangle  \le c_i(\BIx,\BIz) \quad \forall (\BIx,\BIz)\in\CK^{(r)}_i,\; \forall 1\le i\le N,\\
& \langle\BIu^{(j)},\Bchi\rangle= 0  \quad \forall 1\le j\le m_0, \qquad  \Bchi\in\R^{q}.
\end{split}
\label{eqn:cp-tf-properties-proof-lp}
\end{align}
Moreover, let us define 
$s_i(\Bchi,\BIx,\BIz):=c_i(\BIx,\BIz)-\langle\BIa_i(\BIx,\BIz),\Bchi\rangle$ $\forall\Bchi\in\R^{q}$, 
$\forall (\BIx,\BIz)\in\nobreak\CX_i\times\nobreak\CX_0$, 
$\forall 1\le\nobreak i\le\nobreak N$,
and denote $S^{\star}:=\big\{\Bchi\in\R^{q}:\langle\BIa_i(\BIx,\BIz),\Bchi\rangle\le c_i(\BIx,\BIz)$ 
$\forall (\BIx,\BIz)\in\CK^{(0)}_i,$ 
$\forall 1\le\nobreak i\le\nobreak N,$ $\langle\BIu^{(l)},\Bchi\rangle=\nobreak0$ 
$\forall 1\le\nobreak l\le\nobreak m_0,$ 
$\langle\bar{\BIg},\Bchi\rangle\ge\alpha_{\sspar}^\star\big\}$.
Since the set of optimizers of $\CPDLPDUAL{0}$ is bounded, 
it follows from \citep[Theorem~27.1(f)]{ECrockafellar1970convex}
that $S^{\star}$ is compact.
Let us suppose for the sake of contradiction that Algorithm~\ref{alg:cp-tf} does not terminate and produces an infinite sequence 
$(\Bchi^{(r)})_{r\in\N_0}\subset\R^{q}$,
where we denote 
$\Bchi^{(r)}=\nobreak(y_{1,0}^{(r)},\BIy_{1}^{(r)\TRANSP},\BIw_{i}^{(r)\TRANSP},\ldots,y_{N,0}^{(r)},\BIy_{N}^{(r)\TRANSP},\BIw_{N}^{(r)\TRANSP})^\TRANSP$
for 
$y_{i,0}^{(r)}\in\R$,
$\BIy_{i}^{(r)}\in\R^{m_i}$,
$\BIw_{i}^{(r)}\in\R^{m_0}$
$\forall 1\le\nobreak i\le\nobreak N$,
$\forall r\in\N_{0}$.
Since 
$\CK^{(0)}_i\subseteq \CK^{(r)}_i\subseteq \CX_i\times\CX_0$ 
$\forall 1\le i\le N$, 
$\forall r\in\N$, 
and since the optimal value of 
(\ref{eqn:cp-tf-properties-proof-lp}) 
is at least~$\alpha_{\sspar}^\star$ for all $r\in\N_0$,
it holds that 
$\Bchi^{(r)}\in S^{\star}$ $\forall r\in\N_0$. 
Thus, 
extracting a subsequence if necessary, 
let us assume without loss of generality that 
$\lim_{r\to\infty}\Bchi^{(r)}=\Bchi^{(\infty)}\in\R^{q}$.
Since Line~\ref{alglin:cp-tf-aggregate} implies that 
$(\tilde{\BIx}_i^{(r)},\tilde{\BIz}_i^{(r)})\in\CK^{(t)}_i$ 
$\forall 1\le i\le N$, 
$\forall t>r$, 
it follows that 
$\langle\BIa_i(\tilde{\BIx}^{(r)}_i,\tilde{\BIz}^{(r)}_i),\Bchi^{(t)}\rangle\le c_i(\tilde{\BIx}^{(r)}_i,\tilde{\BIz}^{(r)}_i)$ 
$\forall 1\le i\le N$, 
$\forall t>r$, 
and hence 
\begin{align}
    s_i\big(\Bchi^{(\infty)},\tilde{\BIx}^{(r)}_i,\tilde{\BIz}^{(r)}_i\big) 
    &= c_i(\tilde{\BIx}^{(r)}_i,\tilde{\BIz}^{(r)}_i)
    -\lim_{t\to\infty} \langle\BIa_i(\tilde{\BIx}^{(r)}_i,\tilde{\BIz}^{(r)}_i),\Bchi^{(t)}\rangle 
    \ge 0 
    \qquad \forall 1\le i\le N,\; \forall r\in\N_0.
    \label{eqn:cp-tf-properties-proof-ineq1}
\end{align}
Moreover, it follows from Line~\ref{alglin:cp-tf-global} that
\begin{align}
    \begin{split}
        s_i\big(\Bchi^{(r)},\tilde{\BIx}^{(r)}_i,\tilde{\BIz}^{(r)}_i\big)
        &= c_i(\tilde{\BIx}^{(r)}_i,\tilde{\BIz}^{(r)}_i)
        -y_{i,0}^{(r)}
        -\langle\BIg_i(\tilde{\BIx}^{(r)}_i),\BIy_{i}^{(r)}\rangle
        -\langle\BIg_0(\tilde{\BIz}^{(r)}_i),\BIw_i^{(r)}\rangle \\
        &\le \underline{\beta}_i^{(r)} - y_{i,0}^{(r)} + \tau^{(r)} \hspace{90pt} \forall r\in\N_0,\;\forall 1\le i\le N.
    \end{split}
    \label{eqn:cp-tf-properties-proof-ineq2}
\end{align}
Subsequently, combining (\ref{eqn:cp-tf-properties-proof-ineq1}), (\ref{eqn:cp-tf-properties-proof-ineq2}), the bound
$\|\BIa_i(\BIx,\BIz)\|_{\infty}\le 1$ $\forall (\BIx,\BIz)\in\CX_i\times\CX_0$, $\forall 1\le\nobreak i\le\nobreak N$,
and Line~\ref{alglin:cp-tf-global-tol}
yields
\begin{align*}
    \limsup_{r\to\infty} y_{i,0}^{(r)}-\underline{\beta}^{(r)}_i 
    &\le \limsup_{r\to\infty}  \tau^{(r)} 
    - s_i\big(\Bchi^{(r)},\tilde{\BIx}_i^{(r)},\tilde{\BIz}_i^{(r)}\big) 
    + s_i\big(\Bchi^{(\infty)},\tilde{\BIx}_i^{(r)},\tilde{\BIz}_i^{(r)}\big)\\
    &\le \limsup_{r\to\infty} \tau^{(r)} 
    + \big|s_i\big(\Bchi^{(r)},\tilde{\BIx}_i^{(r)},\tilde{\BIz}_i^{(r)}\big)
    -s_i\big(\Bchi^{(\infty)},\tilde{\BIx}_i^{(r)},\tilde{\BIz}_i^{(r)}\big)\big| \\
    &\le \limsup_{r\to\infty}\tau^{(r)}
    +\big\|\BIa_i\big(\tilde{\BIx}_i^{(r)},\tilde{\BIz}_i^{(r)}\big)\big\|_{\infty}\big\|\Bchi^{(r)}-\Bchi^{(\infty)}\big\|_{1} \\
    &= \limsup_{r\to\infty} \tau^{(r)} 
    < \frac{\epsilon_{\sspar}}{N} 
    \hspace{99pt}\qquad \forall 1\le i\le N.
\end{align*}
Hence, there exists $\overline{r}\in\N_0$ such that 
$\sum_{i=1}^N \big(y_{i,0}^{(\overline{r})}-\underline{\beta}^{(\overline{r})}_i\big) 
< \epsilon_{\sspar}$, 
which implies by Line~\ref{alglin:cp-tf-termination} that Algorithm~\ref{alg:cp-tf} will terminate at iteration~$\overline{r}$. 
This completes the proof of statement~\ref{props:cp-tf-termination}.

Next, to prove statements~\ref{props:cp-tf-bounds}--\ref{props:cp-tf-primal}, 
we will show that 
$(\hat{y}_{i,0},\hat{\BIy}_i,\hat{\BIw}_i)_{i=1:N}$ 
is a feasible solution of \eqref{eqn:mt-tf-lsip} 
whose objective value is equal to 
$\alpha_{\sspar}^{\ssLB}$ and that 
$(\hat{\theta}_i)_{i=1:N}$ is a feasible solution of \eqref{eqn:mt-tf-dual}
whose objective value is equal to 
$\alpha_{\sspar}^{\ssUB}$. 
Subsequently, 
since Lines~\ref{alglin:cp-tf-termination}--\ref{alglin:cp-tf-bounds} guarantee that 
$\alpha_{\sspar}^{\ssUB}-\alpha_{\sspar}^{\ssLB}
=\sum_{i=1}^N\big(y_{i,0}^{(r)}-\underline{\beta}_i^{(r)}\big)
\le\epsilon_{\sspar}$, statements~\ref{props:cp-tf-bounds}--\ref{props:cp-tf-primal} will follow from the strong duality in Theorem~\ref{thm:mt-tf-duality}. 
On the one hand, it holds by Line~\ref{alglin:cp-tf-global} and Line~\ref{alglin:cp-tf-primal-dual} that
\begin{align*}
    \hspace{20pt}&\hspace{-20pt}c_i(\BIx,\BIz)-\hat{y}_{i,0}
    -\langle\BIg_i(\BIx),\hat{\BIy}_i\rangle
    -\langle\BIg_0(\BIz),\hat{\BIw}_i\rangle\\
    &=c_i(\BIx,\BIz)
    -\underline{\beta}_i^{(r)}
    -\langle\BIg_i(\BIx),\BIy_i^{(r)}\rangle
    -\langle\BIg_0(\BIz),\BIw_i^{(r)}\rangle\\
    &\ge c_i(\BIx,\BIz)
    -\langle\BIg_i(\BIx),\BIy_i^{(r)}\rangle
    -\langle\BIg_0(\BIz),\BIw_i^{(r)}\rangle
    -\min_{(\bar{\BIx},\bar{\BIz})\in\CX_i\times\CX_0}\big\{c_i(\bar{\BIx},\bar{\BIz})
    -\langle\BIg_i(\bar{\BIx}),\BIy_i^{(r)}\rangle
    -\langle\BIg_0(\bar{\BIz}),\BIw_i^{(r)}\rangle\big\}\\
    &\ge 0 
    \hspace{285pt} \forall (\BIx,\BIz)\in\CX_i\times\CX_0,\; \forall 1\le i\le N.
\end{align*}
Moreover, since 
$\big(y_{i,0}^{(r)},\BIy_i^{(r)},\BIw_i^{(r)}\big)_{i=1:N}$ 
is feasible for \eqref{eqn:mt-tf-algo-lp} by Line~\ref{alglin:cp-tf-lp}, 
it holds by Line~\ref{alglin:cp-tf-primal-dual} that 
$\sum_{i=1}^N\hat{\BIw}_i=\sum_{i=1}^N\BIw_i^{(r)}=\veczero_{m_0}$.
Furthermore, it follows from 
Line~\ref{alglin:cp-tf-lp}, 
Line~\ref{alglin:cp-tf-bounds}, and 
Line~\ref{alglin:cp-tf-primal-dual} that
\begin{align*}
    \sum_{i=1}^N\hat{y}_{i,0}+\langle\bar{\BIg}_i,\hat{\BIy}_i\rangle
    &=\sum_{i=1}^N\underline{\beta}_i^{(r)}
    +\langle\bar{\BIg}_i,\BIy_i^{(r)}\rangle
    =\Bigg(\sum_{i=1}^Ny_{i,0}^{(r)}+\langle\bar{\BIg}_i,\BIy_i^{(r)}\rangle\Bigg)
    -\Bigg(\sum_{i=1}^Ny_{i,0}^{(r)}-\underline{\beta}_i^{(r)}\Bigg)
    =\alpha_{\sspar}^{\ssLB}.
\end{align*}
This shows that $(\hat{y}_{i,0},\hat{\BIy}_i,\hat{\BIw}_i)_{i=1:N}$ 
is a feasible solution of \eqref{eqn:mt-tf-lsip} with objective value 
$\alpha_{\sspar}^{\ssLB}$. 
On the other hand, by Line~\ref{alglin:cp-tf-lp}, 
$\big(\theta^{(r)}_{i,\BIx,\BIz}\big)_{(\BIx,\BIz)\in\CK_i^{(r)},\,i=1:N}$, $\Bxi^{(r)}$ is an optimizer of \eqref{eqn:mt-tf-algo-lp-dual}. 
Let us denote $\Bxi^{(r)}=\big(\xi^{(r)}_1,\ldots,\xi^{(r)}_{m_0}\big)^\TRANSP$. 
Consequently, it holds by Line~\ref{alglin:cp-tf-primal-dual} and the constraints of \eqref{eqn:mt-tf-algo-lp-dual} that, 
for $i=1,\ldots,N$, 
$\hat{\theta}_i$ is a positive Borel measure on $\CX_i\times\CX_0$ with finite support which satisfies
\begin{align}
    \hat{\theta}_i(\CX_i\times\CX_0)
    &=\sum_{(\BIx,\BIz)\in\CK_i^{(r)}}\theta^{(r)}_{i,\BIx,\BIz}=1,\nonumber\\
    \int_{\CX_i\times\CX_0}g_{i,j}(\BIx)\DIFFM{\hat{\theta}_i}{\DIFF \BIx,\DIFF \BIz}
    &=\sum_{(\BIx,\BIz)\in\CK_i^{(r)}}\theta^{(r)}_{i,\BIx,\BIz}g_{i,j}(\BIx)
    =\int_{\CX_i}g_{i,j}\DIFFX{\mu_i} & \forall 1\le j\le m_i,
    \label{eqn:cp-tf-proof-primal}\\
    \int_{\CX_i\times\CX_0}g_{0,l}(\BIz)\DIFFM{\hat{\theta}_i}{\DIFF \BIx,\DIFF \BIz}
    &=\sum_{(\BIx,\BIz)\in\CK_i^{(r)}}\theta^{(r)}_{i,\BIx,\BIz}g_{0,l}(\BIz)
    =\xi^{(r)}_l & \forall 1\le l\le m_0.\nonumber
\end{align}
Thus, $\hat{\theta}_i\in\CP(\CX_i\times\CX_0)$ $\forall 1\le i\le N$. 
For $i=1,\ldots,N$, 
let $\hat{\mu}_i$ and $\hat{\nu}_i$ denote the marginals of $\hat{\theta}_i$ on $\CX_i$ and $\CX_0$, respectively. 
It hence follows from (\ref{eqn:cp-tf-proof-primal}) that $\hat{\mu}_i\overset{\CG_i}{\sim}\mu_i$, $\hat{\nu}_i\overset{\CG_0}{\sim}\hat{\nu}_1$ $\forall 1\le i\le N$. 
Moreover, 
it follows from 
Line~\ref{alglin:cp-tf-lp}, 
Line~\ref{alglin:cp-tf-bounds}, 
Line~\ref{alglin:cp-tf-primal-dual}, and 
the strong duality of LP problems that
$\sum_{i=1}^N\int_{\CX_i\times\CX_0}c_i\DIFFX{\hat{\theta}_i}
=\sum_{i=1}^N\sum_{(\BIx,\BIz)\in\CK_i^{(r)}}\theta^{(r)}_{i,\BIx,\BIz}c_i(\BIx,\BIz)
=\alpha^{(r)}
=\alpha_{\sspar}^{\ssUB}$.
Therefore, $(\hat{\theta}_i)_{i=1:N}$ is a feasible solution of \eqref{eqn:mt-tf-dual} with objective value $\alpha_{\sspar}^{\ssUB}$. 
The proof is now complete. 
\endproof

\proof{Proof of Theorem~\ref{thm:mt-tf}.}
    To begin,
    it follows from 
    Line~\ref{alglin:mt-tf-cpalgo} and
    Proposition~\ref{prop:cp-tf-properties} that 
    $(\hat{\theta}_i)_{i=1:N}$,
    $(\hat{y}_{i,0},\hat{\BIy}_i,\hat{\BIw}_i)_{i=1:N}$
    constitute a pair of 
    $(\epsilon_{\sspar},0)$-optimizers of 
    \eqref{eqn:mt-tf-dual} and \eqref{eqn:mt-tf-lsip}.
    Line~\ref{alglin:mt-tf-marginals} then guarantees that 
    $\hat{\mu}_i\overset{\CG_i}{\sim}\mu_i$,
    $\hat{\nu}_i\overset{\CG_0}{\sim}\hat{\nu}_1$
    $\forall 1\le i\le N$.
    Moreover, 
    Proposition~\ref{prop:cp-tf-properties}\ref{props:cp-tf-primal} implies that
    $(\hat{\theta}_i,\hat{\mu}_i,\hat{\nu}_i)_{i=1:N}$
    all have finite supports.
    Next,
    for $i=\nobreak 1,\ldots,N$,
    Lines~\ref{alglin:mt-tf-binding1}--\ref{alglin:mt-tf-reassembly}
    ensure that 
    the law of 
    $(\BIX_i,\BIZ_i):\Omega\to\CX_i\times\CX_0$
    is $\hat{\theta}_i$,
    the law of $\overbar{\BIX}_i:\Omega\to\CX_i$
    is $\mu_i$,
    and
    Theorem~\ref{thm:error-control}\ref{thms:error-control-W1} implies that
    $\EXP\big[\|\BIZ-\BIZ_i\|_2\big]
    =W_1(\hat{\nu},\hat{\nu}_i)
    \le \overline{\epsilon}_{0}(0)$
    and that
    $\EXP\big[\|\BIX_i-\overbar{\BIX}_i\|_2\big]
    =W_1(\hat{\mu}_i,\mu_i)
    \le \overline{\epsilon}_{i}(0)$.
    Combining these properties with 
    Line~\ref{alglin:mt-tf-dualsol},
    Lines~\ref{alglin:mt-tf-primal2}--\ref{alglin:mt-tf-primalcoup},
    and the definition of $\epsilon_{\sstheo}$,
    it follows from Theorem~\ref{thm:mt-tf-approx} that 
    $(\hat{\varphi}_i,\hat{\gamma}_i)_{i=1:N},\hat{\nu}$
    and 
    $(\hat{\varphi}_i,\tilde{\gamma}_i)_{i=1:N},\tilde{\nu}$
    are both $\epsilon_{\sstheo}$-matching equilibria.

    On the other hand,
    since Line~\ref{alglin:mt-tf-dualsol} and the constraints of \eqref{eqn:mt-tf-lsip} guarantee that
    \begin{align*}
        \hat{\varphi}_i^{c_i}(\BIx)
        &=\inf_{\BIz\in\CX_0}\big\{c_i(\BIx,\BIz)-\langle\BIg_0(\BIz),\hat{\BIw}_i\rangle\big\}
        \ge \hat{y}_{i,0}+\langle\BIg_i(\BIx),\hat{\BIy}_i\rangle 
        \qquad \forall \BIx\in\CX_i,\; \forall 1\le i\le N,
    \end{align*}
    we have by Lines~\ref{alglin:mt-tf-primal2}--\ref{alglin:mt-tf-bounds} 
    and
    Proposition~\ref{prop:cp-tf-properties}\ref{props:cp-tf-dual}
    that
    \begin{align*}
        \hspace{20pt}&\hspace{-20pt}\sum_{i=1}^{N}\int_{\CX_i\times\CX_0}\big(c_i(\BIx,\BIz)-\hat{\varphi}_i(\BIz)-\hat{\varphi}_i^{c_i}(\BIx)\big) \DIFFM{\hat{\gamma}_i}{\DIFF\BIx,\DIFF\BIz}\\
        &\le \Bigg(\sum_{i=1}^{N}\int_{\CX_i\times\CX_0}c_i\DIFFX{\hat{\gamma}_i}\Bigg)
        - \int_{\CX_0}\big({\textstyle\sum_{i=1}^{N}\hat{\varphi}_i}\big)\DIFFX{\hat{\nu}} 
        - \Bigg(\sum_{i=1}^{N}\hat{y}_{i,0}+\int_{\CX_i}\langle\BIg_i(\BIx),\hat{\BIy}_i\rangle\DIFFM{\mu_i}{\DIFF\BIx}\Bigg)\\
        &=\Bigg(\sum_{i=1}^{N}\EXP\big[c_i(\overbar{\BIX}_i,\BIZ)\big]\Bigg)
        - \Bigg(\sum_{i=1}^{N}\hat{y}_{i,0}+\langle\bar{\BIg}_i,\hat{\BIy}_i\rangle\Bigg)
        =\hat{\alpha}_{\MT}^{\ssUB}-\hat{\alpha}_{\MT}^{\ssLB}
        =\hat{\epsilon}_{\sssub},\\
        \hspace{20pt}&\hspace{-20pt}\sum_{i=1}^{N}\int_{\CX_i\times\CX_0}\big(c_i(\BIx,\BIz)-\hat{\varphi}_i(\BIz)-\hat{\varphi}_i^{c_i}(\BIx)\big) \DIFFM{\tilde{\gamma}_i}{\DIFF\BIx,\DIFF\BIz}\\
        &\le \Bigg(\sum_{i=1}^{N}\int_{\CX_i\times\CX_0}c_i\DIFFX{\tilde{\gamma}_i}\Bigg)
        - \int_{\CX_0}\big({\textstyle\sum_{i=1}^{N}\hat{\varphi}_i}\big)\DIFFX{\tilde{\nu}} 
        - \Bigg(\sum_{i=1}^{N}\hat{y}_{i,0}+\int_{\CX_i}\langle\BIg_i(\BIx),\hat{\BIy}_i\rangle\DIFFM{\mu_i}{\DIFF\BIx}\Bigg)\\
        &=\Bigg(\sum_{i=1}^{N}\EXP\big[c_i(\overbar{\BIX}_i,\overbar{\BIZ})\big]\Bigg)
        - \Bigg(\sum_{i=1}^{N}\hat{y}_{i,0}+\langle\bar{\BIg}_i,\hat{\BIy}_i\rangle\Bigg)
        =\tilde{\alpha}_{\MT}^{\ssUB}-\hat{\alpha}_{\MT}^{\ssLB}
        =\tilde{\epsilon}_{\sssub},
    \end{align*}
    as well as 
    \begin{align*}
        \tilde{\alpha}_{\MT}^{\ssUB}
        &=\EXP\Big[{\textstyle\sum_{i=1}^{N}c_i(\overbar{\BIX}_i,\overbar{\BIZ})}\Big]
        \le \EXP\Big[{\textstyle\sum_{i=1}^{N}c_i(\overbar{\BIX}_i,\BIZ)}\Big]
        =\hat{\alpha}_{\MT}^{\ssUB}.
    \end{align*}
    The proofs of statements~\ref{thms:mt-tf-algo-equilibrium1}--\ref{thms:mt-tf-algo-bound} are now complete.

    Lastly, to prove
    statement~\ref{thms:mt-tf-algo-control},
    let us fix arbitrary $\epsilon>0$, $\epsilon_{\sspar}\in(0,\epsilon)$,
    and use the procedure
    $\mathtt{InitTestFuncs}\big(\CX_i,\frac{(\epsilon-\epsilon_{\sspar})\wedge 1}{4N(L_c\vee 1)}\big)$ 
    as well as Lemma~\ref{lem:proc-init-test-funcs}
    to construct $\FS_i$
    for $i=0,1,\ldots,N$.
    It thus holds by Lemma~\ref{lem:proc-init-test-funcs}\ref{lems:proc-init-test-funcs-subdivision}
    that $\eta(\FS_i)\le \frac{\epsilon-\epsilon_{\sspar}}{4NL_c}$
    for $i=0,1,\ldots,N$.
    Then, 
    statement~\ref{thms:mt-tf-algo-bound} yields 
    $\tilde{\epsilon}_{\sssub}\le\nobreak \hat{\epsilon}_{\sssub}\le \epsilon_{\sstheo}=\epsilon_{\sspar}+2L_c\big(N\eta(\FS_0)+\sum_{i=1}^{N}\eta(\FS_i)\big)
    \le \epsilon_{\sspar}+\epsilon-\epsilon_{\sspar}=\epsilon$.
    The proof is now complete.
\endproof

\subsection{Proofs of results in Section~\ref{sapx:complexity-init}}\label{sapx:proof-complexity-init}

\proof{Proof of Lemma~\ref{lem:proc-edgewise-divide}.}
    The procedure \EdgewiseDivide translates the $k$-fold edge-wise subdivision procedure described in Section~2 of \citep{ECedelsbrunner1999edgewise} and the proof of the Counting Lemma in \citep[p.27]{ECedelsbrunner1999edgewise}.
    Consequently, the Main Theorem of \citep[p.24]{ECedelsbrunner1999edgewise} ensures that statements~\ref{lems:proc-edgewise-divide-simpcomp} and \ref{lems:proc-edgewise-divide-faces} hold.
    Statement~\ref{lems:proc-edgewise-divide-complexity} follows directly from combining the number of arithmetic operations in each line, as presented in the comments.
    Moreover, in statement~\ref{lems:proc-edgewise-divide-counts},
    $|\FM(\FS)|=k^n\le k^d$ follows from the Counting Lemma \citep[p.27]{ECedelsbrunner1999edgewise},
    and we subsequently get 
    $|V(\FS)|\le (d+1)k^d$ since each $S\in\FS$ has at most $d+1$ extreme points.

    To prove statement~\ref{lems:proc-edgewise-divide-barycoord},
    let us fix an arbitrary $\BIv_{s}\in V(C)$ where $s\in\{0,1,\ldots,n\}$,
    and  
    let us fix an arbitrary $S_j\in\FM(\FS)$ constructed in 
    Lines~\ref{alglin:edgewise-divide-construction-begin}--\ref{alglin:edgewise-divide-construction-end}.
    Observe that $r$ is first initialized to~0 in Line~\ref{alglin:edgewise-divide-construction-begin},
    remains~0 in the first iteration of the nested for-loops in Lines~\ref{alglin:edgewise-divide-row-loop}--\ref{alglin:edgewise-divide-column-loop},
    and is then incremented exactly $n$ times 
    throughout the nested for-loop in Lines~\ref{alglin:edgewise-divide-row-loop}--\ref{alglin:edgewise-divide-column-loop}.
    Hence, within the nested for-loop in Lines~\ref{alglin:edgewise-divide-row-loop}--\ref{alglin:edgewise-divide-column-loop},
    there is at least one execution of
    Line~\ref{alglin:edgewise-divide-check-end} 
    in which $r$ has value~$s$.
    This shows that 
    there exists $l\in\{0,1,\ldots,n\}$ such that 
    $\lambda_{\FF(C),\BIv_s}(\BIu_l)\ge \frac{1}{k}$,
    where $(\BIu_{l})_{l=0:n}$ are the final values of the extreme points of $S_j$
    in Line~\ref{alglin:edgewise-divide-construction-end}.
    This proves statement~\ref{lems:proc-edgewise-divide-barycoord}.

    Lastly, the definition of $\eta(\FS)$ implies that 
    $\eta(\FS)=\max_{S\in\FM(\FS)}\Big\{\max_{\BIv,\BIv'\in V(S)}\big\{\|\BIv-\BIv'\|_2\big\}\Big\}$.
    Thus, 
    to prove statement~\ref{lems:proc-edgewise-divide-meshsize},
    let us fix an arbitrary $S\in\FM(\FS)$, 
    let $(\BIu_l)_{l=0:n}$ denote their final values of the extreme points of $\widehat{S}$ in Line~\ref{alglin:edgewise-divide-construction-end} after the nested for-loops,
    and show that 
    $\|\BIu_{l}-\nobreak\BIu_{l'}\|_2\le \frac{d}{k}\max_{\BIv,\BIv'\in V(C)}\big\{\|\BIv-\BIv'\|_2\big\}$ for any $l,l'\in\{0,1,\ldots,n\}$.
    For every $l\in\{1,\ldots,n\}$,
    observe from \mbox{Lines~\ref{alglin:edgewise-divide-check-begin}--\ref{alglin:edgewise-divide-check-end}} that 
    the condition in Line~\ref{alglin:edgewise-divide-check-begin} is only true for exactly one $\hat{i}\in\{1,\ldots,k\}$,
    and hence $\BIu_{l}-\BIu_{l-1}=\frac{1}{k}(\BIv_{\hat{r}+1}-\BIv_{\hat{r}})$ for some $\hat{r}\in\{0,1,\ldots,n-1\}$,
    where $\BIv_{\hat{r}+1},\BIv_{\hat{r}}\in V(C)$.
    It thus holds that 
    $\|\BIu_{l}-\BIu_{l-1}\|_2\le \frac{1}{k}\max_{\BIv,\BIv'\in V(C)}\big\{\|\BIv-\BIv'\|_2\big\}$ for $l\in\{1,\ldots,n\}$,
    and it follows that
    $\|\BIu_{l}-\BIu_{l'}\|_2\le \sum_{s=l+1}^{l'}\|\BIu_{s}-\BIu_{s-1}\|_2\le \frac{d}{k}\max_{\BIv,\BIv'\in V(C)}\big\{\|\BIv-\BIv'\|_2\big\}$
    for all $l,l'\in\{0,1,\ldots,n\}$ with $l<l'$.
    The proof is now complete.
\endproof

Before proving Lemma~\ref{lem:proc-init-test-funcs}, 
we first establish the following auxiliary results.
In particular, Lemma~\ref{lem:simplex-integrals} below will be used again in the proof of 
Lemma~\ref{lem:proc-int-simplex}.

\begin{lemma}\label{lem:simplicial-complex-and-face}
    Let $C$ be an $n$-simplex in $\R^d$ where $n\in\{0,1,\ldots,d\}$,
    and let $\FK$ be a simplicial complex with $\bigcup_{K\in\FK}K=\nobreak C$.
    Let $F\in\FF(C)$ be an arbitrary non-empty face of $C$, 
    and define $\FK_{F}:=\{K\in\FK:K\subseteq F\}$.
    Then, the following statements hold.
    \begin{enumerate}[label=(\roman*),beginpenalty=10000]
        \item\label{lems:simplicial-complex-and-face-subset}
        It holds for all non-empty $K\in\FK$ that 
        $K\cap F\in \FF(K)$ and $\FF(K\cap F)\subseteq\FF(K)$.

        \item\label{lems:simplicial-complex-and-face-equality}
        $\FK_{F}=\{K\cap F: K\in\FK\}$.
        
    \end{enumerate}
\end{lemma}

\proof{Proof of Lemma~\ref{lem:simplicial-complex-and-face}.}
Let us first prove statement~\ref{lems:simplicial-complex-and-face-subset}.
Let $(\BIv_i)_{i=0:n}\subset \R^d$ be an arbitrary enumeration of $V(C)$.
Since $F\in\FF(C)$, there exists a non-empty $I\subseteq \{0,1,\ldots,n\}$ such that 
$F=\conv\big(\{\BIv_i:i\in I\}\big)$.
Let $I^{c}:=\{0,1,\ldots,n\}\setminus I$.
It thus holds that every $\BIx\in C$ can be uniquely expressed as 
$\sum_{i=0}^{n}\lambda_i \BIv_i$ where $(\lambda_{i})_{i=0:n}\subset [0,1]$ 
satisfy $\sum_{i=0}^{n}\lambda_i=\nobreak 1$,
and that 
$\BIx \in F$ if and only if 
$\lambda_i=\nobreak 0$ $\forall i\in I^c$ in this expression. 
Let us fix an arbitrary non-empty $K\in \FK$,
let 
$m:=|V(K)|-1$,
and let $(\BIu_{j})_{j=0:m}$ be an arbitrary enumeration of $V(K)$.
Moreover, for $j=0,1,\ldots,m$,
let $\BIu_{j}$ be uniquely expressed as 
$\BIu_{j}=\sum_{i=0}^{n}\zeta_{i,j}\BIv_{i}$ 
where $(\zeta_{i,j})_{i=0:n}\subset[0,1]$ satisfy 
$\sum_{i=0}^{n}\zeta_{i,j}=\nobreak 1$.
Subsequently, every $\BIx\in K$ can be uniquely expressed as 
$\BIx=\sum_{j=0}^{m}\xi_j\BIu_j=\sum_{i=0}^{n}\big(\sum_{j=0}^{m}\zeta_{i,j}\xi_j\big)\BIv_i$
where 
$(\xi_{j})_{j=0:m}\subset[0,1]$ satisfy $\sum_{j=0}^{m}\xi_j=\nobreak 1$,
and one checks that 
$\BIx\in K\cap F$ if and only if
$\sum_{j=0}^{m}\zeta_{i,j}\xi_j=\nobreak0$ $\forall i\in I^{c}$ in the aforementioned expression,
which holds if and only if
$\BIx\in \conv\big(\{\BIu_j\}_{j=0:m}\cap F\big)$.
Since $\conv\big(\{\BIu_j\}_{j=0:m}\cap F\big)$ is a face of $K$,
we get $K\cap F\in\FF(K)$,
as well as $\FF(K\cap F)\subseteq \FF(K)$.
The proof of statement~\ref{lems:simplicial-complex-and-face-subset} is now complete.

To prove statement~\ref{lems:simplicial-complex-and-face-equality},
let us denote $\widetilde{\FK}_{F}:=\{K\cap F: K\in\FK\}$ for simplicity.
It is clear that every $K\in\FK_{F}$ satisfies 
$K=K\cap F\in\widetilde{\FK}_{F}$.
Conversely, let us fix an arbitrary non-empty $R\in\widetilde{\FK}_{F}$.
Thus, there exists a non-empty $K\in\FK$ with $K\cap F=R$,
and it follows from statement~\ref{lems:simplicial-complex-and-face-subset}
that 
$R\in\FF(K)\subseteq \FK$ and hence $R\in\FK_{F}$.
The proof is now complete.
\endproof

\begin{lemma}\label{lem:simplex-integrals}
    Let $C$ be an $n$-simplex in $\R^d$ with 
    $V(C)=\{\BIv_j\}_{j=0:n}\subset\R^d$, 
    $n\in\{0,1,\ldots,d\}$.
    Then, the following statements hold.
    \begin{enumerate}[label=(\roman*),beginpenalty=10000]
        \item\label{lems:simplex-integrals-affine}
        It holds for any affine function $f:C\to\R$ that 
        $\int_{C}f\DIFFX{\Lebesgue_{C}}=\frac{\volume_n(C)}{n+1}\big(\sum_{j=0}^{n}f(\BIv_{j})\big)$.
        
        \item\label{lems:simplex-integrals-quadratic}
        It holds for any affine functions $f,f':C\to\R$ that 
        \begin{align*}
            \int_{C}ff'\DIFFX{\Lebesgue_{C}}&=\frac{\volume_n(C)}{(n+2)(n+1)}\Bigg[\sum_{j=0}^{n} f(\BIv_{j})\Bigg(f'(\BIv_{j})+\sum_{j'=0}^{n}f'(\BIv_{j'})\Bigg)\Bigg].
        \end{align*}
    \end{enumerate}
\end{lemma}

\proof{Proof of Lemma~\ref{lem:simplex-integrals}.}
    One can directly verify that statements~\ref{lems:simplex-integrals-affine} and \ref{lems:simplex-integrals-quadratic} hold in the case where $n=0$.
    Thus, we will assume throughout the proof that $n\in\N$.
    We let $\BI_n\in\R^{n\times n}$ denote the $n$-by-$n$ identity matrix,
    let $\BIe_{0}:=\veczero_n\in\R^{n}$, 
    let $\BIe_l\in\R^{n}$ denotes the $l$-th standard basis vector in $\R^n$ for $l=1,\ldots,n$,
    and define 
    \begin{align*}
        \BV_{C}:=\left(\begin{tabular}{cccc}
            $|$ & $|$ & & $|$ \\
            $(\BIv_{1}-\BIv_{0})$ & $(\BIv_{2}-\BIv_{0})$ & $\cdots$ & $(\BIv_{n}-\BIv_{0})$ \\
            $|$ & $|$ & & $|$
        \end{tabular}\right)\in\R^{d\times n}.
    \end{align*}
    To begin, 
    note that $\FF(C)$ is a simplicial complex,
    and it holds for any affine function $f:C\to\R$ that 
    $f(\BIx)=\sum_{j=0}^{n}f(\BIv_{j})\lambda_{\FF(C),\BIv_{j}}(\BIx)$ $\forall \BIx\in C$.
    Consequently, due to the linearity of integration,
    it suffices to show that 
    \begin{align}
        \int_{C}\lambda_{\FF(C),\BIv_j}\DIFFX{\Lebesgue_C}&= \frac{\volume_n(C)}{n+1} 
        & \forall 0\le j\le n,\label{eqn:simplex-integrals-proof-affine}\\
        \int_{C}\lambda_{\FF(C),\BIv_{j}}^2\DIFFX{\Lebesgue_{C}}&=\frac{2\volume_n(C)}{(n+2)(n+1)}
        & \forall 0\le j\le n,\label{eqn:simplex-integrals-proof-quadsame}\\
        \int_{C}\lambda_{\FF(C),\BIv_{j}}\lambda_{\FF(C),\BIv_{j'}}\DIFFX{\Lebesgue_{C}}&=\frac{\volume_n(C)}{(n+2)(n+1)}
        & \forall 0\le j< j'\le n.\label{eqn:simplex-integrals-proof-quaddiff}
    \end{align}

    To prove (\ref{eqn:simplex-integrals-proof-affine}),
    let us fix an arbitrary $j\in\{0,1,\ldots,n\}$,
    define 
    $\tilde{\BIv}_l:=(\BIv_{l}^\TRANSP,0)^\TRANSP\in\R^{d+1}$
    $\forall 0\le l\le n$,
    $\tilde{\BIv}_{n+1}:=(\BIv_{j}^\TRANSP,1)^\TRANSP\in\R^{d+1}$,
    and define 
    $C^{\dagger}_j:=\conv\big(\{\tilde{\BIv}_{0},\tilde{\BIv}_{1},\ldots,\tilde{\BIv}_{n+1}\}\big)\subset \R^{d+1}$.
    One checks that $C^{\dagger}_{j}$ is an $(n+1)$-simplex in $\R^{d+1}$,
    that $\aff(C^{\dagger}_{j})=\aff(C)\times \R$,
    and that 
    \begin{align}
        \INDI_{C^{\dagger}_{j}}\big((\BIx^\TRANSP,t)^\TRANSP\big)
        &=\INDI_{C}(\BIx)\INDI_{[0,\lambda_{\FF(C),\BIv_{j}}(\BIx)]}(t) \qquad \forall \BIx\in\aff(C),\; \forall t\in \R.
        \label{eqn:simplex-integrals-proof-affine-decomposition}
    \end{align}
    Next, let $\Lebesgue_{\aff(C)}\otimes \Lebesgue_{\R}$
    denote the product measure of $\Lebesgue_{\aff(C)}$ and $\Lebesgue_{\R}$,
    which coincides with $\Lebesgue_{\aff(C^\dagger_{j})}$.
    We thus get from (\ref{eqn:simplex-integrals-proof-affine-decomposition}) 
    and Fubini's theorem that 
    \begin{align}
        \begin{split}
            \volume_{n+1}(C^{\dagger}_{j})
            &= \Lebesgue_{\aff(C^\dagger_{j})}(C^\dagger_{j})\\
            &=\int_{\aff(C)\times \R}\INDI_{C^{\dagger}_{j}}\big((\BIx^\TRANSP,t)^\TRANSP\big)\DIFFM{\Lebesgue_{\aff(C)}\otimes \Lebesgue_{\R}}{\DIFF\BIx,\DIFF t}\\
            &=\int_{\aff(C)}\INDI_{C}(\BIx)\bigg(\int_{\R}\INDI_{[0,\lambda_{\FF(C),\BIv_j}(\BIx)]}(t)\DIFFM{\Lebesgue_{\R}}{\DIFF t}\bigg)\DIFFM{\Lebesgue_{\aff(C)}}{\DIFF\BIx}\\
            &=\int_{C}\lambda_{\FF(C),\BIv_{j}}\DIFFX{\Lebesgue_{C}}.
        \end{split}
        \label{eqn:simplex-integrals-proof-affine-step1}
    \end{align}
    On the other hand, we also have
    \begin{align}
        \volume_{n+1}(C^\dagger_{j})=\frac{1}{(n+1)!}\det\big(\BV_{C^\dagger_{j}}^\TRANSP\BV_{C^\dagger_{j}}\big)^{\frac{1}{2}},
        \label{eqn:simplex-integrals-proof-affine-step2}
    \end{align}
    where
    \begin{align*}
        \BV_{C^\dagger_{j}}:=\left(\begin{tabular}{cccc}
            $|$ & $|$ & & $|$ \\
            $(\tilde{\BIv}_{1}-\tilde{\BIv}_{0})$ & $(\tilde{\BIv}_{2}-\tilde{\BIv}_{0})$ & $\cdots$ & $(\tilde{\BIv}_{n+1}-\tilde{\BIv}_{0})$ \\
            $|$ & $|$ & & $|$
        \end{tabular}\right)\in\R^{(d+1)\times (n+1)}.
    \end{align*}
    Now, let us define 
    $\BW_{C^{\dagger}_{j}}:=\Big(\begin{smallmatrix}
        \BI_n & -\BIe_j \\ \veczero_n^\TRANSP & 1
    \end{smallmatrix}\Big)\in \R^{(n+1)\times (n+1)}$.
    One checks that
    $\det\big(\BW_{C^{\dagger}_{j}}\big)=1$, and that 
    $\BV_{C^{\dagger}_{j}}\BW_{C^{\dagger}_{j}}=\Big(\begin{smallmatrix}
        \BV_C & \veczero_n \\ \veczero_n^\TRANSP & 1
    \end{smallmatrix}\Big)$,
    which yields 
    $\det\big(\BV_{C^{\dagger}_{j}}^\TRANSP\BV_{C^{\dagger}_{j}}\big)
    =\det\Big(\big(\BV_{C^{\dagger}_{j}}\BW_{C^{\dagger}_{j}}\big)^\TRANSP\big(\BV_{C^{\dagger}_{j}}\BW_{C^{\dagger}_{j}}\big)\Big)
    =\det(\BV_C^\TRANSP\BV_C)$.
    Subsequently, combining 
    this
    with (\ref{eqn:simplex-integrals-proof-affine-step1}) and (\ref{eqn:simplex-integrals-proof-affine-step2}) leads to 
    $\int_{C}\lambda_{\FF(C),\BIv_{j}}\DIFFX{\Lebesgue_{C}}=\frac{1}{(n+1)!}\det(\BV_C^\TRANSP\BV_C)^{\frac{1}{2}}=\frac{\volume_{n}(C)}{n+1}$,
    which proves (\ref{eqn:simplex-integrals-proof-affine}).

    To prove (\ref{eqn:simplex-integrals-proof-quadsame}),
    let us fix an arbitrary $j\in\{0,1,\ldots,n\}$,
    define 
    $\hat{\BIv}_l:=(\BIv_{l}^\TRANSP,0,0)^\TRANSP\in\R^{d+2}$
    $\forall 0\le\nobreak l\le\nobreak n$,
    define 
    $\hat{\BIv}_{n+1}:=(\BIv_{j}^\TRANSP,1,0)^\TRANSP\in\R^{d+2}$,
    $\hat{\BIv}_{n+2}:=(\BIv_{j}^\TRANSP,0,1)^\TRANSP\in\R^{d+2}$,
    and define 
    $C^{\ddagger}_{j,j}:=\conv\big(\{\hat{\BIv}_{0},\hat{\BIv}_{1},\ldots,\hat{\BIv}_{n+2}\}\big)\subset \R^{d+2}$.
    One checks that $C^{\ddagger}_{j,j}$ is an $(n+2)$-simplex in $\R^{d+2}$,
    that $\aff(C^{\ddagger}_{j,j})=\aff(C)\times \R^2$,
    and that 
    \begin{align}
        \INDI_{C^{\ddagger}_{j,j}}\big((\BIx^\TRANSP,t,r)^\TRANSP\big)
        &=\INDI_{C}(\BIx)\INDI_{R_{\lambda_{\FF(C),\BIv_j}(\BIx)}}(t,r) \qquad \forall \BIx\in\aff(C),\; \forall (t,r)^\TRANSP\in \R^2,
        \label{eqn:simplex-integrals-proof-quadsame-decomposition}
    \end{align}
    where $R_{\lambda_{\FF(C),\BIv_j}(\BIx)}:=\big\{(t,r)^\TRANSP\in\R_{+}^2:t+r\le \lambda_{\FF(C),\BIv_j}(\BIx)\big\}$ $\forall \BIx\in C$.
    Subsequently, let $\Lebesgue_{\aff(C)}\otimes \Lebesgue_{\R^2}$
    denote the product measure of $\Lebesgue_{\aff(C)}$ and $\Lebesgue_{\R^2}$,
    which coincides with $\Lebesgue_{\aff(C^\ddagger_{j,j})}$.
    Thus, (\ref{eqn:simplex-integrals-proof-quadsame-decomposition}) and Fubini's theorem imply that 
    \begin{align}
        \begin{split}
            \volume_{n+2}(C^{\ddagger}_{j,j})
            &= \Lebesgue_{\aff(C^{\ddagger}_{j,j})}(C^\ddagger_{j,j})\\
            &=\int_{\aff(C)\times \R^2}\INDI_{C^{\ddagger}_{j,j}}\big((\BIx^\TRANSP,t,r)^\TRANSP\big)\DIFFM{\Lebesgue_{\aff(C)}\otimes \Lebesgue_{\R^2}}{\DIFF\BIx,\DIFF t,\DIFF r}\\
            &=\int_{\aff(C)}\INDI_{C}(\BIx)\bigg(\int_{\R^2}\INDI_{R_{\lambda_{\FF(C),\BIv_j}(\BIx)}}(t,r)\DIFFM{\Lebesgue_{\R^2}}{\DIFF t,\DIFF r}\bigg)\DIFFM{\Lebesgue_{\aff(C)}}{\DIFF\BIx}\\
            &=\frac{1}{2}\int_{C}\lambda_{\FF(C),\BIv_{j}}^2\DIFFX{\Lebesgue_{C}}.
        \end{split}
        \label{eqn:simplex-integrals-proof-quadsame-step1}
    \end{align}
    Similar to the proof of (\ref{eqn:simplex-integrals-proof-affine}), we use the property that
    \begin{align}
        \volume_{n+2}(C^\ddagger_{j,j})=\frac{1}{(n+2)!}\det\big(\BV_{C^\ddagger_{j,j}}^\TRANSP\BV_{C^\ddagger_{j,j}}\big)^{\frac{1}{2}},
        \label{eqn:simplex-integrals-proof-quadsame-step2}
    \end{align}
    where
    \begin{align*}
        \BV_{C^\ddagger_{j,j}}:=\left(\begin{tabular}{cccc}
            $|$ & $|$ & & $|$ \\
            $(\hat{\BIv}_{1}-\hat{\BIv}_{0})$ & $(\hat{\BIv}_{2}-\hat{\BIv}_{0})$ & $\cdots$ & $(\hat{\BIv}_{n+2}-\hat{\BIv}_{0})$ \\
            $|$ & $|$ & & $|$
        \end{tabular}\right)\in\R^{(d+2)\times (n+2)}.
    \end{align*}
    Let us define 
    $\BW_{C^\ddagger_{j,j}}:=\Bigg(\begin{smallmatrix}
        \BI_n & -\BIe_j & -\BIe_j \\ 
        \veczero_n^\TRANSP & 1 & 0 \\
        \veczero_n^\TRANSP & 0 & 1
    \end{smallmatrix}\Bigg)\in \R^{(n+2)\times (n+2)}$,
    which satisfies
    $\det\big(\BW_{C^\ddagger_{j,j}}\big)=1$,
    $\BV_{C^\ddagger_{j,j}}\BW_{C^\ddagger_{j,j}}=\Bigg(\begin{smallmatrix}
        \BV_C & \veczero_n & \veczero_n \\ 
        \veczero_n^\TRANSP & 1 & 0 \\
        \veczero_n^\TRANSP & 0 & 1 \\
    \end{smallmatrix}\Bigg)$.
    Hence,
    $\det\big(\BV_{C^\ddagger_{j,j}}^\TRANSP\BV_{C^\ddagger_{j,j}}\big)
    =\det\Big(\big(\BV_{C^\ddagger_{j,j}}\BW_{C^\ddagger_{j,j}}\big)^\TRANSP\big(\BV_{C^\ddagger_{j,j}}\BW_{C^\ddagger_{j,j}}\big)\Big)
    =\det(\BV_C^\TRANSP\BV_C)$.
    Now, combining 
    this
    with (\ref{eqn:simplex-integrals-proof-quadsame-step1}) and (\ref{eqn:simplex-integrals-proof-quadsame-step2}) leads to 
    $\frac{1}{2}\int_{C}\lambda_{\FF(C),\BIv_{j}}^2\DIFFX{\Lebesgue_{C}}=\frac{1}{(n+2)!}\det(\BV_C^\TRANSP\BV_C)^{\frac{1}{2}}=\frac{\volume_{n}(C)}{(n+2)(n+1)}$,
    which proves (\ref{eqn:simplex-integrals-proof-quadsame}).

    Lastly, to prove (\ref{eqn:simplex-integrals-proof-quaddiff}),
    let us fix arbitrary $j,j'\in\{0,1,\ldots,n\}$ with $j<j'$,
    define 
    $\check{\BIv}_l:=(\BIv_{l}^\TRANSP,0,0)^\TRANSP\in\R^{d+2}$
    $\forall 0\le\nobreak l\le\nobreak n$,
    define 
    $\check{\BIv}_{n+1}:=(\BIv_{j}^\TRANSP,1,0)^\TRANSP\in\R^{d+2}$,
    $\check{\BIv}_{n+2}:=(\BIv_{j'}^\TRANSP,0,1)^\TRANSP\in\R^{d+2}$,
    and define 
    $C^{\ddagger}_{j,j'}:=\conv\big(\{\check{\BIv}_{0},\check{\BIv}_{1},\ldots,\check{\BIv}_{n+2}\}\big)\subset \R^{d+2}$.
    One checks that $C^{\ddagger}_{j,j'}$ is an $(n+2)$-simplex in $\R^{d+2}$,
    that $\aff(C^{\ddagger}_{j,j'})=\aff(C)\times \R^2$,
    and that 
    \begin{align*}
        \INDI_{C^{\ddagger}_{j,j'}}\big((\BIx^\TRANSP,t,r)^\TRANSP\big)
        &=\INDI_{C}(\BIx)\INDI_{[0,\lambda_{\FF(C),\BIv_j}(\BIx)]}(t)\INDI_{[0,\lambda_{\FF(C),\BIv_j'}(\BIx)]}(r) \quad\;\; \forall \BIx\in\aff(C),\; \forall (t,r)^\TRANSP\in \R^2.
    \end{align*}
    Fubini's theorem then implies that 
    \begin{align}
        \volume_{n+2}(C^{\ddagger}_{j,j'})
        &= \Lebesgue_{\aff(C^\ddagger_{j,j'})}(C^\ddagger_{j,j'})\label{eqn:simplex-integrals-proof-quaddiff-step1}\\
        &=\int_{\aff(C)\times \R^2}\INDI_{C^{\ddagger}_{j,j'}}\big((\BIx^\TRANSP,t,r)^\TRANSP\big)\DIFFM{\Lebesgue_{\aff(C)}\otimes \Lebesgue_{\R^2}}{\DIFF\BIx,\DIFF t,\DIFF r}\nonumber\\
        &=\int_{\aff(C)}\INDI_{C}(\BIx)\bigg(\int_{\R^2}\INDI_{[0,\lambda_{\FF(C),\BIv_j}(\BIx)]}(t)\INDI_{[0,\lambda_{\FF(C),\BIv_j'}(\BIx)]}(r)\DIFFM{\Lebesgue_{\R^2}}{\DIFF t,\DIFF r}\bigg)\DIFFM{\Lebesgue_{\aff(C)}}{\DIFF\BIx}\nonumber\\
        &=\int_{C}\lambda_{\FF(C),\BIv_{j}}\lambda_{\FF(C),\BIv_{j'}}\DIFFX{\Lebesgue_{C}}.\nonumber
    \end{align}
    Again, we use the property that
    \begin{align}
        \volume_{n+2}(C^\ddagger_{j,j'})=\frac{1}{(n+2)!}\det\big(\BV_{C^\ddagger_{j,j'}}^\TRANSP\BV_{C^\ddagger_{j,j'}}\big)^{\frac{1}{2}},
        \label{eqn:simplex-integrals-proof-quaddiff-step2}
    \end{align}
    where
    \begin{align*}
        \BV_{C^\ddagger_{j,j'}}:=\left(\begin{tabular}{cccc}
            $|$ & $|$ & & $|$ \\
            $(\check{\BIv}_{1}-\check{\BIv}_{0})$ & $(\check{\BIv}_{2}-\check{\BIv}_{0})$ & $\cdots$ & $(\check{\BIv}_{n+2}-\check{\BIv}_{0})$ \\
            $|$ & $|$ & & $|$
        \end{tabular}\right)\in\R^{(d+2)\times (n+2)}.
    \end{align*}
    We let
    $\BW_{C^\ddagger_{j,j'}}:=\Bigg(\begin{smallmatrix}
        \BI_n & -\BIe_j & -\BIe_{j'} \\ 
        \veczero_n^\TRANSP & 1 & 0 \\
        \veczero_n^\TRANSP & 0 & 1
    \end{smallmatrix}\Bigg)\in \R^{(n+2)\times (n+2)}$,
    which satisfies
    $\det\big(\BW_{C^\ddagger_{j,j'}}\big)=\nobreak1$ as well as
    $\BV_{C^\ddagger_{j,j'}}\BW_{C^\ddagger_{j,j'}}=\Bigg(\begin{smallmatrix}
        \BV_C & \veczero_n & \veczero_n \\ 
        \veczero_n^\TRANSP & 1 & 0 \\
        \veczero_n^\TRANSP & 0 & 1 \\
    \end{smallmatrix}\Bigg)$. 
    Hence,
    $\det\big(\BV_{C^\ddagger_{j,j'}}^\TRANSP\BV_{C^\ddagger_{j,j'}}\big)
    =\det\Big(\big(\BV_{C^\ddagger_{j,j'}}\BW_{C^\ddagger_{j,j'}}\big)^\TRANSP\big(\BV_{C^\ddagger_{j,j'}}\BW_{C^\ddagger_{j,j'}}\big)\Big)
    =\det(\BV_C^\TRANSP\BV_C)$.
    Combining this
    with (\ref{eqn:simplex-integrals-proof-quaddiff-step1}) and (\ref{eqn:simplex-integrals-proof-quaddiff-step2}) leads to 
    $\int_{C}\lambda_{\FF(C),\BIv_{j}}\lambda_{\FF(C),\BIv_{j'}}\DIFFX{\Lebesgue_{C}}=\frac{1}{(n+2)!}\det(\BV_C^\TRANSP\BV_C)^{\frac{1}{2}}=\frac{\volume_{n}(C)}{(n+2)(n+1)}$,
    which proves (\ref{eqn:simplex-integrals-proof-quaddiff}).
    The proof is now complete.
\endproof

\proof{Proof of Lemma~\ref{lem:proc-init-test-funcs}.}
    To prove statement~\ref{lems:proc-init-test-funcs-complexity},
    let us explain in the following list the number of arithmetic operations in some noteworthy lines.
    \begin{itemize}
        \item \textbf{Lines~\ref{alglin:init-test-funcs-search-overlineD}--\ref{alglin:init-test-funcs-search-k}.}
        Since $Q_i= D_i^2$ by a convex maximization argument,
        Line~\ref{alglin:init-test-funcs-search-overlineD} guarantees that 
        $D_i=Q_i^{\frac{1}{2}}\le \overline{D}_i< 2D_i$
        and that Line~\ref{alglin:init-test-funcs-search-overlineD}
        incurs $O\big(\log(D)\big)$ arithmetic operations.
        Similarly, Line~\ref{alglin:init-test-funcs-search-k}
        guarantees that 
        $\frac{D_id_i}{\eta_{\ssUB}}
        \le \frac{\overline{D}_id_i}{\eta_{\ssUB}}
        \le  k_i
        < \frac{2\overline{D}_id_i}{\eta_{\ssUB}}
        < \frac{4D_id_i}{\eta_{\ssUB}}
        \le \frac{4Dd}{\eta_{\ssUB}}$,
        and that Line~\ref{alglin:init-test-funcs-search-k}
        incurs $O\big(\log\big(\frac{Dd}{\eta_{\ssUB}}\big)\big)$ arithmetic operations.

        \item \textbf{Line~\ref{alglin:init-test-funcs-divide}.} 
        Since $k_i< \frac{4D_id_i}{\eta_{\ssUB}}$, 
        it follows from Lemma~\ref{lem:proc-edgewise-divide}\ref{lems:proc-edgewise-divide-complexity} that
        $\EdgewiseDivide(K_{i,j},k_i)$ incurs 
        $O(k_i^{d_i+1}d_i^2)=O\big(\big(\frac{4Dd}{\eta_{\ssUB}}\big)^{d+1}d^2\big)$ arithmetic operations.
        Additionally, 
        Lemma~\ref{lem:proc-edgewise-divide}\ref{lems:proc-edgewise-divide-counts} implies that
        $t_{i,j}< \big(\frac{4D_id_i}{\eta_{\ssUB}}\big)^{d_i}$ 
        $\forall 1\le j\le J_i$, $\forall 0\le i\le N$.
        
        \item \textbf{Line~\ref{alglin:init-test-funcs-search-existing}.}
        Since
        $t_{i,j}< \big(\frac{4D_id_i}{\eta_{\ssUB}}\big)^{d_i}$,
        and 
        $\varrho_{i,\tilde{j}}\le d_i+1$ 
        $\forall 1\le\nobreak l\le\nobreak t_{i,j}$, 
        $\forall 1\le\nobreak j\le \nobreak J_i$,
        $\forall 1\le\nobreak \tilde{j}\le\nobreak \widetilde{J}_i$, 
        $\forall 0\le\nobreak i\le\nobreak N$,
        Lines~\ref{alglin:init-test-funcs-add-point-begin}--\ref{alglin:init-test-funcs-add-point-end} are executed no more than 
        $\big(\frac{4D_id_i}{\eta_{\ssUB}}\big)^{d_i}J_i(d_i+1)$ times,
        and hence 
        $m_i
        \le \big(\frac{4D_id_i}{\eta_{\ssUB}}\big)^{d_i}J_i(d_i+1)$ 
        throughout the entire procedure.
        Consequently, Line~\ref{alglin:init-test-funcs-search-existing} involves comparisons between no more than
        $\big(\frac{4D_id_i}{\eta_{\ssUB}}\big)^{d_i}J_i(d_i+1)$ pairs of vectors in $\R^{d_i}$, 
        which can be carried out with
        $O\big(\big(\frac{4Dd}{\eta_{\ssUB}}\big)^{d}Jd^2\big)$ arithmetic operations.

    \end{itemize}
    Statement~\ref{lems:proc-init-test-funcs-complexity} then
    follows from combining the number of arithmetic operations in each line, as shown in the comments. 

    Next, let us prove statement~\ref{lems:proc-init-test-funcs-subdivision}.
    Note that
    the analyses of Lines~\ref{alglin:init-test-funcs-search-overlineD}--\ref{alglin:init-test-funcs-search-k} above 
    have already shown that 
    $\overline{D}_i\ge D_i$ as well as
    $\overline{D}_i=O(D)$, and
    the analysis of Line~\ref{alglin:init-test-funcs-search-existing} above has already shown that
    $m_i
    \le \big(\frac{4D_id_i}{\eta_{\ssUB}}\big)^{d_i}J_i(d_i+1)
    =O\big(\big(\frac{4Dd}{\eta_{\ssUB}}\big)^{d}Jd\big)$
    as well as 
    $\widetilde{J}_i=\sum_{j=1}^{J_i}t_{i,j}=O\big(\big(\frac{4Dd}{\eta_{\ssUB}}\big)^{d}J\big)$.
    We will now prove the properties of $\FS_i$.
    Let $(\FC_{i,j})_{j=1:J_i}$ denote the 
    simplicial complexes constructed in Line~\ref{alglin:init-test-funcs-divide}.
    It follows from the construction of $\big(\varrho_{i,\tilde{j}},(\varkappa_{i,\tilde{j},\iota})_{\iota=0:\varrho_{i,\tilde{j}}}\big)_{\tilde{j}=1:\widetilde{J}_i}$
    in \mbox{Lines~\ref{alglin:init-test-funcs-add-simplex-begin}--\ref{alglin:init-test-funcs-add-simplex-end}} and the definition of 
    $\FS_i$ that 
    $\{S_{i,\tilde{j}}\}_{\tilde{j}=1:\widetilde{J}_i}=\bigcup_{j=1}^{J_i}\FM(\FC_{i,j})$,
    and thus
    \begin{align}
        \FS_i&=\bigcup_{\tilde{j}=1}^{\widetilde{J}_i}\FF(S_{i,\tilde{j}})=\bigcup_{j=1}^{J_i}\FC_{i,j}.
        \label{eqn:proc-init-test-funcs-proof-relation}
    \end{align}

    Now, to show that $\FS_i$ is a simplicial complex, 
    let us fix arbitrary non-empty $A,A'\in\FS_i$ where 
    $A\in\FC_{i,j}$ and $A'\in\FC_{i,j'}$ for $j,j'\in\{1,\ldots,J_i\}$,
    and let $K_{i,j},K_{i,j'}\in\FM(\FK_i)$ denote the simplices that generated $\FC_{i,j},\FC_{i,j'}$ in Line~\ref{alglin:init-test-funcs-divide}, respectively.
    It thus holds that $A$ is an $n$-simplex in $\R^{d_i}$ for some $n\in\{0,1,\ldots,d_i\}$.
    Moreover, since $\FC_{i,j}$ is a simplicial complex, 
    it holds that 
    $\FF(A)\subseteq \FC_{i,j}\subseteq \FS_{i}$. 
    Furthermore, since $\FK_i$ is a simplicial complex by assumption,
    it holds that     
    $F:=K_{i,j}\cap K_{i,j'}\in\FF(K_{i,j})\cap \FF(K_{i,j'})$.
    In the case where $F=\emptyset$,
    $A\cap A'=\emptyset \in \FF(A)\cap \FF(A')$ holds trivially.
    In the case where $F\neq \emptyset$,
    let $\FD_{i,j}:=\{S\cap F: S\in\FC_{i,j}\}$
    and $\FD_{i,j'}:=\{S\cap F: S\in\FC_{i,j'}\}$.
    Because of Lemma~\ref{lem:simplicial-complex-and-face}\ref{lems:simplicial-complex-and-face-equality} and 
    Lemma~\ref{lem:proc-edgewise-divide}\ref{lems:proc-edgewise-divide-faces},
    both $\FD_{i,j}$ and $\FD_{i,j'}$ are equal to 
    the output of $\EdgewiseDivide(F,k)$,
    which we denote by $\FD$ for simplicity.
    It follows from Lemma~\ref{lem:proc-edgewise-divide}\ref{lems:proc-edgewise-divide-simpcomp} that $\FD$ is a simplicial complex.
    Since $A\cap A'\subseteq F$,
    $A\cap F \in \FD$, and $A'\cap F\in \FD$,
    we get from Lemma~\ref{lem:simplicial-complex-and-face}\ref{lems:simplicial-complex-and-face-subset} that
    $A\cap A'=(A\cap F) \cap (A'\cap F)\in\FF(A\cap F)\cap \FF(A'\cap\nobreak F)\subseteq \FF(A) \cap \FF(A')$.
    The above results have shown that $\FS_i$ is a simplicial complex.

    To complete the proof of statement~\ref{lems:proc-init-test-funcs-subdivision}, 
    observe that 
    (\ref{eqn:proc-init-test-funcs-proof-relation}) implies that 
    $\FS_i$ is a subdivision of $\FK_i$,
    the constructions of 
    $(\BIv_{i,l})_{l=0:m_i}$, 
    $\big(\varrho_{i,\tilde{j}},(\varkappa_{i,\tilde{j},\iota})_{\iota=0:\varrho_{i,\tilde{j}}},\varpi_{i,\tilde{j}}\big)_{\tilde{j}=1:\widetilde{J}_i}$
    in Lines~\mbox{\ref{alglin:init-test-funcs-divide}--\ref{alglin:init-test-funcs-add-simplex-end}}
    yield $V(\FS_i)=\{\BIv_{i,l}\}_{l=0:m_i}$
    as well as 
    $\conv\big(\{\BIv_{i,\varkappa_{i,\tilde{j},\iota}}\}_{\iota=0:\varrho_{i,\tilde{j}}}\big)\subseteq K_{i,\varpi_{i,\tilde{j}}}$ 
    $\forall 1\le \tilde{j}\le \widetilde{J}_i$.
    Lastly, 
    (\ref{eqn:proc-init-test-funcs-proof-relation}), 
    Lemma~\ref{lem:proc-edgewise-divide}\ref{lems:proc-edgewise-divide-meshsize},
    and Line~\ref{alglin:init-test-funcs-search-k} guarantee that
    $\eta(\FS_i)\le \frac{D_id_i}{k_i}\le \eta_{\ssUB}=\overline{\eta}_i$.
    The proof of statement~\ref{lems:proc-init-test-funcs-subdivision} is now complete.

    Lastly, 
    to prove statement~\ref{lems:proc-init-test-funcs-LB},
    let us fix an arbitrary
    $l\in\{0,1,\ldots,m_i\}$, 
    fix an arbitrary 
    $\widehat{S}\in\nobreak \big\{S\in\nobreak \FM(\FS_i):\BIv_{i,l}\in V(S)\big\}$,
    and show that 
    $\int_{\CX_i}g_{i,l}\DIFFX{\mu_i}\ge\frac{1}{d+2}\big(\frac{\eta_{\ssUB}}{4Dd}\big)^{d+1}\min_{K\in\FM(\FK_i)}\big\{\int_{K}p_{i,K}\DIFFX{\Lebesgue_K}\big\}$.
    Let us
    denote $\hat{n}:=\dim(\widehat{S})$,
    and let $(\hat{\BIv}_\iota)_{\iota=0:\hat{n}}$
    be an arbitrary enumeration of $V(\widehat{S})$.
    Then, there exists $j\in\{1,\ldots,J_i\}$ such that 
    $\widehat{S}\in \FC_{i,j}$, 
    where $\FC_{i,j}$ is constructed in Line~\ref{alglin:init-test-funcs-divide}
    from $K_{i,j}$ with $\dim(K_{i,j})=\hat{n}$.
    Let $(\hat{\BIu}_{\iota})_{\iota=0:\hat{n}}$ be an arbitrary enumeration of 
    $V(K_{i,j})$.
    We hence get from Lemma~\ref{lem:proc-edgewise-divide}\ref{lems:proc-edgewise-divide-simpcomp} that 
    \begin{align}
        \volume_{\hat{n}}(\widehat{S})=\frac{\volume_{\hat{n}}(K_{i,j})}{k_i^{\hat{n}}}.
        \label{eqn:proc-init-test-funcs-proof-proof-integral-ineq1}
    \end{align}
    Moreover,
    combining the property that $p_{i,K_{i,j}}:K_{i,j}\to\R_+$ is affine on $K_{i,j}$ and
    Lemma~\ref{lem:proc-edgewise-divide}\ref{lems:proc-edgewise-divide-barycoord} yields
    \begin{align}
        \begin{split}
            \sum_{\iota=0}^{\hat{n}}p_{i,K_{i,j}}(\hat{\BIv}_\iota)&=\sum_{\iota=0}^{\hat{n}}\sum_{\iota'=0}^{\hat{n}}p_{i,K_{i,j}}(\hat{\BIu}_{\iota'})\lambda_{\FF(K_{i,j}),\hat{\BIu}_{\iota'}}(\hat{\BIv}_{\iota})
            \ge \frac{1}{k_i}\Bigg(\sum_{\iota'=0}^{\hat{n}}p_{i,K_{i,j}}(\hat{\BIu}_{\iota'})\Bigg).
        \end{split}
        \label{eqn:proc-init-test-funcs-proof-proof-integral-ineq2}
    \end{align}
    Furthermore,
    let $\Lebesgue_{\widehat{S}}$ denote the $\hat{n}$-dimensional Lebesgue measure on $\widehat{S}$.
    It thus follows from Lemma~\ref{lem:simplex-integrals}\ref{lems:simplex-integrals-affine} with respect to 
    $C\leftarrow K_{i,j}$,
    $(\BIv_{j})_{j=0:n}\leftarrow (\hat{\BIu}_\iota)_{\iota=0:\hat{n}}$,
    $f\leftarrow p_{i,K_{i,j}}$
    that 
    \begin{align}
        \int_{K_{i,j}}p_{i,K_{i,j}}\DIFFX{\Lebesgue_{K_{i,j}}}
        =\frac{\volume_{\hat{n}}(K_{i,j})}{\hat{n}+1}\Bigg(\sum_{\iota=0}^{\hat{n}}p_{i,K_{i,j}}(\hat{\BIu}_{\iota})\Bigg).
        \label{eqn:proc-init-test-funcs-proof-proof-integral-eq}
    \end{align}
    Finally, since 
    $\lambda_{\FS_i,\BIv_{i,l}}$ is affine on $\widehat{S}$ and satisfies 
    $\lambda_{\FS_i,\BIv_{i,l}}(\hat{\BIv}_{\iota})$\,\,$=\INDI_{\{\BIv_{i,l}=\hat{\BIv}_{\iota}\}}$ 
    for $\iota\in\{0,1,\ldots,\hat{n}\}$,
    applying 
    Lemma~\ref{lem:simplex-integrals}\ref{lems:simplex-integrals-quadratic} 
    with respect to $C\leftarrow \widehat{S}$,
    $(\BIv_{j})_{j=0:n}\leftarrow (\hat{\BIv}_{\iota})_{\iota=0:\hat{n}}$,
    $f\leftarrow \lambda_{\FS_i,\BIv_{i,l}}$,
    $f'\leftarrow p_{i,K_{i,j}}$,
    as well as applying
    (\ref{eqn:proc-init-test-funcs-proof-proof-integral-ineq1}),
    (\ref{eqn:proc-init-test-funcs-proof-proof-integral-ineq2}),
    (\ref{eqn:proc-init-test-funcs-proof-proof-integral-eq}),
    and Line~\ref{alglin:init-test-funcs-search-k}
    leads to
    \begin{align*}
        \int_{\CX_i}g_{i,l}\DIFFX{\mu_i}
        &\ge \int_{\CX_i}\lambda_{\FS_i,\BIv_{i,l}}\INDI_{\widehat{S}}\DIFFX{\mu_i}
        \ge \int_{\widehat{S}}\lambda_{\FS_i,\BIv_{i,l}} p_{i,K_{i,j}}\DIFFX{\Lebesgue_{\widehat{S}}}\\
        &= \frac{\volume_{\hat{n}}(\widehat{S})}{(\hat{n}+2)(\hat{n}+1)}\Bigg[\sum_{\iota=0}^{\hat{n}}\lambda_{\FS_i,\BIv_{i,l}}(\hat{\BIv}_{\iota})\Bigg(p_{i,K_{i,j}}(\hat{\BIv}_{\iota})+\sum_{\iota'=0}^{\hat{n}}p_{i,K_{i,j}}(\hat{\BIv}_{\iota'})\Bigg)\Bigg]\\
        &\ge \frac{\volume_{\hat{n}}(\widehat{S})}{(\hat{n}+2)(\hat{n}+1)}\Bigg(\sum_{\iota=0}^{\hat{n}}p_{i,K_{i,j}}(\hat{\BIv}_{\iota})\Bigg)\\
        &\ge \frac{\volume_{\hat{n}}(K_{i,j})}{(\hat{n}+2)(\hat{n}+1)k_i^{\hat{n}+1}}\Bigg(\sum_{\iota=0}^{\hat{n}}p_{i,K_{i,j}}(\hat{\BIu}_{\iota})\Bigg)
        =\frac{1}{(\hat{n}+2)k_i^{\hat{n}+1}}\int_{K_{i,j}}p_{i,K_{i,j}}\DIFFX{\Lebesgue_{K_{i,j}}}\\
        &\ge \frac{1}{d+2}\bigg(\frac{\eta_{\ssUB}}{4Dd}\bigg)^{d+1}\min_{K\in\FM(\FK_i)}\bigg\{\int_{K}p_{i,K}\DIFFX{\Lebesgue_K}\bigg\}.
    \end{align*}
    The proof is now complete.
\endproof

\proof{Proof of Lemma~\ref{lem:proc-init-type-measure}.}
    To prove statement~\ref{lems:proc-init-type-measure-complexity},
    let us
    recall from Lemma~\ref{lem:proc-init-test-funcs}\ref{lems:proc-init-test-funcs-subdivision} that $\max_{0\le i\le N}\{m_i\}=O(\widetilde{J}d)$, and
    explain in the following list the number of arithmetic operations in some noteworthy lines.
    \begin{itemize}
    \item \textbf{Line~\ref{alglin:init-type-measure-search-face} and
        Line~\ref{alglin:init-type-measure-insert-face}.}
        Observe that for $\tilde{j}=1,\ldots,\widetilde{J}_i$
        and for any non-empty $I\subseteq \{0,1,\ldots,\varrho_{i,\tilde{j}}\}$,
        $\conv\big(\{\BIv_{i,\varkappa_{i,\tilde{j},\iota}}:\iota\in I\}\big)$ is a non-empty set in $\FS_i$.
        Moreover, $\mathtt{FaceSet}$ is a data structure which stores non-empty sets in $\FS_i$ that have been enumerated in the algorithm,
        and the check in Line~\ref{alglin:init-type-measure-search-face} 
        ensures that each non-empty set in $\FS_i$ is visited exactly once. 
        Concretely, 
        $\mathtt{FaceSet}$ is
        implemented as a red-black tree \citep[Chapter~13]{ECcormen2009introduction}
        storing integer-valued keys.
        To encode any
        $\{\varkappa_{i,\tilde{j},\iota}\}_{\iota\in I}$ with an integer-valued key, 
        the elements of
        $\{\varkappa_{i,\tilde{j},\iota}\}_{\iota\in I}$
        are first sorted into a list 
        $(\varkappa_{i,\tilde{j},(r)})_{r=0:|I|-1}$ in ascending order, 
        that is, 
        $\{\varkappa_{i,\tilde{j},(r)}\}_{r=0:|I|-1}
        =\{\varkappa_{i,\tilde{j},\iota}\}_{\iota\in I}$
        and 
        $\varkappa_{i,\tilde{j},(r)}\le \varkappa_{i,\tilde{j},(r+1)}$
        $\forall r\in\{0,1,\ldots,|I|-2\}$.
        Then, the value
        $\sum_{r=0}^{|I|-1}(m_i+1)^r(\varkappa_{i,\tilde{j},(r)}+1)$ is computed and used as the unique key encoding 
        $\{\varkappa_{i,\tilde{j},\iota}\}_{\iota\in I}$.
        Note that 
        the computation of the key
        takes $O(d\log(d))$ arithmetic operations.
        Subsequently, 
        since $|\FS_i|
        =O(2^d\widetilde{J})$ by Lemma~\ref{lem:proc-init-test-funcs}\ref{lems:proc-init-test-funcs-subdivision},
        the properties of the red-black tree ensure
        that the number of arithmetic operations when
        inserting a new key 
        or when searching for an existing key in $\mathtt{FaceSet}$
        is logarithmic in the total number of keys.
        Hence,
        Line~\ref{alglin:init-type-measure-search-face} and
        Line~\ref{alglin:init-type-measure-insert-face} each incurs
        $O\big(d\log(d)+\log(|\FS_i|)\big)=O\big(d\log(d)+\log(\widetilde{J})\big)$ arithmetic operations.

        \item \textbf{Line~\ref{alglin:init-type-measure-enumerate-face}.} Observe that 
        the number of arithmetic operations in Line~\ref{alglin:init-type-measure-compute-integral} is no fewer than the number of arithmetic operations in Line~\ref{alglin:init-type-measure-encode-face},
        since the procedure $\mathtt{Int}$ needs to process its input $F_{i,\tilde{r}}$ encoded in Line~\ref{alglin:init-type-measure-encode-face}.
        Consequently,
        the total number of arithmetic operations in the for-loop in 
        Line~\ref{alglin:init-type-measure-enumerate-face}
        is $O\big(2^d(T_{\mathtt{Int}}+\widetilde{J}d)\big)$.
    
        \item \textbf{Lines~\ref{alglin:init-type-measure-probmat}--\ref{alglin:init-type-measure-aggregate}.}
        Since the analyses of Line~\ref{alglin:init-type-measure-search-face} and
        Line~\ref{alglin:init-type-measure-insert-face} above 
        have established that $\widetilde{R}_i=|\FS_i|-1=O(2^d\widetilde{J})$,
        Line~\ref{alglin:init-type-measure-probmat} and Line~\ref{alglin:init-type-measure-aggregate} both incur
        $O(\widetilde{R}_im_i)=O(2^d\widetilde{J}^2d)$ arithmetic operations.
        
    \end{itemize}
    Statement~\ref{lems:proc-init-type-measure-complexity} then
    follows from combining the number of arithmetic operations in each line, as shown in the comments. 

    Statement~\ref{lems:proc-init-type-measure-auxiliary} follows from 
    the observation that $V(\FS_i)\subseteq \FS_i$,
    the above analyses of Line~\ref{alglin:init-type-measure-search-face} and
    Line~\ref{alglin:init-type-measure-insert-face},
    the encoding of $\CG_i$,
    and the definitions of $F_{i,\tilde{r}}$ and $\tilde{\varpi}_{i,\tilde{r}}$
    in Lines~\ref{alglin:init-type-measure-encode-face}--\ref{alglin:init-type-measure-simplex-index}.

    To prove statement~\ref{lems:proc-init-type-measure-integrals},
    notice that 
    for any $S\in \FS_i$,
    $\lambda_{\FS_i,\BIv}(\BIx)=\nobreak0$ $\forall \BIx\in \relint(S)$ holds whenever $\BIv\notin\nobreak V(S)$.
    This implies that
    $\int_{\CX_i}\lambda_{\FS_i,\BIv}\INDI_{\relint(S)}\DIFFX{\mu_i}=\nobreak0$ $\forall \BIv\in V(\FS_i)\setminus V(S)$,
    $\forall S\in\nobreak\FS_i$.
    It hence follows from Line~\ref{alglin:init-type-measure-compute-integral} 
    and the assumption about the oracle procedure $\mathtt{Int}$ in \ref{enc:typemeasures}
    that 
    $p_{l,\tilde{r}}=\int_{\CX_i}\lambda_{\FF(F_{i,\tilde{r}}),\BIv_{i,l}}$~\!\!\!$\INDI_{\relint(F_{i,\tilde{r}})}\DIFFX{\mu_i}=$~\!\!$\int_{\CX_i}\lambda_{\FS_i,\BIv_{i,l}}$~\!\!\!$\INDI_{\relint(F_{i,\tilde{r}})}\DIFFX{\mu_i}$ 
    $\forall l\in\{0,1,\ldots,m_i\}$, 
    $\forall \tilde{r}\in\nobreak\{1,\ldots,\widetilde{R}_i\}$
    when the for-loop in Line~\ref{alglin:init-type-measure-enumerate-face} terminates.
    Subsequently, since $\{\relint(S):S\in\FS_i\}$ is a disjoint partition of $\CX_i$ (see, e.g., \citep[p.10]{ECmunkres1984elements}),
    Line~\ref{alglin:init-type-measure-probmat} and statement~\ref{lems:proc-init-type-measure-auxiliary} show that
    $\sum_{\tilde{r}=1}^{\widetilde{R}_i}p_{l,\tilde{r}}=
    \int_{\CX_i}\lambda_{\FS_i,\BIv_{i,l}}\big(\sum_{\tilde{r}=1}^{\widetilde{R}_i}$~\!$\INDI_{\relint(F_{i,\tilde{r}})}\big)\DIFFX{\mu_i}=\int_{\CX_i}g_{i,l}\DIFFX{\mu_i}$ 
    $\forall 0\le l\le m_i$.
    Line~\ref{alglin:init-type-measure-aggregate} then yields 
    $\bar{\BIg}_i=\nobreak\big(\int_{\CX_i}g_{i,1}\DIFFX{\mu_i},\ldots,\int_{\CX_i}g_{i,m_i}\DIFFX{\mu_i}\big)^\TRANSP$
    as well as 
    $\BP_i\vecone_{\widetilde{R}_i}=\big(\int_{\CX_i}g_{i,0}\DIFFX{\mu_i},\allowbreak\int_{\CX_i}g_{i,1}\DIFFX{\mu_i},\ldots,\allowbreak\int_{\CX_i}g_{i,m_i}\DIFFX{\mu_i}\big)^\TRANSP$,
    and proves statement~\ref{lems:proc-init-type-measure-integrals}.

    Lastly, let us prove statement~\ref{lems:proc-init-type-measure-error}.
    For $\tilde{r}=1,\ldots,\widetilde{R}_i$,
    one observes that
    $\PROB[S=\tilde{r}]
    =\sum_{l=0}^{m_i} p_{l,\tilde{r}}
    =\int_{\CX_i}\!\big(\sum_{l=0}^{m_i}\lambda_{\FS_i,\BIv_{i,l}}\big)$~\!\!$\INDI_{\relint(F_{i,\tilde{r}})}\DIFFX{\mu_i}
    =\mu_i\big(\relint(F_{i,\tilde{r}})\big)$.
    Subsequently, for any $E\in\CB(\CX_i)$,
    it holds by the law of total probability that 
    \begin{align*}
        \PROB[\overbar{\BIX}_i\in E]
        &=\sum_{\tilde{r}\in\{\tilde{r}':\PROB[S=\tilde{r}']>0\}}\PROB[\overbar{\BIX}_i\in E|S=\tilde{r}]\PROB[S=\tilde{r}]
        =\sum_{\tilde{r}=1}^{\widetilde{R}_i}\mu_i\big(E\cap \relint(F_{i,\tilde{r}})\big)=\mu_i(E).
    \end{align*}
    Thus, the law of the random variable $\overbar{\BIX}_i$ is equal to $\mu_i$.
    Moreover, Lines~\ref{alglin:init-type-measure-initzero}--\ref{alglin:init-type-measure-compute-integral} 
    imply that 
    for all $l\in\{0,1,\ldots,m_i\}$, $\tilde{r}\in\{1,\ldots,\widetilde{R}_i\}$,
    $\PROB[Q=\nobreak l,\;S=\nobreak \tilde{r}]=p_{l,\tilde{r}}=\nobreak0$
    if $\BIv_{i,l}\notin V(F_{i,\tilde{r}})$.
    Consequently, since 
    it holds that 
    \begin{align*}
        \PROB[\BIX_i=\BIv_{i,l},\;\overbar{\BIX}_i\in F_{i,\tilde{r}}|Q=\nobreak l,\;S=\nobreak \tilde{r}]=\nobreak 1 
        \qquad \forall (l,\tilde{r})\in\big\{(l',\tilde{r}'):\PROB[Q=l',\; S=\tilde{r}']>0\big\},
    \end{align*}
    we get by the law of total expectation that 
    \begin{align*}
        \EXP\big[\|\BIX_i-\overbar{\BIX}_i\|_2\big]
        &=\sum_{(l,\tilde{r})\in\{(l',\tilde{r}'):\PROB[Q=l',\; S=\tilde{r}']>0\}}\EXP\big[\|\BIX_i-\overbar{\BIX}_i\|_2\big| Q=l,\; S=\tilde{r}\big]\PROB[Q=l,\; S=\tilde{r}]
        \le \eta(\FS_i).
    \end{align*}
    The proof is now complete.
\endproof

\subsection{Proofs of results in Section~\ref{sapx:complexity-lsip}}\label{sapx:proof-complexity-lsip}

Let us introduce the following notations that are used throughout this subsection.
We define
\begin{align*}
    c_{\ssmin}&:= \min_{1\le i\le N,\,0\le j\le m_i,\,0\le l\le m_0}\big\{c_i(\BIv_{i,j},\BIv_{0,l})\big\},\\
    c_{\ssmax}&:= \max_{1\le i\le N,\,0\le j\le m_i,\,0\le l\le m_0}\big\{c_i(\BIv_{i,j},\BIv_{0,l})\big\},\\
    \bar{g}_{\ssmin}&:=\min_{1\le i\le N,\,0\le j\le m_i}\bigg\{\int_{\CX_i}g_{i,j}\DIFFX{\mu_i}\bigg\},\\
    q&:=\sum_{i=1}^{N}(1+m_i+m_0),\\
    \bar{\BIg}&:=(1,\bar{\BIg}_1^\TRANSP,\veczero_{m_0}^\TRANSP,\ldots,1,\bar{\BIg}_N^\TRANSP,\veczero_{m_0}^\TRANSP)^\TRANSP\in\R^{q},\\
    \CV&:=\big\{(i,j,0):j\in\{0,1,\ldots,m_i\},\; i\in\{1,\ldots,N\}\big\} \cup \big\{(i,0,l): l\in\{1,\ldots,m_0\},\; i\in\{1,\ldots,N\}\big\}, \\
    \overline{\CI}&:=\big\{(i,j,l):j\in\{0,1,\ldots,m_i\},\; l\in\{0,1,\ldots,m_0\},\; i\in\{1,\ldots,N\}\big\}.
\end{align*}
For any $\CI\subseteq \overline{\CI}$,
let \LPTAG{$\CI$} denote the following LP problem:
\begin{align}
    \begin{split}
        \maximize_{(y_{i,0},\BIy_i,\BIw_i)}\quad & \sum_{i=1}^{N} y_{i,0} + \langle\bar{\BIg}_i,\BIy_i\rangle\\
        \mathrm{subject~to}\quad & y_{i,0}+\langle\BIg_i(\BIv_{i,j}),\BIy_i\rangle + \langle\BIg_0(\BIv_{0,l}),\BIw_i\rangle \le c_i(\BIv_{i,j},\BIv_{0,l})-c_{\ssmax} \quad \forall (i,j,l)\in\CI,\\
        & \sum_{i=1}^{N}\BIw_i\ge\veczero_{m_0},\qquad y_{i,0}\in\R, \; \BIy_i\in\R^{m_i},\; \BIw_i\in\R^{m_0} \quad \forall 1\le i\le N.
    \end{split}
    \tag{$\LP\normalfont{\text{-}}\CI$}
    \label{eqn:mt-tf-lp-any-index}
\end{align}
Observe that \LPTAG{$\CI$} is always feasible.
Moreover, 
for any $\CI\subseteq \overline{\CI}$,
let $\tau^{\star}_{\CI}\in\nobreak\R\cup\{\infty\}$ denote the optimal value of 
\LPTAG{$\CI$},
and define the set $S_{\CI,\tau}\subset \R^{q}$
for any $\tau\in\R\cup\{-\infty\}$ as follows:
\begin{align*}
    S_{\CI,\tau}&:=\left\{(y_{1,0},\BIy_1^\TRANSP,\BIw_1^\TRANSP,\ldots,y_{N,0},\BIy_N^\TRANSP,\BIw_N^\TRANSP)^\TRANSP \in\R^{q} : 
    \begin{tabular}{l}
        $(y_{i,0},\BIy_i,\BIw_i)_{i=1:N}$ is feasible for \LPTAG{$\CI$}\\
        and satisfies $\sum_{i=1}^Ny_{i,0}+\langle\bar{\BIg}_i,\BIy_i\rangle\ge \tau$
    \end{tabular}\right\}.
\end{align*}
In particular, 
$S_{\CI,-\infty}$ corresponds to the feasible set of \LPTAG{$\CI$},
and for every $\tau\in\R$,
$S_{\CI,\tau}$ is the $\tau$-superlevel set of \LPTAG{$\CI$}.
Furthermore, let us define 
$B_{\infty}(R):=\big\{\Bchi\in\R^{q}:\|\Bchi\|_{\infty}\le\nobreak R\big\}$
for any $R>0$.

The proof of Theorem~\ref{thm:proc-solve-LSIP} relies on the following properties of \LPTAG{$\CI$} which are presented as a lemma.
\begin{lemma}\label{lem:lp-any-index}
    Under Setting~\ref{sett:Euclidean}, 
    the following statements hold.
    \begin{enumerate}[label=(\roman*),beginpenalty=10000]
        \item\label{lems:lp-any-index-minint}
        It holds that $\bar{g}_{\ssmin}>0$.

        \item\label{lems:lp-any-index-Linfinity-bound}
        It holds for any $\Delta\ge 0$ 
        and for $M_{\Delta}:=(1+\bar{g}_{\ssmin}^{-1})N^2(c_{\ssmax}-c_{\ssmin}+\nobreak\Delta)$
        that $S_{\CV,\tau^{\star}_{\overline{\CI}}-\Delta}\subseteq B_{\infty}(M_{\Delta})$.
        Consequently,
        for any $\CI$ satisfying $\CV\subseteq\CI\subseteq\overline{\CI}$,
        for any $\tau\ge \tau^{\star}_{\overline{\CI}}-1$, 
        and for $M:=\nobreak(1+\nobreak\bar{g}_{\ssmin}^{-1})N^2(c_{\ssmax}-c_{\ssmin}+1)$,
        it holds that 
        $S_{\CI,\tau}\subseteq B_{\infty}(M)$ and hence 
        $\tau^\star_{\CI}<\infty$.

        \item\label{lems:lp-any-index-Linfinity-ball}
        For any $\CI\subseteq\overline{\CI}$ satisfying 
        $\tau^\star_{\CI}<\infty$ and for any $\epsilon>0$,
        $S_{\CI,\tau^\star_{\CI}-\epsilon}$ contains a closed
        $\|\cdot\|_{\infty}$-ball with radius~$\frac{\epsilon}{6N}$.

        \item\label{lems:lp-any-index-transform}
        Let
        $(y_{i,0},\BIy_i,\BIw_i)_{i=1:N}$
        be feasible for \LPTAG{$\overline{\CI}$},
        let
        $\hat{y}_{i,0}:=y_{i,0}-L_c\big(\eta(\FS_i)+\eta(\FS_0)\big)+c_{\ssmax}$,
        $\hat{\BIy}_i:=\BIy_i$ 
        for $i=1,\ldots,N$,
        and let 
        $\hat{\BIw}_1:=-\sum_{i=2}^{N}\BIw_i$,
        $\hat{\BIw}_i:=\BIw_i$ 
        for $i=2,\ldots,N$.
        Then, $(\hat{y}_{i,0},\hat{\BIy}_i,\hat{\BIw}_i)_{i=1:N}$
        is feasible for \eqref{eqn:mt-tf-lsip}.

    \end{enumerate}
\end{lemma}

\proof{Proof of Lemma~\ref{lem:lp-any-index}.}
    Let us first prove statement~\ref{lems:lp-any-index-minint}.
    Let $i\in\{1,\ldots,N\}$, 
    $j\in\{0,1,\ldots,m_i\}$ be arbitrary,
    and let $K\in\FM(\FK_i)$ satisfy $\BIv_{i,j}\in K$.
    Moreover, let $p_{i,K}:K\to\R_+$ be the affine function in Setting~\ref{sett:Euclidean}.
    Thus, $g_{i,j}=\lambda_{\FS_i,\BIv_{i,j}}$ and 
    $p_{i,K}$ are both continuous on $K$,
    where $\lambda_{\FS_i,\BIv_{i,j}}$ is not equal to the zero function on $\relint(K)$,
    and $p_{i,K}$ is positive on $\relint(K)$.
    Subsequently, it holds by the property of $\mu_i$ 
    in Setting~\ref{sett:Euclidean} that 
    $\int_{\CX_i}g_{i,j}\DIFFX{\mu_i}
    \ge \int_{K}\lambda_{\FS_i,\BIv_{i,j}}p_{i,K}\DIFFX{\Lebesgue_K}>\nobreak0$.
    This proves statement~\ref{lems:lp-any-index-minint}.

    Now, let us fix an arbitrary $\Delta\ge 0$ and prove statement~\ref{lems:lp-any-index-Linfinity-bound}.
    To begin, 
    one checks that the dual LP problem of \LPTAG{$\CV$} is feasible,
    which shows that $\tau^\star_{\overline{\CI}}\le \tau^\star_{\CV}<\nobreak\infty$.
    Thus,
    applying duality arguments to \LPTAG{$\CV$} and \mbox{\LPTAG{$\overline{\CI}$}}
    yields $\tau^\star_{\CV}\le\nobreak 0$
    and $\tau^\star_{\overline{\CI}}\ge\nobreak N(c_{\ssmin}-c_{\ssmax})$.
    Let us fix an arbitrary 
    $\Bchi=\nobreak(y_{1,0},y_{1,1},\ldots,y_{1,m_1},\allowbreak w_{1,1},\ldots,w_{1,m_0},\ldots, y_{N,0},\allowbreak y_{N,1},\ldots, y_{N,m_N},\allowbreak w_{N,1},\ldots,w_{N,m_0})^\TRANSP\in S_{\CV,\tau^\star_{\overline{\CI}}-\Delta}$.
    Moreover, let us denote
    $\bar{g}_{i,j}:=\nobreak\int_{\CX_i}g_{i,j}\DIFFX{\mu_i}$,
    $\tilde{c}_{i,j,l}:=\nobreak c_i(\BIv_{i,j},\BIv_{0,l})$
    $\forall 0\le\nobreak j\le\nobreak m_i$,
    $\forall 0\le\nobreak l\le\nobreak m_0$,
    $\forall 1\le i\le N$
    for notational simplicity.
    Notice that
    $\bar{g}_{i,j}\ge \bar{g}_{\ssmin}$ 
    $\forall 0\le\nobreak j\le\nobreak m_i$,
    $\forall 1\le i\le N$,
    and that 
    $c_{\ssmin}-c_{\ssmax}\le \tilde{c}_{i,j,l}\le 0$
    $\forall 0\le j\le m_i$,
    $\forall 0\le\nobreak l\le\nobreak m_0$,
    $\forall 1\le i\le N$.
    It thus holds by the definitions of 
    $(\bar{\BIg}_i)_{i=1:N}$, $(\BIg_i(\,\cdot\,))_{i=0:N}$,
    and $\CV$
    that
    \begin{align}
        y_{i,0}&\le \tilde{c}_{i,0,0}  & \forall 1\le i\le N,\label{eqn:lp-any-index-proof-ineq-yi0}\\
        y_{i,0}+y_{i,j}&\le \tilde{c}_{i,j,0} & \forall 1\le j\le m_i,\; \forall 1\le i\le N,\label{eqn:lp-any-index-proof-ineq-yij}\\
        y_{i,0}+w_{i,l}&\le \tilde{c}_{i,0,l} & \forall 1\le l\le m_0,\; \forall 1 \le i\le N,\label{eqn:lp-any-index-proof-ineq-wil}\\
        \sum_{i=1}^{N}w_{i,l}&\ge0 & \forall 1\le l\le m_0,\label{eqn:lp-any-index-proof-eq-wsum}\\
        \sum_{i=1}^{N}\Bigg(y_{i,0}+\sum_{j=1}^{m_i}\bar{g}_{i,j}y_{i,j}\Bigg)&\ge \tau^\star_{\overline{\CI}}-\Delta\ge N(c_{\ssmin}-c_{\ssmax}-\Delta).\label{eqn:lp-any-index-proof-ineq-obj}
    \end{align}
    We already get from (\ref{eqn:lp-any-index-proof-ineq-yi0}) that $y_{i,0}\le 0$ $\forall 1\le i\le N$.
    Let us establish the bound $\|\Bchi\|_{\infty}\le\nobreak M_{\Delta}$ through the following 4 steps.
    \begin{itemize}
        \item Step~1: showing that 
        $y_{i,j}\ge -(1+\bar{g}_{\ssmin}^{-1})N(c_{\ssmax}-c_{\ssmin}+\Delta)$ 
        $\forall 1\le j\le m_i$, $\forall 1\le i\le N$.

        \item Step~2: showing that 
        $y_{i,0}\ge -(1+\bar{g}_{\ssmin}^{-1})N(c_{\ssmax}-c_{\ssmin}+\Delta)$ 
        $\forall 1\le i\le N$.
        
        \item Step~3: showing that 
        $y_{i,j}\le (1+\bar{g}_{\ssmin}^{-1})N(c_{\ssmax}-c_{\ssmin}+\Delta)$
        $\forall 1\le j\le m_i$, $\forall 1\le i\le N$.
        
        \item Step~4: showing that 
        $-(1+\bar{g}_{\ssmin}^{-1})N^2(c_{\ssmax}-c_{\ssmin}+\Delta)\le w_{i,l}\le (1+\bar{g}_{\ssmin}^{-1})N(c_{\ssmax}-c_{\ssmin}+\Delta)$
        $\forall 1\le l\le m_0$, $\forall 1\le i\le N$.
        
    \end{itemize}

    \underline{Step~1.}
    By the symmetry among $i=1,\ldots,N$ and the symmetry among $j=1,\ldots,m_i$, 
    let us suppose for the sake of contradiction that 
    $y_{1,1}<-(1+\bar{g}_{\ssmin}^{-1})N(c_{\ssmax}-c_{\ssmin}+\Delta)$,
    and let us define
    $\hat{y}_{1,1}:=\nobreak c_{\ssmin}-c_{\ssmax}$.
    Since 
    $y_{1,0}+\hat{y}_{1,1}\le c_{\ssmin}-c_{\ssmax}\le \tilde{c}_{1,1,0}$,
    one checks that 
    $\hat{\Bchi}:=(y_{1,0},\hat{y}_{1,1},y_{1,2},\ldots,y_{1,m_1},\allowbreak w_{1,1},\ldots,w_{1,m_0},\ldots, y_{N,0}, \allowbreak y_{N,1},\ldots,y_{N,m_N},\allowbreak w_{N,1},\ldots,w_{N,m_0})^\TRANSP\in S_{\CV,-\infty}$,
    and (\ref{eqn:lp-any-index-proof-ineq-obj}) implies that the objective value of $\hat{\Bchi}$ is given by
    \begin{align*}
        \bar{g}_{1,1}(\hat{y}_{1,1}-y_{1,1}) + \sum_{i=1}^{N}\Bigg(y_{i,0}+\sum_{j=1}^{m_i}\bar{g}_{i,j}y_{i,j}\Bigg)
        &> \bar{g}_{1,1}\big[c_{\ssmin} - c_{\ssmax} + (1+\bar{g}_{\ssmin}^{-1})N(c_{\ssmax}-c_{\ssmin} + \Delta)\big] \\
        &\qquad + N(c_{\ssmin}-c_{\ssmax}-\Delta)\\
        &> 0,
    \end{align*}
    which contradicts $\tau^\star_{\CV}\le 0$.
    Therefore, it holds that 
    $y_{i,j}\ge -(1+\bar{g}_{\ssmin}^{-1})N(c_{\ssmax}-c_{\ssmin}+\Delta)$ 
    $\forall 1\le j\le m_i$, $\forall 1\le i\le N$,
    and Step~1 is now complete.

    \underline{Step~2.}
    By the symmetry among $i=1,\ldots,N$,
    let us suppose for the sake of contradiction that 
    $y_{1,0}<-(1+\bar{g}_{\ssmin}^{-1})N(c_{\ssmax}-c_{\ssmin}+\Delta)$,
    and let us define 
    $\hat{y}_{1,0}:=N(c_{\ssmin}-c_{\ssmax})$,
    $\hat{y}_{1,j}:=y_{1,0}+y_{1,j}-\hat{y}_{1,0}$
    $\forall 1\le\nobreak j\le\nobreak m_1$,
    $\hat{w}_{1,l}:=(N-1)(c_{\ssmax}-c_{\ssmin})$ $\forall 1\le\nobreak l\le\nobreak m_0$,
    $\hat{y}_{i,0}:=y_{i,0}$ $\forall 2\le\nobreak i\le\nobreak N$,
    $\hat{y}_{i,j}:=y_{i,j}$ $\forall 1\le\nobreak j\le\nobreak m_i$, $\forall 2\le\nobreak i\le\nobreak N$,
    $\hat{w}_{i,l}:=c_{\ssmin}-c_{\ssmax}$ $\forall 1\le\nobreak l\le\nobreak m_0$, $\forall 2\le\nobreak i\le\nobreak N$.
    One checks from (\ref{eqn:lp-any-index-proof-ineq-yij}) that 
    \begin{align*}
        \hat{y}_{1,0}&\le c_{\ssmin} - c_{\ssmax} \le \tilde{c}_{1,0,0}, \\
        \hat{y}_{1,0} + \hat{y}_{1,j}&=y_{1,0} + y_{1,j}\le \tilde{c}_{1,j,0} & \forall 1\le j\le m_1,\\
        \hat{y}_{1,0} + \hat{w}_{1,l}&= c_{\ssmin} - c_{\ssmax} \le \tilde{c}_{1,0,l} & \forall 1\le l \le m_0, \\
        \hat{y}_{i,0}+\hat{w}_{i,l}&= y_{i,0}+c_{\ssmin} - c_{\ssmax}\le c_{\ssmin} - c_{\ssmax}\le \tilde{c}_{i,0,l} & \forall 1\le l\le m_0,\; \forall 2\le i\le N, \\
        \sum_{i=1}^{N}\hat{w}_{i,l}&= (N-1)(c_{\ssmax}-c_{\ssmin}) + (N-1)(c_{\ssmin}-c_{\ssmax})\ge 0 & \forall 1\le l\le m_0.
    \end{align*}
    Hence,
    $\hat{\Bchi}:=(\hat{y}_{1,0},\hat{y}_{1,1},\ldots,\hat{y}_{1,m_1},\allowbreak \hat{w}_{1,1},\ldots,\hat{w}_{1,m_0},\ldots,\allowbreak \hat{y}_{N,0},\hat{y}_{N,1},\ldots,\hat{y}_{N,m_N},\allowbreak \hat{w}_{N,1},\ldots,\hat{w}_{N,m_0})^\TRANSP\in S_{\CV,-\infty}$.
    Moreover, using (\ref{eqn:lp-any-index-proof-ineq-obj})
    and the identity $\sum_{j=1}^{m_1}\bar{g}_{1,j}=1-\bar{g}_{1,0}$,
    we derive the objective value of $\hat{\Bchi}$ to be 
    \begin{align*}
        \sum_{i=1}^{N}\Bigg(\hat{y}_{i,0}+\sum_{j=1}^{m_i}\bar{g}_{i,j}\hat{y}_{i,j}\Bigg)
        &=\bar{g}_{1,0}(\hat{y}_{1,0}-y_{1,0}) + \sum_{i=1}^{N}\Bigg(y_{i,0}+\sum_{j=1}^{m_i}\bar{g}_{i,j}y_{i,j}\Bigg)\\
        &> \bar{g}_{1,0}\big[N(c_{\ssmin} - c_{\ssmax}) 
        + (1+\bar{g}_{\ssmin}^{-1})N(c_{\ssmax}-c_{\ssmin}+\Delta)\big] + N(c_{\ssmin}-c_{\ssmax}-\Delta)\\
        &> 0,
    \end{align*}
    which contradicts $\tau^\star_{\CV}\le 0$.
    Therefore, it holds that $y_{i,0}\ge -(1+\bar{g}_{\ssmin}^{-1})N(c_{\ssmax}-c_{\ssmin}+\Delta)$ $\forall 1\le i\le N$,
    and Step~2 is now complete.

    \underline{Step~3.}
    Now, combining the conclusion of Step~2 and 
    (\ref{eqn:lp-any-index-proof-ineq-yij}) leads to 
    \begin{align*}
        y_{i,j}&\le \tilde{c}_{i,j,0}-y_{i,0} \le (1+\bar{g}_{\ssmin}^{-1})N(c_{\ssmax}-c_{\ssmin}+\Delta) &  \forall 1\le j\le m_i,\; \forall 1\le i\le N.
    \end{align*}
    This completes Step~3.

    \underline{Step~4.}
    Similar to Step~3,
    combining the conclusion of Step~2 and 
    (\ref{eqn:lp-any-index-proof-ineq-wil}) leads to 
    \begin{align*}
        w_{i,l}&\le \tilde{c}_{i,0,l}-y_{i,0} \le (1+\bar{g}_{\ssmin}^{-1})N(c_{\ssmax}-c_{\ssmin}+\Delta) &  \forall 1\le l\le m_0,\; \forall 1\le i\le N.
    \end{align*}
    Subsequently, (\ref{eqn:lp-any-index-proof-eq-wsum}) implies that 
    \begin{align*}
        w_{i,l}&\ge-\Bigg(\sum_{i'\in\{1,\ldots,N\}\setminus\{i\}}w_{i',l}\Bigg)\ge -(1+\bar{g}_{\ssmin}^{-1})N^2(c_{\ssmax}-c_{\ssmin}+\Delta) &  \forall 1\le l\le m_0,\; \forall 1\le i\le N.
    \end{align*}
    The proof of statement~\ref{lems:lp-any-index-Linfinity-bound} is now complete.

    To prove statement~\ref{lems:lp-any-index-Linfinity-ball},
    let us fix an arbitrary $\CI\subseteq \overline{\CI}$ with 
    $\tau^{\star}_{\CI}<\infty$, 
    as well as an arbitrary $\epsilon>0$.
    Since \LPTAG{$\CI$} is a feasible LP problem with finite optimal value,
    it admits an optimizer.
    Hence, let us fix an arbitrary optimizer  
    $\Bchi^\star=(y^\star_{1,0},\BIy_1^{\star\TRANSP},\BIw_1^{\star\TRANSP},\ldots,\allowbreak y^\star_{N,0},\BIy_N^{\star\TRANSP},\BIw_N^{\star\TRANSP})^\TRANSP\in\R^{q}$
    of \LPTAG{$\CI$},
    where $y^\star_{i,0}\in\nobreak\R$, 
    $\BIy^\star_{i}\in\nobreak\R^{m_i}$,
    $\BIw^\star_{i}\in\nobreak\R^{m_0}$ $\forall 1\le i\le N$.
    Next, let us define 
    $\hat{y}_{i,0}:={y^\star_{i,0}-\frac{2\epsilon}{3N}}$, 
    $\hat{\BIy}_i:=\BIy_i^\star$, 
    $\hat{\BIw}_i:=\nobreak\BIw_i^\star+\frac{\epsilon}{6N}\vecone_{m_0}$
    $\forall 1\le i\le N$,
    and define $\hat{\Bchi}:=(\hat{y}_{1,0},\hat{\BIy}_{1}^\TRANSP,\hat{\BIw}_1^\TRANSP,\ldots,\hat{y}_{N,0},\hat{\BIy}_N^\TRANSP,\hat{\BIw}_N^\TRANSP)^\TRANSP\in\R^{q}$.
    In the following, we will show that 
    $S_{\CI,\tau^\star_{\CI}-\epsilon}$ contains the closed $\|\cdot\|_{\infty}$-ball with radius $\frac{\epsilon}{6N}$ centered at~$\hat{\Bchi}$.
    To that end, let us fix an arbitrary 
    $\Brho=(u_{1,0},\BIu_1^\TRANSP,\BIr_1^\TRANSP,\ldots,u_{N,0},\BIu_N^\TRANSP,\BIr_N^\TRANSP)^\TRANSP$
    with $\|\Brho\|_{\infty}\le 1$, 
    where $u_{i,0}\in\R$, $\BIu_{i}\in\R^{m_i}$, $\BIr_i\in\R^{m_0}$ $\forall 1\le\nobreak i\le\nobreak N$.
    Moreover, let us define 
    $y^{\circ}_{i,0}:=\hat{y}_{i,0}+\frac{\epsilon}{6N}u_{i,0}$,
    $\BIy^{\circ}_{i}:=\hat{\BIy}_{i}+\frac{\epsilon}{6N}\BIu_{i}$,
    $\BIw^{\circ}_{i}:=\hat{\BIw}_{i}+\frac{\epsilon}{6N}\BIr_{i}$
    $\forall 1\le i\le N$,
    and denote 
    $\Bchi^{\circ}:=(y^{\circ}_{1,0},\BIy_1^{\circ\TRANSP},\BIw_1^{\circ\TRANSP},\ldots,\allowbreak y^\circ_{N,0},\BIy_N^{\circ\TRANSP},\BIw_N^{\circ\TRANSP})^\TRANSP$.
    Since 
    $\Bchi^{\circ}=\hat{\Bchi}+\frac{\epsilon}{6N}\Brho$,
    $\Bchi^{\circ}$ is an arbitrary point in the closed $\|\cdot\|_{\infty}$-ball with radius $\frac{\epsilon}{6N}$ centered at~$\hat{\Bchi}$,
    and it thus suffices to prove that $\Bchi^{\circ}\in S_{\CI,\tau^\star_{\CI}-\epsilon}$.

    First, the definitions of $(\BIg_i(\,\cdot\,))_{i=0:N}$ guarantee that
    \begin{align}
        y^{\circ}_{i,0}+\langle\BIg_i(\BIv_{i,j}),\BIy^{\circ}_i\rangle+\langle\BIg_0(\BIv_{0,l}),\BIw^{\circ}_i\rangle
        &=\big(y^{\star}_{i,0}+\langle\BIg_i(\BIv_{i,j}),\BIy^{\star}_i\rangle+\langle\BIg_0(\BIv_{0,l}),\BIw^{\star}_i\rangle\big)\nonumber\\
        &\qquad -\frac{2\epsilon}{3N}+\frac{\epsilon}{6N}\langle\BIg_0(\BIv_{0,l}),\vecone_{m_0}\rangle\nonumber\\
        &\qquad+\frac{\epsilon}{6N}\big(u_{i,0}+\langle\BIg_i(\BIv_{i,j}),\BIu_i\rangle+\langle\BIg_0(\BIv_{0,l}),\BIr_i\rangle\big)\label{eqn:lp-any-index-proof-Linfinity-ball-1} \\
        &\le c_i(\BIv_{i,j},\BIv_{0,l})-c_{\ssmax}-\frac{2\epsilon}{3N}+\frac{\epsilon}{6N}\nonumber\\
        &\qquad + \frac{\epsilon}{6N}\big(|u_{i,0}|+\|\BIg_i(\BIv_{i,j})\|_{1}\|\BIu_{i}\|_{\infty}+\|\BIg_0(\BIv_{0,l})\|_{1}\|\BIr_{i}\|_{\infty}\big)\nonumber
        \\
        &\le c_i(\BIv_{i,j},\BIv_{0,l})-c_{\ssmax}\qquad\qquad\qquad\qquad\qquad\qquad\quad \forall (i,j,l)\in\CI.\nonumber
    \end{align}
    Second, it holds that
    \begin{align}
        \begin{split}
        \sum_{i=1}^N \BIw^{\circ}_{i}&=\sum_{i=1}^N\bigg(\BIw^{\star}_{i}+\frac{\epsilon}{6N}\vecone_{m_0}+\frac{\epsilon}{6N}\BIr_{i}\bigg)\ge\frac{\epsilon}{6}\Bigg(\vecone_{m_0}+\frac{1}{N}\sum_{i=1}^{N}\BIr_{i}\Bigg) \ge \veczero_{m_0}.
        \end{split}
        \label{eqn:lp-any-index-proof-Linfinity-ball-2}
    \end{align}
    Third, the definitions of $(\bar{\BIg}_i)_{i=1:N}$ show that
    \begin{align}
        \begin{split}
        \sum_{i=1}^N y^{\circ}_{i,0}+\langle\bar{\BIg}_i,\BIy^{\circ}_{i}\rangle 
        &= \Bigg(\sum_{i=1}^Ny^{\star}_{i,0}+\langle\bar{\BIg}_i,\BIy_i^{\star}\rangle\Bigg)-\frac{2\epsilon}{3}+\frac{\epsilon}{6N}\Bigg(\sum_{i=1}^Nu_{i,0}+\langle\bar{\BIg}_i,\BIu_i\rangle\Bigg)\\
        &\ge \tau^\star_{\CI}-\frac{2\epsilon}{3}-\frac{\epsilon}{6N} \Bigg(\sum_{i=1}^N |u_{i,0}|+\|\bar{\BIg}_i\|_1\|\BIu_i\|_{\infty}\Bigg)\\
        &\ge \tau^\star_{\CI}-\frac{2\epsilon}{3}-\frac{\epsilon}{3}
        =\tau^\star_{\CI}-\epsilon.
        \end{split}
        \label{eqn:lp-any-index-proof-Linfinity-ball-3}
    \end{align}
    Combining (\ref{eqn:lp-any-index-proof-Linfinity-ball-1})--(\ref{eqn:lp-any-index-proof-Linfinity-ball-3}),
    we conclude that  $\Bchi^{\circ}\in S_{\CI,\tau^\star_{\CI}-\epsilon}$.
    The proof of statement~\ref{lems:lp-any-index-Linfinity-ball} is now complete.

    Lastly, to prove statement~\ref{lems:lp-any-index-transform},
    let us fix an arbitrary $(y_{i,0},\BIy_i,\BIw_i)_{i=1:N}$
    that is feasible for \LPTAG{$\overline{\CI}$},
    and define
    $\hat{y}_{i,0}:=y_{i,0}-L_c\big(\eta(\FS_i)+\eta(\FS_0)\big)+c_{\ssmax}$,
    $\hat{\BIy}_i:=\BIy_i$ $\forall 1\le\nobreak i\le\nobreak N$,
    $\hat{\BIw}_1:=-\sum_{i=2}^{N}\BIw_i$,
    $\hat{\BIw}_i:=\BIw_i$ $\forall 2\le\nobreak i\le\nobreak N$.
    Subsequently, let us 
    fix an arbitrary $i\in\{1,\ldots,N\}$ and
    arbitrary $S_{i}\in \FS_i$, $S_{0}\in \FS_0$.
    Since $c_i$ is $L_{c}$-Lipschitz continuous,
    it follows from the definitions of $\eta(\FS_i)$ and $\eta(\FS_0)$ that 
    \begin{align*}
        \begin{split}
            c_{i}(\BIx,\BIz)&\ge c_{i}(\BIu,\BIv)-L_{c}\|\BIx-\BIu\|_{2}-L_{c}\|\BIz-\BIv\|_{2}\\
            &\ge c_{i}(\BIu,\BIv)-L_{c}\big(\eta(\FS_i)+\eta(\FS_0)\big) \qquad \forall \BIx\in S_i,\; \forall \BIz\in S_{0},\; \forall \BIu\in V(S_i),\; \forall \BIv\in V(S_0).
        \end{split}
    \end{align*}
    Combining this with the facts that 
    $\sum_{\BIu\in V(S_i)}\lambda_{\FS_i,\BIu}(\BIx)=\sum_{\BIv\in V(S_0)}\lambda_{\FS_0,\BIv}(\BIz)=1$ $\forall \BIx\in S_i$, $\forall \BIz\in S_0$,
    we get 
    \begin{align*}
        c_{i}(\BIx,\BIz)&\ge \Bigg(\sum_{\BIu\in V(S_i)}\sum_{\BIv\in V(S_0)}\lambda_{\FS_i,\BIu}(\BIx)\lambda_{\FS_0,\BIv}(\BIz)c_{i}(\BIu,\BIv)\Bigg) - L_{c}\big(\eta(\FS_i)+\eta(\FS_0)\big)
        \qquad \forall \BIx\in S_i,\; \forall \BIz\in S_{0}.
    \end{align*}
    Let us define 
    $\hat{c}_{i}(\BIx,\BIz):=\sum_{\BIu\in V(S_i)}\sum_{\BIv\in V(S_0)}\lambda_{\FS_i,\BIu}(\BIx)\lambda_{\FS_0,\BIv}(\BIz)c_{i}(\BIu,\BIv)$
    $\forall \BIx\in S_i$, $\forall \BIz\in S_0$.
    Observe that 
    the function $S_0\ni \BIz\mapsto \hat{c}_i(\BIx,\BIz)\in\nobreak\R$ is affine for each $\BIx\in S_i$,
    the function $S_i\ni \BIx\mapsto \hat{c}_i(\BIx,\BIz)\in\nobreak\R$ is affine for each $\BIz\in S_0$,
    that $\hat{c}_i(\BIu,\BIv)=c_i(\BIu,\BIv)$ $\forall \BIu\in V(S_i)$, $\forall \BIv\in V(S_0)$,
    and that 
    $S_i\ni\BIx\mapsto\langle\BIg_i(\BIx),\BIy_i\rangle\in\nobreak\R$,
    $S_0\ni\BIz\mapsto\langle\BIg_0(\BIz),\BIw_i\rangle\in\nobreak\R$ 
    are both affine.
    Consequently, 
    since every affine function attains it minimum over a simplex at an extreme point of the simplex,
    we get 
    \begin{align*}
            \hspace{20pt}&\hspace{-20pt}\min_{\BIx\in S_i,\, \BIz\in S_0}\big\{c_i(\BIx,\BIz)-\langle\BIg_i(\BIx),\BIy_i\rangle-\langle\BIg_0(\BIz),\BIw_i\rangle\big\}\\
            &\ge \min_{\BIx\in S_i,\, \BIz\in S_0}\big\{\hat{c}_i(\BIx,\BIz)-\langle\BIg_i(\BIx),\BIy_i\rangle-\langle\BIg_0(\BIz),\BIw_i\rangle\big\} - L_{c}\big(\eta(\FS_i)+\eta(\FS_0)\big)\\
            &= \min_{\BIx\in S_i,\, \BIv\in V(S_0)}\big\{\hat{c}_i(\BIx,\BIv)-\langle\BIg_i(\BIx),\BIy_i\rangle-\langle\BIg_0(\BIv),\BIw_i\rangle\big\} - L_{c}\big(\eta(\FS_i)+\eta(\FS_0)\big)\\
            &= \min_{\BIu\in V(S_i),\, \BIv\in V(S_0)}\big\{\hat{c}_i(\BIu,\BIv)-\langle\BIg_i(\BIu),\BIy_i\rangle-\langle\BIg_0(\BIv),\BIw_i\rangle\big\} - L_{c}\big(\eta(\FS_i)+\eta(\FS_0)\big)\\
            &= \min_{\BIu\in V(S_i),\, \BIv\in V(S_0)}\big\{c_i(\BIu,\BIv)-\langle\BIg_i(\BIu),\BIy_i\rangle-\langle\BIg_0(\BIv),\BIw_i\rangle\big\} - L_{c}\big(\eta(\FS_i)+\eta(\FS_0)\big)\\
            &\ge y_{i,0} + c_{\ssmax} - L_{c}\big(\eta(\FS_i)+\eta(\FS_0)\big)=\hat{y}_{i,0}.
    \end{align*}
    Since $S_i\in \FS_i$ and $S_0\in \FS_0$ are arbitrary and 
    since $\bigcup_{S\in \FS_i}S=\CX_i$, $\bigcup_{S\in \FS_0}S=\CX_0$,
    it thus follows that 
    \begin{align}
        \hat{y}_{i,0}+\langle\BIg_i(\BIx),\BIy_i\rangle+\langle\BIg_0(\BIz),\BIw_i\rangle\le c_i(\BIx,\BIz) 
        \qquad \forall (\BIx,\BIz)\in\CX_i\times\CX_0,\; \forall 1\le i\le N.
        \label{eqn:lp-any-index-proof-transform-ineq}
    \end{align}
    Moreover, it holds by definition that  
    $\BIg_0(\BIz)\in\R^{m_0}_{+}$ $\forall \BIz\in\CX_0$,
    and that
    $\hat{\BIy}_i=\BIy_i$,
    $\hat{\BIw}_i\le \BIw_i$ $\forall 1\le i\le N$.
    Combining this with (\ref{eqn:lp-any-index-proof-transform-ineq})
    leads to
    \begin{align*}
        \hat{y}_{i,0}+\langle\BIg_i(\BIx),\hat{\BIy}_i\rangle+\langle\BIg_0(\BIz),\hat{\BIw}_i\rangle
        &\le \hat{y}_{i,0}+\langle\BIg_i(\BIx),\BIy_i\rangle+\langle\BIg_0(\BIz),\BIw_i\rangle\le c_i(\BIx,\BIz)  \\*
        & \hspace{60pt}\forall (\BIx,\BIz)\in\CX_i\times\CX_0,\; \forall 1\le i\le N.
    \end{align*}
    Finally, one checks that 
    $\sum_{i=1}^{N}\hat{\BIw}_i=\veczero_{m_0}$,
    which shows that 
    $(\hat{y}_{i,0},\hat{\BIy}_i,\hat{\BIw}_i)_{i=1:N}$
    is feasible for \eqref{eqn:mt-tf-lsip}.
    The proof is now complete.
\endproof

Next, let us establish some crucial properties of $\GetSepHyp$ in the following lemma, which show that it satisfies the requirements for the separation oracle in the volumetric center algorithm of \citet{ECvaidya1996new}.

\begin{lemma}\label{lem:proc-get-sep-hyp}%
Under Setting~\ref{sett:Euclidean},
let us denote $m:=\max_{0\le i\le N}\{m_i\}$,
and let $T_{\mathtt{Eval}}$
denote the number of arithmetic operations incurred when evaluating $c_i(\BIx,\BIz)$ at any $(\BIx,\BIz)\in\CX_i\times\CX_0$.
Let $\Bchi\in\R^{q}$ be arbitrary,
and let 
$(\mathtt{Inside},\mathtt{Tag},\BIa)$ be the output of 
$\GetSepHyp\big(N,\allowbreak(\CG_i)_{i=0:N},\allowbreak (\bar{\BIg}_i,c_i)_{i=1:N},c_{\ssmax},\allowbreak \Bchi\big)$.
Then, the following statements hold.
\begin{enumerate}[label=(\roman*),beginpenalty=10000]
    \item\label{lems:proc-get-sep-hyp-complexity}
    $\GetSepHyp\big(N,\allowbreak(\CG_i)_{i=0:N},\allowbreak (\bar{\BIg}_i,c_i)_{i=1:N},c_{\ssmax},\allowbreak \Bchi\big)$
    incurs $O(Nm^2T_{\mathtt{Eval}})$ arithmetic operations.
    
    \item\label{lems:proc-get-sep-hyp-separation}
    If $\mathtt{Inside}=\mathtt{False}$, then it holds that 
    $\Bchi\notin S_{\overline{\CI},-\infty}$,
    and that 
    $\langle\BIa,\Bomega\rangle\ge \langle\BIa,\Bchi\rangle$ $\forall \Bomega\in S_{\overline{\CI},-\infty}$.
    If~$\mathtt{Inside}=\mathtt{True}$, then it holds that
    $\Bchi\in S_{\overline{\CI},-\infty}$,
    and that 
    $\langle\BIa,\Bomega\rangle\ge \langle\BIa,\Bchi\rangle$ $\forall \Bomega\in S_{\overline{\CI},\langle\bar{\BIg},\Bchi\rangle}$.
\end{enumerate}
\end{lemma}

\proof{Proof of Lemma~\ref{lem:proc-get-sep-hyp}.}
The number of arithmetic operations in $\GetSepHyp$ is self-evident.
Note that Line~\ref{alglin:get-sep-hyp-balancecut} is executed at most once throughout the for-loop in Line~\ref{alglin:get-sep-hyp-balancecut-forloop}.
Similarly, Lines~\ref{alglin:get-sep-hyp-feascutstart}--\ref{alglin:get-sep-hyp-feascutend} are executed at most once throughout the nested for-loops in Lines~\ref{alglin:get-sep-hyp-feascut-forloop}--\ref{alglin:get-sep-hyp-feascut-forloop-innermost}.
This proves statement~\ref{lems:proc-get-sep-hyp-complexity}.

In the following, let us
denote $\Bchi=(y_{1,0},y_{1,1},\ldots,y_{1,m_1},w_{1,1},\ldots,w_{1,m_0},\ldots,y_{N,0},y_{N,1},\ldots,y_{N,m_N},\allowbreak w_{N,1},\ldots,w_{N,m_0})^\TRANSP$ and
prove statement~\ref{lems:proc-get-sep-hyp-separation}.
First, in the case where $\GetSepHyp$ terminates in Line~\ref{alglin:get-sep-hyp-balancecut-return} and returns $\mathtt{Inside}=\mathtt{False}$,
it holds that 
$\langle\BIa,\Bchi\rangle=\sum_{i=1}^{N}w_{i,l}<0$,
which implies that $\Bchi\notin S_{\overline{\CI},-\infty}$.
Moreover, the constraints of \LPTAG{$\overline{\CI}$} 
guarantee that $\langle\BIa,\Bomega\rangle\ge 0 > \langle\BIa,\Bchi\rangle$ 
$\forall\Bomega\in S_{\overline{\CI},-\infty}$.

Second, in the case where $\GetSepHyp$ terminates in Line~\ref{alglin:get-sep-hyp-feascut-return} and returns $\mathtt{Inside}=\mathtt{False}$,
it holds that 
$\langle\BIa,\Bchi\rangle=-(y_{i,0}+\INDI_{\{j\ne 0\}}y_{i,j}+\INDI_{\{l\ne 0\}}w_{i,l})=-y_{i,0}-\big(\sum_{\bar{j}=1}^{m_i}g_{i,\bar{j}}(\BIv_{i,j})y_{i,\bar{j}}\big)-\big(\sum_{\bar{l}=1}^{m_0}g_{0,\bar{l}}(\BIv_{0,l})w_{i,\bar{l}}\big)<c_{\ssmax}-c_i(\BIv_{i,j},\BIv_{0,l})$,
which implies that $\Bchi\notin S_{\overline{\CI},-\infty}$.
Moreover, the constraints of \LPTAG{$\overline{\CI}$} 
guarantee that $\langle\BIa,\Bomega\rangle\ge c_{\ssmax}-c_i(\BIv_{i,j},\BIv_{0,l}) > \langle\BIa,\Bchi\rangle$ 
$\forall\Bomega\in S_{\overline{\CI},-\infty}$.

Third, in the case where $\GetSepHyp$ terminates in Line~\ref{alglin:get-sep-hyp-objcut-return} and returns $\mathtt{Inside}=\mathtt{True}$,
then Line~\ref{alglin:get-sep-hyp-balancecut-check} and Line~\ref{alglin:get-sep-hyp-feascut-check} guarantee that 
$\sum_{i=1}^{N}w_{i,l}\ge 0$ $\forall 1\le l\le m_0$,
and that
$y_{i,0}+\big(\sum_{\bar{j}=1}^{m_i}g_{i,\bar{j}}(\BIv_{i,j})y_{i,\bar{j}}\big)+\big(\sum_{\bar{l}=1}^{m_0}g_{0,\bar{l}}(\BIv_{0,l})w_{i,\bar{l}}\big)
=y_{i,0}+\INDI_{\{j\ne0\}}y_{i,j}+\INDI_{\{l\ne 0\}}w_{i,l}
\le c_i(\BIv_{i,j},\BIv_{0,l})-c_{\ssmax}$
$\forall (i,j,l)\in\overline{\CI}$.
This shows that 
$\Bchi\in\nobreak S_{\overline{\CI},-\infty}$.
Furthermore, since $\BIa=\bar{\BIg}$,
it follows from the definition of $S_{\overline{\CI},\langle\bar{\BIg},\Bchi\rangle}$
that 
$\Bomega\in S_{\overline{\CI},\langle\bar{\BIg},\Bchi\rangle}$
only if 
$\langle\BIa,\Bomega\rangle=\langle\bar{\BIg},\Bomega\rangle\ge \langle\bar{\BIg},\Bchi\rangle = \langle\BIa,\Bchi\rangle$.
The proof is now complete.
\endproof

We are now ready to prove Theorem~\ref{thm:proc-solve-LSIP}.

\proof{Proof of Theorem~\ref{thm:proc-solve-LSIP}.}
Let us first present an overview of $\SolveLSIPDual$, 
which is a direct application of the volumetric center algorithm of \citet{ECvaidya1996new} for solving \LPTAG{$\overline{\CI}$}.
Lines~\ref{alglin:solve-LSIP-dual-init-start}--\ref{alglin:solve-LSIP-dual-init-iternum} carry out 
the initialization steps before the main iterations of 
the volumetric center algorithm.
Specifically, 
Line~\ref{alglin:solve-LSIP-dual-init-box}
constructs 
$P_{0}:=\big\{\Bchi\in\R^{q}:\langle\BIa_{j},\Bchi\rangle\ge b_{j}\; \forall 1\le j\le 2q\big\}$
as the initial polytope.
Observe that $P_{0}=\nobreak B_{\infty}(M)$,
and that Line~\ref{alglin:solve-LSIP-dual-init-center} sets 
$\Bchi_{0}$ to be the volumetric center of $P_{0}$;
see \citep[p.294--295]{ECvaidya1996new} for the definition and properties of the volumetric center.
Moreover, 
Line~\ref{alglin:solve-LSIP-dual-init-iternum} computes the number of iterations $\overline{r}\in\N$ of 
the volumetric center algorithm.
We would like to remark that the constant terms in 
Lines~\ref{alglin:solve-LSIP-dual-init-iternum}, 
\ref{alglin:solve-LSIP-dual-init-check-addorrem},
\ref{alglin:solve-LSIP-dual-RHS},
\ref{alglin:solve-LSIP-dual-add-Newtonloop},
\ref{alglin:solve-LSIP-dual-rem-Newtonloop}
are calculated by choosing $\delta\leftarrow 10^{-4}$, $\varepsilon\leftarrow 10^{-7}$ in \citep{ECvaidya1996new}.
\citet[p.295]{ECvaidya1996new} allows $\delta$ and $\varepsilon$ 
to be arbitrary constants that satisfy 
$\delta\le 10^{-4}$ and $\varepsilon\le 10^{-3}\delta$.
We fix these constants in our presentation in order to avoid notational clutter.

Subsequently, Lines~\ref{alglin:solve-LSIP-dual-forloop}--\ref{alglin:solve-LSIP-dual-rem-end} carry out the main iterations of the volumetric center algorithm.
The algorithm maintains a polytope throughout the iterations via the list $\big\{(\BIa_j,b_j,\iota_j,t_j)\big\}_{j=1:\overline{k}}$,
where each
$(\BIa_j,b_j)\in\nobreak\R^{q}\times\R$ describes a closed half-space 
$\big\{\Bchi\in\R^{q}:\langle\BIa_j,\Bchi\rangle\ge b_j\big\}$,
and 
$t_j\in\{0,1\}$ indicates whether this closed half-space is present in the characterization of the current polytope.
In subsequent analyses, let us denote the polytope maintained at the end of iteration~$r$ by
$P_{r}:=\big\{\Bchi\in\R^{q}:\langle\BIa_j,\Bchi\rangle\ge b_j\; \forall j\in\{1,\ldots,\overline{k}\},\;t_{j}=1\}$,
where $\overline{k}$ and $(t_j)_{j=1:\overline{k}}$ 
correspond to their values at the end of iteration~$r$.
At any point of the volumetric center algorithm,
we say a closed half-space 
$\big\{\Bchi\in\R^{q}:\langle\BIa_j,\Bchi\rangle\ge b_j\big\}$
is in the characterization of the current polytope if $t_j=1$.
Moreover, 
each $\iota_j$ is a tuple which is either $(0,0,0)$ or an element of $\overline{\CI}$ where the corresponding inequality constraint had been violated 
at some point of the algorithm.
These tuples are subsequently gathered into a set $\widehat{\CI}$ in Line~\ref{alglin:solve-LSIP-dual-postproc-indexset} 
which is then used in $\SolveLSIPPrimal$ to approximately solve \eqref{eqn:mt-tf-dual}.
Throughout the iterations,
the procedure $\GetVolumetricQuantities$ computes quantities related to the volumetric center in \citep[p.294--295]{ECvaidya1996new}.
Line~\ref{alglin:solve-LSIP-dual-init-check-addorrem} decides 
whether to add a new closed half-space to the characterization of the current polytope, 
or to remove an existing closed half-space from the characterization of the current polytope.
Concretely, if Line~\ref{alglin:solve-LSIP-dual-init-check-addorrem} decides 
to add a new closed half-space,
then Lines~\ref{alglin:solve-LSIP-dual-add-start}--\ref{alglin:solve-LSIP-dual-add-end}
add the new closed half-space $\big\{\Bchi\in\R^{q}:\langle\hat{\BIa},\Bchi\rangle\ge \hat{b}\big\}$ to the characterization of the current polytope,
i.e., $P_{r-1}$ is updated to $P_{r}$ where $P_{r}=P_{r-1}\cap \big\{\Bchi\in\R^{q}:\langle\hat{\BIa},\Bchi\rangle\ge \hat{b}\big\}$,
as reflected by
Line~\ref{alglin:solve-LSIP-dual-add}.
The for-loop in
Line~\ref{alglin:solve-LSIP-dual-add-Newtonloop}
then begins at an approximate volumetric center $\Bchi_{r-1}$ of $P_{r-1}$
and performs a constant number of Newton-type updates to
compute an approximate volumetric center $\Bchi_{r}$ of $P_{r}$.
On the other hand, 
if Line~\ref{alglin:solve-LSIP-dual-init-check-addorrem} decides 
to remove an existing closed half-space,
then 
Lines~\ref{alglin:solve-LSIP-dual-rem-start}--\ref{alglin:solve-LSIP-dual-rem-end}
remove the closed half-space 
$\big\{\Bchi\in\R^{q}:\langle\BIa_{j_{\hat{l}}},\Bchi\rangle\ge b_{j_{\hat{l}}}\big\}$ from the characterization of the current polytope,
i.e., $P_{r-1}$ is updated to $P_{r}$ where $P_{r-1}=P_{r}\cap \big\{\Bchi\in\R^{q}:\langle\BIa_{j_{\hat{l}}},\Bchi\rangle\ge b_{j_{\hat{l}}}\big\}$,
as reflected by
Line~\ref{alglin:solve-LSIP-dual-rem-start}.
The for-loop in
Line~\ref{alglin:solve-LSIP-dual-rem-Newtonloop}
then starts at an approximate volumetric center $\Bchi_{r-1}$ of $P_{r-1}$
and performs a constant number of Newton-type updates to
compute an approximate volumetric center $\Bchi_{r}$ of $P_{r}$.
Additionally,
we maintain $\hat{\Bchi}\in\R^{q}$ and $\hat{\tau}\in\R\cup\{-\infty\}$
throughout the main iterations by Lines~\ref{alglin:solve-LSIP-dual-init-center}, 
\ref{alglin:solve-LSIP-dual-inside-check}, 
\ref{alglin:solve-LSIP-dual-inside-update}.
One observes that, after
Lines~\ref{alglin:solve-LSIP-dual-inside-check}--\ref{alglin:solve-LSIP-dual-inside-update}
have been executed at least once,
$\hat{\Bchi}$ and $\hat{\tau}$
satisfy 
$\hat{\Bchi}\in S_{\overline{\CI},-\infty}$
and 
$\hat{\tau}=\langle\bar{\BIg},\hat{\Bchi}\rangle$.
Therefore, when the for-loop in Line~\ref{alglin:solve-LSIP-dual-forloop} terminates,
$\hat{\Bchi}$ is the incumbent feasible solution of \LPTAG{$\overline{\CI}$},
whereas 
$\hat{\tau}$ is the best objective value encountered in the algorithm.

After the main iterations,
Lines~\ref{alglin:solve-LSIP-dual-postproc-start}--\ref{alglin:solve-LSIP-dual-postproc-indexset}
carry out the post-processing steps.
Specifically, 
Line~\ref{alglin:solve-LSIP-dual-postproc-solution}
modifies $\hat{\Bchi}$ into a feasible solution 
$(\hat{y}_{i,0},\hat{\BIy}_i,\hat{\BIw}_i)_{i=1:N}$ 
of \eqref{eqn:mt-tf-lsip}, as guaranteed by statement~\ref{thms:proc-solve-LSIP-dual-feasibility},
whereas 
Line~\ref{alglin:solve-LSIP-dual-postproc-indexset}
produces an index set $\widehat{\CI}\subseteq\overline{\CI}$.
Subsequently, in $\SolveLSIPPrimal$,
the algorithm of \citet{ECvandenbrand2020deterministic} is used to solve 
the dual LP problem of 
\LPTAG{$\widehat{\CI}_{+}$} with respect to $\widehat{\CI}_{+}:=\widehat{\CI}\cup \CV$,
in order to compute a possibly infeasible approximate optimizer of \eqref{eqn:mt-tf-dual}.

The proof of
statement~\ref{thms:proc-solve-LSIP-dual-complexity}
follows from the rate $q=O(Nm)$,
Lemma~\ref{lem:proc-get-sep-hyp}\ref{lems:proc-get-sep-hyp-complexity},
and the properties of the volumetric center algorithm established by \citet{ECvaidya1996new}.
We explain in the following list some crucial properties that are used for deriving the number of arithmetic operations in $\SolveLSIPDual$.
\begin{itemize}
    \item In Line~\ref{alglin:solve-LSIP-dual-init-iternum},
    observe that $\overline{r}=O\big(Nm\log\big(\frac{Nm(c_{\ssmax}-c_{\ssmin}+1)}{\epsilon_0\bar{g}_{\ssmin}}\big)\big)$;
    note that $\bar{g}_{\ssmin}>0$ due to Lemma~\ref{lem:lp-any-index}\ref{lems:lp-any-index-minint}.
    Hence, the for-loop in Line~\ref{alglin:solve-LSIP-dual-forloop} is carried out $O\big(Nm\log\big(\frac{Nm(c_{\ssmax}-c_{\ssmin}+1)}{\epsilon_0\bar{g}_{\ssmin}}\big)\big)$ times,
    and the value of $\overline{k}$ when $\SolveLSIPDual$ terminates
    is thus equal to $2q+\overline{r}=O\big(Nm\log\big(\frac{Nm(c_{\ssmax}-c_{\ssmin}+1)}{\epsilon_0\bar{g}_{\ssmin}}\big)\big)$.
    
    \item Line~\ref{alglin:solve-LSIP-dual-init-check-addorrem} and Line~\ref{alglin:solve-LSIP-dual-rem-start}
    guarantee that the number of closed half-spaces in the characterization of the current polytope $P_{r}$ never exceeds $10^7q=O(Nm)$,
    as discussed in \citep[p.297]{ECvaidya1996new}.
    As a consequence, one can check that $|\widehat{\CI}|=O(Nm)$.
    Furthermore, 
    Line~\ref{alglin:solve-LSIP-dual-activecuts} and Line~\ref{alglin:solve-LSIP-dual-postproc-indexset} 
    can be implemented by keeping track of 
    $\big\{j\in\{1,\ldots,\overline{k}\}:t_j=1\big\}$ throughout the algorithm,
    whose cardinality is at most $10^7q=O(Nm)$.

    \item The number of arithmetic operations incurred when computing the quantities related to the volumetric center in $\GetVolumetricQuantities$ and Line~\ref{alglin:solve-LSIP-dual-RHS} of $\SolveLSIPDual$ are derived in \citep[p.297--298]{ECvaidya1996new}.
\end{itemize}
The proof of statement~\ref{thms:proc-solve-LSIP-dual-complexity} is now complete.

To prove statement~\ref{thms:proc-solve-LSIP-primal-complexity},
observe that 
statement~\ref{thms:proc-solve-LSIP-dual-complexity} and Line~\ref{alglin:solve-LSIP-primal-aggindex} 
of $\SolveLSIPPrimal$
imply that 
$p=|\widehat{\CI}_{+}|\le |\widehat{\CI}|+|\CV|=O(Nm)$
when the for-loop in Line~\ref{alglin:solve-LSIP-primal-buildconstr} terminates,
and thus $p+m_0=O(Nm)$.
Moreover, since $\widehat{\CI}_{+}\supseteq \CV$, 
one checks that the matrix $\BA$ has rank $q$.
Subsequently, the number of arithmetic operations incurred by 
Line~\ref{alglin:solve-LSIP-primal-solver}
of $\SolveLSIPPrimal$ follows from \citep[Theorem~1.1]{ECvandenbrand2020deterministic} with respect to 
$A\leftarrow \BA$, 
$b\leftarrow \bar{\BIg}$,
$c\leftarrow\tilde{\BIc}$,
$n\leftarrow p+m_0$,
$R\leftarrow N+1$,
$\delta\leftarrow \frac{\varsigma}{4(N+1)(p+1)}$.
Note that in the notation of
\citep[Theorem~1.1]{ECvandenbrand2020deterministic},
$R>0$ corresponds to an upper bound for 
$\sup\big\{\|\Bpsi\|_1:\Bpsi\in\R^{p+m_0}_{+},$ $\BA\Bpsi=\nobreak\bar{\BIg}\big\}$.
Hence, we are using here the property that
\begin{align}
    \sup\big\{\|\Bpsi\|_1:\Bpsi\in\R^{p+m_0}_{+},\;\BA\Bpsi=\bar{\BIg}\big\}\le N+1,
    \label{eqn:proc-solve-LSIP-proof-primal-L1bound}
\end{align} 
which we will prove later in the proof of statement~\ref{thms:proc-solve-LSIP-primal-feasibility}.
The number of arithmetic operations incurred by the remaining lines of 
$\SolveLSIPPrimal$ is evident from the property that $|\widehat{\CI}_{+}|=\nobreak O(Nm)$.

Next, let $\overline{\tau}$ denote the value of $\hat{\tau}$ when $\SolveLSIPDual$ terminates, and
let us divide the remainder of this proof into the following 5 steps.
\begin{itemize}
    \item Step~1: showing that $S_{\widehat{\CI},\overline{\tau}} \cap B_{\infty}(M)\subseteq P_{\overline{r}}$.
    
    \item Step~2: showing that $\Lebesgue_{\R^q}(P_{\overline{r}})< \Lebesgue_{\R^{q}}\big(B_{\infty}\big(\frac{\epsilon_0}{6N}\big)\big)$.

    \item Step~3: proving statement~\ref{thms:proc-solve-LSIP-dual-feasibility}
    and showing that 
    $\sum_{i=1}^{N}\hat{y}_{i,0}+\langle\bar{\BIg}_i,\hat{\BIy}_i\rangle
    =\overline{\tau}-L_{c}\big(N\overline{\eta}_0+\sum_{i=1}^{N}\overline{\eta}_i\big)$.
    
    \item Step~4: proving (\ref{eqn:proc-solve-LSIP-proof-primal-L1bound}) and statement~\ref{thms:proc-solve-LSIP-primal-feasibility},
    as well as showing that 
    $\sum_{i=1}^{N}\int_{\CX_i\times\CX_0}c_i\DIFFX{\hat{\theta}_i}\le \tau^\star_{\widehat{\CI}_{+}}+(c_{\ssmax}-c_{\ssmin})\varsigma$.
    
    \item Step~5: showing that $\tau^\star_{\widehat{\CI}_{+}}\le \overline{\tau}+\epsilon_0$ and proving statement~\ref{thms:proc-solve-LSIP-gap}.
\end{itemize}

\underline{Step~1.}
Let us fix an arbitrary closed half-space 
$\widehat{H}:=\big\{\Bchi\in\R^{q}:\langle\hat{\BIa},\Bchi\rangle\ge \hat{b}\big\}$
characterizing $P_{\overline{r}}$.
In the case where $(\hat{\BIa},\hat{b})=(\BIa_{j},b_{j})$ for some $j\in\{1,\ldots,2q\}$, Line~\ref{alglin:solve-LSIP-dual-init-box} of $\SolveLSIPDual$ guarantees that 
$S_{\widehat{\CI},\overline{\tau}} \cap B_{\infty}(M)\subseteq B_{\infty}(M)\subseteq \widehat{H}$.
Otherwise, there exists $r\in\{1,\ldots,\overline{r}\}$ 
such that 
$(\BIa_{2q+r},b_{2q+r})=(\hat{\BIa},\hat{b})$ and 
$t_{2q+r}=1$ are generated in iteration~$r$ of the for-loop in Line~\ref{alglin:solve-LSIP-dual-forloop} of $\SolveLSIPDual$.
Let $(\mathtt{Inside}_{r},\mathtt{Tag}_{r},\BIa_{2q+r})$
denote the output of $\GetSepHyp$ in Line~\ref{alglin:solve-LSIP-dual-add-start} in iteration~$r$,
and let us analyze the three cases below.

\emph{Case~1: $\mathtt{Inside}_{r}=\mathtt{False}$, $\mathtt{Tag}_{r}=(0,0,0)$.}
In this case, observe from Line~\ref{alglin:get-sep-hyp-balancecut} of $\GetSepHyp$ that 
$\BIa_{2q+r}=(0,\veczero_{m_1}^\TRANSP,\BIe_{0,l}^\TRANSP,0,\veczero_{m_2}^\TRANSP,\BIe_{0,l}^\TRANSP,\ldots,0,\veczero_{m_N}^\TRANSP,\BIe_{0,l}^\TRANSP)^\TRANSP$ for some $l\in\{1,\ldots,m_0\}$,
and thus the constraints of \LPTAG{$\widehat{\CI}$},
Line~\ref{alglin:get-sep-hyp-balancecut-check} of $\GetSepHyp$,
and Line~\ref{alglin:solve-LSIP-dual-RHS} of $\SolveLSIPDual$
guarantee that 
$\langle\BIa_{2q+r},\Bchi\rangle\ge 0> \langle\BIa_{2q+r},\Bchi_{r-1}\rangle>b_{2q+r}$ $\forall \Bchi\in S_{\widehat{\CI},-\infty}$.
This implies that 
$S_{\widehat{\CI},\overline{\tau}} \cap B_{\infty}(M)\subseteq S_{\widehat{\CI},-\infty}\subseteq \widehat{H}$.

\emph{Case~2: $\mathtt{Inside}_{r}=\mathtt{False}$, $\mathtt{Tag}_{r}\ne(0,0,0)$.}
In this case, let us denote $\mathtt{Tag}_{r}=(i_{r},j_{r},l_{r})$.
It thus holds by definition that $(i_{r},j_{r},l_{r})\in\widehat{\CI}$.
Subsequently,
the constraints of \LPTAG{$\widehat{\CI}$},
Lines~\ref{alglin:get-sep-hyp-feascut-check}--\ref{alglin:get-sep-hyp-feascutend} of $\GetSepHyp$,
and
Line~\ref{alglin:solve-LSIP-dual-RHS} of $\SolveLSIPDual$
guarantee that
$\langle\BIa_{2q+r},\Bchi\rangle\ge c_{\ssmax}-c_{i_{r}}(\BIv_{i_{r},j_{r}},\BIv_{0,l_{r}})> \langle\BIa_{2q+r},\Bchi_{r-1}\rangle>b_{2q+r}$
$\forall\Bchi\in S_{\widehat{\CI},-\infty}$.
This implies that 
$S_{\widehat{\CI},\overline{\tau}} \cap B_{\infty}(M)\subseteq S_{\widehat{\CI},-\infty}\subseteq \widehat{H}$.

\emph{Case~3: $\mathtt{Inside}_{r}=\mathtt{True}$.}
In this case, observe from Line~\ref{alglin:get-sep-hyp-objcut-return} of $\GetSepHyp$
that $\BIa_{2q+r}=\nobreak\bar{\BIg}$.
\mbox{Lines~\ref{alglin:solve-LSIP-dual-inside-check}--\ref{alglin:solve-LSIP-dual-RHS}} 
of $\SolveLSIPDual$ then guarantee that 
$\overline{\tau}>-\infty$ and
$\langle\BIa_{2q+r},\Bchi\rangle\ge \overline{\tau}\ge \langle\BIa_{2q+r},\Bchi_{r-1}\rangle>b_{2q+r}$ $\forall \Bchi\in S_{\widehat{\CI},\overline{\tau}}$.
This implies that 
$S_{\widehat{\CI},\overline{\tau}} \cap B_{\infty}(M)\subseteq S_{\widehat{\CI},\overline{\tau}}\subseteq \widehat{H}$.

Summarizing the cases above, we conclude that $S_{\widehat{\CI},\overline{\tau}} \cap B_{\infty}(M)\subseteq P_{\overline{r}}$.
This completes Step~1. 

\underline{Step~2.}
Applying an inequality in \citep[p.297]{ECvaidya1996new} for bounding the number of iterations (with respect to 
$\pi^k\leftarrow \Lebesgue_{\R^{q}}(P_{\overline{r}})$,
$n\leftarrow q$,
$\varepsilon\leftarrow 10^{-7}$,
$k\leftarrow \overline{r}$,
$\rho^0\leftarrow \frac{1}{2}\log\det\Big(\sum_{j=1}^{2q}\frac{\BIa_{j}\BIa_{j}^\TRANSP}{(\langle\BIa_j,\Bchi_{0}-b_{j}\rangle)^2}\Big)=q\big(\frac{\log(2)}{2}-\log(M)\big)$
in the notation of \citet{ECvaidya1996new}),
we get 
\begin{align*}
    \log\big(\Lebesgue_{\R^{q}}(P_{\overline{r}})\big)
    &\le q\log(2\cdot10^{7}q)-q\bigg(\frac{\log(2)}{2}-\log(M)\bigg) - \frac{10^{-7}}{2}\overline{r}.
\end{align*}
Since Line~\ref{alglin:solve-LSIP-dual-init-iternum} 
of $\SolveLSIPDual$
guarantees 
$\overline{r}\ge 2\cdot 10^{7}q\big(18+\log\big(\frac{qNM}{\epsilon_0}\big)\big)$,
it subsequently follows that
\begin{align*}
    \log\big(\Lebesgue_{\R^{q}}(P_{\overline{r}})\big)
    &\le q\bigg(\log(3\sqrt{2}\cdot 10^7)-18 + \log\bigg(\frac{\epsilon_0}{3N}\bigg)\bigg)
    < q\log\bigg(\frac{\epsilon_0}{3N}\bigg)
    = \log\bigg(\Lebesgue_{\R^{q}}\bigg(B_{\infty}\bigg(\frac{\epsilon_0}{6N}\bigg)\bigg)\bigg).
\end{align*}
Step~2 is now complete.

\underline{Step~3.}
Let us first show that $\overline{\tau}>-\infty$.
To that end, suppose for the sake of contradiction that 
$\overline{\tau}=\nobreak-\infty$.
On the one hand, since $\epsilon_0<1$,
Lemma~\ref{lem:lp-any-index}\ref{lems:lp-any-index-Linfinity-bound}
implies that 
$S_{\overline{\CI},\tau^\star_{\overline{\CI}}-\epsilon_0}
=S_{\overline{\CI},\tau^\star_{\overline{\CI}}-\epsilon_0} \cap B_{\infty}(M)\subseteq S_{\widehat{\CI},\overline{\tau}}\cap B_{\infty}(M)$.
The conclusion of Step~1 then shows that 
$S_{\overline{\CI},\tau^\star_{\overline{\CI}}-\epsilon_0}\subseteq P_{\overline{r}}$,
which yields 
$\Lebesgue_{\R^{q}}(S_{\overline{\CI},\tau^\star_{\overline{\CI}}-\epsilon_0})\le \Lebesgue_{\R^{q}}(P_{\overline{r}})$.
On the other hand,
Lemma~\ref{lem:lp-any-index}\ref{lems:lp-any-index-Linfinity-ball}
implies that 
$S_{\overline{\CI},\tau^\star_{\overline{\CI}}-\epsilon_0}$ 
contains a closed $\|\cdot\|_{\infty}$-ball with radius~$\frac{\epsilon_0}{6N}$,
and we thus get
$\Lebesgue_{\R^{q}}\big(B_{\infty}\big(\frac{\epsilon_0}{6N}\big)\big)
\le \Lebesgue_{\R^{q}}(S_{\overline{\CI},\tau^\star_{\overline{\CI}}-\epsilon_0})
\le \Lebesgue_{\R^{q}}(P_{\overline{r}})$,
contradicting the conclusion of Step~2.
Thus, $\overline{\tau}>-\infty$,
and therefore 
Line~\ref{alglin:solve-LSIP-dual-inside-update} of $\SolveLSIPDual$
is executed at least once.
Lemma~\ref{lem:proc-get-sep-hyp}\ref{lems:proc-get-sep-hyp-separation} then 
ensures that $\hat{\Bchi}\in S_{\overline{\CI},-\infty}$
when the for-loop in Line~\ref{alglin:solve-LSIP-dual-forloop} of $\SolveLSIPDual$ terminates.
Consequently, 
the properties that $\overline{\eta}_i\ge \eta(\FS_i)$ $\forall 0\le i\le N$,
Lines~\ref{alglin:solve-LSIP-dual-postproc-start}--\ref{alglin:solve-LSIP-dual-postproc-solution} of $\SolveLSIPDual$
and Lemma~\ref{lem:lp-any-index}\ref{lems:lp-any-index-transform}
yield that 
$(\hat{y}_{i,0},\hat{\BIy}_i,\hat{\BIw}_i)_{i=1:N}$ is feasible for \eqref{eqn:mt-tf-lsip}.
This completes the proof of statement~\ref{thms:proc-solve-LSIP-dual-feasibility}.
Furthermore, it holds that 
\begin{align}
    \sum_{i=1}^{N}\hat{y}_{i,0}+\langle\bar{\BIg}_i,\hat{\BIy}_i\rangle
    &=\langle\bar{\BIg},\hat{\Bchi}\rangle-L_{c}\Bigg(N\overline{\eta}_0+\sum_{i=1}^{N}\overline{\eta}_i\Bigg)
    =\overline{\tau}-L_{c}\Bigg(N\overline{\eta}_0+\sum_{i=1}^{N}\overline{\eta}_i\Bigg).
    \label{eqn:proc-solve-LSIP-proof-dualobj}
\end{align}
Step~3 is now complete.

\underline{Step~4.}
One observes from Lines~\ref{alglin:solve-LSIP-primal-buildconstr}--\ref{alglin:solve-LSIP-primal-aggregate} of $\SolveLSIPPrimal$
that the LP problem:
\begin{align*}
    \minimize_{\Bpsi}\quad & \langle\tilde{\BIc},\Bpsi\rangle \\
    \mathrm{subject~to}\quad & \BA\Bpsi=\bar{\BIg},\quad \Bpsi\in\R^{p+m_0}_{+}
\end{align*}
is equivalent to the following LP problem:
\begin{align}
    \begin{split}
        \minimize_{(\theta_{i,j,l}),\,(\xi_{l})}\quad & \sum_{i=1}^{N} \sum_{(j,l)\in\widehat{\CI}_i} \big(c_i(\BIv_{i,j},\BIv_{0,l})-c_{\ssmax}\big)\theta_{i,j,l}\\
        \mathrm{subject~to}\quad & \sum_{(j,l)\in\widehat{\CI}_i}\theta_{i,j,l}=1 
        \hspace{144.6pt} \quad \forall 1\le i\le N, \\
        & \sum_{(j,l)\in\widehat{\CI}_i}g_{i,\bar{j}}(\BIv_{i,j})\theta_{i,j,l}= \int_{\CX_i}g_{i,\bar{j}}\DIFFX{\mu_i} 
        \quad \forall 1\le \bar{j}\le m_i,\; \forall 1\le i\le N,\\
        & \sum_{(j,l)\in\widehat{\CI}_i}g_{0,\bar{l}}(\BIv_{0,l})\theta_{i,j,l}=\xi_{\bar{l}} 
        \hspace{41.9pt}\quad \forall 1\le \bar{l}\le m_0,\; \forall 1\le i\le N,\\
        & \theta_{i,j,l}\ge 0 \quad \forall (i,j,l)\in\widehat{\CI}_{+}, \hspace{47.3pt}\qquad 
        \xi_{l}\ge 0 \quad \forall 1\le l\le m_0,
    \end{split}
    \label{eqn:mt-tf-lp-any-index-primal}
\end{align}
which is the dual LP problem of \LPTAG{$\widehat{\CI}_{+}$}.
Before we apply \citep[Theorem~1.1]{ECvandenbrand2020deterministic}, 
let us first derive some bounds related to 
$\BA$, $\bar{\BIg}$, and $\tilde{\BIc}$.
First, observe that $\|\tilde{\BIc}\|_{\infty}\le c_{\ssmax}-c_{\ssmin}$
and that 
$\|\bar{\BIg}\|_{1}\le 2N$.
Second, observe from
Lines~\ref{alglin:solve-LSIP-primal-ineqconstr}--\ref{alglin:solve-LSIP-primal-ineqconstr-end} of $\SolveLSIPPrimal$
that each of
$\BIa_1,\ldots,\BIa_{p}\in\R^{q}$ corresponds to an index $(i,j,l)\in\widehat{\CI}_{+}$.
Since $\CV\subseteq \widehat{\CI}_{+}$,
it holds that 
$\sum_{t=1}^{p}\|\BIa_{t}\|_{1}
=\big(\sum_{i=1}^{N}(1+2m_i+2m_0)\big)+3(p-q)
=3p-q-N$.
Moreover, Line~\ref{alglin:solve-LSIP-primal-balanceconstr} implies that 
$\sum_{l=1}^{m_0}\|\BIa_{p+l}\|_{1}=Nm_0$.
Thus, the sum of all entries in the matrix $\BA$
is equal to 
$\sum_{t=1}^{p+m_0}\|\BIa_{t}\|_{1}=3p-q-N+Nm_0< 3p$.
Third, for any $\big((\theta_{i,j,l})_{(i,j,l)\in\widehat{\CI}_{+}},(\xi_{l})_{l=1:m_0}\big)$ that is feasible for 
(\ref{eqn:mt-tf-lp-any-index-primal}),
it holds that 
\begin{align*}
    \sum_{(i,j,l)\in\widetilde{\CI}_{+}}\theta_{i,j,l}
    &=\sum_{i=1}^{N}\sum_{(j,l)\in\widehat{\CI}_{i}}\theta_{i,j,l}=N,\\
    \sum_{\bar{l}=1}^{m_0}\xi_{\bar{l}}
    &=\sum_{\bar{l}=1}^{m_0}\Bigg(\sum_{(j,l)\in\widehat{\CI}_{1}}\theta_{1,j,l}\INDI_{\{l=\bar{l}\}}\Bigg)=\sum_{(j,l)\in\widehat{\CI}_1}\theta_{1,j,l}=1.
\end{align*}
This proves (\ref{eqn:proc-solve-LSIP-proof-primal-L1bound}).

Now, 
let us denote by $\Bpsi:=(\Btheta^\TRANSP,\Bxi^\TRANSP)\in\R^{p+m_0}_{+}$
the output of \citep[Algorithm~2]{ECvandenbrand2020deterministic}
with input $\big(\BA,\bar{\BIg},\tilde{\BIc},\frac{\varsigma}{4(N+1)(p+1)}\big)$ in Line~\ref{alglin:solve-LSIP-primal-solver} of $\SolveLSIPPrimal$,
where $\Btheta\in\R_{+}^{p}$, $\Bxi\in\R_{+}^{m_0}$.
Since the optimal value of (\ref{eqn:mt-tf-lp-any-index-primal}) 
is equal to $\tau^\star_{\widehat{\CI}_{+}}$ by the strong duality of LP problems,
\citep[Theorem~1.1]{ECvandenbrand2020deterministic} with respect to 
$A\leftarrow \BA$, 
$b\leftarrow \bar{\BIg}$,
$c\leftarrow\tilde{\BIc}$,
$n\leftarrow p+m_0$,
$R\leftarrow N+1$,
$\delta\leftarrow \frac{\varsigma}{4(N+1)(p+1)}$
implies that 
\begin{align}
    \langle\tilde{\BIc},\Bpsi\rangle&\le \tau^\star_{\widehat{\CI}_{+}} + \frac{(N+1) \|\tilde{\BIc}\|_{\infty} \varsigma}{4(N+1)(p+1)}< \tau^\star_{\widehat{\CI}_{+}} + \frac{\varsigma}{4}(c_{\ssmax}-c_{\ssmin}), \label{eqn:proc-solve-LSIP-proof-primal-objerr}\\
    \begin{split}
    {\|\BA\Bpsi-\bar{\BIg}\|_1}&\le \frac{\varsigma}{4(N+1)(p+1)} \Bigg((N+1)\Bigg(\sum_{t=1}^{p+m_0}\|\BIa_{t}\|_{1}\Bigg)+\|\bar{\BIg}\|_1\Bigg)\\
    &< \frac{\varsigma}{4(N+1)(p+1)}\big(3p(N+1)+2N\big)< \frac{3\varsigma}{4}.
    \end{split}\label{eqn:proc-solve-LSIP-proof-primal-constrerr}
\end{align}
In the following, let us denote 
$\BA\Bpsi-\bar{\BIg}=(\zeta_{1,0},\zeta_{1,1},\ldots,\zeta_{1,m_1},\omega_{1,1},\ldots,\omega_{1,m_0},\ldots,\zeta_{N,0},\zeta_{N,1},\ldots,\zeta_{N,m_N},\allowbreak\omega_{N,1},\ldots,\omega_{N,m_0})^\TRANSP$ and
denote $\Btheta=(\theta_{i,j,l})^\TRANSP_{(i,j,l)\in\widehat{\CI}_{+}}$,
$\Bxi=(\xi_{1},\ldots,\xi_{m_0})^\TRANSP$.
Subsequently, the definitions of $\BA$ and $\bar{\BIg}$ yield
\begin{align}
    \zeta_{i,0}&= \Bigg(\sum_{(j,l)\in\widehat{\CI}_i}\theta_{i,j,l}\Bigg) - 1 &  \forall 1\le i\le N, \label{eqn:proc-solve-LSIP-proof-primal-constrerr-const}\\
    \zeta_{i,\bar{j}}&=\Bigg(\sum_{(j,l)\widehat{\CI}_i}\theta_{i,j,l} g_{i,\bar{j}}(\BIv_{i,j})\Bigg) - \int_{\CX_i}g_{i,\bar{j}}\DIFFX{\mu_i} & \forall 1\le \bar{j}\le m_i,\; \forall 1\le i\le N, \label{eqn:proc-solve-LSIP-proof-primal-constrerr-gi}\\
    \omega_{i,\bar{l}}&=\Bigg(\sum_{(j,l)\in\widehat{\CI}_i}\theta_{i,j,l} g_{0,\bar{l}}(\BIv_{0,l})\Bigg)- \xi_{\bar{l}} & \forall 1\le \bar{l}\le m_0,\; \forall 1\le i\le N.\label{eqn:proc-solve-LSIP-proof-primal-constrerr-g0}
\end{align}
It follows from (\ref{eqn:proc-solve-LSIP-proof-primal-constrerr}) 
and the assumption $\varsigma\in(0,1)$ that 
\begin{align}
    \sum_{i=1}^N\Bigg(|\zeta_{i,0}|+\Bigg(\sum_{\bar{j}=1}^{m_i}|\zeta_{i,\bar{j}}|\Bigg)+\Bigg(\sum_{\bar{l}=1}^{m_0}|\omega_{i,\bar{l}}|\Bigg)\Bigg) &=\|\BA\Bpsi-\bar{\BIg}\|_{1}<\frac{3\varsigma}{4}<1.\label{eqn:proc-solve-LSIP-proof-primal-constrerr-sum}
\end{align}
In particular, (\ref{eqn:proc-solve-LSIP-proof-primal-constrerr-const}) 
and (\ref{eqn:proc-solve-LSIP-proof-primal-constrerr-sum})
show that 
$\sum_{(j,l)\in\widehat{\CI}_i}\theta_{i,j,l}=1+\zeta_{i,0}>0$ $\forall 1\le i\le N$.

Next, 
it holds by Lines~\ref{alglin:solve-LSIP-primal-normalize}--\ref{alglin:solve-LSIP-primal-buildmeas} of $\SolveLSIPPrimal$ 
that $\hat{\theta}_i\in\CP(\CX_i\times\CX_0)$,
$\support(\hat{\theta}_{i})\subseteq V(\FS_i)\times V(\FS_{0})$ 
$\forall 1\le i\le N$.
For $i=1,\ldots,N$, let $\hat{\mu}_i$ and $\hat{\nu}_i$ denote the marginals of $\hat{\theta}_i$ on $\CX_i$ and $\CX_0$, respectively.
Using 
(\ref{eqn:proc-solve-LSIP-proof-primal-constrerr-const}), 
(\ref{eqn:proc-solve-LSIP-proof-primal-constrerr-gi}), and 
(\ref{eqn:proc-solve-LSIP-proof-primal-constrerr-sum}),
we get 
\begin{align*}
    \sum_{\bar{j}=1}^{m_i}\bigg|\int_{\CX_i}g_{i,\bar{j}}\DIFFX{\hat{\mu}_i} - \int_{\CX_i}g_{i,\bar{j}}\DIFFX{\mu}_i\bigg|
    &\le \sum_{\bar{j}=1}^{m_i}\left(\Bigg|\int_{\CX_i}g_{i,\bar{j}}\DIFFX{\hat{\mu}_i} - \Bigg(\sum_{(j,l)\in \widehat{\CI}_i}\theta_{i,j,l}g_{i,\bar{j}}(\BIv_{i,j})\Bigg)\Bigg| + |\zeta_{i,\bar{j}}|\right)\\
    &=\sum_{\bar{j}=1}^{m_i}\left(\frac{|\zeta_{i,0}|}{1+\zeta_{i,0}}\Bigg(\sum_{(j,l)\in \widehat{\CI}_i}\theta_{i,j,l}g_{i,\bar{j}}(\BIv_{i,j})\Bigg)+|\zeta_{i,\bar{j}}|\right)\\
    &= \frac{|\zeta_{i,0}|}{1+\zeta_{i,0}}\Bigg(\sum_{(j,l)\in \widehat{\CI}_i} \sum_{\bar{j}=1}^{m_i} \theta_{i,j,l}g_{i,\bar{j}}(\BIv_{i,j})\Bigg)
    + \Bigg(\sum_{\bar{j}=1}^{m_i}|\zeta_{i,\bar{j}}|\Bigg)\\
    &\le |\zeta_{i,0}|+ \Bigg(\sum_{\bar{j}=1}^{m_i}|\zeta_{i,\bar{j}}|\Bigg)<\varsigma \qquad\qquad\qquad\qquad \forall 1\le i\le N.
\end{align*}
This shows that 
$\bar{\mu}_i\abovebelowset{\CG_i}{\varsigma}{\sim}\mu_i$ $\forall 1\le i\le N$. 
Similarly,
(\ref{eqn:proc-solve-LSIP-proof-primal-constrerr-const}) and
(\ref{eqn:proc-solve-LSIP-proof-primal-constrerr-g0})
imply that 
\begin{align*}
    \sum_{\bar{l}=1}^{m_0}\bigg|\int_{\CX_0}g_{0,\bar{l}}\DIFFX{\hat{\nu}_i} - \xi_{\bar{l}} \bigg|
    &\le \sum_{\bar{l}=1}^{m_0}\left(\Bigg|\int_{\CX_0}g_{0,\bar{l}}\DIFFX{\hat{\nu}_i} - \Bigg(\sum_{(j,l)\in \widehat{\CI}_i}\theta_{i,j,l}g_{0,\bar{l}}(\BIv_{0,l})\Bigg)\Bigg| + |\omega_{i,\bar{l}}|\right)\\
    &=\sum_{\bar{l}=1}^{m_0}\left(\frac{|\zeta_{i,0}|}{1+\zeta_{i,0}}\Bigg(\sum_{(j,l)\in \widehat{\CI}_i}\theta_{i,j,l}g_{0,\bar{l}}(\BIv_{0,l})\Bigg)+|\omega_{i,\bar{l}}|\right)\\
    &= \frac{|\zeta_{i,0}|}{1+\zeta_{i,0}}\Bigg(\sum_{(j,l)\in \widehat{\CI}_i} \sum_{\bar{l}=1}^{m_0} \theta_{i,j,l}g_{0,\bar{l}}(\BIv_{0,l})\Bigg)
    + \Bigg(\sum_{\bar{l}=1}^{m_0}|\omega_{i,\bar{l}}|\Bigg)\\
    &\le |\zeta_{i,0}|+ \Bigg(\sum_{\bar{l}=1}^{m_0}|\omega_{i,\bar{l}}|\Bigg)  \qquad\qquad\qquad\qquad\qquad \forall 1\le i\le N.
\end{align*}
This and (\ref{eqn:proc-solve-LSIP-proof-primal-constrerr-sum}) yield
\begin{align*}
    \sum_{\bar{l}=1}^{m_0}\bigg|\int_{\CX_0}g_{0,\bar{l}}\DIFFX{\hat{\nu}_i} - \int_{\CX_0}g_{0,\bar{l}}\DIFFX{\hat{\nu}_1}\bigg|
    &\le \Bigg(\sum_{\bar{l}=1}^{m_0}\bigg|\int_{\CX_0}g_{0,\bar{l}}\DIFFX{\hat{\nu}_i} - \xi_{\bar{l}}\bigg|\Bigg)
    + \Bigg(\sum_{\bar{l}=1}^{m_0}\bigg|\int_{\CX_0}g_{0,\bar{l}}\DIFFX{\hat{\nu}_1} - \xi_{\bar{l}}\bigg|\Bigg)\\
    &\le |\zeta_{i,0}|+ \Bigg(\sum_{\bar{l}=1}^{m_0}|\omega_{i,\bar{l}}|\Bigg) + |\zeta_{1,0}|+ \Bigg(\sum_{\bar{l}=1}^{m_0}|\omega_{1,\bar{l}}|\Bigg)
    < \varsigma \qquad \forall 2\le i\le N,
\end{align*}
which shows that 
$\bar{\nu}_i\abovebelowset{\CG_0}{\varsigma}{\sim}\bar{\nu}_1$ 
$\forall 1\le i\le N$.
The proof of statement~\ref{thms:proc-solve-LSIP-primal-feasibility} is now complete.
Moreover, we derive from
(\ref{eqn:proc-solve-LSIP-proof-primal-objerr}),
(\ref{eqn:proc-solve-LSIP-proof-primal-constrerr-const}), and
(\ref{eqn:proc-solve-LSIP-proof-primal-constrerr-sum})
that
\begin{align*}
    \sum_{i=1}^N\int_{\CX_i\times\CX_0}c_i\DIFFX{\hat{\theta}_i}&=\sum_{i=1}^N \frac{1}{1+\zeta_{i,0}}\Bigg(\sum_{(j,l)\in\widehat{\CI}_i}\theta_{i,j,l}c_i(\BIv_{i,j},\BIv_{0,l})\Bigg)\\
    &=\Bigg(\sum_{(i,j,l)\in\widehat{\CI}_{+}}\theta_{i,j,l}c_i(\BIv_{i,j},\BIv_{0,l})\Bigg) - \Bigg(\sum_{i=1}^{N}\frac{\zeta_{i,0}}{1+\zeta_{i,0}}\Bigg(\sum_{(j,l)\in\widehat{\CI}_i}\theta_{i,j,l}c_i(\BIv_{i,j},\BIv_{0,l})\Bigg)\Bigg)\\
    &\le \langle\tilde{\BIc},\Bpsi\rangle + (c_{\ssmax}-c_{\ssmin})\sum_{i=1}^{N}\frac{|\zeta_{i,0}|}{1+\zeta_{i,0}}\Bigg(\sum_{(j,l)\in\widehat{\CI}_i}\theta_{i,j,l}\Bigg)\\
    &< \tau^\star_{\widehat{\CI}_{+}} + \frac{\varsigma}{4}(c_{\ssmax}-c_{\ssmin}) + (c_{\ssmax}-c_{\ssmin})\sum_{i=1}^N|\zeta_{i,0}|
    < \tau^\star_{\widehat{\CI}_{+}} + (c_{\ssmax}-c_{\ssmin})\varsigma.
\end{align*}
Step~4 is now complete.

\underline{Step~5.}
Let us suppose for the sake of contradiction that 
$\tau^\star_{\widehat{\CI}_{+}}>\overline{\tau}+\epsilon_0$.
On the one hand, since 
Lemma~\ref{lem:lp-any-index}\ref{lems:lp-any-index-Linfinity-bound} ensures
$\tau^\star_{\widehat{\CI}_{+}}<\infty$,
Lemma~\ref{lem:lp-any-index}\ref{lems:lp-any-index-Linfinity-ball} implies that 
$S_{\widehat{\CI}_{+},\tau^\star_{\widehat{\CI}_{+}}-\epsilon_0}$
contains a closed $\|\cdot\|_{\infty}$-ball with radius $\frac{\epsilon_0}{6N}$,
and hence 
$\Lebesgue_{\R^q}\big(S_{\widehat{\CI}_{+},\tau^\star_{\widehat{\CI}_{+}}-\epsilon_0}\big)\ge \Lebesgue_{\R^{q}}\big(B_{\infty}\big(\frac{\epsilon_0}{6N}\big)\big)$.
On the other hand, since $\epsilon_0<1$,
Lemma~\ref{lem:lp-any-index}\ref{lems:lp-any-index-Linfinity-bound}
shows that 
$S_{\widehat{\CI}_{+},\tau^\star_{\widehat{\CI}_{+}}-\epsilon_0}\subseteq B_{\infty}(M)$.
Combining this with the conclusion of Step~1 results in 
$S_{\widehat{\CI}_{+},\tau^\star_{\widehat{\CI}_{+}}-\epsilon_0} 
= S_{\widehat{\CI}_{+},\tau^\star_{\widehat{\CI}_{+}}-\epsilon_0} \cap B_{\infty}(M)
\subseteq S_{\widehat{\CI},\overline{\tau}} \cap B_{\infty}(M)
\subseteq P_{\overline{r}}$.
This implies that 
$\Lebesgue_{\R^{q}}\big(B_{\infty}\big(\frac{\epsilon_0}{6N}\big)\big)
\le \Lebesgue_{\R^q}\big(S_{\widehat{\CI}_{+},\tau^\star_{\widehat{\CI}_{+}}-\epsilon_0}\big)
\le \Lebesgue_{\R^q}(P_{\overline{r}})$,
contradicting the conclusion of Step~2.
Therefore, we conclude that 
$\tau^\star_{\widehat{\CI}_{+}}\le \overline{\tau}+\epsilon_0$.
Lastly, combining this with the conclusions of Step~3 and Step~4
completes the proof of statement~\ref{thms:proc-solve-LSIP-gap}.
The proof is now complete.
\endproof

\subsection{Proofs of results in Section~\ref{sapx:complexity-AME}}\label{sapx:proof-complexity-AME}

Let us first present in the following lemma the properties of the auxiliary sub-procedure 
$\GetCoords$.

\begin{lemma}\label{lem:proc-get-coords}
    For any $n$-simplex $C$ in $\R^d$ with $V(C)=\{\BIv_l\}_{l=0:n}$, $n\in\{0,1,\ldots,d\}$,
    and for any $\BIx\in\nobreak \R^d$,
    $\GetCoords(C,\BIx)$ incurs 
    $O(d^\omega)$ arithmetic operations and the output
    $\big(\mathtt{Inside},\allowbreak(h_l)_{l=0:n}\big)$ satisfy the following properties.
    If $\mathtt{Inside}=\mathtt{True}$, then it holds that  
    $h_l=\lambda_{\FF(C),\BIv_l}(\BIx)$ $\forall l\in\{0,1,\ldots,n\}$.
    If $\mathtt{Inside}=\mathtt{False}$, then it holds that  
    $h_l=\nobreak0$ $\forall l\in\{0,1,\ldots,n\}$.
\end{lemma}

\proof{Proof of Lemma~\ref{lem:proc-get-coords}.}
    The number of arithmetic operations in 
    $\GetCoords(C,\BIx)$ is self-evident.
    To prove the properties of the output,
    observe that $\BIx\in C$ if and only if there exists $(\lambda_l)_{l=0:n}\subset [0,1]$ with $\sum_{l=0}^{n}\lambda_l=\nobreak1$ 
    and $\sum_{l=0}^{n}\lambda_l\BIv_l=\nobreak\BIx$,
    which translates into 
    the system of linear equations $\BV\Blambda=\big(\begin{smallmatrix}
        \BIx \\ 1
    \end{smallmatrix}\big)$
    that 
    Line~\ref{alglin:get-coords-linear-system} attempts to solve.
    Subsequently, in the case where such $(\lambda_l)_{l=0:n}$ exists, 
    it holds by definition that $\lambda_l=\lambda_{\FF(C),\BIv_l}(\BIx)$ $\forall l\in\{0,1,\ldots,n\}$.
    The proof is now complete.
\endproof

\proof{Proof of Lemma~\ref{lem:proc-construct-trans-funcs}.}
    The number of arithmetic operations in $\ConstructTransFuncs$ is self-evident.
    Since $\bigcup_{\tilde{j}=1}^{\widetilde{J}_0}S_{0,\tilde{j}}=\CX_0$ (see~\ref{enc:testfuncs}), Lemma~\ref{lem:proc-get-coords}
    ensures that 
    $\mathtt{Inside}$ in Line~\ref{alglin:eval-trans-func-get-coord} of $\EvalTransFunc$ will be 
    $\mathtt{True}$ in at least one iteration.
    The number of arithmetic operations in $\EvalTransFunc$ then follows directly from Lemma~\ref{lem:proc-get-coords}.
    Next, let us fix an arbitrary
    $i\in\{1,\ldots,N\}$,
    an arbitrary $\BIz\in\CX_0$, and an arbitrary
    $\tilde{j}\in\{1,\ldots,\widetilde{J}_0\}$ with $\BIz\in\nobreak S_{0,\tilde{j}}$.
    Then, it holds that 
    $\lambda_{\FS_0,\BIv_{0,l}}(\BIz)=\lambda_{\FF(S_{0,\tilde{j}}),\BIv_{0,l}}(\BIz)$ $\forall l\in\{\varkappa_{0,\tilde{j},\iota}\}_{\iota=0:\varrho_{0,\tilde{j}}}$, and that
    $\lambda_{\FS_0,\BIv_{0,l}}(\BIz)=\nobreak0$ $\forall l\in\nobreak\{0,1,\ldots,m_0\}\setminus \{\varkappa_{0,\tilde{j},\iota}\}_{\iota=0:\varrho_{0,\tilde{j}}}$.
    Consequently, with $w_{i,0}:=0$,
    we get 
    $\hat{\varphi}_i(\BIz)=\sum_{l=0}^{m_0}w_{i,l}\lambda_{\FS_i,\BIv_{0,l}}(\BIz)=\sum_{\iota=0}^{\varrho_{0,\tilde{j}}}w_{i,\varkappa_{0,\tilde{j},\iota}}\lambda_{\FF(S_{0,\tilde{j}}),\BIv_{0,\varkappa_{0,\tilde{j},\iota}}}(\BIz)$,
    and Lemma~\ref{lem:proc-get-coords} guarantees 
    that $\varphi=\hat{\varphi}_i(\BIz)$.
    The proof is now complete.
\endproof

The following lemma presents the properties of the auxiliary sub-procedure $\NormalizeRows$. 
\begin{lemma}\label{lem:proc-normalize-rows}
    Let $n,n'\in\N$, let $\BP\in\R_+^{n\times n'}$ satisfy 
    $\langle\vecone_n,\BP\vecone_{n'}\rangle=\nobreak 1$,
    and let $\overrightarrow{\BP}$ be the output of $\NormalizeRows(\BP)$.
    Let us denote the entries of $\BP$ by $(p_{i,j})_{i=0:n-1,\,j=0:n'-1}$,
    denote the entries of $\overrightarrow{\BP}$ by 
    $(q_{i,j})_{i=0:n-1,\,j=0:n'-1}$,
    and define $(u_0,u_1,\ldots,u_{n-1})^\TRANSP:=\BP\vecone_{n'}$.
    Then, the following statements hold.
    \begin{enumerate}[label=(\roman*),beginpenalty=10000]
        \item\label{lems:proc-normalize-rows-complexity}
        $\NormalizeRows(\BP)$ incurs $O(nn')$ arithmetic operations.

        \item\label{lems:proc-normalize-rows-condprob}
        It holds that $\sum_{j=0}^{n'-1}q_{i,j}=1$ $\forall i\in\{0,1,\ldots,n-1\}$.
        Let $(\Omega,\CF,\PROB)$ be a probability space,
        let $Q:\Omega\to\{0,1,\ldots,n-1\}$ be a random variable such that 
        $\PROB[Q=i]=u_i$ $\forall i\in\{0,1,\ldots,n-1\}$,
        and let 
        $S:\Omega\to\{0,1,\ldots,n'-1\}$ be a random variable such that 
        $\PROB[S=j|Q=i]=q_{i,j}$ $\forall j\in\{0,1,\ldots,n'-1\}$, $\forall i\in\big\{i'\in\{0,1,\ldots,n-1\}: u_{i'}> 0\big\}$.
        Then, it holds that 
        $\PROB[Q=i,\;S=j]=p_{i,j}$ $\forall i\in\{0,1,\ldots,n-\nobreak1\}$, 
        $\forall j\in\{0,1,\ldots,n'-1\}$.
        
    \end{enumerate}
\end{lemma}

\proof{Proof of Lemma~\ref{lem:proc-normalize-rows}.}
    Statement~\ref{lems:proc-normalize-rows-complexity} is evident from the procedure itself.
    In statement~\ref{lems:proc-normalize-rows-condprob}, 
    $\sum_{j=0}^{n'-1}q_{i,j}=1$ $\forall i\in\{0,1,\ldots,n-1\}$ follows from Line~\ref{alglin:normalize-rows-normalize} and Line~\ref{alglin:normalize-rows-degenerate}.
    The remainder of statement~\ref{lems:proc-normalize-rows-condprob} follows from elementary properties of conditional probability.
    The proof is complete.
\endproof

\proof{Proof of Lemma~\ref{lem:proc-construct-allocs}.}
    To begin, statement~\ref{lems:proc-construct-allocs-complexity} follows from combining 
    the number of arithmetic operations incurred by $\NormalizeRows$ in
    Lemma~\ref{lem:proc-normalize-rows}\ref{lems:proc-normalize-rows-complexity}
    and the number of arithmetic operations incurred by the rest of the lines
    of $\ConstructAllocs$,
    as presented in the comments.

    In the following, let us divide the proofs of statements~\ref{lems:proc-construct-allocs-samp-qual}--\ref{lems:proc-construct-allocs-error}
    into 4 steps.
    \begin{itemize}
        \item Step~1: explicitly constructing the random variables 
        $\BIZ:\Omega\to\CX_0$, 
        $(\BIZ_i:\Omega\to\CX_0)_{i=1:N}$,
        $(\BIX_i:\Omega\to\nobreak\CX_i)_{i=1:N}$, 
        $(\overbar{\BIX}_i:\Omega\to\CX_i)_{i=1:N}$,
        as well as the auxiliary random variables 
        $(\BIX^\dagger_i:\Omega\to\CX_i)_{i=1:N}$.

        \item Step~2: showing that for $i=1,\ldots,N$, 
        the law of $(\BIX_i,\BIZ_i):\Omega\to \CX_i\times\CX_0$
        is $\hat{\theta}_i$,
        the law of $\overbar{\BIX}_i$ is~$\mu_i$,
        and that 
        $\EXP\big[\|\BIZ-\BIZ_i\|_2\big]\le \overline{\epsilon}_0(\varsigma)$,
        $\EXP\big[\|\BIX_i-\overbar{\BIX}_i\|_2\big]\le \overline{\epsilon}_i(\varsigma)$.

        \item Step~3: showing that the law of $\BIZ$ and the law of $\widetilde{\BIZ}$ in statement~\ref{lems:proc-construct-allocs-samp-qual} are both identical to $\hat{\nu}$ and proving statement~\ref{lems:proc-construct-allocs-samp-qual}.
        
        \item Step~4: showing that for $i=1,\ldots,N$, the law of $(\overbar{\BIX}_i,\BIZ):\Omega\to\CX_i\times\CX_0$ and the law of $(\widetilde{\BIX}_i,\widetilde{\BIZ}):\Omega\to\CX_i\times\CX_0$ in
        statement~\ref{lems:proc-construct-allocs-samp-alloc} are identical,
        and proving statements~\ref{lems:proc-construct-allocs-samp-alloc} and \ref{lems:proc-construct-allocs-error}.
    \end{itemize}

    \underline{Step~1.}
    For $i=1,\ldots,N$, 
    let 
    $\BTheta_i\in\R_{+}^{(m_i+1)\times(m_0+1)}$,
    $\Bmu_i\in\R_{+}^{m_i+1}$,
    $\hat{\Bmu}_i\in\R_{+}^{m_i+1}$,
    $\hat{\Bmu}_{i,\ssmin}\in\R_{+}^{m_i+1}$,
    $\hat{\Bnu}_i\in\nobreak\R_{+}^{m_0+1}$,
    $\hat{\Bnu}_{i,\ssmin}\in\R_{+}^{m_0+1}$,
    $\BA_i,\overrightarrow{\BA}_i\in\nobreak\R_{+}^{(m_i+1)\times(m_i+1)}$,
    $\BB_i,\overrightarrow{\BB}_i\in\nobreak\R_{+}^{(m_0+1)\times(m_0+1)}$,
    $\overleftarrow{\BTheta}_i\in\nobreak\R_{+}^{(m_0+1)\times(m_i+1)}$,
    $\BP_i,\overrightarrow{\BP}_i\in\nobreak\R_{+}^{(m_i+1)\times \widetilde{R}_i}$, and
    $\BT_i\in\nobreak\R_{+}^{(m_0+1)\times \widetilde{R}_i}$
    be computed in 
    Lines~\ref{alglin:construct-allocs-begin}--\ref{alglin:construct-allocs-matprod} of $\ConstructAllocs$,
    denote 
    $\Bmu_i=(\mu_{i,0},\mu_{i,1},\ldots,\mu_{i,m_i})^\TRANSP$,
    $\hat{\Bmu}_i=(\hat{\mu}_{i,0},\hat{\mu}_{i,1},\ldots,\hat{\mu}_{i,m_i})^\TRANSP$,
    $\hat{\Bnu}_i=(\hat{\nu}_{i,0},\hat{\nu}_{i,1},\ldots,\hat{\nu}_{i,m_0})^\TRANSP$,
    and denote the entries of the matrices 
    $\BTheta_i$, $\overleftarrow{\BTheta}_i$,
    $\BA_i$, $\overrightarrow{\BA}_i$,
    $\BB_i$, $\overrightarrow{\BB}_i$,
    $\BP_i$, $\overrightarrow{\BP}_i$,
    $\BT_i$ by
    $\BTheta_i=(\tilde{\theta}_{i,j,l})_{j=0:m_i,\,l=0:m_0}$,
    $\overleftarrow{\BTheta}_i=(\cev{\theta}_{i,l,j})_{l=0:m_0,\,j=0:m_i}$,
    $\BA_i=(a_{i,j,j'})_{j=0:m_i,\,j'=0:m_i}$,
    $\overrightarrow{\BA}_i=(\vec{a}_{i,j,j'})_{j=0:m_i,\,j'=0:m_i}$,
    $\BB_i=(b_{i,l,l'})_{l=0:m_0,\,l'=0:m_0}$,
    $\overrightarrow{\BB}_i=(\vec{b}_{i,l,l'})_{l=0:m_0,\,l'=0:m_0}$,
    $\BP_i=(p_{i,j,\tilde{r}})_{j=0:m_i,\,\tilde{r}=1:\widetilde{R}_i}$,
    $\overrightarrow{\BP}_i=(\vec{p}_{i,j,\tilde{r}})_{j=0:m_i,\,\tilde{r}=1:\widetilde{R}_i}$,
    $\BT_i=(t_{i,l,\tilde{r}})_{l=0:m_0,\,\tilde{r}=1:\widetilde{R}_i}$.
    One checks that $\langle\Bmu_i,\vecone_{m_i+1}\rangle=\langle\hat{\Bmu}_i,\vecone_{m_i+1}\rangle=\langle\hat{\Bnu}_i,\vecone_{m_0+1}\rangle=\nobreak 1$
    $\forall i\in\{1,\ldots,N\}$,
    and it thus holds by 
    Lines~\ref{alglin:construct-allocs-coup-start}--\ref{alglin:construct-allocs-coup-end} of $\ConstructAllocs$ that
    \begin{align}
        \BA_i\vecone_{m_i+1}&=\hat{\Bmu}_i=\BTheta_i\vecone_{m_0+1}, 
        & 
        \BA_i^\TRANSP\vecone_{m_i+1}&=\Bmu_i=\BP_i\vecone_{\widetilde{R}_i}
        & \forall i\in\{1,\ldots,N\},\label{eqn:proc-construct-allocs-proof-marginals1}\\  
        \BB_i\vecone_{m_0+1}&=\hat{\Bnu}_1, 
        &
        \BB_i^\TRANSP\vecone_{m_0+1}&=\hat{\Bnu}_{i}=\BTheta_i^\TRANSP\vecone_{m_0+1}
        & \forall i\in\{1,\ldots,N\}.\label{eqn:proc-construct-allocs-proof-marginals2}
    \end{align}

    Subsequently, 
    let the random variables
    $Q:\Omega\to\{0,1,\ldots,m_0\}$,
    $\big(Q_i:\Omega\to\{0,1,\ldots,m_0\}\big)_{i=1:N}$,
    $\big(S_i:\Omega\to\{0,1,\ldots,m_i\}\big)_{i=1:N}$,
    $\big(S_i^\dagger:\Omega\to\{0,1,\ldots,m_i\}\big)_{i=1:N}$,
    $\big(\bar{S}_i:\Omega\to\{1,\ldots,\widetilde{R}_i\}\big)_{i=1:N}$
    be defined as follows:
    \begin{align*}
        \PROB[Q=l]&=\hat{\nu}_{1,l} 
            & \forall l\in\{0,1,\ldots,m_0\}&,\\
        \PROB[Q_i=l'|Q=l]&= \vec{b}_{i,l,l'} 
            & \forall l'\in\{0,1,\ldots,m_0\}, \; \forall l\in\big\{\bar{l}:\PROB[Q=\bar{l}]>0\big\},\; \forall i\in\{1,\ldots,N\}&, \\
        \PROB[S_i=j|Q_i=l]&= \cev{\theta}_{i,l,j} 
            & \forall j\in\{0,1,\ldots,m_i\}, \; \forall l\in\big\{\bar{l}:\PROB[Q_i=\bar{l}]>0\big\},\; \forall i\in\{1,\ldots,N\}&, \\
        \PROB[S_i^\dagger=j'|S_i=j]&= \vec{a}_{i,j,j'} 
            & \forall j'\in\{0,1,\ldots,m_i\}, \; \forall j\in\big\{\bar{j}:\PROB[S_i=\bar{j}]>0\big\},\; \forall i\in\{1,\ldots,N\}&, \\
        \PROB[\bar{S}_i=\tilde{r}|S_i^\dagger=j]&= \vec{p}_{i,j,\tilde{r}} 
            & \forall \tilde{r}\in\{1,\ldots,\widetilde{R}_i\}, \; \forall j\in\big\{\bar{j}:\PROB[S_i^\dagger=\bar{j}]>0\big\},\; \forall i\in\{1,\ldots,N\}&.
    \end{align*}
    It thus follows from 
    (\ref{eqn:proc-construct-allocs-proof-marginals1}),
    (\ref{eqn:proc-construct-allocs-proof-marginals2}), and
    repeated applications of Lemma~\ref{lem:proc-normalize-rows}\ref{lems:proc-normalize-rows-condprob} 
    that 
    \begin{align}
        \PROB[Q_i=l',\;Q=l]&= b_{i,l,l'} 
            & \forall l'\in\{0,1,\ldots,m_0\}, \; \forall l\in\{0,1,\ldots,m_0\},\; \forall i\in\{1,\ldots,N\},
            \label{eqn:proc-construct-allocs-proof-Qi-Q} \\
        \PROB[S_i=j,\;Q_i=l]&= \tilde{\theta}_{i,j,l} 
            & \forall j\in\{0,1,\ldots,m_i\}, \; \forall l\in\{0,1,\ldots,m_0\},\; \forall i\in\{1,\ldots,N\},
            \label{eqn:proc-construct-allocs-proof-Si-Qi}\\
        \PROB[S_i^\dagger=j',\;S_i=j]&= a_{i,j,j'} 
            & \forall j'\in\{0,1,\ldots,m_i\}, \; \forall j\in\{0,1,\ldots,m_i\},\; \forall i\in\{1,\ldots,N\},
            \label{eqn:proc-construct-allocs-proof-Si-Sidag}\\
        \PROB[\bar{S}_i=\tilde{r},\;S_i^\dagger=j]&= p_{i,j,\tilde{r}} 
            & \forall \tilde{r}\in\{1,\ldots,\widetilde{R}_i\}, \; \forall j\in\{0,1,\ldots,m_i\},\; \forall i\in\{1,\ldots,N\}.
            \label{eqn:proc-construct-allocs-proof-barSi-Sidag}
    \end{align}
    Moreover, 
    note that 
    (\ref{eqn:proc-construct-allocs-proof-barSi-Sidag}) and Lemma~\ref{lem:proc-init-type-measure}\ref{lems:proc-init-type-measure-error}
    with respect to $Q\leftarrow S_i^\dagger$,
    $S\leftarrow \bar{S}_i$ 
    guarantee that 
    \begin{align}
        \PROB[\bar{S}_i=\tilde{r}]&=\mu_i\big(\relint(F_{i,\tilde{r}})\big)
        & \forall \tilde{r}\in\{1,\ldots,\widetilde{R}_i\},\; \forall i\in\{1,\ldots,N\}.
        \label{eqn:proc-construct-allocs-proof-barSi}
    \end{align}
    
    Now, let us define the random variables
    $\BIZ:\Omega\to\CX_0$, 
    $(\BIZ_i:\Omega\to\CX_0)_{i=1:N}$,
    $(\BIX_i:\Omega\to\CX_i)_{i=1:N}$,
    $(\BIX_i^\dagger:\Omega\to\CX_i)_{i=1:N}$
    as follows:
    \begin{align*}
        \BIZ&:=\sum_{l=0}^{m_0}\BIv_{0,l}\INDI_{\{Q=l\}},
        & \BIZ_i&:=\sum_{l=0}^{m_0}\BIv_{0,l}\INDI_{\{Q_i=l\}},
        & \forall i\in\{1,\ldots,N\},\\
        \BIX_i&:=\sum_{j=0}^{m_i}\BIv_{i,j}\INDI_{\{S_i=j\}}, 
        & \BIX_i^\dagger&:=\sum_{j=0}^{m_i}\BIv_{i,j}\INDI_{\{S_i^\dagger=j\}}
        &\forall i\in\{1,\ldots,N\}.
    \end{align*}
    Lastly, let 
    $(\overbar{\BIX}_i:\Omega\to\CX_i)_{i=1:N}$
    be random variables where for $i=1,\ldots,N$,
    $\overbar{\BIX}_i$ and $\BIZ$
    are conditionally independent given $\bar{S}_i$,
    such that 
    \begin{align*}
        \PROB[\overbar{\BIX}_i\in E|\bar{S}_i=\tilde{r}] &= \frac{\mu_i\big(E\cap \relint(F_{i,\tilde{r}})\big)}{\mu_i\big(\relint(F_{i,\tilde{r}})\big)} 
        \qquad \forall E\in\CB(\CX_i),\; \forall \tilde{r}\in\big\{\tilde{r}':\PROB[\bar{S}_i=\tilde{r}']>0\big\},\; \forall i\in\{1,\ldots,N\}.
    \end{align*}
    Observe that $(\overbar{\BIX}_i)_{i=1:N}$ are well-defined because of 
    (\ref{eqn:proc-construct-allocs-proof-barSi}).
    This completes Step~1.

    \underline{Step~2.}
    Let us fix an arbitrary $i\in\{1,\ldots,N\}$ in this step.
    Observe that Lines~\ref{alglin:construct-allocs-begin}--\ref{alglin:construct-allocs-assignmatrix} of $\ConstructAllocs$
    and the encoding of $\hat{\theta}_i$ in \ref{enc:par-primal}
    imply that
    \begin{align*}
        \PROB[\BIX_i=\BIx,\;\BIZ_i=\BIz]
        &=\sum_{j=0}^{m_i}\sum_{l=0}^{m_0}\tilde{\theta}_{i,j,l}\INDI_{\{\BIx=\BIv_{i,j},\,\BIz=\BIv_{0,l}\}}
        =\sum_{\iota=1}^{n_i}\hat{\theta}_{i,\iota}\INDI_{\{\BIx=\BIv_{i,j_{i,\iota}},\,\BIz=\BIv_{0,l_{i,\iota}}\}}\\
        &=\hat{\theta}_i\big(\{(\BIx,\BIz)\}\big)
        \qquad\qquad\qquad\qquad\qquad \forall (\BIx,\BIz)\in\CX_i\times\CX_0,
    \end{align*}
    and thus the law of 
    $(\BIX_i,\BIZ_i):\Omega\to\CX_i\times\CX_0$ is equal to 
    $\hat{\theta}_i$.
    Moreover,
    in view of (\ref{eqn:proc-construct-allocs-proof-barSi-Sidag}),
    applying Lemma~\ref{lem:proc-init-type-measure}\ref{lems:proc-init-type-measure-error}
    with respect to $Q\leftarrow S_i^\dagger$,
    $S\leftarrow \bar{S}_i$ 
    yields that the law of $\overbar{\BIX}_i$ is $\mu_i$,
    and that 
    \begin{align}
        \EXP\big[\|\BIX_i^\dagger-\overbar{\BIX}_i\|_2\big]\le \eta(\FS_{i}).
        \label{eqn:proc-construct-allocs-proof-bound-Xidag-barXi}
    \end{align}

    Next, 
    to prove $\EXP\big[\|\BIZ-\BIZ_i\|_2\big]\le \overline{\epsilon}_0(\varsigma)$ and
    $\EXP\big[\|\BIX_i-\overbar{\BIX}_i\|_2\big]\le \overline{\epsilon}_i(\varsigma)$,
    let us define 
    \begin{align*}
        \BD_i&:=\big(\|\BIv_{i,j}-\BIv_{i,j'}\|_2\big)_{j=0:m_i,\,j'=0:m_i}\in\R_{+}^{(m_i+1)\times (m_i+1)},\\
        \BD_0&:=\big(\|\BIv_{0,l}-\BIv_{0,l'}\|_2\big)_{l=0:m_0,\,l'=0:m_0}\in\R_{+}^{(m_0+1)\times (m_0+1)}, \\
        \overline{D}_i&:=\max_{0\le j<j'\le m_i}\big\{\|\BIv_{i,j}-\BIv_{i,j'}\|_2\big\}=\max_{\BIx,\BIx'\in\CX_i}\big\{\|\BIx-\BIx'\|_2\big\},\\
        \overline{D}_0&:=\max_{0\le l<l'\le m_0}\big\{\|\BIv_{0,l}-\BIv_{0,l'}\|_2\big\}=\max_{\BIx,\BIx'\in\CX_0}\big\{\|\BIx-\BIx'\|_2\big\}.
    \end{align*}
    On the one hand, notice that in the case where $\hat{\Bmu}_i=\Bmu_i$,
    (\ref{eqn:proc-construct-allocs-proof-Si-Sidag})
    and Line~\ref{alglin:construct-allocs-coup-same-first} of $\ConstructAllocs$ guarantee that
    $\EXP\big[\|\BIX_i-\nobreak\BIX_i^\dagger\|_2\big]=0=\frac{\overline{D}_i}{2}\|\hat{\Bmu}_i-\Bmu_i\|_1$.
    In the case where $\hat{\Bmu}_i\ne\Bmu_i$,
    since 
    $\|\hat{\Bmu}_i-\nobreak\hat{\Bmu}_{i,\ssmin}\|_1
    =\|\Bmu_i-\nobreak\hat{\Bmu}_{i,\ssmin}\|_1
    =1 - \langle\hat{\Bmu}_{i,\ssmin},\vecone_{m_i+1}\rangle
    =\frac{1}{2}\|\hat{\Bmu}_i-\Bmu_i\|_1$,
    we have by 
    (\ref{eqn:proc-construct-allocs-proof-Si-Sidag})
    and Lines~\ref{alglin:construct-allocs-coup-diff-first-1}--\ref{alglin:construct-allocs-coup-diff-first-2} of $\ConstructAllocs$ that
    \begin{align*}
        \EXP\big[\|\BIX_i-\BIX_i^\dagger\|_2\big]
        &=\sum_{j=0}^{m_i}\sum_{j'=0}^{m_i}a_{i,j,j'}\|\BIv_{i,j}-\BIv_{i,j'}\|_2
        =\tr(\BA_i^\TRANSP\BD_i)\\
        &= \tr\big(\diag(\hat{\Bmu}_{i,\ssmin})\BD_i\big)+\big(1-\langle\hat{\Bmu}_{i,\ssmin},\vecone_{m_i+1}\rangle\big)^{-1}\langle\hat{\Bmu}_i-\hat{\Bmu}_{i,\ssmin},\BD_i(\Bmu_i-\hat{\Bmu}_{i,\ssmin})\rangle\nonumber\\
        &\le \big(1-\langle\hat{\Bmu}_{i,\ssmin},\vecone_{m_i+1}\rangle\big)^{-1} \|\hat{\Bmu}_i-\hat{\Bmu}_{i,\ssmin}\|_{1}\|\BD_i(\Bmu_i-\hat{\Bmu}_{i,\ssmin})\|_{\infty}\nonumber\\
        &\le \overline{D}_i\|\Bmu_i-\hat{\Bmu}_{i,\ssmin}\|_1
        =\frac{\overline{D}_i}{2}\|\hat{\Bmu}_i-\Bmu_i\|_1,\nonumber
    \end{align*}
    Similarly,
    in the case where $\hat{\Bnu}_1=\hat{\Bnu}_i$,
    (\ref{eqn:proc-construct-allocs-proof-Qi-Q}) and Line~\ref{alglin:construct-allocs-coup-same-second} guarantee that 
    $\EXP\big[\|\BIZ-\BIZ_i\|_2\big]=0=\frac{\overline{D}_0}{2}\|\hat{\Bnu}_1-\hat{\Bnu}_i\|_1$.
    In the case where $\hat{\Bnu}_1\ne\hat{\Bnu}_i$,
    since
    $\|\hat{\Bnu}_1-\nobreak\hat{\Bnu}_{i,\ssmin}\|_1
    =\|\hat{\Bnu}_i-\hat{\Bnu}_{i,\ssmin}\|_1
    =1-\langle\hat{\Bnu}_{i,\ssmin},\vecone_{m_0+1}\rangle
    =\frac{1}{2}\|\hat{\Bnu}_1-\hat{\Bnu}_i\|_1$,
    we have by (\ref{eqn:proc-construct-allocs-proof-Qi-Q})
    and Lines~\ref{alglin:construct-allocs-coup-diff-second-1}--\ref{alglin:construct-allocs-coup-end}
    of $\ConstructAllocs$ that
    \begin{align*}
        \EXP\big[\|\BIZ-\BIZ_i\|_2\big]
        &=\sum_{l=0}^{m_0}\sum_{l'=0}^{m_0}b_{i,l,l'}\|\BIv_{0,l}-\BIv_{0,l'}\|_2
        =\tr(\BB_i^\TRANSP\BD_0)\\
        &= \tr\big(\diag(\hat{\Bnu}_{i,\ssmin})\BD_0\big)+\big(1-\langle\hat{\Bnu}_{i,\ssmin},\vecone_{m_0+1}\rangle\big)^{-1}\langle\hat{\Bnu}_1-\hat{\Bnu}_{i,\ssmin},\BD_0(\hat{\Bnu}_i-\hat{\Bnu}_{i,\ssmin})\rangle\nonumber\\
        &\le \big(1-\langle\hat{\Bnu}_{i,\ssmin},\vecone_{m_0+1}\rangle\big)^{-1} \|\hat{\Bnu}_1-\hat{\Bnu}_{i,\ssmin}\|_{1}\|\BD_0(\hat{\Bnu}_i-\hat{\Bnu}_{i,\ssmin})\|_{\infty}\nonumber\\
        &\le \overline{D}_0\|\hat{\Bnu}_i-\hat{\Bnu}_{i,\ssmin}\|_1
        =\frac{\overline{D}_0}{2}\|\hat{\Bnu}_1-\hat{\Bnu}_i\|_1.\nonumber
    \end{align*}
    On the other hand, 
    recalling that the marginals of $\hat{\theta}_i$ on $\CX_i$ and $\CX_0$ are $\hat{\mu}_i$ and $\hat{\nu}_i$, respectively,
    and that the marginals of $\hat{\theta}_1$ on $\CX_1$ and $\CX_0$ are $\hat{\mu}_1$ and $\hat{\nu}_1$, respectively,
    we have by Lemma~\ref{lem:proc-init-type-measure}\ref{lems:proc-init-type-measure-integrals} and 
    Lines~\ref{alglin:construct-allocs-begin}--\ref{alglin:construct-allocs-margvec} of $\ConstructAllocs$ that
    \begin{align*}
        \mu_{i,j}&=\int_{\CX_i}g_{i,j}\DIFFX{\mu_i} 
        & \forall j\in\{0,1,\ldots,m_i\}, \\
        \hat{\mu}_{i,j}&=\sum_{l=0}^{m_0}\tilde{\theta}_{i,j,l}
        =\sum_{\iota=1}^{n_i}\hat{\theta}_{i,\iota}\INDI_{\{j_{i,\iota}=j\}}
        =\sum_{\iota=1}^{n_i}\hat{\theta}_{i,\iota}g_{i,j}(\BIv_{i,j_{i,\iota}})
        =\int_{\CX_i}g_{i,j}\DIFFX{\hat{\mu}_i}  
        & \forall j\in\{0,1,\ldots,m_i\}, \\
        \hat{\nu}_{1,l}&=\sum_{j=0}^{m_1}\tilde{\theta}_{1,j,l}
        =\sum_{\iota=1}^{n_1}\hat{\theta}_{1,\iota}\INDI_{\{l_{1,\iota}=l\}}
        =\sum_{\iota=1}^{n_1}\hat{\theta}_{1,\iota}g_{0,l}(\BIv_{0,l_{1,\iota}})
        =\int_{\CX_0}g_{0,l}\DIFFX{\hat{\nu}_1}  
        & \forall l\in\{0,1,\ldots,m_0\}, \\
        \hat{\nu}_{i,l}&=\sum_{j=0}^{m_i}\tilde{\theta}_{i,j,l}
        =\sum_{\iota=1}^{n_i}\hat{\theta}_{i,\iota}\INDI_{\{l_{i,\iota}=l\}}
        =\sum_{\iota=1}^{n_i}\hat{\theta}_{i,\iota}g_{0,l}(\BIv_{0,l_{i,\iota}})
        =\int_{\CX_0}g_{0,l}\DIFFX{\hat{\nu}_i}  
        & \forall l\in\{0,1,\ldots,m_0\}.
    \end{align*}
    Subsequently,
    since 
    $\sum_{j=0}^{m_i}g_{i,j}(\BIx)=\nobreak 1$ $\forall \BIx\in\CX_i$
    and 
    $\sum_{l=0}^{m_0}g_{0,l}(\BIz)=\nobreak 1$ $\forall \BIz\in\CX_0$,
    the assumptions
    $\hat{\mu}_i\abovebelowset{\CG_i}{\varsigma}{\sim}\mu_i$, 
    $\hat{\nu}_i\abovebelowset{\CG_0}{\varsigma}{\sim}\hat{\nu}_1$
    imply that
    \begin{align*}
        \|\hat{\Bmu}_i-\Bmu_i\|_1
        &=\sum_{j=0}^{m_i}\bigg|\int_{\CX_i}g_{i,j}\DIFFX{\hat{\mu}_i}-\int_{\CX_i}g_{i,j}\DIFFX{\mu_i}\bigg|
        \le 2\sum_{j=1}^{m_i}\bigg|\int_{\CX_i}g_{i,j}\DIFFX{\hat{\mu}_i}-\int_{\CX_i}g_{i,j}\DIFFX{\mu_i}\bigg|
        \le 2\varsigma,\\
        \|\hat{\Bnu}_1-\hat{\Bnu}_i\|_1
        &\le \sum_{l=0}^{m_0}\bigg|\int_{\CX_0}g_{0,l}\DIFFX{\hat{\nu}_1}-\int_{\CX_0}g_{0,l}\DIFFX{\hat{\nu}_i}\bigg|
        \le 2\sum_{l=1}^{m_0}\bigg|\int_{\CX_0}g_{0,l}\DIFFX{\hat{\nu}_1}-\int_{\CX_0}g_{0,l}\DIFFX{\hat{\nu}_i}\bigg|
        \le 2\varsigma.
    \end{align*}
    Combining these inequalities with 
    (\ref{eqn:proc-construct-allocs-proof-bound-Xidag-barXi}),
    $\EXP\big[\|\BIX_i-\BIX_i^\dagger\|_2\big]\le \frac{\overline{D}_i}{2}\|\hat{\Bmu}_i-\Bmu_i\|_1$, and
    $\EXP\big[\|\BIZ-\BIZ_i\|_2\big]\le \frac{\overline{D}_0}{2}\|\hat{\Bnu}_1-\hat{\Bnu}_i\|_1$
    yields 
    \begin{align*}
        \EXP\big[\|\BIX_i-\overbar{\BIX}_i\|_2\big]
        &\le \EXP\big[\|\BIX_i-\BIX^\dagger_i\|_2\big] + \EXP\big[\|\BIX_i^\dagger-\overbar{\BIX}_i\|_2\big]
        \le \overline{D}_i\varsigma + \eta(\FS_i)
        \le \overline{\epsilon}_i(\varsigma),\\
        \EXP\big[\|\BIZ-\BIZ_i\|_2\big]
        &\le \overline{D}_0\varsigma
        \le \overline{\epsilon}_0(\varsigma).
    \end{align*}
    This completes Step~2.

    \underline{Step~3.}
    Observe that 
    Line~\ref{alglin:construct-allocs-qualitymeasure} of $\ConstructAllocs$
    sets the encoding of $\hat{\nu}$ to
    $\big(m_0+\nobreak1,\allowbreak d_0,\allowbreak(\BIv_{0,l},\hat{\nu}_{1,l})_{l=0:m_0}\big)$,
    we thus get 
    $\hat{\nu}=\sum_{l=0}^{m_0}\hat{\nu}_{1,l}\delta_{\BIv_{0,l}}$.
    Moreover, we have
    $n_0=m_0+1$ and $(\BIz_{0,l},\hat{\nu}_{0,l})_{l=0:n_0-1}=(\BIv_{0,l},\hat{\nu}_{1,l})_{l=0:m_0}$ in $\SampQual$.
    In this step, we let 
    $\widetilde{Q}:=\min\big\{l\in\{0,1,\ldots,m_0\}:\sum_{s=0}^{l}\hat{\nu}_{1,s}\ge U_1\big\}$
    denote the random variable defined in $\SampQual$,
    and let 
    $\widetilde{\BIZ}:=\BIz_{0,\widetilde{Q}}
    =\BIv_{0,\widetilde{Q}}
    =\sum_{l=0}^{m_0}\BIv_{0,l}\INDI_{\{\widetilde{Q}=l\}}$
    denote the random variable in the output of 
    $\SampQual\big(\hat{\nu},(U_t)_{t\in\N}\big)$.
    We hence get
    $\PROB[\widetilde{Q}=\nobreak l]=\PROB\big[\sum_{s=0}^{l-1}\hat{\nu}_{1,s}< U_1 \le \sum_{s=0}^{l}\hat{\nu}_{1,s}\big]=\hat{\nu}_{1,l}$ $\forall l\in\nobreak\{0,1,\ldots,m_0\}$.
    Hence, the law of $Q$ and the law of $\widetilde{Q}$ are identical,
    and we conclude that 
    the law of $\BIZ$ defined in Step~1 and the law of $\widetilde{\BIZ}$ in statement~\ref{lems:proc-construct-allocs-samp-qual} are both identical
    to $\hat{\nu}$.
    Lastly, the number of arithmetic operations in $\SampQual\big(\hat{\nu},(U_t)_{t\in\N}\big)$ is $O(m_0+1)=O(m)$.
    This completes the proof of statement~\ref{lems:proc-construct-allocs-samp-qual} and Step~3.

    \underline{Step~4.}
    In this step, let us fix an arbitrary $i\in\{1,\ldots,N\}$.
    Again, we have
    $n_0=m_0+1$ and $(\BIz_{0,l},\hat{\nu}_{0,l})_{l=0:n_0-1}=(\BIv_{0,l},\hat{\nu}_{1,l})_{l=0:m_0}$ in $\SampAlloc$.
    Subsequently, we let
    $\widetilde{Q}:=\min\big\{l\in\{0,1,\ldots,m_0\}:\sum_{s=0}^{l-1}\hat{\nu}_{1,s}\ge U_1\big\}$,
    $\widetilde{S}:=\min\big\{\tilde{r}\in\{1,\ldots,\widetilde{R}_i\}:\sum_{\tilde{l}=1}^{\tilde{r}}t_{\widetilde{Q},\tilde{l}}\ge\nobreak U_2\big\}$
    denote the random variables defined in $\SampAlloc$,
    and let $\widetilde{\BIX}_i$ and
    $\widetilde{\BIZ}:=\BIz_{0,\widetilde{Q}}
    =\BIv_{0,\widetilde{Q}}
    =\sum_{l=0}^{m_0}\BIv_{0,l}\INDI_{\{\widetilde{Q}=l\}}$ 
    denote the random variables in the output of 
    $\SampAlloc\big(\hat{\gamma}_i,(U_t)_{t\in\N}\big)$.
    One can use the same argument in Step~3 to show that 
    $\PROB[\widetilde{Q}=\nobreak l]=\nobreak\hat{\nu}_{1,l}$ $\forall l\in\nobreak\{0,1,\ldots,m_0\}$.
    Similarly, for each 
    $l\in\big\{l':\PROB[\widetilde{Q}=l']>0\big\}$,
    it holds that 
    $\PROB[\widetilde{S}=\nobreak\tilde{r}|\widetilde{Q}=\nobreak l]
    =\PROB\big[\sum_{\tilde{s}=1}^{\tilde{r}-1}t_{l,\tilde{s}}< U_2 \le \sum_{\tilde{s}=1}^{\tilde{r}}t_{l,\tilde{s}}\big]=t_{l,\tilde{r}}$
    $\forall \tilde{r}\in\{1,\ldots,\widetilde{R}_i\}$.
    Moreover, since $\BT_i=\nobreak\overrightarrow{\BB}_{i}\overleftarrow{\BTheta}_i\overrightarrow{\BA}_i\overrightarrow{\BP}_i$
    by Line~\ref{alglin:construct-allocs-matprod} of $\ConstructAllocs$,
    it follows from properties of conditional probabilities that 
    the law of $(\widetilde{Q},\widetilde{S}):\Omega\to\{0,1,\ldots,m_0\}\times\{1,\ldots,\widetilde{R}_i\}$
    is identical to the law of 
    $(Q,\bar{S}_i):\Omega\to\{0,1,\ldots,m_0\}\times\{1,\ldots,\widetilde{R}_i\}$.
    Thus, (\ref{eqn:proc-construct-allocs-proof-barSi}) guarantees that 
    $\PROB[\widetilde{S}=\nobreak\tilde{r}]=\PROB[\bar{S}_i=\nobreak\tilde{r}]=\mu_i\big(\relint(F_{i,\tilde{r}})\big)$ 
    $\forall \tilde{r}\in\nobreak\{1,\ldots,\widetilde{R}_i\}$,
    which then yields 
    $\PROB\big[\mu_i\big(\relint(F_{i,\widetilde{S}})\big)>0\big]=\nobreak 1$.
    Furthermore, 
    Line~\ref{alglin:samp-alloc-condsamp} of $\SampAlloc$
    and the assumption about the oracle procedure $\mathtt{Samp}$
    in \ref{enc:typemeasures}
    guarantee that 
    $\widetilde{\BIX}_i$ and $\widetilde{\BIZ}$ are conditionally independent given $\widetilde{S}$,
    and that 
    \begin{align*}
        \PROB[\widetilde{\BIX}_i\in E|\widetilde{S}=\tilde{r}]
        =\frac{\mu_i(E \cap \relint(F_{i,\tilde{r}}))}{\mu_i(\relint(F_{i,\tilde{r}}))}
        =\PROB[\overbar{\BIX}_i\in E|\bar{S}_i=\tilde{r}]
        \qquad\forall E\in\CB(\CX_i),\; \forall \tilde{r}\in\big\{\tilde{r}':\PROB[\widetilde{S}=\tilde{r}']>0\big\}.
    \end{align*}
    Consequently, it follows that 
    \begin{align*}
        \PROB[\widetilde{\BIX}_i\in E_1,\; \widetilde{\BIZ}\in E_2]
        =\PROB[\overbar{\BIX}_i\in E_1,\; \BIZ \in E_2]
        \qquad \forall E_1\in\CB(\CX_i),\; \forall E_2\in\CB(\CX_0).
    \end{align*}
    We thus conclude that 
    the law of $(\overbar{\BIX}_i,\BIZ):\Omega\to\CX_i\times\CX_0$ defined in Step~1
    and the law of $(\widetilde{\BIX}_i,\widetilde{\BIZ}):\Omega\to\CX_i\times\CX_0$ in
    statement~\ref{lems:proc-construct-allocs-samp-alloc} are identical.
    Lastly, 
    the number of arithmetic operations in 
    $\SampAlloc\big(\hat{\gamma}_i,(U_t)_{t\in\N}\big)$
    is $O(m_0+1+\widetilde{R}_i+T_{\mathtt{Samp}})=O(\widetilde{R}+m+T_{\mathtt{Samp}})$.
    This completes the proof of statement~\ref{lems:proc-construct-allocs-samp-alloc} and Step~4.
    The proof is now complete.
\endproof

\subsection{Proofs of results in Section~\ref{sapx:complexity-particular-case}}\label{sapx:proof-complexity-particular-case}

\begin{lemma}\label{lem:proc-get-vol}
    For any $n$-simplex $C=\conv\big(\{\BIv_l\}_{l=0:n}\big)$ in $\R^d$ where $n\in\{0,1,\ldots,d\}$,
    $\GetVol(C)$ incurs 
    $O(d^3)$ arithmetic operations and returns 
    $w=\volume_n(C)$.
\end{lemma}

\proof{Proof of Lemma~\ref{lem:proc-get-vol}.}
    The number of arithmetic operations is evident.
    It follows from the 
    Gram determinant formula that
    $\volume_n(C)=\frac{1}{n!}\det(\BV^\TRANSP\BV)^{\frac{1}{2}}$.
    Since Lines~\ref{alglin:get-vol-QR}--\ref{alglin:get-vol-det} result in
    $\det(\BV^\TRANSP\BV)=\det(\BR^\TRANSP\BQ^\TRANSP\BQ\BR)=\det(\BR^\TRANSP\BR)=\det(\BR_1^\TRANSP\BR_1)=\big(\prod_{i=1}^{n}r_{i,i}\big)^2$,
    it holds that 
    $w=\volume_n(C)$.
    The proof is complete.
\endproof

\proof{Proof of Lemma~\ref{lem:proc-int-simplex}.}
    Let us first prove statement~\ref{lems:proc-int-simplex-int}.
    The number of arithmetic operations incurred by $\IntSimplex$ follows directly from 
    the number of arithmetic operations incurred by $\GetCoords$ and $\GetVol$ in
    Lemma~\ref{lem:proc-get-coords} and Lemma~\ref{lem:proc-get-vol},
    as well as the number of arithmetic operations incurred when accessing 
    the associative array $\mathtt{Map}_{\mu_i}$ discussed in \ref{enc:typemeasures-particular}.
    Let us denote 
    $\hat{\BIu}_{\iota}:=\BIu_{i,\kappa_{i,\hat{j},\iota}}$ 
    $\forall \iota\in\nobreak\{0,1,\ldots,\rho_{i,\hat{j}}\}$ for notational simplicity, which yields
    $V(K_{i,\hat{j}})=\{\hat{\BIu}_{\iota}\}_{\iota=0:\rho_{i,\hat{j}}}$, 
    and
    observe that 
    Line~\ref{alglin:int-simplex-get-coords} 
    of $\IntSimplex$
    and Lemma~\ref{lem:proc-get-coords} imply that 
    $\zeta_{s,\iota}=\lambda_{\FF(K_{i,\hat{j}}),\hat{\BIu}_{\iota}}(\BIv_{s})$ 
    $\forall s\in\nobreak\{0,1,\ldots,n\}$, 
    $\forall \iota\in\nobreak\{0,1,\ldots,\rho_{i,\hat{j}}\}$.

    Next, let us consider an arbitrary non-empty $D\in\FF(K_{i,\hat{j}})$,
    where $V(D)=\{\hat{\BIu}_{\iota_r}\}_{r=0:q}$, 
    $q\in\nobreak\{0,1,\ldots,\rho_{i,\hat{j}}\}$,
    $\{\iota_{r}\}_{r=0:q}\subseteq\{0,1,\ldots,\rho_{i,\hat{j}}\}$.
    Subsequently, one can observe that 
    $C\subseteq D$ if and only if 
    $\sum_{s=0}^{n}\lambda_{\FF(K_{i,\hat{j}}),\hat{\BIu}_{\iota}}(\BIv_s)=\nobreak 0$ $\forall \iota\in\{0,1,\ldots,\rho_{i,\hat{j}}\}\setminus \{\iota_r\}_{r=0:q}$.
    From this and Line~\ref{alglin:int-simplex-search-face},
    we can conclude that 
    $F:=K_{i,\hat{j},(\hat{\iota}_r)}=\conv\big(\{\hat{\BIu}_{\hat{\iota}_r}\}_{r=0:\hat{n}}\big)$ 
    is equal to the intersection of all $D\in \FF(K_{i,\hat{j}})$ that is a superset of $C$.
    This means that $F$ is the ``smallest'' face in $\FF(K_{i,\hat{j}})$ containing~$C$.
    Moreover, \citep[Theorem~18.2]{ECrockafellar1970convex} and its proof
    imply that $F$ is also the unique face in $\FF(K_{i,\hat{j}})$ which satisfies 
    $\relint(C)\subseteq \relint(F)$.
    Thus, the assumption about $\mu_i$ in Setting~\ref{sett:Euclidean}\ref{settc:Euclidean-particular} yields
    $\int_{\CX_i}\lambda_{\FF(C),\BIv_l}\INDI_{\relint(C)}\DIFFX{\mu_i}
    =\int_{F}\lambda_{\FF(C),\BIv_l}p_{i,F}\INDI_{\relint(C)}\DIFFX{\Lebesgue_{F}}$
    $\forall l\in\nobreak\{0,1,\ldots,n\}$.
    Since $\dim(C)=n\le \hat{n}= \dim(F)$, 
    if $n<\hat{n}$ then we have 
    $\Lebesgue_{F}\big(\relint(C)\big)=\nobreak0$ and hence
    $\int_{\CX_i}\lambda_{\FF(C),\BIv_l}\INDI_{\relint(C)}\DIFFX{\mu_i}=\nobreak0$
    $\forall l\in\nobreak\{0,1,\ldots,n\}$,
    which is handled by Line~\ref{alglin:int-simplex-degenerate}.
    Otherwise, since $\Lebesgue_F$ restricted to $\CB(C)$ is equal to $\Lebesgue_{C}$, we get 
    $\int_{\CX_i}\lambda_{\FF(C),\BIv_l}\INDI_{\relint(C)}\DIFFX{\mu_i}
    =\int_{C}\lambda_{\FF(C),\BIv_l}p_{i,F}\DIFFX{\Lebesgue_{C}}$ 
    $\forall l\in\nobreak\{0,1,\ldots,n\}$.
    Subsequently, 
    applying Lemma~\ref{lem:simplex-integrals}\ref{lems:simplex-integrals-quadratic} with respect to 
    $C\leftarrow C$,
    $(\BIv_j)_{j=0:n}\leftarrow (\BIv_s)_{s=0:n}$,
    $f\leftarrow \lambda_{\FF(C),\BIv_{l}}$,
    $f'\leftarrow\nobreak p_{i,F}$
    for each $l\in\nobreak\{0,1,\ldots,n\}$
    leads to 
    $\int_{\CX_i}\lambda_{\FF(C),\BIv_l}\INDI_{\relint(C)}\DIFFX{\mu_i}
    = \frac{\volume_n(C)}{(n+2)(n+1)}\big(p_{i,F}(\BIv_{l})+\sum_{s=0}^{n}p_{i,F}(\BIv_s)\big)$
    $\forall l\in\nobreak\{0,1,\ldots,n\}$.
    Since 
    $p_{i,F}$ is affine on $F$, it holds that
    \begin{align}
        p_{i,F}(\BIv_{s})
        =\sum_{r=0}^{\hat{n}}\lambda_{\FF(K_{i,\hat{j}}),\hat{\BIu}_{\hat{\iota}_r}}(\BIv_s)p_{i,F}(\hat{\BIu}_{\hat{\iota}_r})
        =\sum_{r=0}^{\hat{n}}\zeta_{s,\hat{\iota}_r}p_{i,F}(\hat{\BIu}_{\hat{\iota}_r})
        \qquad\forall s\in\{0,1,\ldots,n\},
        \label{eqn:proc-int-simplex-proof-int-densidentity}
    \end{align}
    and thus Lines~\ref{alglin:int-simplex-get-dens}--\ref{alglin:int-simplex-calc-int} guarantee 
    that $h_l=\int_{\CX_i}\lambda_{\FF(C),\BIv_l}\INDI_{\relint(C)}\DIFFX{\mu_i}$
    $\forall l\in\nobreak\{0,1,\ldots,n\}$.
    The proof of statement~\ref{lems:proc-int-simplex-int} is now complete.

    Now, let us prove statement~\ref{lems:proc-int-simplex-samp}.
    Same as $\IntSimplex$,
    the number of arithmetic operations incurred by $\SampSimplex$ follows directly from 
    the number of arithmetic operations incurred by $\GetCoords$ in
    Lemma~\ref{lem:proc-get-coords},
    as well as the number of arithmetic operations incurred when accessing 
    the associative array $\mathtt{Map}_{\mu_i}$ discussed in \ref{enc:typemeasures-particular}.
    Next, 
    let us denote 
    $\hat{\BIu}_{\iota}:=\BIu_{i,\kappa_{i,\hat{j},\iota}}$ $\forall \iota\in\nobreak\{0,1,\ldots,\rho_{i,\hat{j}}\}$ for notational simplicity,
    and
    observe that 
    \mbox{Lines~\ref{alglin:samp-simplex-encode-outer-simplex}--\ref{alglin:samp-simplex-search-face}} of $\SampSimplex$
    are identical to 
    \mbox{Lines~\ref{alglin:int-simplex-encode-outer-simplex}--\ref{alglin:int-simplex-search-face}} of $\IntSimplex$.
    Consequently, we obtain from the proof of statement~\ref{lems:proc-int-simplex-int} that 
    $\zeta_{s,\iota}=\lambda_{\FF(K_{i,\hat{j}}),\hat{\BIu}_{\iota}}(\BIv_{s})$ 
    $\forall s\in\nobreak\{0,1,\ldots,n\}$,
    $\forall \iota\in\nobreak\{0,1,\ldots,\rho_{i,\hat{j}}\}$,
    and that
    $F:=K_{i,\hat{j},(\hat{\iota}_r)}=\conv\big(\{\hat{\BIu}_{\hat{\iota}_r}\}_{r=0:\hat{n}}\big)$
    is the unique face in $\FF(K_{i,\hat{j}})$ which satisfies 
    $\relint(C)\subseteq \relint(F)$.

    In the following, 
    let $S:\Omega\to\{0,1,\ldots,n\}$ be the random variable defined in Line~\ref{alglin:samp-simplex-samp-component},
    and let $Q_l:=\nobreak U_{(l+1)}-U_{(l)}$ $\forall l\in\{0,1,\ldots,n\}$.
    Moreover, let us define $\tilde{\BIv}_{l}:=(\BIv_{l}^\TRANSP,0)^\TRANSP\in\R^{d_i+1}$ 
    $\forall l\in\nobreak\{0,1,\ldots,n\}$,
    and for $s=0,1,\ldots,n$, let us define
    $\tilde{\BIv}_{n+1,s}:=(\BIv_{s}^\TRANSP,1)^\TRANSP\in\R^{d_i+1}$, 
    $\widetilde{C}_{s}:=\conv\big(\{\tilde{\BIv}_0,\tilde{\BIv}_1,\ldots,\tilde{\BIv}_{n},\tilde{\BIv}_{n+1,s}\}\big)\subset \R^{d_i+1}$.
    Moreover, let us define the random variables 
    $\BIX_s:=\BIv_{s}+\sum_{l=0}^{n}Q_{l}(\BIv_l-\BIv_{s})$,
    $\widetilde{\BIX}_s:=\tilde{\BIv}_{n+1,s}+\sum_{l=0}^{n}Q_{l}(\tilde{\BIv}_{l}-\tilde{\BIv}_{n+1,s})$.
    Observe that 
    $U_{(n+1)}=\sum_{l=0}^{n}Q_l$
    and thus
    $\BIX=\sum_{s=0}^{n}\BIX_s\INDI_{\{S=s\}}$.
    Let us divide the remainder of this proof into the following three steps.
    \begin{itemize}
        \item Step~1: showing that $\hat{n}=n$ and that
        $\PROB[S=s]=\frac{\volume_n(C)p_{i,F}(\BIv_{s})}{(n+1)\mu_i(\relint(C))}$ 
        $\forall s\in\{0,1,\ldots,n\}$.

        \item Step~2: showing that 
        $\PROB[\BIX_s\in E]=\frac{n+1}{\volume_n(C)}\int_{C}\lambda_{\FF(C),\BIv_s}\INDI_{E}\DIFFX{\Lebesgue_C}$
        $\forall E\in\CB(\CX_i)$, $\forall s\in\nobreak\{0,1,\ldots,n\}$.

        \item Step~3: showing that 
        $\PROB[\BIX\in E]=\tilde{\mu}_{i,C}(E)$ $\forall E\in\CB(\CX_i)$.
    \end{itemize}

    \underline{Step~1.}
    Since $\dim(C)=n\le \hat{n}= \dim(F)$,
    and since 
    the assumption about $\mu_i$ in Setting~\ref{sett:Euclidean}\ref{settc:Euclidean-particular}
    implies that 
    $\mu_i\big(\relint(C)\big)=\int_{F}p_{i,F}\INDI_{\relint(C)}\DIFFX{\Lebesgue_{F}}$,
    the assumption $\mu_i\big(\relint(C)\big)>0$ necessarily guarantees that $\hat{n}=n$.
    We hence get that
    $\aff(F)=\aff(C)$, and that
    $\Lebesgue_{F}$ restricted to $\CB(C)$ is equal to $\Lebesgue_{C}$, which yield
    \begin{align}
        \mu_i\big(\relint(C)\big)=\int_{F}p_{i,F}\INDI_{\relint(C)}\DIFFX{\Lebesgue_{F}}=\int_{C}p_{i,F}\DIFFX{\Lebesgue_{C}}.
        \label{eqn:proc-int-simplex-proof-samp-normalizationidentity}
    \end{align}
    On the other hand, 
    Line~\ref{alglin:samp-simplex-calc-int} and (\ref{eqn:proc-int-simplex-proof-int-densidentity}) imply that
    \begin{align}
        q_l=p_{i,F}(\BIv_l) 
        \qquad\forall l\in\{0,1,\ldots,n\}.
        \label{eqn:proc-int-simplex-proof-samp-compweightidentity}
    \end{align}
    Applying Lemma~\ref{lem:simplex-integrals}\ref{lems:simplex-integrals-affine} 
    with respect to 
    $C\leftarrow C$,
    $(\BIv_{j})_{j=0:n}\leftarrow(\BIv_{s})_{s=0:n}$,
    $f\leftarrow p_{i,F}$
    in (\ref{eqn:proc-int-simplex-proof-samp-normalizationidentity})
    and then combining the result with (\ref{eqn:proc-int-simplex-proof-samp-compweightidentity})
    lead to 
    \begin{align*}
        \mu_i\big(\relint(C)\big)
        &=\frac{\volume_n(C)}{n+1}\sum_{l=0}^{n}p_{i,F}(\BIv_l)
        =\frac{\volume_n(C)}{n+1}\sum_{l=0}^{n}q_l.
    \end{align*}
    It hence follows from (\ref{eqn:proc-int-simplex-proof-samp-compweightidentity}) that
    \begin{align}
        \begin{split}
            \PROB[S=s]&=\PROB\big[\big({\textstyle\sum_{l=0}^{s-1}q_l}\big)/\big({\textstyle\sum_{l=0}^{n}q_l}\big)< U_{n+2}\le \big({\textstyle\sum_{l=0}^{s}q_l}\big)/\big({\textstyle\sum_{l=0}^{n}q_l}\big)\big] \\
            &=\frac{q_s}{\sum_{l=0}^{n}q_l} 
            = \frac{\volume_n(C)p_{i,F}(\BIv_{s})}{(n+1)\mu_i\big(\relint(C)\big)}
            \hspace{40pt}\qquad \forall s\in\{0,1,\ldots,n\}.
        \end{split}
        \label{eqn:proc-samp-simplex-proof-samp-compweight}
    \end{align}
    Step~1 is now complete.

    \underline{Step~2.}
    Let us fix an arbitrary $s\in\{0,1,\ldots,n\}$ in this step.
    One observes that 
    $\tilde{\BIv}_{0},\tilde{\BIv}_1,\ldots,\tilde{\BIv}_{n},\tilde{\BIv}_{n+1,s}$ are $n+\nobreak1$ affinely independent points in $\R^{d_i+1}$,
    and hence $\widetilde{C}_{s}$ is an $(n+1)$-simplex.
    We thus have
    \begin{align*}
        \volume_{n+1}(\widetilde{C}_s)&=\frac{1}{(n+1)!}\det\big(\BV_{\widetilde{C}_s}^\TRANSP\BV_{\widetilde{C}_s}\big)^{\frac{1}{2}},\\
        \volume_{n}(C)&=\frac{1}{n!}\det\big(\BV_{C}^\TRANSP\BV_{C}\big)^{\frac{1}{2}},
    \end{align*}
    where
    \begin{align*}
        \BV_{\widetilde{C}_s}&:=\left(\begin{tabular}{cccc}
            $|$ & $|$ & & $|$ \\
            $(\tilde{\BIv}_{1}-\tilde{\BIv}_{0})$ & $(\tilde{\BIv}_{2}-\tilde{\BIv}_{0})$ & $\cdots$ & $(\tilde{\BIv}_{n+1,s}-\tilde{\BIv}_{0})$ \\
            $|$ & $|$ & & $|$
        \end{tabular}\right)\in\R^{(d_i+1)\times (n+1)},\\
        \BV_{C}&:=\left(\begin{tabular}{cccc}
            $|$ & $|$ & & $|$ \\
            $(\BIv_{1}-\BIv_{0})$ & $(\BIv_{2}-\BIv_{0})$ & $\cdots$ & $(\BIv_{n}-\BIv_{0})$ \\
            $|$ & $|$ & & $|$
        \end{tabular}\right)\in\R^{d_i\times n}.
    \end{align*}
    Now, let us denote 
    by $\BI_n\in\R^{n\times n}$ the $n$-by-$n$ identity matrix,
    denote $\BIe_0:=\veczero_{n}$,
    denote the $l$-th standard basis vector of $\R^{n}$ by $\BIe_l$ for $l=1,\ldots,n$,
    and define
    $\BW_{\widetilde{C}_s}:=\Big(\begin{smallmatrix}
        \BI_n & -\BIe_s \\ \veczero_n^\TRANSP & 1
    \end{smallmatrix}\Big)\in \R^{(n+1)\times (n+1)}$.
    One checks that
    $\det\big(\BW_{\widetilde{C}_s}\big)=1$, and that 
    $\BV_{\widetilde{C}_s}\BW_{\widetilde{C}_s}=\Big(\begin{smallmatrix}
        \BV_C & \veczero_n \\ \veczero_n^\TRANSP & 1
    \end{smallmatrix}\Big)$,
    which yields 
    $\det\big(\BV_{\widetilde{C}_s}^\TRANSP\BV_{\widetilde{C}_s}\big)
    =\det\Big(\big(\BV_{\widetilde{C}_s}\BW_{\widetilde{C}_s}\big)^\TRANSP\big(\BV_{\widetilde{C}_s}\BW_{\widetilde{C}_s}\big)\Big)
    =\det(\BV_C^\TRANSP\BV_C)$.
    We subsequently get 
    \begin{align}
        \volume_{n+1}(\widetilde{C}_s)=\frac{\volume_n(C)}{n+1}.
        \label{eqn:proc-samp-simplex-proof-samp-volumeidentity}
    \end{align}

    On the other hand, it follows from
    \citep[Theorem~2.1]{ECdevroye1986nonuniform} 
    that 
    the random variable 
    $(Q_0,Q_1,\ldots,Q_{n})^\TRANSP:\Omega\to\R^{n+1}$
    is uniformly distributed over the $(n+1)$-simplex 
    $\Delta_{n+1}:=\big\{(h_0,h_1,\ldots,h_n)^\TRANSP:\sum_{l=0}^{n}h_l\le 1\big\}\subset\R^{n+1}$.
    Since the map 
    $\Delta_{n+1}\ni(h_0,h_1,\ldots,h_n)^\TRANSP\mapsto 
    \tilde{\BIv}_{n+1,s}+\sum_{l=0}^{n}h_l(\tilde{\BIv}_{l}-\tilde{\BIv}_{n+1,s})\in \widetilde{C}_s$ 
    is affine and bijective, 
    it follows that 
    $\widetilde{\BIX}_s$ is uniformly distributed over $\widetilde{C}_s$,
    in the sense that 
    \begin{align*}
        \PROB[\widetilde{\BIX}_s\in \widetilde{E}]
        &=\frac{\Lebesgue_{\aff(\widetilde{C}_{s})}(\widetilde{E}\cap \widetilde{C}_{s})}{\volume_{n+1}(\widetilde{C}_s)}
        \qquad \forall \widetilde{E}\in \CB(\CX_i\times\R).
    \end{align*}
    In particular, 
    we can use (\ref{eqn:proc-samp-simplex-proof-samp-volumeidentity})
    to get 
    \begin{align}
        \PROB[\BIX_s\in E]
        &= \PROB[\widetilde{\BIX}_s\in (E\times \R)]
        = \frac{n+1}{\volume_n(C)}\int_{\aff(\widetilde{C}_s)}\INDI_{E\times \R}\INDI_{\widetilde{C}_s} \DIFFX{\Lebesgue_{\aff(\widetilde{C}_s)}}
        \qquad \forall E\in\CB(\CX_i).
        \label{eqn:roc-samp-simplex-proof-samp-integralidentity}
    \end{align}
    Let $\Lebesgue_{\aff(C)}\otimes \Lebesgue_{\R}$
    denote the product measure of $\Lebesgue_{\aff(C)}$ and $\Lebesgue_{\R}$.
    Since one may check that 
    $\aff(\widetilde{C}_{s})=\nobreak\aff(C)\times \R$,
    $\Lebesgue_{\aff(C)}\otimes\nobreak \Lebesgue_{\R}=\Lebesgue_{\aff(\widetilde{C}_s)}$,
    and that
    $\widetilde{C}_s$ possesses the following property:
    \begin{align*}
        \INDI_{\widetilde{C}_s}\big((\BIx^\TRANSP,r)^\TRANSP\big)
        &= \INDI_{C}(\BIx)\INDI_{[0,\lambda_{\FF(C),\BIv_s}(\BIx)]}(r) 
        \qquad \forall \BIx\in \aff(C),\; \forall r\in\R,
    \end{align*}
    applying Fubini's theorem to (\ref{eqn:roc-samp-simplex-proof-samp-integralidentity}) yields 
    \begin{align}
        \begin{split}
            \PROB[\BIX_s\in E]
            &=\frac{n+1}{\volume_n(C)} \int_{\aff(C)\times \R}\INDI_{E}(\BIx)\INDI_{\widetilde{C}_s}\big((\BIx^\TRANSP,r)^\TRANSP\big)\DIFFM{\Lebesgue_{\aff(C)}\otimes \Lebesgue_{\R}}{\DIFF\BIx,\DIFF r}\\
            &=\frac{n+1}{\volume_n(C)} \int_{\aff(C)}\INDI_{E\hspace{1pt}\cap\hspace{1pt} C}(\BIx)\bigg(\int_{\R}\INDI_{[0,\lambda_{\FF(C),\BIv_s}(\BIx)]}(r)\DIFFM{\Lebesgue_{\R}}{\DIFF r}\bigg)\DIFFM{\Lebesgue_{\aff(C)}}{\DIFF\BIx}\\
            &=\frac{n+1}{\volume_n(C)} \int_{C}\lambda_{\FF(C),\BIv_{s}}\INDI_{E}\DIFFX{\Lebesgue_{C}}
            \hspace{110pt}\qquad \forall E\in\CB(\CX_i).
        \end{split}
        \label{eqn:roc-samp-simplex-proof-samp-condprob}
    \end{align}
    Step~2 is now complete.

    \underline{Step~3.}
    Lastly,
    since $p_{i,F}$ is affine on $F$,
    we get 
    \begin{align}
        p_{i,F}(\BIx)
        &= \sum_{l=0}^{n}p_{i,F}(\BIv_l)\lambda_{\FF(C),\BIv_l}(\BIx) 
        \qquad \forall \BIx\in C.
        \label{eqn:proc-int-simplex-proof-int-densproperty}
    \end{align}
    Moreover, since $\BIX=\sum_{s=0}^{n}\BIX_s\INDI_{\{S=s\}}$,
    we get $\PROB[\BIX\in E|S=s]=\PROB[\BIX_s\in E]$ 
    $\forall E\in\CB(\CX_i)$,
    $\forall s\in\nobreak\{0,1,\ldots,n\}$.
    Combining this with 
    (\ref{eqn:proc-samp-simplex-proof-samp-compweight}), 
    (\ref{eqn:roc-samp-simplex-proof-samp-condprob}),
    and (\ref{eqn:proc-int-simplex-proof-int-densproperty}) 
    leads to 
    \begin{align*}
        \PROB[\BIX\in E]
        &=\sum_{s=0}^{n}\PROB[\BIX_s\in E]\PROB[S=s]
        =\frac{1}{\mu_i\big(\relint(C)\big)}\sum_{s=0}^{n}p_{i,F}(\BIv_s)\int_{C}\lambda_{\FF(C),\BIv_{s}}\INDI_{E}\DIFFX{\Lebesgue_{C}}\\
        &=\frac{1}{\mu_i\big(\relint(C)\big)}\int_{C}\left(\textstyle\sum_{s=0}^{n}p_{i,F}(\BIv_s)\lambda_{\sloppy{\FF(C)},\BIv_{s}}(\BIx)\right)\INDI_{E}(\BIx)\DIFFM{\Lebesgue_{C}}{\DIFF\BIx}\\
        &=\frac{1}{\mu_i\big(\relint(C)\big)}\int_{C}p_{i,F}\INDI_{E}\DIFFX{\Lebesgue_{F}}
        =\frac{1}{\mu_i\big(\relint(C)\big)}\int_{F}p_{i,F}\INDI_{E\hspace{1pt}\cap\hspace{1pt}\relint(C)}\DIFFX{\Lebesgue_{F}}\\
        &=\frac{\mu_i\big(E\cap \relint(C)\cap \relint(F)\big)}{\mu_i\big(\relint(C)\big)}
        =\frac{\mu_i\big(E\cap \relint(C)\big)}{\mu_i\big(\relint(C)\big)}
        =\tilde{\mu}_{i,C}(E)
        \qquad \forall E\in\CB(\CX_i).
    \end{align*}
    The proof is now complete.
\endproof

\subsection{Proofs of results in Section~\ref{apx:efficient-algo}}\label{sapx:proof-efficient-algo}

\proof{Proof of Proposition~\ref{prop:OT-construction}.}
In the discrete-to-discrete case, the LP formulation of the optimal transport problem is well-known (see, e.g., \citep[Section~2.3]{ECpeyre2019computational} and \citep[Section~1.3]{ECbenamou2021optimal}). 
Subsequently, by the definition of the random variable 
$\overbar{\BIX}_i$, 
it holds that 
the law of $\overbar{\BIX}_i$ is $\nu$,
and that
$\EXP\big[\|\BIX_i-\overbar{\BIX}_i\|_2\big]$
is equal to the optimal value of (\ref{eqn:discrete-OT}), which is 
equal to $W_1(\mu,\nu)$.

In the discrete-to-continuous case, it follows from the proof of \citep[Lemma~3.1]{ECneufeld2022v7numerical} and \citep[Proposition~3.4]{ECneufeld2022v7numerical} that 
the law of $\overbar{\BIX}_i$ is $\nu$,
and that
$\EXP\big[\|\BIX_i-\overbar{\BIX}_i\|_2\big]=W_1(\mu,\nu)$. 

In the one-dimensional case,
let us define 
$F_{\mu}(x):=\mu\big(\CX_i\cap (-\infty,x]\big)
=\sum_{j=1}^{n}p_j\INDI_{\{x\ge \BIx_j\}}$
$\forall x\in\CX_i$,
$F_{\mu}^{-1}(u):=\inf\big\{x\in\CX_i:F_{\mu}(x)\ge u\big\}$
$\forall u\in\nobreak[0,1]$,
$\bar{U}:=F_J+p_JU$,
and divide the proof into the following 4 steps.
\begin{itemize}
    \item Step~1: showing that 
    $F_{\sigma(j)}=\sum_{l=1}^{j-1}p_{\sigma(l)}$ for $j=1,\ldots,n$,
    $0=F_{\sigma(1)}<F_{\sigma(2)}<\cdots<F_{\sigma(n)}<1$,
    $F_{\sigma(j+1)}-\nobreak F_{\sigma(j)}=p_{\sigma(j)}$ for $j=\nobreak1,\ldots,n-1$,
    and that 
    $1-F_{\sigma(n)}=p_{\sigma(n)}$.
    
    \item Step~2: showing that $\bar{U}:\Omega\to[0,1]$ has distribution $\mathrm{Uniform}[0,1]$.
    
    \item Step~3: showing that $F_{\mu}^{-1}(\bar{U})=\BIX_i$ $\PROB$-almost surely.
    
    \item Step~4: proving that the law of $\overbar{\BIX}_i$ is $\nu$ and that 
    $\EXP\big[\|\BIX_i-\overbar{\BIX}_i\|_2\big]=W_1(\mu,\nu)$.
\end{itemize}

\underline{Step~1.}
One verifies from the definition that 
$F_{\sigma(j)}=\sum_{l=1}^{j-1}p_{\sigma(l)}$ for $j=1,\ldots,n$.
The rest of Step~1 follows directly from this identity.

\underline{Step~2.}
Using the conclusion of Step~1, one checks that
\begin{align*}
    \PROB[\bar{U}\le u]
    &=\sum_{j=1}^{n}\PROB[\bar{U}\le u| J=\sigma(j)]\PROB[J=\sigma(j)]
    =\sum_{j=1}^{n}\PROB\Big[U\le {\textstyle\frac{u-F_{\sigma(j)}}{p_{\sigma(j)}}}\Big]p_{\sigma(j)}\\
    &=\sum_{j=1}^{n}\bigg(\frac{(u-F_{\sigma(j)})^+}{p_{\sigma(j)}}\wedge 1\bigg) p_{\sigma(j)}
    =\sum_{j=1}^{n}\big((u-F_{\sigma(j)})^{+} \wedge p_{\sigma(j)}\big)
    =u
    \qquad \forall u\in[0,1].
\end{align*}
This completes Step~2.

\underline{Step~3.}
Let us fix an arbitrary $j\in\{1,\ldots,n\}$
as well as an arbitrary $u\in(F_{\sigma(j)},F_{\sigma(j)}+p_{\sigma(j)})$.
On the one hand, observe from the conclusion of Step~1 that 
$F_{\mu}(\BIx_{\sigma(j)})=\sum_{l=1}^{j}p_{\sigma(l)}=F_{\sigma(j)}+p_{\sigma(j)}>u$,
which implies that 
$F_{\mu}^{-1}(u)\le \BIx_{\sigma(j)}$.
On the other hand,
it holds for sufficiently small $\epsilon>0$ that 
$F_{\mu}(\BIx_{\sigma(j)}-\epsilon)
=\sum_{l=1}^{j-1}p_{\sigma(l)}=F_{\sigma(j)}<u$,
which implies that 
$F_{\mu}^{-1}(u)\ge \BIx_{\sigma(j)}-\epsilon$.
We thus conclude that 
$F_{\mu}^{-1}(u)=\BIx_{\sigma(j)}$ $\forall u\in(F_{\sigma(j)},F_{\sigma(j)}+p_{\sigma(j)})$.
Now, observe that on the event $\{\BIX_i=\BIx_{\sigma(j)}\}$,
it holds $\PROB$-almost surely that 
$\bar{U}\in(F_{\sigma(j)},F_{\sigma(j)}+p_{\sigma(j)})$
and hence
$F_{\mu}^{-1}(\bar{U})=\BIx_{\sigma(j)}=\BIX_i$.
Step~3 is now complete.

\underline{Step~4.}
Since $\BIX_i=F_{\mu}^{-1}(\bar{U})$, $\overbar{\BIX}_i=F_{\nu}^{-1}(\bar{U})$ 
$\PROB$-almost surely
where $\bar{U}$ is a $\mathrm{Uniform}[0,1]$ random variable,
it follows from classical results about one-dimensional optimal transport
that the law of $\overbar{\BIX}_i$ is $\nu$ and that 
$\EXP\big[\|\BIX_i-\overbar{\BIX}_i\|_2\big]=W_1(\mu,\nu)$;
see, e.g., \citep[Section~3.1]{ECrachev1998mass}.
The proof is now complete. 
\endproof

\subsection{Proofs of results in Section~\ref{apx:mmot}}\label{sapx:proof-mmot}

\proof{Proof of Corollary~\ref{cor:mmot}}
    It follows from the constructions of the random variables 
    $(\bar{X}_i)_{i=1:N}$ in Theorem~\ref{thm:mt-tf-approx} that 
    the law of $\bar{X}_i$ is $\mu_i$ for $i=1,\ldots,N$,
    and it thus holds that
    $\mu\in\Gamma(\mu_1,\ldots,\mu_N)$.
    Moreover, it follows from the constructions of the random variable 
    $\bar{Z}$ and the probability measures $(\tilde{\gamma}_i)_{i=1:N}$ in Theorem~\ref{thm:mt-tf-approx} that
    \begin{align*}
        \int_{\CX_1\times\cdots\times\CX_N}\bar{c}\DIFFX{\mu}
        &=\EXP\Big[\min_{z\in\CX_0}{\textstyle\sum_{i=1}^{N}}c_i(\bar{X}_i,z)\Big]
        =\sum_{i=1}^{N}\EXP\big[c_i(\bar{X}_i,\bar{Z})\big]
        =\sum_{i=1}^{N}\int_{\CX_i\times\CX_0}c_i\DIFFX{\tilde{\gamma}_i}.
    \end{align*}
    Since Theorem~\ref{thm:mt-tf-approx} guarantees that 
    $(\hat{\varphi}_i,\tilde{\gamma}_i)_{i=1:N},\tilde{\nu}$ is an 
    $\epsilon_{\ssapx}$-matching equilibrium, 
    it holds that 
    $\sum_{i=1}^{N}\int_{\CX_i\times\CX_0}c_i\DIFFX{\tilde{\gamma}_i} 
    \le \alpha_{\MT}^{\star} + \epsilon_{\ssapx}$.
    Subsequently, 
    \citep[Proposition~3]{ECcarlier2010matching} implies that 
    $\mu$ is an $\epsilon_{\ssapx}$-optimizer of (\ref{eqn:mt-mmot}).
    The proof is now complete.
\endproof

\proof{Proof of Corollary~\ref{cor:complexity-mmot}}
    Let $(\hat{\theta}_i)_{i=1:N}$ be constructed in Line~\ref{alglin:complexity-overall-solve-LSIP-primal} of Algorithm~\ref{alg:complexity-overall} and define
    $\varsigma:=\nobreak\frac{\epsilon}{8N(L_c\vee 1)\max_{0\le i\le N}\{\overline{D}_i\}}$,
    let $(\hat{\varphi}_i,\hat{\gamma}_i)_{i=1:N},\hat{\nu}$ 
    be the output of Algorithm~\ref{alg:complexity-overall},
    let $(\Omega,\CF,\PROB)$ be a probability space,
    let $(U_t)_{t\in\N}$ be an infinite sequence of independent 
    $\mathrm{Uniform}[0,1]$ random variables on $\Omega$,
    and let 
    $\big((\widetilde{\BIX}_i)_{i=1:N},\overline{t}\big)$
    be the output of 
    $\SampMMOT\big((\hat{\gamma}_i)_{i=1:N},(U_t)_{t\in\N}\big)$.
    Moreover, let $\hat{\mu}\in\nobreak\CP(\CX_1\times\nobreak\cdots\times\nobreak\CX_N)$ be the law of the random variable 
    $(\widetilde{\BIX}_1,\ldots,\widetilde{\BIX}_N):\Omega\to\CX_1\times\cdots\times\CX_N$.
    Then, it follows from the proofs of Lemma~\ref{lem:proc-construct-allocs}\ref{lems:proc-construct-allocs-samp-alloc} and 
    Lemma~\ref{lem:proc-construct-allocs}\ref{lems:proc-construct-allocs-error} that 
    one can construct random variables 
    $\BIZ:\Omega\to\CX_0$, 
    $(\BIZ_i:\Omega\to\CX_0)_{i=1:N}$,
    $(\BIX_i:\Omega\to\CX_i)_{i=1:N}$, and 
    $(\overbar{\BIX}_i:\Omega\to\CX_i)_{i=1:N}$ 
    such that 
    the law of $\BIZ$ is $\hat{\nu}$,
    the law of $(\overbar{\BIX}_1,\ldots,\overbar{\BIX}_N):\Omega\to\CX_1\times\cdots\times\CX_N$ is~$\hat{\mu}$,
    and for $i=1,\ldots,N$,
    the law of $(\BIX_i,\BIZ_i):\Omega\to \CX_i\times\CX_0$
    is $\hat{\theta}_i$,
    the law of $\overbar{\BIX}_i$ is $\mu_i$,
    the law of $(\overbar{\BIX}_i,\BIZ):\Omega\to\CX_i\times\CX_0$
    is $\hat{\gamma}_i$,
    and it holds that 
    $\EXP\big[\|\BIZ-\BIZ_i\|_2\big]\le \overline{\epsilon}_0(\varsigma)$,
    $\EXP\big[\|\BIX_i-\overbar{\BIX}_i\|_2\big]\le \overline{\epsilon}_{i}(\varsigma)$.
    It subsequently follows from Corollary~\ref{cor:mmot} that 
    $\hat{\mu}$ is an $\epsilon$-optimizer of~(\ref{eqn:mt-mmot}).
    The number of arithmetic operations incurred by $\SampMMOT$
    follows directly from the proofs of Lemma~\ref{lem:proc-construct-allocs}\ref{lems:proc-construct-allocs-samp-alloc} and Theorem~\ref{thm:complexity-explicit}\ref{thms:complexity-explicit-general}.
    The proof is now complete.
\endproof



\setlength{\bibsep}{1pt}
\bibliographystyle{informs2014} 


\end{document}